\documentclass[a4paper]{article}
\usepackage[latin1]{inputenc} 
\usepackage[T1]{fontenc} 
\usepackage{RR,RRthemes}
\usepackage{hyperref}
\usepackage{amssymb}
\usepackage[all]{xy}
\usepackage[latin1]{inputenc}
\usepackage{amsfonts}
\usepackage{amsmath}
\usepackage{amssymb}
\usepackage{url}
\usepackage{hyperref}%
\usepackage{graphicx}
\providecommand{\U}[1]{\protect\rule{.1in}{.1in}}
\usepackage[all]{xy}
\usepackage{color}
\RRdate{July 2011}

\newcommand{\point}{\mbox{\LARGE .}}
\def\blbx{\hbox{\vrule height 5pt width 5pt depth 0pt}\medskip}
\newcommand{\cqfd}{\hfill\blbx \\}
\newcommand{\proof}{\noindent\mbox{\bf Proof:}\\}
\usepackage{dsfont}
\newcommand{\indicator}[1]{\mathds{1}_{ {#1}}}
\def\un{\indicator{}}

\newcommand{\EE}{\mathbb{E}}

\newcommand{\LL}{\mathbb{L}}
\newcommand{\MM}{\mathbb{M}}
\newcommand{\NN}{\mathbb{N}}
\newcommand{\PP}{\mathbb{P}}
\newcommand{\QQ}{\mathbb{Q}}
\newcommand{\RR}{\mathbb{R}}

\newcommand{\Aa}{ {\cal A }}
\newcommand{\Ba}{ {\cal B }}

\newcommand{\Na}{ {\cal N }}

\newcommand{\Ea}{ {\cal E }}
\newcommand{\Sa}{ {\cal S }}

\newcommand{\Fa}{ {\cal F }}
\newcommand{\Ga}{ {\cal G }}

\newcommand{\Xa}{ {\cal X }}
\newcommand{\Ma}{ {\cal M }}
\newcommand{\Ta}{ {\cal T}}

\newcommand{\Pa}{ {\cal P }}
\newcommand{\Za}{ {\cal Z }}

\newcommand{\preuve}{\noindent\mbox{\bf Proof:}\\}
\def\blbx{\hbox{\vrule height 5pt width 5pt depth 0pt}\medskip}

\def \BB{\mathbb{B}}
\def \PP{\mathbb{P}}
\def \RR{\mathbb{R}}
\def \EE{\mathbb{E}}
\def \QQ{\mathbb{Q}}

\def \LL{\mathbb{L}}

\def \e{\epsilon}

\def\d{\delta}

\def\e{\varepsilon}

\RRauthor{
Pierre Del Moral\thanks{ Centre INRIA Bordeaux et Sud-Ouest \& Institut de Math\'ematiques de Bordeaux , Universit\'e de Bordeaux I, 351 cours de la Lib\'eration 33405 Talence cedex, France, Pierre.Del-Moral@inria.fr}, 
 Peng Hu\thanks{ Centre INRIA Bordeaux et Sud-Ouest \& Institut de Math\'ematiques de Bordeaux , Universit\'e de Bordeaux I, 351 cours de la Lib\'eration 33405 Talence cedex, France, Peng.Hu@inria.fr}, Liming Wu\thanks{Institute of Applied Mathematics, AMSS, CAS, Siyuan Building, Beijing, China \& Universit\'e Blaise Pascal, 63177 Aubi\`ere, France, wuliming@amt.ac.cn, Li-Ming.Wu@math.univ-bpclermont.fr}
}
\authorhead{Del Moral, Hu \& Wu}
\RRtitle{Sur les propri\'et\'es concentrales de processus particulaires d'interaction}
\RRetitle{On the concentration properties of Interacting particle processes}
\titlehead{On the concentration properties of Interacting particle processes}
\RRresume{Ces notes de cours pr\'esentent de nouvelles in\'egalit\'es de concentration exponentielles pour les processus empiriques en interaction associ\'es \`a des mod\`eles particulaires de type Feynman-Kac. Nous analysons diff\`erents mod\`eles stochastiques, notamment des mesures d'occupation courante de population g\'en\'etiques, des mod\`eles  historiques bas\'es sur des \'evolutions d'arbres g\'en\'ealogiques, des estimations d'\'energies libres, ainsi que des mod\`eles de cha\^ines de Markov particulaires \`a rebours. Nous illustrons ces r\'esultats avec une s\'erie d'applications li\'ees \`a la physique num\'erique et la biologie, l'optimisation stochastique, le traitement du signal et la statistique bay\'esienne, avec de nombreux algorithmes probabilistes d'apprentissage automatique. Un accent particulier est donn\'e \`a la mod\'elisation stochastique de ces algorithmes de Monte Carlo, et \`a l'analyse quantitative de leurs performances. Nous examinons notamment la convergence de filtres particulaires, des ``Island models'' de type g\'en\'etique, des processus de ponts markoviens, ainsi que diverses m\'ethodes de type MCMC en interaction.
  
  }
\RRabstract{
These lecture notes present  some new concentration inequalities for 
Feynman-Kac  particle  processes. We analyze different types of stochastic particle models, including
particle profile occupation measures,  genealogical tree based evolution models, particle free energies, as well as backward Markov chain particle models. We illustrate these results with a series of topics related to computational physics and biology, stochastic optimization, signal processing and bayesian statistics, and many other probabilistic machine learning algorithms.
 Special emphasis is given to the stochastic modeling and the quantitative performance analysis of a series of advanced Monte Carlo methods, including particle filters, genetic type island models,  Markov bridge models,
 interacting particle Markov chain Monte Carlo methodologies. }
\RRmotcle{Propri\'et\'es de concentration exponentielle, processus empiriques en interaction, interpr\'etations particulaires de formules de Feynman-Kac.}
\RRkeyword{Concentration properties, Feynman-Kac particle processes, stochastic particle models}
\RRprojet{ALEA}
\RRdomaine{1} 
\RRthemeProj{cqfd} 
\RRdomaineProjBis{cqfd} 
\RCBordeaux 

\begin{document}
 \RRNo{7677}

\makeRR   


\newtheorem{theorem}{Theorem}[section]
\newtheorem{definition}{Definition}[section]
\newtheorem{lemma}{Lemma}[section]
\newtheorem{remark}{Remark}[section]
\newtheorem{prop}{Proposition}[section]
\newtheorem{cor}{Corollary}[section]

\section{Stochastic particle methods}
\subsection{Introduction}

Stochastic particle methods have come to play a significant role in applied probability, numerical physics, Bayesian statistics, probabilistic machine learning, and engineering sciences. 

They are increasingly used to solve a variety of problems. To name a few, 
nonlinear filtering equations, data assimilation problems, rare event sampling, hidden Markov chain parameter estimation, stochastic control  problems, financial mathematics. There are also used in computational physics for  free energy computations, and Schr\"odinger operator's ground states estimation problems, as well as in computational chemistry for sampling the conformation of polymers in a given solvent.

Understanding rigorously these new particle Monte Carlo methodologies leads to fascinating mathematics related to Feynman-Kac path integral theory and their interacting particle interpretations~\cite{dmsp,Delm04,dd-2012}. In the last two decades, this line of research has been developed 
by using methods from stochastic analysis of interacting particle systems and nonlinear semigroup models in distribution spaces,
but it has also generated difficult questions that cannot be addressed without developing new mathematical tools.

Let us survey some
of the important challenges that arise. 

For numerical  applications, 
it is essential to obtain non asymptotic quantitative information on the convergence of the algorithms. Asymptotic theory, including central limit theorems, moderate deviations, and large deviations principles have clearly limited practical values.
An overview of these asymptotic results in the context of mean field and Feynman-Kac particle models can be found in the series of articles~\cite{dawson-dm,DeGu:la,DeGu:ce,dm-jacod,dm-zajic,douc}.

Furthermore, when solving a given concrete problem, it is important to obtain explicit non asymptotic 
error bounds estimates to ensure that the stochastic algorithm is provably correct. While non asymptotic
 propagation of chaos results provide some insights on the bias properties of these models, they rarely provide useful
 effective convergence rates.

Last but not least, it is essential to analyze the robustness properties, and more particularly 
the uniform performance of particle algorithms w.r.t. the time horizon.
By construction, these important questions are intimately related to the
 stability properties of complex nonlinear Markov chain semigroups associated with the limiting measure valued process.
 This line of ideas has been further developed in the articles~\cite{dawson-dm,DeGu:ont,dmsp,rio-2009}, and in the books~\cite{Delm04,dd-2012}.

Without any doubt, one of the most powerful mathematical tools to analyze the deviations of Monte Carlo based approximations
is the theory of empirical processes and measure concentration theory. In the last two decades, these new tools have become
one of the most important step forward in infinite dimensional stochastic analysis, 
advanced machine learning techniques, as well as in the development of a statistical non asymptotic theory. 

In recent years, a lot of effort has been devoted to describing the behavior of the supremum norm of
empirical functionals around the mean value of the norm.
For an overview of these subjects, we refer the reader to the seminal books of D. Pollard~\cite{pollard},
and the one of A.N. Van der Vaart and J.A. Wellner~\cite{wellner}, and the remarkable articles by E. Gin\'e~\cite{gine}, M. Ledoux~\cite{ledoux-c}, and M. Talagrand~\cite{talagrand-1,talagrand-2,talagrand-3}, and the more recent article by R. Adamczak~\cite{adamczak}. 
 The best constants in Talagrand's concentration
inequalities were obtained by Th. Klein, and E. Rio~\cite{klein}. In this article, the authors proved the functional
version of Bennett's and Bernstein's inequalities for sums of independent random variables.

To main difficulties we encountered in applying these concentration inequalities to interacting particle models
are of different order:

Firstly,  all of the concentration inequalities developed in the literature on empirical processes
still involve the mean value of the supremum norm empirical functionals. 
In practical situation, these tail style inequalities can only be used
if we have some precise information on the magnitude of the mean value of the supremum norm 
of the functionals.

On the other hand, 
the range of application the theory of empirical processes and measure concentration theory  is restricted to independent random samples, or equivalently product measures, and more recently to mixing
Markov chain models. In the reverse angle, stochastic particle techniques are not based on fully independent sequences, nor on Markov chain Monte Carlo principles, but on interacting particle samples combined with complex nonlinear Markov chain semigroups.
More precisely, besides the fact that particle models are built sequentially using conditionally independent random samples,
their respective conditional distributions are still random. In addition, they 
strongly  depend in a nonlinear way on the occupation measure of the current population.

In summary, the concentration analysis of 
interacting particle processes require the development of new stochastic perturbation style techniques to control the interaction
propagation and the independence
degree between the samples.

The first article  extending empirical processes theory to particle models is a joint work of the first author with M. Ledoux~\cite{dm-ledoux}. In this work, we proved Glivenko-Cantelli and Donsker theorems under entropy conditions, as well as non asymptotic exponential bounds for Vapnik-Cervonenkis classes of sets or functions. 
Nevertheless, in practical situations these non asymptotic  results tend to be a little disappointing, with very poor 
constants that degenerate w.r.t. the time horizon.

The second most important result on the concentration properties of mean field particle model is the recent article of the first author  with E. Rio~\cite{rio-2009}.  This article is only concerned with  finite marginal model. The authors generalize the classical Hoeffding, Bernstein and Bennett inequalities for independent random sequences to interacting particle systems. 
 
In these notes, we survey some of these results, and we provide new concentration inequalities for
interacting empirical processes. We emphasize that these lectures don't give a comprehension treatment of the theory of interacting empirical processes. To name a few missing topics, we do not discuss
large deviation principles w.r.t. the strong $\tau$-topology, Donsker type fluctuation theorems, moderate deviation principles, and continuous time models. The first two topics and developed in the monograph~\cite{Delm04},  the third one is developed in~\cite{liming-shulan}, the last one is still an open research subject. 

These notes emphasize a single stochastic perturbation method, with second order expansion entering the stability properties of the limiting Feynman-Kac semigroups. The concentration results attained are probably not the best possible of their kind. We have chosen to strive for just enough generality to derive useful and {\em uniform concentration inequalities w.r.t. the time horizon}, without having to impose complex and often unnatural regularity conditions to squeeze them into the general theory of empirical processes.

Some of the results are borrowed from the recent article~\cite{rio-2009}, and many others are new. These notes should be complemented with the book~\cite{doucet}, the books~\cite{Delm04,dd-2012}, and the article~\cite{DeGu:ont}. A very 
basic knowledge in statistics and machine learning theory will be useful, but not necessary. Good backgrounds
in Markov chain theory and in stochastic semigroup analysis are necessary.

We have made our best to give a self-contained presentation, with detailed proofs assuming however some familiarity with Feynman-Kac models, and basic facts on the theory of Markov chains on abstract state spaces. Only in section~\ref{sconvex}, we have skipped the proof of some tools from convex analysis. We hope that the essential ideas are still accessible to the readers. 

It is clearly not the scope of these lecture notes to  give an exhaustive list of references to articles 
in computational physics, engineering sciences, and machine learning 
presenting heuristic like particle algorithms to solve a specific estimation problem. Up to a few exceptions, we have only
provided references to articles with rigorous and well founded mathematical treatments on particle models.
We already
apologize for possible errors, or for references that have been omitted
due to the lack of accurate information.\\

These notes grew from series of lectures the first author gave 
 in the
Computer Science and Communications Research Unit, of the
University of Luxembourg in February and March 2011. They were reworked,
with the addition of new material on the concentration of empirical processes  
 for a course given at the Sino-French Summer Institute in Stochastic Modeling and Applications (CNRS-NSFC Joint institute of Mathematics), held at the Academy of Mathematics and System Science, Beijing, on June 2011. The Summer Institute was ably organized by Fuzhou Gong, Ying Jiao, Gilles Pag\`es, and Mingyu Xu, and the members of the scientific committee, including
 Nicole El Karoui, Zhiming Ma, Shige Peng, Liming Wu, Jia-An Yan, and Nizar Touzi. The first author
 is grateful to them for giving to him the opportunity to experiment on a receptive audience with material not entirely polished.\\
 
 In reworking the lectures, we have tried to resist the urge to push the 
 analysis to general classes of mean field particle models, in the spirit of the recent joint article with E. Rio~\cite{rio-2009}. Our principal objective has been to develop just enough analysis to handle four types of Feynman-Kac interacting particle processes; namely, genetic dynamic population models, genealogical tree based algorithms, particle free energies, as well as backward Markov chain particle models. These application models do not exhaust the possible uses of the theory developed in these lectures.

\subsection{A brief review on particle algorithms}
Stochastic particle methods belong to the class of Monte Carlo methods. They
can be thought as an universal
 particle methodology for sampling 
complex distributions in highly dimensional state spaces.

We can distinguish two different classes
of models; namely, diffusion type interacting processes, and interacting jump particle models.
Feynman-Kac particle methods belongs to the second class of models, with rejection-recycling jump
type interaction mechanisms. In contrast to conventional
acceptance-rejection type techniques, Feynman-Kac particle methods are
equipped with an adaptive and interacting recycling strategy.

The common central feature of all the Monte Carlo particle methodologies developed so far is to solve discrete generation, or continuous time integro-differential
equations in distribution spaces. The first heuristic like description of these probabilistic  techniques in mathematical physics goes back to the Los Alamos report~\cite{everett}, and the article~\cite{everett2}  by C.J. Everett and S. Ulam in 1948, and the 
short article by N. Metropolis and 
S. Ulam~\cite{ulam}, published in 1949.

In some instances, the flow of measures is dictated by the problem at hand.
In advanced signal processing, the conditional distributions of the signal given partial and noisy observations
are given by the so-called nonlinear filtering equation in distribution space (see for instance~\cite{dm96,dm98,dmsp,Delm04,dd-2012}, and references therein).  

Free energies and 
Schr\"odinger operator's ground states are 
given by the quasi-invariant distribution of a Feynman-Kac conditional distribution flow of non absorbed particles
in absorbing media. We refer the reader to the articles by E. Canc\`es, B. Jourdain and T. Leli\`evre~\cite{benj2},
M. El Makrini, B. Jourdain and T. Leli\`evre~\cite{benjamin}, M. Rousset~\cite{rousset}, the pair of articles
of the first author with L. Miclo~\cite{dm3,dmsp}, the one with A. Doucet~\cite{soft}, 
and the monograph~\cite{Delm04}, and the references therein.

In mathematical biology, branching processes and
infinite population models are also expressed by nonlinear parabolic type
integro-differential equations. Further details on this subject can be found in the articles by D.A. Dawson and his co-authors~\cite{dawson,dawson2,dawson3}, the works of E.B. Dynkin~\cite{dynkin}, and J.F. Le Gall~\cite{legall}, and more particularly the seminal book of S.N. Ethier and T.G. Kurtz~\cite{kurtz}, and the pioneering article by W. Feller~\cite{feller}.

 In other instances, we formulate a given estimation problem in terms  
 a sequence of distributions with increasing complexity on state space models 
with increasing dimension.  These stochastic evolutions can be related
to decreasing temperature schedules in Boltzmann-Gibbs measures, multilevel decompositions 
for rare event excursion models on critical level sets, decreasing subsets strategies for sampling tail style
distributions, and many other sequential importance sampling plan. For a more thorough discussion on these models
we refer the reader to~\cite{delmoraldoucetjasra2006}.

From the pure probabilistic point of view,
any flow of probability measures can be interpreted as the evolution of the laws of the random states
of a Markov process.  In contrast to conventional Markov chain models, the Markov transitions
of these chains may depend on the distribution of the current random state. The mathematical foundations of these discrete generation models have been started in 1996 in~\cite{dm96} in the context of nonlinear filtering problems. Further analysis was developed in a joint work~\cite{dmsp} of the first author with L. Miclo published in 2000. For a more thorough discussion on the origin and the performance  analysis of these discrete generation models, we also refer to the monograph~\cite{Delm04}, and the joint articles of the first author with A. Guionnet~\cite{DeGu:la,DeGu:ce,dmg-cras,DeGu:ont}, and  M. Kouritzin~\cite{dkm}.

The continuous time version of these nonlinear type Markov chain models  take their origins from the 1960s, with the development of fluid 
mechanisms and statistical physics. We refer the reader to the pioneering
works of H.P. McKean~\cite{MR36:4647,MR38:1759}, see also the more recent 
treatments by N.~Bellomo and M.~Pulvirenti~\cite{MR2001g:82084,MR2001a:82003}, the series of articles by C. Graham and S. M\'el\'eard on interacting jump models~\cite{GrMe:sto,MR2003a:82062,MR98k:35214}, the articles by
S. M\'el\'eard on {B}oltzmann equations~\cite{Me:asympt,MR99g:60103,MR2002a:60164,MR98m:82064}, and the 
lecture notes of A.S. Sznitman~\cite{Sz:topi},
and references therein.

In contrast to conventional Markov chain Monte Carlo techniques, 
these McKean type nonlinear Markov chain models can be thought as  perfect importance 
sampling strategies, in the sense that the desired {\em target measures coincide} at any time step with the law of the random
states of a Markov chain. Unfortunately, as we mentioned above, the transitions of these chains
depend on the distributions of its random states. Thus, they cannot be sampled without an additional level of approximation.
One natural solution is to use a mean field particle interpretation model.
 These stochastic techniques belong 
 to the class of stochastic population models, with free evolutions mechanisms, coupled with branching and/or 
 adaptive interacting  jumps. At any time step, the occupation measure of the population of individuals 
 approximate the solution of the nonlinear equation, when the size of the system tends to $\infty$.

 In genetic algorithms and sequential Monte Carlo  literature, the reference free evolution model is interpreted
 as a reference sequence of twisted Markov chain
samplers. These chains are used  to perform the mutation/proposal transitions. As in conventional Markov chain Monte Carlo methods, the interacting jumps are interpreted as
 an acceptance-rejection transition, equipped with sophisticated interacting and adaptive recycling
 mechanism.  In Bayesian statistics and engineering sciences, 
the resulting adaptive particle sampling model is often coined as 
 a sequential Monte Carlo algorithm, genetic procedures, or simply Sampling Importance Resampling methods,
 mainly because it is based on importance sampling plans and online approximations of a flow of probability measures.
 
Since the 1960s, the adaptive particle recycling strategy
has also been associated in biology 
and engineering science with several heuristic-like paradigms, with a proliferation of botanical 
names, depending the application area they are thought: bootstrapping, switching, replenishing, pruning, enrichment, cloning, reconfigurations, 
resampling, rejuvenation, acceptance/rejection, spawning. 

 Of course, the idea of duplicating online
better-fitted individuals and moving them one step forward to
explore state-space regions is the basis of various stochastic search
algorithms. To name a few:

{\em Particle and bootstrap filters, Rao-Blackwell particle filters, sequential Monte Carlo methods,
sequentially Interacting Markov chain Monte Carlo, genetic type search algorithms, 
Gibbs cloning search techniques, interacting simulated annealing algorithms, sampling-importance resampling methods, quantum Monte Carlo walkers, adaptive population Monte Carlo sampling models, and many others evolutionary type Monte Carlo methods}.
 
For a more detailed discussion on these models, with precise references we refer the reader to the three books~\cite{Delm04,dd-2012,doucet}.

\subsection{Feynman-Kac path integrals}

Feynman-Kac measures represent the distribution of the
paths of a Markov process, weighted by a collection of potential functions. These functional models are natural mathematical extensions of the traditional changes of probability
measures, commonly used in importance sampling
technologies, Bayesian inference, and in nonlinear filtering modeling. 

These stochastic models are defined in terms of only two
 ingredients:

A Markov chain  $X_{n}$, with Markov transition $M_{n}$ on some measurable state spaces
$(E_{n},\Ea_n)$ with initial distribution $\eta_0$, and a  sequence of $(0,1]$-valued potential functions $%
G_{n}$ on the set $E_{n}$.

The Feynman-Kac path measure associated with the
pairs $(M_{n},G_{n})$ is the probability measure $\mathbb{Q}_{n}$ on the
product state space $${\bf E}_n:=\left( E_{0}\times \ldots \times
E_{n}\right) $$ defined by the following formula 
\begin{equation}
d\mathbb{Q}_{n}:=\frac{1}{\mathcal{Z}_{n}}~\left\{ \prod_{0\leq
p<n}G_{p}(X_{p})\right\} ~d\mathbb{P}_{n}  \label{defQn}
\end{equation}%
where $\mathcal{Z}_{n}$ is a normalizing constant and $\mathbb{P}_{n}$ is
the distribution of the random paths $${\bf X_n}=(X_0,\ldots,X_n)\in {\bf E}_n$$ of the Markov
process $X_{p}$ from the origin $p=0$  with initial distribution $\eta_0$, up to the current time $p=n$. We also
denote by 
\begin{equation}\label{def-Gamma}
\Gamma _{n}=\mathcal{Z}_{n}~\mathbb{Q}_{n}
\end{equation} 
its unnormalized
version. 

The prototype model we have in head is the traditional particle absorbed Markov chain model
\begin{equation}\label{absorption-intro}  
X_{n}^c\in E^c_n:=E_n\cup\{c\}\stackrel{\tiny { absorption}~\sim (1-G_n)}{  
-\!\!\!\!-\!\!\!\!-\!\!\!\!-\!\!\!\!-\!\!\!\!-\!\!\!\!-\!\!\!\!-\!\!\!\!-\!\!\!\!-\!\!\!\!-\!\!\!\!-\!\!\!\!-\!\!\!\!-\!\!\!\!-\!\!\!\!-\!\!\!\!-\!\!\!\!  
\longrightarrow}  
\widehat{X}_n^c\stackrel{\tiny exploration ~\sim M_{n+1}}{  
-\!\!\!\!-\!\!\!\!-\!\!\!\!-\!\!\!\!-\!\!\!\!-\!\!\!\!-\!\!\!\!-\!\!\!\!-\!\!\!\!-\!\!\!\!-\!\!\!\!-\!\!\!\!-\!\!\!\!-\!\!\!\!-\!\!\!\!-\!\!\!\!-\!\!\!\!  
\longrightarrow} X_{n+1}^c 
\end{equation}

The chain $X^c_n$ starts at some initial state $X^c_0$ randomly chosen with distribution $\eta_0$. During the absorption stage, we set
$\widehat{X}_n^c=X_n^c$ with probability $G_n(X_n)$, otherwise we put the particle  in an auxiliary cemetery state $\widehat{X}_n^c=c$. When the particle $\widehat{X}_n^c$ is still alive (that is, if we have $\widehat{X}_n^c\in E_n$), it performs an elementary move $\widehat{X}_n^c\leadsto X_{n+1}^c $ according to the Markov transition $M_{n+1}$. Otherwise, the particle is absorbed and we set $X_{p}^c=\widehat{X}_{p}^c=c$, for any time $p>n$.

If we let $T$ be the first time $\widehat{X}_n^c=c$, then we have the Feynman-Kac representation formulae
$$
\QQ_n=\mbox{\rm Law}((X_0^c,\ldots,X_n^c)~|~T\geq n)
\quad\mbox{\rm
and}\quad
\Za_n=\mbox{\rm Proba}\left(T\geq n\right)
$$
For a more thorough discussion on the variety of application domains of Feynman-Kac models,
we refer the reader to chapter~\ref{FK-applications-sec}.

We also
denote by $\eta_n$ and $\gamma_n$, the $n$-th time marginal of $\QQ_n$ and $\Gamma_n$.
It is a simple exercise to check that
\begin{equation}\label{eq-flot-fk}
\gamma _{n}=\gamma _{n-1}Q_{n}\quad \mbox{\rm and}\quad \eta _{n+1}=\Phi
_{n+1}(\eta _{n}):=\Psi _{G_{n}}(\eta _{n})M_{n+1}
\end{equation}
with the positive integral operator
$$
Q_{n}(x,dy)=G_{n-1}(x)~ M_{n}(x,dy)
$$
and the Boltzmann-Gibbs transformation
\begin{equation}\label{def-BG-ref-intro}
\Psi _{G_{n}}(\eta _{n})(dx)=\frac{1}{\eta_n(G_n)}~G_n(x)~\eta_n(dx)
\end{equation}
In addition, the normalizing constants $\Za_n$ 
can be expressed in terms of the flow of marginal measures $\eta_p$, from the origin $p=0$ up to the current time
$n$, with the following multiplicative formulae:
\begin{equation}\label{multiplicative}
\Za_n:=\gamma_n(\un)=\EE\left(\prod_{0\leq
p<n}G_{p}(X_{p})\right)=\prod_{0\leq
p<n}\eta_p(G_{p})
\end{equation}
This multiplicative formula is easily checked
using the induction
 $$\gamma_{n+1}(1)=\gamma_{n}(G_{n})=\eta_n(G_n)~%
\gamma_n(1)$$

The abstract formulae discussed above are more general than it may appear. For instance,  they can be used
to analyze without further work path spaces models, including historical processes or transition space models, as well as finite excursion models. These functional models also encapsulated  quenched Feynman-Kac models, Brownian type bridges and 
linear Gaussian Markov chains conditioned on starting and end points.
 
For a more thorough discussion on these path space models, we refer the reader to section 2.4, section 2.6, chapters 11-12 in the monograph~\cite{Delm04}, as
well as to the section~\ref{FK-applications-sec}, in the former lecture notes.

When the Markov transitions $M_n$ are
absolutely continuous with respect to some measures $\lambda_n$ on $%
E_n $, and for any $(x,y)\in \left(E_{n-1}\times E_n\right)$ we have 
\begin{equation}
H_n(x,y):=\frac{dM_n(x,\mbox{\LARGE .})}{d\lambda_n}(y)>0\label{reg-H}
\end{equation}
we also have the following backward formula
\begin{equation}  \label{backward}
\mathbb{Q}_n(d(x_0,\ldots,x_n))=\eta_n(dx_n)~\prod_{q=1}^{n}\MM_{q,%
\eta_{q-1}}(x_q,dx_{q-1})
\end{equation}
with the the collection of Markov transitions
defined  by 
\begin{equation}  \label{backwardt}
\MM_{n+1,\eta_n}(x,dy)~\propto~G_{n}(y)~H_{n+1}(y,x)~\eta_n(dy)
\end{equation}
The proof of this formula is housed in section~\ref{historical-sec}.

Before launching into the description of the particle approximation of these 
models, we end this section with some connexions between discrete generation 
Feynman-Kac models and more conventional continuous time models arising in physics
and scientific computing.

The Feynman-Kac
models presented above  play a central role in the numerical analysis of certain partial
differential equations, offering a natural way to solve these functional
integral models by simulating random paths of stochastic processes. These
Feynman-Kac models were originally presented by Mark Kac in 1949~\cite{kac}
for continuous time processes.

These continuous time models are used in
molecular chemistry and computational physics to calculate the ground state
energy of some Hamiltonian operators associated with some potential function
$V$ describing the energy of a molecular configuration (see for instance~%
\cite{benj2,dm3,Delm04,benjamin,rousset}, and references therein).
To better connect these partial differential equation models with (\ref%
{defQn}), let us assume that $M_{n}(x_{n-1},dx_{n})$ is the Markov
probability transition $X_{n}=x_{n}\leadsto X_{n+1}=x_{n+1}$ coming from a
discretization in time $X_{n}=X_{t_{n}}^{\prime }$ of a continuous time $E$%
-valued Markov process $X_{t}^{\prime }$ on a given time mesh $%
(t_{n})_{n\geq 0}$ with a given time step $(t_{n}-t_{n-1})=\Delta t$. For
potential functions of the form $G_{n}=e^{-V\Delta t}$, the measures $%
\mathbb{Q}_{n}\simeq _{\Delta t\rightarrow 0}\mathbb{Q}_{t_{n}}$ represents
the time discretization of the following distribution:
\begin{equation*}
d\mathbb{Q}_{t}=\frac{1}{\mathcal{Z}_{t}}~\exp {\left(
-\int_{0}^{t}~V(X_{s}^{\prime })~ds\right) }~d\mathbb{P}_{t}^{X^{\prime }}
\end{equation*}%
where $\mathbb{P}_{t}^{X^{\prime }}$ stands for the distribution of the
random paths $(X_{s}^{\prime })_{0\leq s\leq t}$ with a given infinitesimal
generator $L$. The marginal distributions $\gamma _{t}$ at time $t$ of the
unnormalized measures $\mathcal{Z}_{t}~d\mathbb{Q}_{t}$ are the solution of
the so-called imaginary time Schroedinger equation, given in weak
formulation on sufficiently regular function $f$ by the following intregro--differential
equation
\begin{equation*}
\frac{d}{dt}~\gamma _{t}(f):=\gamma _{t}(L^{V}(f))\quad \mbox{\rm with}\quad
L^{V}=L-V
\end{equation*}%
The errors introduced by the discretization of the time are well understood
for regular models, we refer the interested reader to~\cite{djjp,dimasi,korez2,picard} in the context of nonlinear filtering.

\subsection{Interacting particle systems}\label{sec-ips-fk-intro}

The stochastic particle interpretation of the Feynman-Kac measures (\ref{defQn}) starts with a population of $N$ 
candidate possible solutions $(\xi^1_0,\ldots,\xi^N_0)$ randomly chosen w.r.t. some distribution $\eta_0$.

The coordinates
$\xi^i_0$ also called individuals or phenotypes, with $1\leq N$.  The random evolution of the particles is
decomposed into two main steps : the free exploration and the adaptive selection transition.

During the updating-selection stage,\index{Selection transition} multiple individuals in the current population  $(\xi^1_n,\ldots,\xi^N_n)$ at time $n\in\NN$ 
are stochastically selected based
on the fitness function $G_n$.
In practice, we choose a random proportion $B^i_n$ of an
existing solution $\xi^i_n$ in the current population with a mean value $\propto G_n(\xi^i_n)$ to breed a brand new generation of "improved"
solutions 
$(\widehat{\xi}^1_n,\ldots,\widehat{\xi}^N_n)$. For instance,  for every
index $i$, with a probability $\epsilon _{n}G_{n}(\xi _{n}^{i})$, we set $%
\widehat{\xi }_{n}^{i}=\xi _{n}^{i}$, otherwise we replace $\xi _{n}^{i}$
with a new individual $\widehat{\xi }_{n}^{i}=\xi _{n}^{j}$ randomly chosen
from the whole population with a probability proportional to $G_{n}(\xi
_{n}^{j})$. The parameter $\epsilon _{n}\geq 0$ is a tuning parameter that
must satisfy the constraint $\epsilon _{n}G_{n}(\xi _{n}^{i})\leq 1$, for
every $1\leq i\leq N$.
During the prediction-mutation stage, every selected individual $\widehat{\xi}^i_n$ moves to a new solution
$\xi^i_{n+1}=x$ randomly chosen in $E_{n+1}$, with a distribution $M_{n+1}(\widehat{\xi}^i_n,dx)$.

If we interpret the updating-selection transition as a birth and death process, then
arises the important notion of the ancestral line of a current individual.
More precisely, when a particle $\widehat{\xi }_{n-1}^{i}\longrightarrow \xi
_{n}^{i}$ evolves to a new location $\xi _{n}^{i}$, we can interpret $%
\widehat{\xi }_{n-1}^{i}$ as the parent of $\xi _{n}^{i}$. Looking backwards
in time and recalling that the particle $\widehat{\xi }_{n-1}^{i}$ has
selected a site $\xi _{n-1}^{j}$ in the configuration at time $(n-1)$, we
can interpret this site $\xi _{n-1}^{j}$ as the parent of $\widehat{\xi }%
_{n-1}^{i}$ and therefore as the ancestor denoted $\xi _{n-1,n}^{i}$ at
level $(n-1)$ of $\xi _{n}^{i}$. Running backwards in time we may trace the
whole ancestral line 
\begin{equation}
\xi _{0,n}^{i}\longleftarrow \xi _{1,n}^{i}\longleftarrow \ldots
\longleftarrow \xi _{n-1,n}^{i}\longleftarrow \xi _{n,n}^{i}=\xi _{n}^{i}
\label{ancestraline}
\end{equation}%

 Most of the terminology we have used is drawn from filtering and genetic evolution theories.
  
 In filtering, the  former particle model is dictated by the two steps 
 prediction-updating learning equations of the conditional distributions of a signal process, 
 given some  noisy and partial observations. In this setting, the potential functions represent
 the likelihood function of the current observation, while the free exploration transitions are related to the
 Markov transitions of the signal process. 
 
  In biology, the mutation-selection particle model presented above  is used to mimic genetic evolutions of biological organisms and more generally natural evolution processes. For instance, in gene analysis, each population of individuals represents a chromosome and each individual particle is called a gene. In this setting the fitness potential function is usually time-homogeneous and it represents the quality and the adaptation potential value  of the set of genes in a chromosome~\cite{holland}. These particle algorithms are also used in population analysis to model changes in the structure of population in time and in space.

The different types of particle approximation measures associated with the
genetic type particle model described above are summarized
 in the following synthetic picture
corresponding to the case $N=3$:

  \begin{center}     
\hskip.3cm
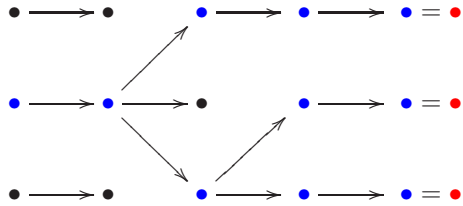
\xymatrix{
&{\bullet} \ar[r]&{\bullet}&\textcolor{blue}{\bullet}\ar[r]&~\textcolor{blue}{\bullet}\ar[r]&~\textcolor{blue}{\bullet}=\textcolor{red}{ \bullet}\\
&\textcolor{blue}{\bullet} \ar[r]&\textcolor{blue}{\bullet}\ar[dr]\ar[ur]\ar[r]&{\bullet}&~\textcolor{blue}{\bullet}\ar[r]&~\textcolor{blue}{\bullet}=\textcolor{red}{ \bullet}\\
&{\bullet}\ar[r]&{\bullet}&\textcolor{blue}{\bullet}\ar[ur]\ar[r]&~\textcolor{blue}{\bullet}\ar[r]&~\textcolor{blue}{\bullet}=\textcolor{red}{ \bullet}
}  
  \end{center}

In the next 4 sections we give an overview of the 4 particle approximation measures can be be extracted from
the interacting population evolution model described above. We also provide some basic formulation of the concentration inequalities
  that will be treated in greater detail later. As a service to the reader we also provide precise pointers to their location within the 
  following chapters. We already mention that the 
  proofs of these results are quite subtle.

  The precise form of the constants in
  these exponential inequalities  depends on the contraction properties of Feynman-Kac flows. Our stochastic analysis requires to combine the stability properties of the nonlinear semigroup
  of the Feynman-Kac distribution flow $\eta_n$, with deep convergence results of empirical processes theory associated with interacting random samples. 
 
\subsubsection{Current population models}

The occupation measures of the current population, represented by the red dots
in the above figure
$$\eta _{n}^{N}:=\frac{1}{N}\sum_{i=1}^{N}\delta_{\xi _{n}^{i}}$$ converge to the $n$-th time marginals $\eta_n$ of the Feynman-Kac
measures $\QQ_n$. We shall measure the performance of these particle estimates through several concentration inequalities, with a
special emphasis on uniform inequalities w.r.t. the time parameter. Our results will basically be stated as follows.

1) For any time horizon $n\geq 0$, any bounded function $f$,  any $N\geq 1$,  and for any $x\geq 0$, the probability of the event
$$
\left[\eta^N_{n}-\eta_n\right](f)\leq \frac{c_1}{N}~\left(1+x+\sqrt{x}\right)+\frac{c_2}{\sqrt{N}}~\sqrt{x}
$$
is greater than $1-e^{-x}$. In the above display,   $c_1$ stands for a finite constant related to {\em the bias} of the particle model,
while $c_2$ is related to {\em the variance}
of the scheme. The values of $c_1$ and $c_2$ don't depend on the time parameter. 

We already mention one important consequence of these uniform concentration inequalities for time homogeneous Feynman-Kac models. Under some regularity conditions, the flow of measures $\eta_n$ tends to some fixed point distribution $\eta_{\infty}$, in the sense that
\begin{equation}\label{stab-prop-intro-ref}
\|\eta_n-\eta_{\infty}\|_{\rm\tiny tv}\leq c_3~e^{-\delta n}
\end{equation}
for some finite positive constants $c_3$ and $\delta$. 
The connexions between these limiting measures and
the top of the spectrum of Schr\"odinger  operators is discussed in section~\ref{doob-h-sec}. We also refer the reader to section
~\ref{yaglom-sec} for a discussion on these quasi-invariant measures and Yaglom limits. 
Quantitative contraction theorems for Feynman-Kac semigroups 
are developed in the section~\ref{quantitative-contract-sec}. As a  direct consequence of the above inequalities,  we find that
for any $x\geq 0$, the probability of the following events is is greater than $1-e^{-x}$
$$
\left[\eta^N_{n}-\eta_{\infty}\right](f)\leq \frac{c_1}{N}~\left(1+x+\sqrt{x}\right)+\frac{c_2}{\sqrt{N}}~\sqrt{x}+c_3~e^{-\delta n}
$$

2) For any $x=(x_i)_{1\leq i\leq d}\in E_n=\RR^d$, we set $(-\infty,x]=\prod_{i=1}^d(-\infty,x_i]$ and
we consider the repartition functions
$$
F_n(x)=\eta_n\left(1_{(-\infty,x]}\right)\quad\mbox{\rm and}\quad
F^N_n(x)=\eta^N_n\left(1_{(-\infty,x]}\right)
$$
 The probability of the following event
$$
\sqrt{N}~\left\|F^N_n-F_n\right\|\leq c~\sqrt{d~(x+1)}
$$
is greater than $1-e^{-x}$, for any $x\geq 0$, for some universal constant $c<\infty$ that doesn't depend on the dimension, nor on the time parameter. Furthermore, under the stability properties (\ref{stab-prop-intro-ref}), 
if we set
$$
F_{\infty}(x)=\eta_{\infty}\left(1_{(-\infty,x]}\right)
$$
 then, the probability of the following event
$$
\left\|F^N_n-F_{\infty}\right\|\leq \frac{c}{\sqrt{N}}~\sqrt{d~(x+1)}+c_3~e^{-\delta n}
$$
is greater than $1-e^{-x}$, for any $x\geq 0$, for some universal constant $c<\infty$ that doesn't depend on the dimension.

For more precise statements, we refer the reader to corollary~\ref{cor-eta-1-ref},  and respectively to corollary~\ref{cor-cells-cv}.

The concentration properties of the particle measures $\eta^N_n$ around their limiting values
are developed in chapter~\ref{fk-particle-chap}. In section~\ref{stoch-perturbation-sec}, we design a stochastic perturbation analysis
that allows to enter the stability properties of the limiting Feynman-Kac semigroup. Finite marginal models are discussed in section~\ref{finite-marg-ips-concentration}. Section~\ref{interacting-ips-sec} is concerned with the concentration inequalities of interacting particle processes w.r.t. some collection of functions.

\subsubsection{Particle free energy models}\label{sec-intro-pfree}

Mimicking the multiplicative formula (\ref{multiplicative}), we set
and
\begin{equation}\label{free-energy}
\Za^N_n=\prod_{0\leq p<n}\eta^N_p(G_p) \quad\mbox{\rm and}\quad
\gamma^N_n(dx)=\Za^N_n\times \eta^N_n(dx)
\end{equation}

We already mention that these rather complex particle models provide an unbiased estimate of the
unnormalized measures. That is, we have that
\begin{equation}\label{unbias-prop}
\EE\left(\eta^N_n(f_n)~\prod_{0\leq p<n}\eta^N_p(G_p)\right)=\EE\left(f_n(X_n)~\prod_{0\leq p<n}G_p(X_p)\right)
\end{equation}

The concentration properties of the {\em unbiased}  particle free energies $\Za^N_n$ around their limiting values $\Za_n$
are developed in section~\ref{section-pfe}.
Our results  will basically be stated as follows. 
 
For any $N\geq 1$,  and any $\epsilon\in \{+1,-1\}$, the probability of each of the following events
$$
\frac{\epsilon}{n}\log{\frac{\Za^N_n}{\Za_n}}
\leq  \frac{c_1}{N}~\left(1+x+\sqrt{x}\right)+\frac{c_2}{\sqrt{N}}~\sqrt{x}
$$
is greater than $1-e^{-x}$. In the above display,   $c_1$ stands for a finite constant related to {\em the bias} of the particle model,
while $c_2$ is related to {\em the variance}
of the scheme. Here again, the values of $c_1$ and $c_2$ don't depend on the time parameter. A more precise statement is provided in corollary~\ref{cor-gammaN-ref-intro}.

\subsubsection{Genealogical tree model}
The occupation measure  of the $N$-genealogical tree model represented by the lines linking the blue dots converges as $N\rightarrow \infty $ to the
 distribution $\mathbb{Q}_{n}$
 \begin{equation}\label{genealogy-sec}
\lim_{N\rightarrow \infty }\frac{1}{N}\sum_{i=1}^{N}\delta_{(\xi _{0,n}^{i},\xi
_{1,n}^{i},\ldots ,\xi _{n,n}^{i})}=\QQ_n
\end{equation}
Our concentration inequalities will basically be stated as follows. A more precise statement is provided in corollary~\ref{ref-cor-genealogical-tree-intro}.
 
For any $n\geq 0$, any bounded function ${\bf f_n}$ on the path space ${\bf E_n}$,  s.t. $\|{\bf f_n}\|\leq 1$,
and any $N\geq 1$, the probability of each of the following events
$$
\begin{array}{l}
\left[
\frac{1}{N}\sum_{i=1}^{N}{\bf f_n}(\xi _{0,n}^{i},\xi
_{1,n}^{i},\ldots ,\xi _{n,n}^{i})-\QQ_n({\bf f_n})\right]
\\
\\
\leq  \displaystyle c_1~\frac{n+1}{N}~\left(1+x+\sqrt{x}\right)+c_2~\sqrt{\frac{(n+1)}{N}}~\sqrt{x}
\end{array}
$$
is greater than $1-e^{-x}$. In the above display,   $c_1$ stands for a finite constant related to {\em the bias} of the particle model,
while $c_2$ is related to {\em the variance}
of the scheme. Here again, the values of $c_1$ and $c_2$ don't depend on the time parameter. 

The concentration properties of genealogical tree occupation measures can be derived more or less directly from
the ones of the current population models. This rather surprising assertion comes from the fact that the $n$-th time
 marginal $\eta_n$ of a Feynman-Kac measure associated with a reference historical Markov process has the same form as in 
 the measure (\ref{defQn}). This equivalence principle between  $\QQ_n$ and the marginal measures
 are developed  in section~\ref{historical-sec}, dedicated to in historical Feynman-Kac models. 
 
Using these properties, we prove concentration properties for interacting empirical processes associated with genealogical tree models. Our concentration inequalities will basically be stated as follows. A more precise statement is provided in section~\ref{interacting-ips-sec}. We let $\Fa_n$ be the set of product functions of  indicator of cells  in the path space
${\bf E_n}=\left(\RR^{d_0}\times\ldots,\times\RR^{d_n}\right)$,  for some $d_p\geq 1$, $p\geq 0$. We also denote by $\eta^N_n$
the occupation measure of the genealogical tree model.In this notation, the probability of the following event
$$
\sup_{{\bf f_n}\in\Fa_n}{\left|\eta^N_n({\bf f_n})-\QQ_n({\bf f_n})\right|}\leq c~(n+1)~\sqrt{\frac{\sum_{0\leq p\leq n}d_p}{N}~(x+1)}
$$
is greater than $1-e^{-x}$, for any $x\geq 0$, for some universal constant $c<\infty$ that doesn't depend on the dimension. 

\subsubsection{Complete genealogical tree models}

Mimicking the backward model (\ref{backward}) and the above formulae, we set 
\begin{equation}\label{backNnU}
\Gamma^N_n=\Za^N_n\times \mathbb{Q}^N_n
\end{equation}
with
$$
\mathbb{Q}^N_n(d(x_0,\ldots,x_n))=\eta_n^N(dx_n)~\prod_{q=1}^{n}\MM_{q,%
\eta^N_{q-1}}(x_q,dx_{q-1})
$$

Notice that the computation of sums w.r.t. these particle measures 
are reduced to summations over the particles locations $\xi^i_n$. It is therefore natural to identify
a population of individual $(\xi^1_n,\ldots,\xi^N_n)$ at time $n$ to the ordered set of indexes $\{1,\ldots,N\}$. In this case, the occupation measures and the functions are identified with the following
line and column vectors\index{Occupation measures}
$$
\eta_n^N:=\left[\frac{1}{N},\ldots,\frac{1}{N}\right]\quad\mbox{\rm and}\quad\mbox{f}_n:=
\left(
\begin{array}
[c]{c}%
f_n(\xi^1_n)\\
\vdots\\
f_n(\xi^N_n)
\end{array}
\right) 
$$
and the matrices  $\MM_{n,\eta^N_{n-1}}$ by the $(N\times N)$ matrices
\begin{equation}\label{random-matrix-ref}
\qquad \MM_{n,\eta^N_{n-1}}:=
\left(
\begin{array}
[c]{ccc}%
\MM_{n,\eta^N_{n-1}}(\xi^1_n,\xi^1_{n-1}) & \cdots & \MM_{n,\eta^N_{n-1}}(\xi^1_n,\xi^N_{n-1}) \\
\vdots & \vdots & \vdots\\
\MM_{n,\eta^N_{n-1}}(\xi^N_n,\xi^1_{n-1}) & \cdots & \MM_{n,\eta^N_{n-1}}(\xi^N_n,\xi^N_{n-1}) 
\end{array}
\right) \end{equation}
with the $(i,j)$-entries
$$
\MM_{n,\eta^N_{n-1}}(\xi^i_{n},\xi^j_{n-1})=
\frac{G_{n-1}(\xi^j_{n-1})
H_n(\xi^j_{n-1},\xi^i_n)}{\sum_{k=1}^N G_{n-1}(\xi^k_{n-1})
H_n(\xi^k_{n-1},\xi^i_n)}
$$
For instance, the $\QQ_n$-integration of normalized additive linear functionals 
of the form
\begin{equation}\label{ref-additive-norm-intro}
{\bf f}_n(x_0,\ldots,x_n)=\frac{1}{n+1}\sum_{0\leq p\leq n}f_p(x_p)
\end{equation}
is given the particle matrix approximation model
$$
\QQ^N_n({\bf f}_n)=\frac{1}{n+1}\sum_{0\leq p\leq n} 
 \eta_n^N \MM_{n,\eta^N_{n-1}}\MM_{2,\eta^N_{1}}\MM_{p+1,\eta^N_{p}}(f_p)
$$
These type of additive functionals arise in the calculation of the sensitivity measures
discussed in section~\ref{sensitivity-meas}.

The concentration properties of the particle measures $\QQ ^N_n$ around the Feynman-Kac measures  $\QQ_n$
are developed in section~\ref{section-backward-ref}. Special emphasis is given to the additive functional models (\ref{ref-additive-norm-intro}).
In section~\ref{stoch-perturbation-sec-back}, we extend the stochastic perturbation methodology developed in  section~\ref{stoch-perturbation-sec} for time marginal model to the particle backward Markov chain associated with the random stochastic matrices (\ref{random-matrix-ref}).
This technique allows to enter not only the stability properties of the limiting Feynman-Kac semigroup, but also
the ones of the particle backward Markov chain model.

Our concentration inequalities will basically be stated as follows. A more precise statement is provided in corollary~\ref{ref-cor-final-back-intro} and in corollary~\ref{cor-fin-lecture}.

For any $n\geq 0$, any normalized additive functional of the form (\ref{ref-additive-norm-intro}), with
$\max_{0\leq p\leq n}{\|f_p\|}\leq 1$,
and any $N\geq 1$, the probability of each of the following events
$$
\left[\QQ ^N_n-\QQ_n\right](\overline{\bf f}_n)
\leq \displaystyle c_1~\frac{1}{N}~(1+(x+\sqrt{x}))+
c_2~\sqrt{\frac{x}{N(n+1)}}
$$
is greater than $1-e^{-x}$. In the above display,   $c_1$ stands for a finite constant related to {\em the bias} of the particle model,
while $c_2$ is related to {\em the variance}
of the scheme. Here again, the values of $c_1$ and $c_2$ don't depend on the time parameter.

For any $a=(a_i)_{1\leq i\leq d}\in E_n=\RR^d$, we denote by $C_a$ the cell 
 $$
C_a:=(-\infty,a]=\prod_{i=1}^d(-\infty,a_i]
$$
and ${\bf f}_{a,n}$ the additive functional
$$
{\bf f}_{a,n}(x_0,\ldots,x_n)=\frac{1}{n+1}\sum_{0\leq p\leq n}1_{\left(-\infty,a\right]}(x_p)
$$
The probability of the following event
$$
\sup_{a\in\RR^d}{\left|\QQ^N_n({\bf f_{a,n}})-\QQ_n({\bf f_{a,n}})\right|}\leq c~\sqrt{\frac{d}{{N}}(x+1)}
$$
is greater than $1-e^{-x}$, for any $x\geq 0$, for some constant $c<\infty$ that doesn't depend on the dimension, nor on the time horizon.

\begin{remark}
One way to turn all of these inequalities in term of Bernstein style concentration inequalities is as follows.
For any exponential inequality of the form
$$
\forall x\geq 0\qquad \PP\left(X\leq ax+\sqrt{2bx}+c\right)\leq 1-e^{-x}
$$
for some non negative constants $(a,b,c)$, we also have
$$
\forall y\geq 0\qquad \PP\left(X\leq y+c\right)\leq 1-\exp{\left(-\frac{y^2}{2\left(b+ay\right)}\right)}
$$
A proof of this result is provided in lemma~\ref{bernstein-lem}.
\end{remark}

\subsection{Basic notation}

This section provides some background from stochastic analysis and integral operator theory 
we require for our proofs. Most of the results with detailed proofs can be located in the book~\cite{Delm04},
on Feynman-Kac formulae and interacting particle methods. 
Our proofs also contain cross-references to this well rather known material, so the reader may wish to skip this section and enter
 directly to the chapter~\ref{FK-applications-sec} dedicated to some application domains of Feynman-Kac models.

\subsubsection{Integral operators}
We denote respectively by $\mathcal{M}(E)$, $\mathcal{M}_{0}(E)$, $\Pa(E)$,
and $\mathcal{B}(E)$, the set of all finite signed measures
on some measurable space $(E,\mathcal{E})$, the convex subset of measures with
null mass,  the set of all probability measures, and the Banach space of all
bounded and measurable functions $f$  equipped  with
the uniform norm $\Vert f\Vert$. We also denote by
 $\mbox{Osc}_{1}(E)$, and by $\mathcal{B}_1(E)$ the
set of $\mathcal{E}$-measurable functions $f$ with oscillations $\mbox{osc}(f)\leq1$, and respectively with $\|f\|\leq 1$. We let $$\mu(f)=\int~\mu(dx)~f(x)$$ be the Lebesgue integral of a function
$f\in\mathcal{B}(E)$, with respect to a measure $\mu\in\mathcal{M}(E)$.

We recall that the total variation distance on ${\cal M}(E)$
is defined for any $\mu\in{\cal M}(E)$ by
$$
\|\mu\|_{\tiny\rm tv}=\frac{1}{2}~\sup_{(A,B)\in\Ea^2}(\mu(A)-\mu(B))
$$
We recall that a bounded integral operator $M$ from a measurable space
$(E,\mathcal{E})$ into an auxiliary measurable space $(F,\mathcal{F})$ is an
operator $f\mapsto M(f)$ from $\mathcal{B}(F)$ into $\mathcal{B}(E)$ such that
the functions $$
M(f)(x):=\int_{F}M(x,dy)f(y)$$
are $\mathcal{E}$-measurable and bounded, for any $f\in\mathcal{B}(F)$. A Markov kernel 
is a positive and bounded integral operator $M$ with $M(1)=1$. Given a pair of bounded integral operators $(M_1,M_2)$, we let $(M_1M_2)$ the composition operator defined by
$(M_1M_2)(f)=M_1(M_2(f))$. For time homogeneous state spaces, we denote
by $M^m=M^{m-1}M=MM^{m-1}$ the $m$-th composition of a given bounded integral operator $M$, with
$m\geq 1$.
A
bounded integral operator $M$ from a measurable space $(E,\mathcal{E})$ into
an auxiliary measurable space $(F,\mathcal{F})$ also generates a dual operator
$$\mu(dx)\mapsto(\mu M)(dx)=\int\mu(dy)M(y,dx)$$ from $\mathcal{M}(E)$ into $\mathcal{M}(F)$ defined by $(\mu
M)(f):=\mu(M(f))$. We also used the notation
$$
K\left(\left[f-K(f)\right]^2\right)(x):=K\left(\left[f-K(f)(x)\right]^2\right)(x)
$$
for some bounded integral operator $K$ and some bounded function $f$. 

We prefer to avoid unnecessary abstraction and technical assumptions, so we frame the standing assumption that 
all the test functions are in the unit sphere, and the integral operators, and all the random variables
are sufficiently regular that we are justified in computing  integral transport equations, regular versions of conditional expectations, and so forth.

\subsubsection{Contraction coefficients}
When the bounded integral operator  $M$ has a constant mass, that is, when $M(1)\left(  x\right)  =M(1)\left(
y\right)  $ for any $(x,y)\in E^{2}$, the operator $\mu\mapsto\mu M$ maps
$\mathcal{M}_{0}(E)$ into $\mathcal{M}_{0}(F)$. In this situation, we let
$\beta(M)$ be the Dobrushin coefficient of a bounded integral operator $M$
defined by the formula
$$\beta(M):=\sup{\ \{\mbox{\rm osc}(M(f))\;;\;\;f\in\mbox{\rm Osc}(F)\}} $$
Notice that $\beta(M)$ is the operator norm of 
 $M$ on ${\cal M}_0(E)$, and we  have
the equivalent formulations  
\begin{eqnarray*} 
\beta(M)&=& \sup{\{\|M(x,\point)-M(y,\point)\|_{\tiny\rm tv}\;;\; (x,y)\in E^2\}} \\
&=&\sup_{\mu\in {\cal M}_0(E)}{\|\mu M\|_{\tiny\rm tv}}/{\|\mu\|_{\tiny\rm 
    tv}}\label{norm}
\end{eqnarray*}  
A detailed proof of these well known formulae can be found in~\cite{Delm04}.

Given a positive and bounded potential function $G$ on $E$, we also denote by $\Psi_G$
the Boltzmann-Gibbs mapping  from $\Pa(E)$ into itself defined for any $\mu\in\Pa(E)$ by
$$
\Psi_G(\mu)(dx)=\frac{1}{\mu(G)}~G(x)~\mu(dx)
$$
For $]0,1]$-valued potential functions, we also mention that $\Psi_G(\mu)$ can be expressed as a non linear Markov 
transport equation
\begin{equation}\label{nonlin-transport}
\Psi_G(\mu)=\mu S_{\mu,G}
\end{equation}
with the Markov transitions
$$
S_{\mu,G}(x,dy)=G(x)~\delta_x(dy)+\left(1-G(x)\right)~\Psi_G(\mu)(dy)
$$
We notice that
$$
\Psi_G(\mu)-\Psi_G(\nu)=(\mu-\nu)S_{\mu}+\nu(S_{\mu}-S_{\nu})
$$
and
$$
\nu(S_{\mu}-S_{\nu})=(1-\nu(G))\left[\Psi_G(\mu)-\Psi_G(\nu)\right]
$$
from which we find the formula
$$
\Psi_G(\mu)-\Psi_G(\nu)=\frac{1}{\nu(G)}~(\mu-\nu)S_{\mu}
$$
In addition, using the fact that
$$
\forall (x,A)\in (E,\Ea)\quad S_{\mu}(x,A)\geq (1-\|G\|)~\Psi_G(\mu)(A)
$$
we prove that 
$
\beta(S_{\mu})\leq \|G\|$ and $$
\|\Psi_G(\mu)-\Psi_G(\nu)\|_{\tiny\rm tv}\leq \frac{\|G\|}{\mu(G)\vee\nu(G)}~\|\mu-\nu\|_{\tiny\rm tv}
$$

If we set $\Phi(\mu)=\Psi_G(\mu)M$, for some Markov transition  $M$, then  we have the decomposition 
\begin{equation}\label{tv-Phi-equal}
\Phi(\mu)-\Phi(\nu)=\frac{1}{\nu(G)}~(\mu-\nu)S_{\mu}M
\end{equation}
for any couple of measures $\nu,\mu$ on $E$. From the previous discussion, 
we also find the following Lipschitz estimates
\begin{equation}\label{tv-Phi-contract}
\|\Phi(\mu)-\Phi(\nu)\|_{\tiny\rm tv}\leq \frac{\|G\|}{\mu(G)\vee\nu(G)}~\beta(M)~\|\mu-\nu\|_{\tiny\rm tv}
\end{equation}
We end this section with an interesting contraction property of a Markov transition
\begin{equation}\label{def-Mg-ref}
M_G(x,dy)=\frac{M(x,dy)G(y)}{M(G)(x)}=\Psi_G(\delta_xM)(dy)
\end{equation}
associated with a
$]0,1]$-valued potential function $G$, with 
\begin{equation}\label{def-g-ref}
g=\sup_{x,y}G(x)/G(y)<\infty
\end{equation}
It is easily checked that
\begin{eqnarray*}
\left|M_G(f)(x)-M_G(f)(y)\right|&=&\left|\Psi_G(\delta_xM)(f)-\Psi_G(\delta_yM)(f)\right|\\
&\leq& g~\|\delta_xM-\delta_yM\|_{\tiny\rm tv}
\end{eqnarray*}
from which we conclude that
\begin{equation}\label{trickch4}
\beta\left(M_G\right)\leq g~\beta\left(M\right)
\end{equation}

\subsubsection{Orlicz norms and Gaussian moments}

We let $\pi_{\psi}[Y]$ be the Orlicz norm of an $\RR$-valued random variable $Y$
associated with the 
 convex function $\psi(u)=e^{u^2}-1$, and defined by
$$
\pi_{\psi}(Y)=\inf{\{a\in (0,\infty)\;:\;\EE(\psi(|Y|/a))\leq 1\}}
$$ 
with the convention $\inf_{\emptyset}=\infty$. Notice that
$$
\pi_{\psi}(Y)\leq c\Longleftrightarrow \EE\left(\psi(Y/c)\right)\leq 1
$$

For instance, the Orlicz norm  of a Gaussian and centred random variable $U$, s.t. $E(U^2)=1$, is given by $\pi_{\psi}(U)=\sqrt{8/3}$. We also recall that 
\begin{eqnarray}
\EE\left(U^{2m}\right)&=&b(2m)^{2m}:=(2m)_{m}~2^{-m}\nonumber\\
\EE\left(|U|^{2m+1}\right)&\leq& b(2m+1)^{2m+1} :=\frac{(2m+1)_{(m+1)}}{\sqrt{m+1/2}}~2^{-(m+1/2)}\nonumber\\&&\label{collec}
\end{eqnarray}
with $(q+p)_{p}:=(q+p)!/q!$.  The second assertion comes from the fact that
\begin{eqnarray*}
\EE\left(U^{2m+1}\right)^2&\leq& \EE\left(U^{2m}\right)~\EE\left(U^{2(m+1)}\right)
\end{eqnarray*}
and therefore
\begin{eqnarray*}
b(2m+1)^{2(2m+1)}&=& \EE\left(U^{2m}\right)~\EE\left(U^{2(m+1)}\right)\\
&=&2^{-(2m+1)}~(2m)_{m}~(2(m+1))_{(m+1)}
\end{eqnarray*}
This formula is a direct consequence of the following decompositions
\[
(2(m+1))_{(m+1)}=\frac{(2(m+1))!}{(m+1)!}=2~\frac{(2m+1)!}{m!}=2~(2m+1)_{(m+1)}%
\]
and
\[
(2m)_{m}=\frac{1}{2m+1}~\frac{(2m+1)!}{m!}=\frac{1}{2m+1}~(2m+1)_{(m+1)}%
\]
We also mention that
\begin{equation}\label{b2mbm}
b(m)\leq b(2m)
\end{equation}
Indeed, for even numbers $m=2p$ we have
$$
b(m)^{2m}=b(2p)^{4p}=\EE(U^{2p})^2\leq \EE(U^{4p})=b(4p)^{4p}=b(2m)^{2m}
$$
and for odd numbers $m=(2p+1)$, we have
\begin{eqnarray*}
b(m)^{2m}&=&b(2p+1)^{2(2p+1)}=\EE\left(U^{2p}\right)~\EE\left(U^{2(p+1)}\right)\\
&\leq&\EE\left(\left(U^{2p}\right)^{\frac{(2p+1)}{p}}\right)^{\frac{p}{2p+1}}~
\EE\left(\left(U^{2(p+1)}\right)^{\frac{(2p+1)}{p+1}}\right)^{\frac{p+1}{2p+1}}\\
&=&\EE\left(U^{2(2p+1)}\right)=b(2(2p+1))^{2(2p+1)}=b(2m)^{2m}
\end{eqnarray*}

\section{Some application domains}\label{FK-applications-sec}
\subsection{Introduction}
Feynman-Kac particle methods are also termed quantum Monte Carlo methods in computational physics, genetic algorithms in computer sciences, and particle filters and-or
sequential Monte Carlo methods in information theory, as well as in bayesian
statistics. 

The mathematical foundations of these advanced interacting Monte Carlo methodologies are now fifteen years old~\cite{dm96}.
Since this period, so many descriptions and variants of these models have been published in applied probability,
signal processing and Bayesian statistic literature. For a detailed discussion on their application domains with a precise 
bibliography of who first did that when, 
we refer the reader to any of the following references~\cite{dmsp,Delm04}, and ~\cite{dd-2012,doucet}.

In the present section, we merely content ourselves in illustrating the rather abstract models (\ref{defQn})
with the Feynman-Kac representations of 20 more or less well known conditional distributions, including three more recent applications
related to island particle models, functional kinetic parameter derivatives, and gradient analysis of Markov semigroups.

The forthcoming series of examples, combined with their mean field  particle interpretation models described in section~\ref{sec-ips-fk-intro}, also illustrate the ability of the Feynman-Kac particle 
methodology to  solve complex conditional distribution flows
as well as their normalizing constants.

Of course, this selected list of applications does not attempt to be exhaustive. 
The topics selection is largely influenced by the personal taste of the authors. 
A complete description
on how particle methods are applied in each application model area would of course
require
separate volumes, with precise computer simulations and comparisons with
different types of particle models and other existing algorithms.  

We also limit ourselves to describing the key ideas in a simple way, often sacrificing generality. 
Some applications are nowadays 
routine, and in this case we provide precise 
pointers to existing more application-related 
articles in the literature. Reader who wishes to know more about some specific application of these particle algorithms is invited to consult the referenced papers.

One natural path of "easy reading" will probably be 
to choose a familiar or attractive application area
and to explore some selected parts of the lecture notes in terms
of this choice. Nevertheless, this advice must not be taken
too literally. To see the impact of 
particle
methods, it is essential to understand the full force
of Feynman-Kac modeling techniques on
various research domains. Upon doing so, the reader will
have a powerful weapon for the discovery of new particle
interpretation models. 
The principal challenge is to understand the theory
 well enough to reduce them to 
practice.

\subsection{Boltzmann-Gibbs measures}

 \subsubsection{Interacting Markov chain Monte Carlo methods}\label{i-mcmc-intro}
 
Suppose we are given a sequence of target probability measures on some measurable state space
$E$ of the following form
\begin{equation}\label{mcmcintrod}
\mu_n(dx)=\frac{1}{\Za_n}~\left\{\prod_{0\leq p<n} h_p(x)\right\}~\lambda(dx)
\end{equation}
with some sequence of bounded nonnegative potential functions $$h_n~:~x\in E\mapsto h_n(x)\in (0,\infty)$$
and some
reference probability measure $\lambda$ on $E$. In the above displayed formula, $\Za_n$ stands for a normalizing 
constant. We use the convention $\prod_{\emptyset}=1$ and $\mu_0=\lambda$.

We further assume that we have a dedicated Markov chain Monte Carlo 
transition
$\mu_{n}=\mu_{n}{M_n}$ with prescribed target measures $\mu_n$, at any time step. Using the fact that
$$
\mu_{n+1}=\Psi _{h_{n}}(\mu _{n})
$$
with the Boltzmann-Gibbs transformations defined in (\ref{def-BG-ref-intro}), we prove that
$$
\mu_{n+1}=\mu_{n+1}M_{n+1}=\Psi _{h_{n}}(\mu_{n})M_{n+1}
$$
from which we conclude that
$$
\mu_{n}(f)=\EE\left(f(X_n)~\prod_{0\leq p<n}h_p(X_p)\right)~/~\EE\left(\prod_{0\leq p<n}h_p(X_p)\right)
$$
with the reference Markov chain
$$
\PP\left(X_n\in dx~|~X_{n-1}\right)=M_n(X_{n-1},dx)
$$
In addition, we have
$$
\Za_{n+1}=\int \left[\prod_{0\leq p<n}h_p(x)\right] ~h_n(x)~\mu(dx)=  \Za_{n}~\mu_n(h_n)=\prod_{0\leq p\leq n}\mu_p(h_p)
$$

We illustrate these rather abstract models with two applications related respectively to probability restriction
models and stochastic optimization simulated annealing type models. For a more thorough discussion on these interacting MCMC models, and related sequential Monte Carlo methods, we refer the reader to~\cite{Delm04,delmoraldoucetjasra2006}.
\subsubsection{Probability restrictions}\label{proba-restrict-sec}

If we choose Markov chain Monte Carlo type local moves $$\mu_{n}=\mu_{n}M_n$$ with some prescribed 
target Boltzmann-Gibbs measures
$$
\mu_n(dx)\propto 1_{A_n}(x)~\lambda(dx)
$$
associated with a sequence of decreasing subsets  $A_n\downarrow$, and some reference measure $\lambda$, then we find that $\mu_n=\eta_n$ and $\Za_n=\lambda(A_n)$, as soon as the potential functions in (\ref{defQn}) and (\ref{mcmcintrod}) are chosen so that
$$
G_n=h_n=1_{A_{n+1}}
$$

This stochastic model arise in several application domains. In computer science literature, the corresponding particle approximation models
are sometimes called subset methods, sequential sampling plans, randomized algorithms, or level splitting algorithms. They were used to solve complex NP-hard
combinatorial counting problems~\cite{cgrv}, extreme quantile probabilities~\cite{cerou3,ghm2011}, and
uncertainty propagations in numerical codes~\cite{caron4}.

 \subsubsection{ Stochastic optimization}

If we choose Markov chain Monte Carlo type local moves $\mu_{n}=\mu_{n}{M_n}$ with some prescribed target Boltzmann-Gibbs measures
$$
\mu_n(dx)\propto e^{-\beta_n V(x)}~\lambda(dx)
$$
associated with a sequence of increasing inverse temperature parameters  $\beta_n\uparrow$, and some reference measure $\lambda$, then we find that $\mu_n=\eta_n$ and $\Za_n=\lambda(e^{-\beta_n V})$ as soon as the potential functions  in (\ref{defQn})  and (\ref{mcmcintrod}) are chosen  so that
$$G_n=h_n=e^{-(\beta_{n+1}-\beta_{n})V}$$
For instance, we can assume that the Markov transition $M_n=\Ma_{n,\beta_n}^{m_n}$ is the $m_n$-iterate of the following
Metropolis Hasting
transitions 
\begin{equation}\label{i-MHM}
\begin{array}{l}
\Ma_{n,\beta_{n}}(x,dy)
\\
\\
=K_n(x,dy)~\min{\left(  1,e^{-\beta_{n}(V(y)-V(x))}\right)  }\nonumber\\
\\\qquad+\left(  1-\int
_{z}K_n(x,dz)~\min{\left(  1,e^{-\beta_{n}(V(z)-V(x))}\right)  }\right)
~\delta_{x}(dy)
\end{array}
\end{equation}

We finish this section with assorted collection of enriching  comments on interacting Markov chain Monte Carlo algorithms associated with the
Feynman-Kac models described above. 

Conventional Markov chain Monte Carlo methods ({\em abbreviated MCMC methods}) with time varying target measures $\mu_n$
can be seen as a single particle
model with only mutation explorations according to the Markov transitions $M_n=K_n^{m_n}$, where $K^{m_n}_n$ stands for the iteration of an MCMC transition $K_n$ s.t. $\mu_n=\mu_nK_n$. In this situation, we choose a judicious increasing sequence
$m_n$ so that the non homogeneous Markov chain is sufficiently stable, even if the target measures become more and more complex to sample.  When the target measure is fixed, say of the form $\mu_T$ for some large $T$, the MCMC sampler again uses a single particle with behave as a Markov chain with time homogenous transitions  $M_T$. The obvious drawback with
these two conventional MCMC samplers is that the user does not know how many steps 
are really needed to be close to the equilibrium target measure. A wrong choice will return samples with a distribution far from the desired target measure. 

Interacting MCMC methods run a population of MCMC samplers that interact one each other through a recycling-updating 
mechanism so that the occupation measure of the current measure converge to the target measure, when we increase the population sizes. In contrast with conventional MCMC methods, there are no burn-in time questions, nor any 
quantitative  analysis to estimate the convergence to equilibrium of the MCMC chain..

\subsubsection{Island particle models}\label{islands-sec}
In this section, we provide a brief discussion on interacting colonies and island particle models
arising in mathematical biology and evolutionary computing literature~\cite{island1,island2}.
The evolution of these stochastic island models is again defined in terms of a free evolution and a selection transition.
During the free evolution each island evolves separately as
a single mean field particle model with mutation-selection mechanism between the individual in the island population.
The selection
pressure between islands is related to the average fitness of the individuals in  island population. A colony 
with poor fitness is killed,
and replaced by brand new generations of "improved" individuals coming from better fitted islands.

For any measurable function $f_n$ on $E_n$, we set 
$$
X^{(0)}_n=X_n\in E^{(0)}_n:=E_n\quad\mbox{\rm and}\quad f^{(0)}_n=f_n
$$
 and we denote by 
 $$
X^{(1)}_n=\left(X^{(1,i)}_n\right)_{1\leq i\leq N_1}\in E_n^{(1)}:=\left(E^{(0)}_n\right)^{N_1}
 $$
the $N_1$-particle model associated with the reference Markov chain $X^{(0)}_n$, and the potential function $G^{(0)}_n$.

To get one step further, we denote by $f^{(1)}_n$ the empirical mean valued function
on $E_n^{(1)}$ defined by
$$ f^{(1)}_n(X^{(1)}_n)=\frac{1}{N_1}\sum_{i=1}^{N_1}f^{(0)}_n(X^{(1,i)}_n)
$$
In this notation, the potential value of the random state $X^{(1)}_n$ is given by the formula
$$
G_n^{(1)}(X^{(1)}_n)
:=\frac{1}{N_1}\sum_{i=1}^NG^{(0)}_n(X^{(1,i)}_n)
$$

By construction, we have the almost sure property
$$
N_1=1\Longrightarrow X^{(0)}_n=X^{(1)}_n
\quad\mbox{\rm and}\quad
G^{(1)}_n(X^{(1)}_n)=G^{(0)}_n(X^{(0)}_n)
$$
More interestingly, 
by the unbiased properties (\ref{unbias-prop}) we have {\em for any population size $N_1$}
$$
\EE\left(f^{(1)}_n(X^{(1)}_n)~\prod_{0\leq p<n}G_p^{(1)}(X^{(1)}_p)\right)=\EE\left(f_n^{(0)}(X^{(0)}_n)~\prod_{0\leq p<n}G_p^{(0)}(X^{(0)}_p)\right)
$$
Iterating this construction, we let
$$
X^{(2)}_n=\left(X^{(2,i)}_n\right)_{1\leq i\leq N_2}\in E_n^{(2)}:=\left(E^{(1)}_n\right)^{N_2}
 $$
the $N_2$-particle model associated with the reference Markov chain $X^{(1)}_n$, and the potential function $G^{(1)}_n$.
For any function $f_n^{(1)}$ on $E^{(1)}_n$, we denote by $f^{(2)}_n$ the empirical mean valued function
on $E_n^{(2)}$ defined by
$$ f^{(2)}_n(X^{(2)}_n)=\frac{1}{N_2}\sum_{i=1}^{N_2}f^{(1)}_n(X^{(2,i)}_n)
$$
In this notation, the potential value of the random state $X^{(2)}_n$ is given by the formula
$$
G_n^{(2)}(X^{(2)}_n)
:=\frac{1}{N_2}\sum_{i=1}^NG^{(1)}_n(X^{(2,i)}_n)
$$
and for any {\em for any population size $N_2$}
$$
\EE\left(f^{(2)}_n(X^{(2)}_n)~\prod_{0\leq p<n}G_p^{(2)}(X^{(2)}_p)\right)=\EE\left(f_n^{(0)}(X^{(0)}_n)~\prod_{0\leq p<n}G_p^{(0)}(X^{(0)}_p)\right)
$$

\subsubsection{Particle Markov chain Monte Carlo methods}\label{p-mcmc}
In this section, we present an interacting particle version of the particle Markov chain Monte Carlo method developed
in the recent seminal article by C. Andrieu, A. Doucet, and R. Holenstein~\cite{andrieu-doucet}.

We consider a  collection of Markov transition and positive potential functions $(M_{\theta,n},G_{\theta,n})$ that depend on some
random variable $\Theta=\theta$, with distribution $\nu$ on some state space $S$. We let $\eta_{\theta,n}$ be the $n$-time marginal
of the Feynman-Kac measures defined as in (\ref{defQn}), by replacing $(M_{n},G_{n})$ by $(M_{\theta,n},G_{\theta,n})$. 
We also consider  the probability distribution $P(\theta,d\xi)$ of the $N$-particle model 
$$
\xi:=\left(\xi_{\theta,0},\xi_{\theta,1},\ldots,\xi_{\theta,T}\right) 
$$
on the interval $[0,T]$, with mutation transitions $M_{\theta,n}$, and potential selection functions $G_{\theta,n}$, with $n\leq T$.
We fix  a large time horizon $T$, and for any $0\leq n\leq T$, we set
\begin{equation}\label{def-mu-xi-theta}
\mu_n(d(\xi,\theta))=\frac{1}{\Za_n}~\left\{\prod_{0\leq p<n} h_p(\xi,\theta)\right\}~\lambda(d(\xi,\theta))
\end{equation}
with some sequence of bounded nonnegative potential functions $$h_n~:~(\xi,\theta)\in \left(\prod_{0\leq p\leq T}E_p^N\right)\times  S
\mapsto h_n(\xi,\theta)\in (0,\infty)
$$
the reference measure $\lambda$ given by
$$
\lambda(d(\xi,\theta))=\nu(d\theta)~P(\theta,d\xi)
$$
and some normalizing constants ${\Za_n}$. Firstly, we observe that these target measures have the same form as in (\ref{mcmcintrod}).
Thus, they can be sampled using the Interacting Markov chain Monte Carlo methodology presented in section~\ref{i-mcmc-intro}. 

Now, we examine the situation where $h_p$ is given by the empirical mean value of the potential
function $G_{\theta,p}$ w.r.t. the occupation measures $\eta^N_{\theta,p}$ of the $N$-particle
model $\xi_{\theta,p}=\left(\xi^i_{\theta,p}\right)_{1\leq i\leq N}$ associated with the realization $\Theta=\theta$; more formally, we have that 
$$
h_p(\xi,\theta)=\frac{1}{N}\sum_{1\leq i\leq N} G_{\theta,p}(\xi_{\theta,p}^i)=\eta^N_{\theta,p}\left(G_{\theta,p}\right)
$$
Using the unbiased property  of the particle free energy models presented in (\ref{unbias-prop}), we clearly have
\begin{eqnarray*}
\int P(\theta,d\xi)~\left\{\prod_{0\leq p<n} h_p(\xi,\theta)\right\}&=&\EE\left(\prod_{0\leq p<n} \eta^N_{\theta,p}\left(G_{\theta,p}\right)\right)
\\
&=&\prod_{0\leq p<n} \eta_{\theta,p}(G_{\theta,p})
\end{eqnarray*}
from which we conclude that the $\Theta$-marginal of $\mu_n$ is given by the following equation
$$
\left(\mu_n\circ\Theta^{-1}\right)(d\theta)=\frac{1}{\Za_n}~\left\{\prod_{0\leq p<n} \eta_{\theta,p}(G_{\theta,p})\right\}~\nu(d\theta)
$$

We end this section with some comments on these distributions. 

As the initiated reader may have certainly noticed, the marginal analysis derived above coincides with 
one developed in section~\ref{islands-sec} dedicated to island particle models.

We also mention that the measures $\mu_n$
introduced in (\ref{def-mu-xi-theta}) can be approximated using the interacting  Markov chain Monte Carlo methodology
presented in section~\ref{i-mcmc-intro}, or the particle MCMC methods introduced in the article \cite{andrieu-doucet}.

Last but not least, we observe that 
$$
\prod_{0\leq p<n} \eta_{\theta,p}(G_{\theta,p})=\EE\left(\prod_{0\leq p<n}G_{\theta,p}(X_{\theta,p})\right)=\Za_n(\theta)
$$
where $X_{\theta,n}$ stand for the Markov chain with transitions $M_{\theta,n}$, and initial distribution $\eta_{\theta,0}$.
In the r.h.s. of the above displayed formulae, $\Za_n(\theta)$ stands for the normalizing constant of the Feynman-Kac measures defined as in (\ref{defQn}), by replacing $(M_{n},G_{n})$ by $(M_{\theta,n},G_{\theta,n})$. This shows that
$$
\left(\mu_n\circ\Theta^{-1}\right)(d\theta)=\frac{1}{\Za_n}~\Za_n(\theta)~\nu(d\theta)
$$
The goal of some stochastic optimization problem is to extract the parameter $\theta$ that
minimizes some mean value functional of the form
$$
\theta\mapsto  \Za_n(\theta)=\EE\left(\prod_{0\leq p<n}G_{\theta,n}(X_{\theta,p})\right)
$$

For convex functionals, we can use gradient type techniques using the Backward Feynman-Kac derivative interpretation models developed in section~\ref{sensitivity-meas} (see also  the three joint articles of the first author with A. Doucet, and S. S. Singh~\cite{sss,sss2,sss3}). 

When $\nu$ is the uniform measure
over some compact set $S$,
an alternative approach is to estimate
the measures (\ref{def-mu-xi-theta}) by some empirical measure 
$$\frac{1}{N}\sum_{1\leq i\leq N}\delta_{\left(\xi^{(i)}_n,\theta^{(i)}_n\right)}\in \Pa\left(\left(\prod_{0\leq p\leq n}E_p^N\right)\times  S\right)
$$ and to select the sampled state $$
\begin{array}{l}
\left(\xi^{(i)}_n,\theta^{(i)}_n\right)\\
\\
:=
\left(\left(
\left(\xi^{(i,j)}_{0,n}\right)_{1\leq j\leq N},\left(\xi^{(i,j)}_{1,n}\right)_{1\leq j\leq N}, \ldots,\left(\xi^{(i,j)}_{n,n}\right)_{1\leq j\leq N}
\right),\theta^{(i)}_n\right)\\
\\
\in\left( \left(E_0^N\times E_1^N\times\ldots\times E^N_n\right)\times S\right)
\end{array}
$$ that maximizes
the empirical objective functional 
$$
i\in \{1,\ldots,N\}\mapsto 
\prod_{0\leq p<n} \frac{1}{N}\sum_{1\leq j\leq N}G_{\theta^{(i)}_n,p}(\xi^{(i,j)}_{p,n},\theta^{(i)}_n)
$$

\subsubsection{Markov bridges and chains with fixed terminal value}

In many applications, it is important to sample paths of Markov chains with prescribed fixed terminal conditions.

When the left end starting point is distributed w.r.t. to a given regular probability measure $\pi$, we can use the time reversal Feynman-Kac formula presented by the first author and A. Doucet in~\cite{dm-metropolis}.
More precisely, for time homogeneous models $(G_n,M_n)=(G,M)$ in transition spaces, if we consider the Metropolis-Hasting ratio  $$G(x_1,x_2)=\frac{\pi(dx_2)K(x_2,dx_1)}{\pi(dx_1)M(x_1,dx_2)}$$
then we find that
$$
\QQ_n=\mbox{\rm Law}_{\pi}^K((X_0,\ldots,X_n)~|~X_n=x_n)
$$
where $\mbox{\rm Law}_{\pi}^K$ stands for the distribution of the Markov chain starting with an initial condition $\pi$ and 
evolving according some Markov transition $K$. The proof of these formulas are rather technical, we refer the reader the article~\cite{dm-metropolis}n and to the monograph~\cite{Delm04}.

For initial and terminal fixed end-points, we need to consider the paths distribution of Markov bridges. 
As we mentioned in the introduction, on page~\pageref{multiplicative}, these Markov bridges are particular instances of the reference
Markov chains of the abstract  Feynman-Kac model (\ref{defQn}). Depending on the choice of the potential functions in (\ref{defQn}), these Markov bridge models can be associated with several application domains, including  filtering problems or rare event analysis of bridge processes.

We assume that the 
elementary Markov transitions $M_n$ of the chain $X_n$ satisfy the regularity condition (\ref{reg-H}) for some
density functions $H_n$ and some reference measure $\lambda_n$. In this situation, the semigroup Markov 
transitions $M_{p,n+1}=M_{p+1}M_{p+2}\ldots M_{n+1}$ are
absolutely continuous with respect to the measure $\lambda_{n+1}$, for any $0\leq p\leq n $, and we have
$$
M_{p,n+1}(x_p,dx_{n+1})=H_{p,n+1}(x_p,x_{n+1})~\lambda_{n+1}(dx_{n+1})
$$
with the density function
$$
H_{p,n+1}(x_p,x_{n+1})=M_{p,n}\left(H_{n+1}(\point,x_{n+1})\right)(x_p)
$$
Thanks to these regularity conditions, we readily check that the paths distribution of Markov bridge starting at $x_0$ and ending at 
$x_{n+1}$ at the final time horizon $(n+1)$ are given by
$$
\begin{array}{l}
\BB_{(0,x_0),(n+1,x_{n+1})}\left(d(x_1,\ldots,x_n)\right)\\
\\
:=\PP\left((X_1,\ldots,X_n)\in d(x_1,\ldots,x_n)~|~X_0=x_0,~X_{n+1}=x_{n+1}\right)\\
\\
=\prod_{1\leq p\leq n}~M_p(x_{p-1},dx_p)~\frac{dM_{p,n+1}(x_p,\point)}{
dM_{p-1,n+1}(x_{p-1},\point)}(x_{n+1})\\
\\
=\prod_{1\leq p\leq n}~\frac{
M_p(x_{p-1},dx_p)~H_{p,n+1}(x_p,x_{n+1})
}{
M_p\left(H_{p,n+1}(\point,x_{n+1})\right)(x_{p-1})
}
\end{array}
$$

Using some abusive Bayesian notation, we can rewrite these formula as follows
$$
\begin{array}{l}
p((x_1,\ldots,x_n)~|~(x_0,x_{n+1}))
\\
\\=\frac{p(x_{n+1}|x_n)}{p(x_{n+1}|x_{n-1})}~p(x_n|x_{n-1})\ldots \frac{p(x_{n+1}|x_p)}{p(x_{n+1}|x_{p-1})}~p(x_p|x_{p-1})\ldots\\
\\\hskip7cm\ldots\frac{p(x_{n+1}|x_1)}{p(x_{n+1}|x_{0})}~p(x_1|x_{0})
\end{array}
$$
with 
$$
\frac{dM_p(x_{p-1},\point)}{d\lambda_p}(x_p)=p(x_p|x_{p-1})\quad\mbox{\rm and}\quad
H_{p,n+1}(x_p,x_{n+1})=p(x_{n+1}|x_p)
$$
For linear-Gaussian models, the Markov bridge transitions $$
\frac{
M_p(x_{p-1},dx_p)~H_{p,n+1}(x_p,x_{n+1})
}{
M_p\left(H_{p,n+1}(\point,x_{n+1})\right)(x_{p-1})
}
= \frac{p(x_{n+1}|x_p)}{p(x_{n+1}|x_{p-1})}~p(x_p|x_{p-1})\lambda_p(dx_p)
$$
can be explicitly computed 
using the traditional regression formula, or equivalently the updating step of the Kalman filter. 

\subsection{Rare event analysis}\label{rare-event-sec}
\subsubsection{ Importance sampling and twisted measures}
 
Computing the probability
of some events of the form $\{V_n(X_n)\geq a\}$, for some energy like function $V_n$
and some threshold $a$ is often performed using the importance sampling distribution
of the state variable $X_n$ with some multiplicative Boltzmann weight function $e^{\beta V_n(X_n)}$ associated with some temperature parameter $\beta$. These twisted measures
can be described by a Feynman-Kac model in transition space by setting 
$$G_n(X_{n-1},X_n)=e^{\beta [V_n(X_n)-V_{n-1}(X_{n-1})]}$$
For instance, it is easily checked that
 \begin{eqnarray*}
 \PP\left(V_n(X_n)\geq a\right)&=&\EE\left(f_n(X_n)~e^{V_n(X_n)}\right)\\
&=&\EE\left({\bf f_n}({\bf X}_n)~\prod_{0\leq p< n}G_{p}({\bf X}_p)\right)
 \end{eqnarray*}
with
$$
{\bf X}_n=(X_n,X_{n+1})\quad \mbox{\rm and}\quad
 G_n({\bf X}_n)=e^{V_{n+1}(X_{n+1})-V_n(X_n)}
$$
and the test function
$$
 f_n({\bf X}_n)=1_{V_n(X_n)\geq a}~e^{-V_n(X_n)}
$$
In the same vein, we have $$
 \begin{array}{l}
 \EE\left(\varphi_n(X_0,\ldots,X_n)~|~V_n(X_n)\geq a\right)\\
 \\
 =\EE\left(  F_{n,\varphi_n}(X_0,\ldots,X_n)~e^{V_n(X_n)}\right)/\EE\left(  F_{n,1}(X_0,\ldots,X_n)~e^{V_n(X_n)}\right)\\
\\
=\QQ_n( F_{n,\varphi_n})/\QQ_n( F_{n,1})
 \end{array}
$$
with the function
$$
 F_{n,\varphi_n}(X_0,\ldots,X_n)=\varphi_n(X_0,\ldots,X_n)~1_{V_n(X_n)\geq a}~e^{-V_n(X_n)}
$$
We illustrate these rather abstract formula with a Feynman-Kac formulation of European style call options
with exercise price $a$ at time $n$. The prices of these financial contracts are given by formulae of the following
form
$$
\begin{array}{l}
\EE\left((V_n(X_n)-a)_+\right)\\
\\
=\EE\left((V_n(X_n)-a)~1_{V_n(X_n)\geq a}\right)\\
\\
=\PP\left(V_n(X_n)\geq a\right)\times
\EE\left((V_n(X_n)-a)~\left|~V_n(X_n)\geq a\right.\right)
\end{array}
$$
It is now a simple exercise to check that these formulae fit with the Feynman-Kac importance sampling model discussed above. 
Further details on these models, including applications in fiber optics communication
and financial risk analysis
can also be found in the couple of articles~\cite{dm-garnier1,dm-garnier2} and in the article~\cite{carmona-fouque}.
 \subsubsection{ Rare event excursion models}

If we consider Markov excursion type models between a sequence of decreasing subsets $A_n$ or hitting a absorbing level $B$, then choosing an indicator potential function that detects 
if the $n$-th excursion hits $A_n$ before $B$ we find that
$$
\QQ_n=\mbox{\rm Law}(~X~\mbox{\rm hits}~A_n~|~X~\mbox{\rm hits}~A_n~\mbox{\rm before}~B)
$$
and
$$
 \PP(X~\mbox{\rm hits}~A_n~\mbox{before}~B)=\EE\left(\prod_{0\leq p\leq n}1_{A_p}(X_{T_p})\right)=\EE\left(\prod_{0\leq p<n} G_p({\bf X}_{p})\right)
$$
 with the random times
 $$
 T_n:=\inf{\left\{p\geq T_{n-1}~:~X_p\in (A_n\cup B)\right\}}
 $$
and the excursion models
 $$
 {\bf X}_n=(X_p)_{p\in [T_{n},T_{n+1}]}\quad\mbox{\rm \&}\quad
 G_n({\bf X}_n)=1_{A_{n+1}}(X_{T_{n+1}})
 $$
In this notation, it is also easily checked that
 $$
 \EE\left({\bf f_n}(X_{[0,T_{n+1}]})~|~X~\mbox{\rm hits}~A_n~\mbox{before}~B\right)=\QQ_n({\bf f}_n)
 $$
For a more thorough discussion on these excursion particle model, we refer the reader to the series of articles~\cite{cerou,cerou2,lezaud,Delm04,johansen}.
\subsection{Sensitivity measures}

\subsubsection{Kinetic  sensitivity measures}\label{sensitivity-meas}
We let $\theta\in\RR^d$ be some parameter that may represent some kinetic type parameters
related  free evolution model or to the adaptive potential functions. 
We assume that  the free evolution model  $X^{(\theta)}_k$ associated to
some value of the parameter $\theta$, is given by a one-step probability transition of the form
$$
M^{(\theta)}_k(x,dx^{\prime}):=
\PP\left(X^{(\theta)}_k\in dx^{\prime}|X^{(\theta)}_{k-1}=x\right)=
H^{(\theta)}_k(x,x^{\prime})~\lambda_k(dx^{\prime})
$$
for some positive density functions $H_k^{(\theta)}$ and some reference
measures $\lambda_k$. We also consider a collection of 
functions $G_k^{(\theta)}=e^{-V_k^{(\theta)}}$ that depend on $\theta$. We also assume that
the gradient and the Hessian of the logarithms of these functions w.r.t. the parameter $\theta$ are well defined.
We let $ \Gamma^{\theta}_n$ be the Feynman-Kac measure associated with a given value of $\theta$
defined for any function ${\bf f_n}$ on the path space ${\bf E_n}$ by
\begin{equation}\label{greekstheta}
 \Gamma^{\theta}_n({\bf f_n})=\EE\left(
{\bf f_n}(X_0^{(\theta)},\ldots,X_n^{(\theta)})\displaystyle\prod_{0\leq p<n}G_p^{(\theta)}\left(X_p^{(\theta)}\right)
\right)
\end{equation}
We denote by $\Gamma_{n}^{(\theta,N)}$ the $N$-particle approximation measures associated with a given value of the parameter $\theta$ and defined in (\ref{backNnU}).

By using simple derivation calculations, we prove that the first
order derivative of the option value w.r.t. $\theta$ is given by  
\begin{eqnarray*}
 \nabla   \Gamma_{n}^{(\theta)}({\bf f_n})&=&
 \Gamma_{n}^{(\theta)}({\bf f_n}\Lambda^{(\theta)}_n)\\
   \nabla^2    \Gamma_{n}^{(\theta)}({\bf f_n})&=&\Gamma_{n}^{(\theta)}\left[{\bf f_n}
  (\nabla \LL^{(\theta)}_n)^{\prime}  (\nabla \LL^{(\theta)}_n)+ {\bf f_n}  \nabla^2  \LL^{(\theta)}_n  \right]
\end{eqnarray*}
with
$$
\Lambda^{(\theta)}_n:= \nabla \LL^{(\theta)}_n
$$
and the additive functional 
$$ \LL^{(\theta)}_n(x_0,\ldots,x_n):=
\sum_{p=1}^n  \log{\left(G_{p-1}^{(\theta)}(x_{p-1})H^{(\theta)}_{p}(x_{p-1},x_p)\right)}  
  $$
 These quantities are approximated by the unbiased particle models
  \begin{eqnarray*}
 \nabla_N   \Gamma_{n}^{(\theta)}({\bf f_n})&:=&
 \Gamma_{n}^{(\theta,N)}({\bf f_n}\Lambda^{(\theta)}_n)\\
   \nabla^2_N    \Gamma_{n}^{(\theta)}({\bf f_n})&=&\Gamma_{n}^{(\theta,N)}\left[{\bf f_n}
  (\nabla \LL^{(\theta)}_n)^{\prime}  (\nabla \LL^{(\theta)}_n)+ {\bf f_n}  \nabla^2  \LL^{(\theta)}_n  \right]
\end{eqnarray*}
For a more thorough discussion on these Backward Feynman-Kac models, we refer the reader to the three joint articles of the first author with A. Doucet, and S. S. Singh~\cite{sss,sss2,sss3}.
\subsubsection{Gradient estimation of Markov semigroup}

We assume that the underlying stochastic evolution
 is given by an iterated $\RR^{d}$-valued random process given by the following equation 
\begin{equation}\label{defXF}
X_{n+1}:=F_{n}(X_{n})=(F_n\circ F_{n-1}\circ \cdots\circ F_0)(X_0)
\end{equation}
starting at some random state $X_0$, with a sequence of random smooth functions
of the form 
\begin{equation}\label{defXFW}
F_n(x)=\Fa_n(x,W_n)
\end{equation}
with some smooth collection of functions
$$ \Fa_n~:~(x,w)\in\RR^{d+d^{\prime}}\mapsto \Fa_n(x,w)\in \RR^d
$$
and some collection of independent, and independent of $s$, random variables $W_n$  taking values in some $\RR^{d^{\prime}}$, with $d^{\prime}\geq 1$. The semigroup of the Markov chain $X_n$ is the expectation operator defined for any
regular function $f_n$ and any state $x$ by 
$$
P_{n+1}(f_{n+1})(x):=\EE\left(f_{n+1}(X_{n+1})~|~X_0=x\right)=\EE\left(f(X_{n+1}(x))\right)
$$
with the random flows $\left(X_{n}(x)\right)_{n\geq  0}$ defined for any
$n\geq  0$ by the following equation
$$
X_{n+1}(x)=F_{n}(X_{n}(x))
$$
with the initial condition $X_{0}(x)=x$.

By construction, for any $1\leq i,j\leq d$ and any $x\in\RR^d$ we have the first variational equation \index{First variational equation}
\begin{equation}\label{Xpnij}
 \frac{\partial X_{n+1}^i}{\partial x^j}(x)= \sum_{1\leq k\leq d}\frac{\partial F_n^i}{\partial x^k}(X_{n}(x))~
  \frac{\partial X_{n}^k}{\partial x^j}(x)
  \end{equation}
  This clearly implies that
  \begin{equation}\label{Ppnij}
 \frac{\partial P_{n+1}(f)}{\partial x^j}(x)=\EE\left(\sum_{1\leq i\leq d}
  \frac{\partial f}{\partial x^i}(X_{n+1}(x))~\frac{\partial X_{n+1}^i}{\partial x^j}(x)\right)
  \end{equation}

 We denote by  $V_{n}=(V_{n}^{(i,j)})_{1\leq i,j\leq d}$ and  $A_{n}=(A_{n}^{(i,j)})_{1\leq i,j\leq d}$ the random $(d\times d)$ matrices with the $i$-th line and $j$-th column entries
 $$
V_{n}^{(i,j)}(x)=\frac{\partial S_{n}^i}{\partial x^j}(x)$$
and
$$
A_{n}^{(i,j)}(x)=\frac{\partial F^i_n}{\partial x^j}(x)=\frac{\partial \Fa_n^i(\point,W_n)}{\partial x^j}(x):=\Aa_{n}^{(i,j)}(x,W_n)
$$ 

In this notation, the equation (\ref{Xpnij})  can be rewritten in terms of
the following random matrix formulae
\begin{eqnarray}
  V_{n+1}(x)&=&A_{n}(X_n(x))~V_{n}(x)
  :=\prod_{p=0}^nA_p(X_p(x))\label{pro1st}
\end{eqnarray}
 with a product $\prod_{p=0}^nA_p$ of noncommutative 
random elements $A_p$ taken in the order $A_n$, $A_{n-1}$,\ldots, $A_0$.
In the same way, the equation  (\ref{Ppnij}) can be rewritten as
\begin{equation}\label{fkmatrices}
\nabla P_{n+1}(f_{n+1})(x)
=\EE\left(\nabla f_{n+1}(X_{n+1})~V_{n+1}~|~X_0=x\right)
\end{equation}
with
$$
V_{n+1}:=\prod_{0\leq p\leq n}A_p(X_p)
$$
 
We equip the space $\RR^d$ with some norm $\|\point\|$.
We assume that for any state $U_0$ in the unit sphere $\Sa^{d-1}$, we have
$$
\left\|V_{n+1}~U_0
\right\|>0
$$
In this situation, we have
the multiplicative formulae
$$
\nabla f_{n+1}(X_{n+1})~V_{n+1}~U_0=
\left[\nabla f_{n+1}(X_{n+1})~U_{n+1}\right]~\prod_{0\leq p\leq n}\left\|
A_p(X_p)~U_{p}\right\|
$$
with the well defined $\Sa^{d-1}$-valued Markov chain defined by
$$
U_{n+1}={A_n(X_n)U_{n}}/{\left\|
A_n(X_n)U_{n}
\right\|}\left(\Leftrightarrow U_{n+1}=\frac{V_{n+1}~U_0}{\left\|
V_{n+1}~U_0
\right\|}\right)
$$
If we choose $U_0=u_0$, then we obtain the following Feynman-Kac interpretation
of the gradient of a semigroup
\begin{equation}\label{gradient-sg}
\nabla P_{n+1}(f_{n+1})(x)~u_0=\EE\left(F_{n+1}(\Xa_{n+1})~\prod_{0\leq p\leq n}
\Ga_p\left(\Xa_p\right)
\right)
\end{equation}
In the above display, $\Xa_n$ is  the Markov chain sequence
$$
\Xa_n:=\left(X_{n},U_{n},W_n\right)
$$
starting at $(x,u_0,W_0)$, 
and the functions $F_{n+1}$ and $G_n$ are defined by
$$
F_{n+1}(x,u,w):=\nabla f_{n+1}(x)~u
\quad\mbox{\rm
and}
\quad
\Ga_n\left(x,u,w\right):=\left\|
\Aa_n(x,w)~u\right\|
$$

 In computational physics literature, the mean particle approximations of these non commutative Feynman-Kac models are often referred as Resampled Monte Carlo methods~\cite{vanneste}.
 
Roughly speaking, besides the fact that formula (\ref{gradient-sg}) provides a explicit functional Feynman-Kac description
of the the gradient of a Markov semigroup, 
the random evolution model $U_n$ on the unite sphere may be degenerate. More precisely,
the Markov chain $\Xa_n=\left(X_{n},U_{n},W_n\right)$ may not satisfy the regularity properties stated in section~\ref{regularity-conditions}.
We end this section with some rather crude upper bound that can be estimated uniformly w.r.t. the time parameter
under appropriate regularity conditions on the reduced Markov chain model $(X_n,W_n)$. Firstly, we notice that
$$
\Ga_n\left(x,u,w\right):=\left\|
\Aa_n(x,w)~u\right\|\leq \begin{array}[t]{rcl}
G_n(x,w)&:=&
\left\|
\Aa_n(x,w)\right\|\\
&:=&\sup_{u\in \Sa^{d-1}}{\left\|
\Aa_n(x,w)~u\right\|}
\end{array}$$ 
 This implies that
 \begin{eqnarray*}
\left\|\nabla P_{n+1}(f_{n+1})(x)\right\|&:=&\sup_{1\leq i\leq d}{\left|\frac{\partial}{\partial x^i}~P_{n+1}(f_{n+1})(x)\right|}\\
&\leq&\|F_{n+1}\| \times \EE\left(\prod_{0\leq p\leq n}
G_p\left(X_p,W_p\right)
\right)
\end{eqnarray*}
The r.h.s. functional expectation in the above equation can be approximated using the particle approximation (\ref{free-energy})
of the multiplicative  Feyman-Kac formulae (\ref{multiplicative}), with reference Markov chain $(X_n,W_n)$ and potential functions $G_n$.

 \subsection{Partial observation models}

 \subsubsection{Nonlinear filtering models}\label{filtering-sec}
In this section we introduce one of the most important example of estimation problem with partial observation, namely the nonlinear filtering model.
This model has been the starting point of the application of particle models to engineering sciences, and
more particularly  to advanced signal processing. 

The first rigorous subject on the stochastic  modeling, and the rigorous theoretical analysis of
particle filters has been started in the mid 1990's in the article~\cite{dm96}. For a detailed discussion on the application domains of particle filtering, with a precise 
bibliography we refer the reader  to any of the following references~\cite{dmsp,Delm04}, and ~\cite{dd-2012, doucet}.

The  typical model is given by a reference Markov chain model
  $X_n$, and some partial and noisy observation  $Y_n$. The pair process $(X_n,Y_n)$ usually forms a Markov chain on some product space $E^X\times E^Y$ with elementary transitions given
 \begin{equation}\label{Markov-XY}
 \begin{array}{l}
 \PP\left((X_n,Y_n)\in d(x,y)~|~(X_{n-1},Y_{n-1})\right)\\
 \\=M_n(X_{n-1},dx)~
 g_n(x,y)~\lambda_n(dy)
 \end{array}
 \end{equation} 
for some positive likelihood function $g_n$, and some reference probability measure $\lambda_n$ on $E^Y$, and the elementary Markov transitions
$M_n$ of the Markov chain $X_n$. 
If we take 
\begin{equation}\label{likelihood-ref}
G_n(x_n)=p_n(y_n|x_n)=g_n(x_n,y_n)
\end{equation} the likelihood function of a given observation $Y_n=y_n$ and a signal state $X_n=x_n$ associated with a filtering or an hidden Markov chain problem, then
 we find that
$$
\mathbb{Q}_{n}=\mbox{\rm Law}((X_0,\ldots,X_n)~|~\forall 0\leq p<n\quad Y_p= y_p)
$$
and
$$ \Za_{n+1}=p_n(y_0,\ldots,y_{n})
$$
In this context, the optimal one step predictor $\eta_n$  and the optimal filter $\widehat{\eta}_n$ are given by the $n$-th time
marginal distribution
\begin{equation}\label{pred-def}
\eta^{[y_0,\ldots,y_{n-1}]}_n=\eta_n=\mbox{\rm Law}\left(X_n~|~\forall 0\leq p<n\quad Y_p= y_p\right)
\end{equation}
and
\begin{equation}\label{update-def}
\widehat{\eta}^{[y_0,\ldots,y_{n}]}_n=\widehat{\eta}_n=\Psi_{G_n}(\eta_n)=\mbox{\rm Law}\left(X_n~|~\forall 0\leq p\leq n\quad Y_p= y_p\right)
\end{equation}

\begin{remark}
We can combine these filtering models with the probability restriction models discussed in section~\ref{proba-restrict-sec}, or with the 
 rare event analysis presented in section~\ref{rare-event-sec}. For instance, if we replace the potential likelihood function
 $G_n$ defined in (\ref{likelihood-ref}) by the function
 $$
 G_n(x_n)=g_n(x_n,y_n)~1_{A_n}(x_n)
 $$ 
 then we find that
 $$
 \mathbb{Q}_{n}=\mbox{\rm Law}((X_0,\ldots,X_n)~|~\forall 0\leq p<n\quad Y_p= y_p,~X_p\in A_p)
 $$
\end{remark}
 \subsubsection{Approximated filtering models}\label{abc-sec}
 We return to the stochastic  filtering model discussed in section~\ref{filtering-sec}.
In some instance, the likelihood functions $x_n\mapsto g_n(x_n,y_n)$ in (\ref{likelihood-ref}) are computationally 
 intractable, or too expensive to evaluate.

 To solve this problem, a natural solution is to sample pseudo-observations. The central idea is to sample the signal-observation
 Markov chain $${\bf X}_n=\left(X_n,Y_n\right)\in {\bf E^{\bf X}}=(E^X\times E^Y)$$ and compare
 the values of the sampled observations with the real observations.

To describe with some precision these models, we notice that the
 transitions of  ${\bf X}_n$ are given by
 $$
 {\bf M_n}({\bf X_{n-1}}, d(x,y))=M_n(X_{n-1},dx)~
 g_n(x,y)~\lambda_n(dy)
 $$   
 To simplify the presentation, we further assume that $E^Y=\RR^d$, for some $d\geq 1$, and
 we let $g$ be a Borel bounded non negative function such that
 $$
 \int ~g(u) du=1\quad\int u g(u)~du=0\quad\mbox{\rm and}\quad
\int |u|^3 g(u)~du<\infty
 $$
Then, we set for any $\epsilon>0$, and any ${\bf x}=(x,y)\in (E^X\times E^Y)$
$$
g_{\epsilon,n}((x,y),z)=\epsilon^{-d}~g\left((y-z)/\epsilon\right)
$$
 
 Finally, we  let $({\bf X_n,Y^{\epsilon}_n})$ be the Markov chain on the augmented state space $
 \left({\bf E^{\bf X}}\times E^Y\right)=\left((E^X\times E^Y)\times E^Y\right)$ with  transitions given
\begin{equation}\label{pseudo-filtering}
  \begin{array}{l}
 \PP\left(({\bf X_n,Y^{\epsilon}_n})\in d({\bf x},y)~|~({\bf X_{n-1},Y^{\epsilon}_{n-1}})\right)\\
 \\= {\bf M_n}({\bf X_{n-1}},{\bf dx})~g_{\epsilon,n}({\bf x},y)~dy
 \end{array}
 \end{equation}
 This approximated filtering  problem has exactly the same form as the one introduced in (\ref{Markov-XY}). 
Here,  the particle approximation model are defined in terms of signal-observation valued particles, and the selection potential function is given by the pseudo-likelihood functions $g_{\epsilon,n}(\point,y_n)$, where $y_n$ stands for the value of the observation sequence at time $n$.

For a more detailed discussion on these particle models, including the convergence analysis of the 
approximated filtering model,
we refer the reader to the article~\cite{djjp,dm-jacod}. These particle models are sometimes called 
convolution particle filters~\cite{vila}. In Bayesian literature, these approximated filtering  models are  termed as  
Approximate Bayesian Computation (and often abbreviated with {\em the acronym ABC}).
 \subsubsection{Parameter estimation in hidden Markov chain models}

We consider a pair signal-observation filtering model $(X,Y)$ that depend on some
random variable $\Theta$ with distribution $\mu$ on some state space $S$. Arguing as above, if 
we take $$G_{\theta,n}(x_n)=p_n(y_n|x_n,\theta)$$ the likelihood function of a given observation $Y_n=y_n$ and a signal state $X_n=x_n$ and a realization of the parameter  $\Theta=\theta$,
then the  $n$-th time marginal of $\QQ_n$ is given by
$$
\eta_{\theta,n}=\mbox{\rm Law}(X_n~|~\forall 0\leq p<n\quad Y_p= y_p,~\theta)
$$
Using that the multiplicative formula (\ref{multiplicative}), we prove that
$$
\Za_{n+1}(\theta)=p_n(y_0,\ldots,y_{n}|\theta)=\prod_{0\leq p\leq n}\eta_{\theta,p}(G_{\theta,p})
$$
with
\begin{eqnarray*}
\eta_{\theta,p}(G_{\theta,p})&=&p(y_p|y_0,\ldots,y_{p-1},\theta)\\
&=&\int~
p(y_p|x_p,\theta)~dp(x_p|\theta,y_0,\ldots,y_{p-1})\\
&=&\int G_{\theta,p}(x_p)~\eta_{\theta,p}(dx_p)
\end{eqnarray*}
from which we conclude that
$$
\PP(\Theta\in d\theta~|~\forall 0\leq p\leq n\quad Y_p= y_p)=\frac{1}{\Za_n}~\Za_{n}(\theta)~\mu(d\theta)$$
with
$$
\Za_n:=\int~\Za_{n}(\theta)~\mu(d\theta)
$$

In some instance, such as in conditionally linear Gaussian models, the normalizing constants $\Za_{n}(\theta)$ can be computed
explicitly, and we can use a Metropolis-Hasting style Markov chain Monte Carlo method to sample
the target measures $\mu_n$. As in section~\ref{i-mcmc-intro}, we can also turn this scheme into an interacting Markov chain Monte Carlo algorithm.

Indeed, let us choose a Markov chain Monte Carlo type local moves 
$\mu_{n}=\mu_{n}{M_n}$ with prescribed target measures
$$
\mu_{n}(d\theta):=\frac{1}{\Za_n}~\Za_{n}(\theta)~\mu(d\theta)
$$
Notice that
$$
\Za_{n+1}(\theta)=\Za_{n}(\theta)\times \eta_{\theta,n}(G_{\theta,n})\Rightarrow
\mu_{n+1}=\Psi _{G_{n}}(\mu _{n})
$$
with the Boltzmann-Gibbs transformations defined in (\ref{def-BG-ref-intro}) associated with the potential function
$$
G_n(\theta)=\eta_{\theta,n}(G_{\theta,n})
$$
By construction, we have
$$
\mu_{n+1}=\mu_{n+1}M_{n+1}=\Psi _{G_{n}}(\mu_{n})M_{n+1}
$$
from which we conclude that
$$
\mu_{n}(f)=\EE\left(f(\theta_n)~\prod_{0\leq p<n}G_p(\theta_p)\right)~/~\EE\left(\prod_{0\leq p<n}G_p(\theta_p)\right)
$$
with the reference Markov chain
$$
\PP\left(\theta_n\in d\theta|\theta_{n-1}\right)=M_n(\theta_{n-1},d\theta)
$$
In addition, we have
$$
\Za_{n+1}=\int \Za_{n}(\theta)~G_n(\theta)~\mu(d\theta)=  \Za_{n}~\mu_n(G_n)=\prod_{0\leq p\leq n}\mu_p(G_p)
$$
\begin{remark}
For more general models, we can use the particle Markov chain Monte Carlo methodology presented in section~\ref{p-mcmc}. When the likelihood functions are too expensive to evaluate,
we can also combine these particle models with the pseudo-likelihood stochastic models (\ref{pseudo-filtering}) discussed in section~\ref{abc-sec}.
\end{remark}
 \subsubsection{ Interacting Kalman-Bucy filters}

We use the same notation as above, but we assume that $\Theta=(\Theta_n)_{n\geq 0}$ is a random sample of a 
stochastic process $\Theta_n$ taking values in some state spaces $S_n$. If we consider
the Feynman-Kac model associated with the Markov chain $\Xa_n=(\Theta_n,\eta_{\Theta,n})$
and the potential functions $$\Ga_n(\Xa_n)=\eta_{\Theta,n}(G_{\Theta,n})$$ then we find that
$$
\QQ_n=\mbox{\rm Law}(\Theta_0,\ldots,\Theta_n~|~\forall 0\leq p<n\quad Y_p= y_p)
$$
and the $n$-th time marginal are clearly given by
$$\eta_n=\mbox{\rm Law}(\Theta_n~|~\forall 0\leq p<n\quad Y_p= y_p)
$$

Assuming that the pair $(X,Y)$ is a linear and gaussian filtering model given $\Theta$,
the measures $\eta_{\Theta,n}$ coincide with the one step predictor of the Kalman-Bucy filter
and the potential functions $\Ga_n(\Xa_n)$ can be easily computed by gaussian integral calculations. In this situation, the conditional distribution of the parameter $\Theta$ is 
given by a Feynman-Kac model $\QQ_n$ of a the free Markov chain $\Xa_n$ weighted by some  Boltzmann-Gibbs exponential weight function  $$\prod_{0\leq p<n}\Ga_p(\Xa_p)=p_n(y_0,\ldots,y_{n}|\Theta_0,\ldots,\Theta_n)$$ that reflects
 the likelihood of the path sequence $(\Theta_0,\ldots,\Theta_n)$. For a more thorough discussion on these interacting Kalman filters, we refer the reader to section~2.6 and section~12.6 in the monograph~\cite{Delm04}.

 \subsubsection{ Multi-target tracking models}

Multiple-target tracking problems deal with correctly estimating several manoeuvring and interacting targets simultaneously  given a sequence
of noisy and partial observations. At every time $n$, the first moment of the occupation 
measure $\Xa_{n}:=\sum_{i=1}^{N_{n}}\delta_{X^{i}_{n}}$ of some spatial branching signal 
is given for any regular function $f$  by the following formula:
$$
\gamma_{n}(f):=\EE\left(\Xa_{n}(f)\right)
\quad\mbox{\rm with}\quad
\Xa_{n}(f):=\int~f(x)~\Xa_{n}(dx)
$$
For null spontaneous birth measures, these measures coincide with that of an unnormalized  Feynman-Kac
model with some spatial branching potential functions $G_n$ and some free evolution target
model $X_n$.

In more general situations, the approximate  filtering equation is given by
the Malher's multi-objective filtering approximation based on the propagation of the first conditional moments of Poisson approximation models~\cite{maler0,maler}. These evolution equations are rather complex to introduce and notationally consuming. Nevertheless, as  the first moment evolution of any spatial and marked branching process,
they can be abstracted by an unnormalized Feynman-Kac
model with nonlinear potential functions~\cite{caron1,caron2,caron3}.

\subsubsection{Optimal stopping problems with partial observations}

We consider the partially observed Markov chain model discussed in (\ref{Markov-XY}).
 The Snell envelop  associated with an optimal stopping problem 
with finite horizon, payoff style function $f_n(X_n,Y_n)$, and noisy observations $Y_n$ as some Markov process, is given by
$$
U_k:= \sup_{\tau \in \Ta^Y_k} \EE(f_{\tau}(X_{\tau},Y_{\tau})|(Y_0,\ldots,Y_k))
$$   
where $\Ta_k^Y$ stands for the set of all stopping times $\tau$ taking values in $\{k,\ldots, n\}$, whose values are measurable w.r.t. the sigma field generated by the observation sequence
$Y_p$, from $p=0$ up to the current time $k$. We denote by $\eta^{[y_0,\ldots,y_{n-1}]}_{n}$ and $\widehat{\eta}^{[y_0,\ldots,y_{n}]}_{n}$ the 
conditional distributions defined in (\ref{pred-def}) and (\ref{update-def}). In this notation, for any $0\leq k\leq n$ we have that
\begin{equation}\label{ref-snell-us-partial-obs}
\begin{array}{l}
\EE(f_{\tau}(X_{\tau},Y_{\tau})|(Y_0,\ldots,Y_k))
\\
\\=\EE\left(F_{\tau}\left(Y_{\tau},\widehat{\eta}^{[Y_0,\ldots,Y_{\tau}]}_{\tau}\right)~|~(Y_0,\ldots,Y_k)\right)
\end{array}
\end{equation}
with the conditional payoff function
$$
F_p\left(Y_p,\widehat{\eta}^{[Y_0,\ldots,Y_p]}_p\right)=\int~
\widehat{\eta}^{[Y_0,\ldots,Y_p]}_p(dx)~f_{p}(X_{p},Y_{p})
$$

It is rather well known that $$
\Xa_p:=\left(Y_p,\widehat{\eta}^{[Y_0,\ldots,Y_p]}_p\right)
$$
is a Markov chain with elementary transitions defined by
$$
\begin{array}{l}
\EE\left[F_p\left(Y_p,\widehat{\eta}^{[Y_0,\ldots,Y_p]}_p\right)~\left|~\left(Y_{p-1},\widehat{\eta}^{[Y_0,\ldots,Y_{p-1}]}_{p-1}\right)=(y,\mu)\right.\right]
\\
\\
=\displaystyle\int \lambda_p(dy_p)~\mu M_{p}\left(g_p(\point,y_p)\right)~
F_p\left(y_p,\Psi_{g_p(\point,y_p)}\left(\mu M_p\right)
\right)
\end{array}
$$
A detailed proof of this assertion can be found in any textbook on advanced stochastic filtering theory. For instance, the book of W. Runggaldier, L. Stettner~\cite{stettner}
provides a detailed treatment  on  discrete time non linear filtering, and related partially observed control models.
  
Roughly speaking, using some abusive Bayesian notation, we have
\begin{eqnarray*}
\eta^{[y_0,\ldots,y_{p-1}]}_{p}(dx_p)&=&dp_p(x_p~|~(y_0,\ldots,y_{p-1}))\\
&=&\int dp_p(x_p~|~x_{p-1})\times p_n(x_{p-1}~|~(y_0,\ldots,y_{p-1}))\\
&=&\widehat{\eta}^{[y_0,\ldots,y_{p-1}]}_{p-1}M_{p}(dx_p)
\end{eqnarray*}
and
$$
\begin{array}{l}
\Psi_{g_p(\point,y_p)}\left(\widehat{\eta}^{[y_0,\ldots,y_{p-1}]}_{p-1}M_p\right)(dx_p)
\\
\\
=\displaystyle\frac{p(y_p|x_p)}{\int p_p(y_p~|~x^{\prime}_p)~dp_p(x^{\prime}_p~|~(y_0,\ldots,y_{p-1}))}~dp_p(x_p~|~(y_0,\ldots,y_{p-1}))\\
\\=dp_p(x_p~|~(y_0,\ldots,y_{p-1},y_p))
\end{array}
$$
from which we prove that
\begin{eqnarray*}
\mu M_p(g_p(\point,y_p))&=&\int p_p(y_p~|~x_p)~dp_p(x_p~|~(y_0,\ldots,y_{p-1}))\\
&=&p_p(y_p~|~(y_0,\ldots,y_{p-1}))
\end{eqnarray*}
and
$$
\Psi_{g_p(\point,y_p)}\left(\mu M_p\right)=
\widehat{\eta}^{[y_0,\ldots,y_{p}]}_{p}
$$
as soon as $\mu=\widehat{\eta}^{[y_0,\ldots,y_{p-1}]}_{p-1}~\left(
\Rightarrow \mu M_p={\eta}^{[y_0,\ldots,y_{p-1}]}_{p}
\right)$

From the above discussion, we can rewrite
 (\ref{ref-snell-us-partial-obs}) as the Snell envelop of a fully observed
 augmented Markov chain sequence 
 $$
\EE(f_{\tau}(X_{\tau},Y_{\tau})|(Y_0,\ldots,Y_k))=\EE\left(F_{\tau}\left(\Xa_{\tau}\right)~|~(\Xa_0,\ldots,\Xa_{k})\right)
$$

The Markov chain $\Xa_n$ takes values in an infinite dimensional state space, and it can rarely be sampled without some addition level of approximation. 
Using the $N$-particle approximation models, we can replace the chain $\Xa_n$ by the $N$-particle approximation model defined by
$$
\Xa_n^N:=\left(Y_p,\widehat{\eta}^{([Y_0,\ldots,Y_p],N)}_p\right)
$$
where 
$$
\widehat{\eta}^{([Y_0,\ldots,Y_p],N)}_p:=\Psi_{g_p(\point,Y_p)}\left(\eta^{([Y_0,\ldots,Y_{p-1},N)]}_{p-1}\right)
$$ stands for the updated
measure associated  associated with the likelihood
selection functions $g_p(\point,Y_p)$. The $N$-particle approximation of the Snell envelop is now given by
 $$
\EE(f_{\tau}(X_{\tau},Y_{\tau})|(Y_0,\ldots,Y_k))\simeq_{N\uparrow\infty}\EE\left(F_{\tau}\left(\Xa^N_{\tau}\right)~|~(\Xa_0^N,\ldots,\Xa^N_{k})\right)
$$
In this interpretation, the $N$-approximated optimal stopping problem amounts to compute the quantities
$$
U_k^N:= \sup_{\tau \in \Ta^{N}_k} \EE\left(F_{\tau}\left(\Xa^N_{\tau}\right)~|~(\Xa_0^N,\ldots,\Xa^N_{k})\right)
$$   
where $\Ta_k^{N}$ stands for the set of all stopping times $\tau$ taking values in $\{k,\ldots, n\}$, whose values are measurable w.r.t. the sigma field generated by the Markov chain sequence
$\Xa^N_k$, from $p=0$ up to time $k$.

\subsection{Markov chain restriction models}
 \subsubsection{Markov confinement  models}
One of the simplest example of Feynman-Kac conditional distributions is given by 
 choosing indicator functions $G_n=1_{A_n}$ of measurable subsets $A_n\in \Ea_n$
 s.t. $\PP\left(\forall 0\leq p<n~X_p\in A_p\right)>0$. In this situation, it is readily checked that
$$
\mathbb{Q}_{n}=\mbox{\rm Law}((X_0,\ldots,X_n)~|~\forall 0\leq p<n~X_p\in A_p)
$$
and
$$
\Za_{n}=\PP\left(\forall 0\leq p<n~X_p\in A_p\right)
$$
This Markov chain restriction model fits into the particle absorption model (\ref{absorption-intro}) presented in the introduction. 
For a detailed analysis of these stochastic models, and their particle approximations, we refer the reader to the articles~\cite{dm3,dmsp,soft}, and the
monograph~\cite{Delm04}. 
  \subsubsection{Directed polymers and self-avoiding walks}
 
The conformation of polymers in a 
chemical solvent can be seen as the realization of a Feynman-Kac
 distribution of a free Markov chain weighted by some  Boltzmann-Gibbs exponential weight function that reflects
 the attraction or the repulsion forces between the monomers. For instance, if we consider the historical process
$${\bf X_n}=(X_0,\ldots,X_n)$$ and
$${\bf G_n(X_n)}=1_{\not\in\{X_p, ~p<n\}}(X_n)$$ then we find that
\begin{eqnarray*}
\QQ_n&=&\mbox{\rm Law}({\bf X_n}~|~\forall 0\leq p<n\quad {\bf X_p}\in A_p)\\
&=&\mbox{\rm Law}((X_0,\ldots,X_n)~|~\forall 0\leq p<q<n\quad X_p\not=X_q)
\end{eqnarray*}
with the set $A_n=\{{\bf G_n=1}\}$, and  the normalizing constants 
$$ \Za_{n}=\PP(\forall 0\leq p<q<n\quad X_p\not=X_q)
$$

 \subsection{Particle absorption models}
 
We return to the particle absorption model (\ref{absorption-intro}) presented in the introduction. For instance, we can  assume that the potential function $G_n$ and Markov transitions 
$M_n$ are defined by
$
G_n(x)=e^{-V_n(x) h}
$, 
and
\begin{equation}\label{h-ref-gm}
M_n(x,dy)=\left(1-\lambda_n h\right) \delta_x(dy)+\lambda_n h~K_n(x,dy)
\end{equation}
for some non negative and bounded function $V_n$, some positive parameter $\lambda_n\leq 1/h$, $h>0$,
and some Markov transition
$K_n$.

We also mention that the confinement models described above can also be interpreted as a particle absorption
model related to hard obstacles. In branching processes and population dynamics literature, the model $X^c_n$ often represent the number of individuals of a given specie~\cite{ferrari,gosselin,david}. Each individual can die or reproduce.  The state $0\in E_n=\NN$ is interpreted as a trap, or as an hard obstacle, in the sense that the specie disappear as soon as $X^c_n$ hits $0$. For a more thorough discussion on particle motions in absorbing
medium with hard and soft obstacles, we refer the reader to the pair of articles~\cite{dm3,soft}.

\subsubsection{Doob h-processes}\label{doob-h-sec}

We consider a  time homogeneous Feynman-Kac model $(G_n,M_n)=(G,M)$ on some measurable state space $E$, and we set
$$
Q(x,dy)=G(x)M(x,dy)
$$ 
We also assume that $G$ is uniformly bounded above and below by some positive constant, and
the Markov transition
 $M$ is reversible w.r.t. some probability measure $\mu$ on $E$, with $M(x,\point)\simeq \mu$ and
 $dM(x,\point)/d\mu\in\LL_2(\mu)$. 
 We
  denote by $\lambda$ the largest eigenvalue of the integral operator $Q$ on $\LL_2$, and by $h(x)$ a positive eigenvector
$$
Q(h)=\lambda h
$$
The Doob $h$-process corresponding to the ground state eigenfunction $h$ defined above 
is a Markov chain $X^h_n$ with the time homogeneous Markov transition
$$
M^h(x,dy):=\frac{1}{\lambda }\times h^{-1}(x)Q(x,dy)h(y)=\frac{M(x,dy)h(y)}{M(h)(x)}
$$
and initial distribution $\eta^h_0(dx)\propto h(x)~\eta_0(dx)$. 
By construction, we have $G=\lambda h/M(h)$ and therefore
$$
\Gamma_n(d(x_0,\ldots,x_n))
= \lambda^{n}~\eta_0(h)~\PP^h_n(d(x_0,\ldots,x_n))~\frac{1}{h(x_n)}
$$
where $\PP^h_n$ stands for the law of the historical process  $${\bf X^h_n}=(X^h_0,\ldots,X^h_n)$$
We conclude that
$$
d\QQ_n=\frac{1}{\EE(h^{-1}(X^h_n))}~h^{-1}(X^h_n)~~d\PP^h_n
$$
with the normalizing constants
$$
\Za_n=\lambda^{n}~\eta_0(h)~\EE(h^{-1}(X^h_n))
$$

\subsubsection{Yaglom limits and quasi-invariant measures}\label{yaglom-sec}
We return to the time homogeneous Feynman-Kac models introduced in section~\ref{doob-h-sec}.
Using the particle absorption interpretation (\ref{absorption-intro}) we have
$$
\mbox{\rm Law}((X^c_0,\ldots,X^c_n)~|~T^c\geq n)=\frac{1}{\EE(h^{-1}(X^h_n))}~h^{-1}(X^h_n)~~d\PP^h_n
$$
and
\begin{equation}\label{cv-yaglom1}
 \Za_{n}=\PP\left(T^c\geq n\right)=\lambda^{n}~\eta_0(h)~\EE(h^{-1}(X^h_n))\longrightarrow_{n\uparrow\infty}0
\end{equation}
Letting
$\eta^h_n:=\mbox{\rm Law}(X^h_n)$,  we readily prove the following formulae
$$
\eta_n=\Psi_{1/h}(\eta^h_n)\quad\mbox{and}\quad \eta^h_n=\Psi_{h}(\eta_n)
$$
Whenever it exists, the Yaglom limit of the measure $\eta_0$ is  is defined as the limiting of measure 
\begin{equation}\label{cv-yaglom2}
\eta_n\longrightarrow_{n\uparrow\infty}\eta_{\infty}=\Psi_G(\eta_{\infty})M
\end{equation}
of the Feynman-Kac flow $\eta_n$, when n tends to infinity. We also say that $\eta_0$ is quasi-invariant measure is we have
$\eta_0=\eta_n$, for any time step. When the Feynman-Kac flow $\eta_n$ is asymptotically stable, in the sense that
it forgets its initial conditions, we also say that the quasi-invariant measure $\eta_{\infty}$ is the Yaglom measure.
Whenever it exist, we let $\eta^h_{\infty}$ be the invariant measure of the $h$-process $X^h_n$.
Under our assumptions, it is a now simple exercise to check that $$
\eta_{\infty}=\Psi_{M(h)}(\mu)\quad\quad\mbox{\rm and}\quad
\eta_{\infty}^h:=\Psi_{h}(\eta_{\infty})=\Psi_{hM(h)}(\mu)
$$
Quantitative convergence estimates of the limiting formulae (\ref{cv-yaglom1}) and (\ref{cv-yaglom2}) can be  derived
using the stability properties of the Feynman-Kac models developed in chapter~\ref{fk-sg-chap}. For a more thorough discussion
on these particle absorption models, we refer the reader to the articles of the first author with A. Guionnet~\cite{dmg-cras,DeGu:ont},
the ones with L. Miclo~\cite{dmsp,dm3}, the one with A. Doucet~\cite{soft}, and the monographs~\cite{Delm04,dd-2012}.

\section{Feynman-Kac semigroup analysis}\label{fk-sg-chap}

\subsection{Introduction}

As we mentioned in  section~\ref{sec-ips-fk-intro}, the concentration analysis of  particle models 
is intimately related to the regularity  properties of the limiting nonlinear  semigroup.
In this short section, we survey some selected topics on the theory of Feynman-Kac semigroup developed in
the series of articles~\cite{DeGu:ont,dmsp,sss}. For more recent treatments, we also refer the reader to the books~\cite{Delm04,dd-2012}. 

We begin this chapter with a discussion on path space models. Section~\ref{historical-sec}
is concerned with Feynman-Kac historical processes and Backward Markov chain interpretation models
We show that the the $n$-th marginal measures $\eta_n$ of Feynman-Kac model with a reference historical Markov process
coincides with the path space measures $\QQ_n$ introduced in (\ref{defQn}). 

The second part of this section is dedicated to the proof of the Backward Markov chain formulae (\ref{backwardt}). In section~\ref{sec-fk-sg}, we analyze the regularity and the semigroup structure of the normalized and unnormalized Feynman-Kac distribution flows $\eta_n$ and $\gamma_n$.  

Section~\ref{quantitative-contraction} is concerned with the stability properties of the normalized Feynman-Kac distribution flow.
In a first section, section~\ref{regularity-conditions}, we present regularity conditions on the potential functions $G_n$ and on the Markov transitions $M_n$, under which the Feynman-Kac semigroup forgets exponentially fast its initial condition.
Quantitative contraction theorems are provided in section~\ref{quantitative-contract-sec}. 

We illustrate these results with three applications  related respectively to time discretization techniques,
simulated annealing type schemes, and path space models. 

The last two sections of this chapter, section~\ref{mean-field-sec} and section~\ref{loc-sampling-sec}, are concerned with  mean field stochastic particle models and  local sampling random field models.

\subsection{Historical and backward models}\label{historical-sec}

The historical process associated with some reference Markov chain
$X_n$ is defined by the sequence
of random paths 
$$
{\bf X}_n=\left(X_0,\ldots,X_n\right)\in {\bf E_n}:=\left(E_0\times\ldots\times E_n\right)
$$
Notice that the Markov transitions of the chain ${\bf X}_n$ is given for
any ${\bf x}_{n-1}=(x_0,\ldots,x_{n-1})$ and ${\bf y_n}=(y_0,\ldots,y_n)=({\bf y_{n-1}},y_n)$ by
the following formulae
\begin{equation}\label{eqMb-ref}
{\bf M}_n({\bf x}_{n-1},d{\bf y_n})=\delta_{{\bf x}_{n-1}}(d{\bf y_{n-1}})~M_n(y_{n-1},dx_n)
\end{equation}
We consider a sequence of $(0,1]$-valued potential functions ${\bf G}_n$ on ${\bf E_n}$ whose values only depend
on the final state of the paths; that is, we have that 
\begin{equation}\label{eqGb-ref}
{\bf G}_n~:~{\bf x}_n=(x_0,\ldots,x_n)\in {\bf E_n}\mapsto {\bf G}_n({\bf x}_n)=G_n(x_n)\in (0,1]
\end{equation}
with some $(0,1]$-valued potential function $G_n$ on $E_n$.

We let $(\gamma_n,\eta_n)$ the Feynman-Kac model associated with the pair 
 $({\bf G}_n,{\bf M}_n)$ on the path spaces $ {\bf E_n}$. By construction, for any
 function ${\bf f}_n$ on ${\bf E_n}$, we have
 \begin{eqnarray*}
 \gamma_n({\bf f}_n)&=&\EE\left({\bf f}_n({\bf X}_n)~\prod_{0\leq p<n}{\bf G}_p({\bf X}_p)\right)\\
 &=&\EE\left({\bf f}_n(X_0,\ldots,X_n)~\prod_{0\leq p<n}G_p(X_p)\right)
 \end{eqnarray*}
from which we conclude that
\begin{equation}\label{equivalence-principle}
\gamma_n=\Za_n~\QQ_n
\quad\mbox{\rm and}\quad \eta_n=\QQ_n
\end{equation}
 where $\QQ_n$ is the  Feynman-Kac measure on path space associated with the pair 
 $(G_n,M_n)$, and defined in (\ref{defQn}).

We end this section with the proof of the  backward formula (\ref{backwardt}).
Using the decomposition
$$
\QQ_n(d(x_0,\ldots,x_n))=\frac{\Za_{n-1}}{\Za_n}~\QQ_{n-1}(d(x_0,\ldots,x_{n-1}))~Q_{n}(x_{n-1},dx_n)$$
we prove the following formulae
\begin{eqnarray}
\eta_n(dx_n)&=&\frac{\Za_{n-1}}{\Za_n}~\eta_{n-1}Q_{n}(dx_n)\\
\nonumber\\
&=&\frac{\Za_{n-1}}{\Za_n}~\eta_{n-1}\left(G_{n-1}H_n(\point,x_n)\right)~\lambda_n(dx_n)\label{matrixformu}
\end{eqnarray}
and
$$
\frac{\Za_{n}}{\Za_{n-1}}=\eta_{n-1}Q_{n}(\un)= \eta_{n-1}(G_{n-1})
$$

This implies that
$$
\begin{array}{l}
\displaystyle\frac{d\eta_{n-1}Q_n}{d\eta_n}(x_n)\times\frac{d\eta_{n-2}Q_{n-1}}{d\eta_{n-1}}(x_{n-1})\times\cdots\times\frac{d\eta_{0}Q_1}{d\eta_1}(x_1)\\
\\
=\displaystyle\frac{\Za_n}{\Za_{n-1}}\times\frac{\Za_{n-1}}{\Za_{n-2}}\times\cdots\times\frac{\Za_1}{\Za_0}=\Za_n
\end{array}
$$ 

Using these observations, we readily prove the desired backward decomposition formula.

\subsection{Semigroup models}\label{sec-fk-sg}

This section is concerned with the semigroup structure and the weak regularity properties of Feynman-Kac models.

\begin{definition}
We denote by $$\Phi_{p,n}(\eta_p)=\eta_n\quad\mbox{\rm and} \quad\gamma_p Q_{p,n}=\gamma_n$$
with $0\leq p\leq n$, the linear and the nonlinear semigroup associated with the unnormalized and the normalized Feynman-Kac measures. For $p=n$, we use the convention $\Phi_{n,n}=Id$, the identity mapping. 
\end{definition}

Notice that $Q_{p,n}$ has the following functional representation
$$
Q_{p,n}(f_n)(x_p):=\EE\left(f_n(X_n)~\prod_{p\leq q<n}G_q(X_q)~|~X_p=x_p\right)
$$

\begin{definition}
We let  $G_{p,n}$ and $P_{p,n}$ be the potential functions and the Markov transitions 
defined by
$$
Q_{p,n}(\un)(x)=G_{p,n}(x)\quad\mbox{\rm and}\quad P_{p,n}(f)=\frac{Q_{p,n}(f)}{Q_{p,n}(\un)}
$$
we also set 
$$
g_{p,n}:=\sup_{x,y}\frac{G_{p,n}(x)}{G_{p,n}(y)}\quad\mbox{\rm and}\quad\beta(P_{p,n})=\sup \mbox{\rm osc}(P_{p,n}(f))
$$
The r.h.s. supremum is taken the set of  functions $\mbox{\rm Osc}(E)$. To simplify notation, for $n=p+1$ we
have also set $$G_{p,p+1}=Q_{p,p+1}(\un)=G_p$$ and sometimes we write $g_p$ instead of $g_{p,p+1}$. 
\end{definition}

The particle concentration inequalities developed in chapter~\ref{fk-particle-chap} will be expressed in terms of the following parameters.

\begin{definition}\label{deftaukl}
 For any $k,l\geq 0$, we also set
\begin{equation}\label{def-tau}
\tau_{k,l}(n):=\sum_{0\leq p\leq n} g_{p,n}^k~\beta(P_{p,n})^l\quad\mbox{\rm and}\quad
\kappa(n):=\sup_{0\leq p\leq n}{\left(g_{p,n}\beta(P_{p,n})\right)}
\end{equation}

\end{definition}

Using the fact that
\begin{equation}
 \eta _{n}(f_{n}):={%
\eta _{p}Q_{p,n}(f_{n})}/{\eta _{p}Q_{p,n}(1)}\label{etaQpn}
\end{equation}
we prove the following formula
$$
\Phi_{p,n}\left(  \eta_{p}\right) =\Psi_{G_{p,n}}\left(  \eta_{p}\right)P_{p,n}
$$
for  any
$0\leq p\leq n$.

As a direct consequence of (\ref{tv-Phi-equal}) and (\ref{tv-Phi-contract}), we quote the 
following weak regularity property of the Feynman-Kac semigroups.
\begin{prop}\label{propcontract-de}
For $[0,1]$-valued potential function $G_n$, and 
any couple of measures $\nu,\mu$ on the set $E$ s.t. $\mu(G_{p,n})\wedge\nu(G_{p,n})>0$, we have the decomposition 
$$
\Phi_{p,n}(\mu)-\Phi_{p,n}(\nu)=\frac{1}{\nu(G_{p,n})}~(\mu-\nu)S_{G_{p,n},\mu}P_{p,n}
$$
In addition, we have the following Lipschitz estimates
$$
\|\Phi_{p,n}(\mu)-\Phi_{p,n}(\nu)\|_{\tiny\rm tv}\leq \frac{\|G_{p,n}\|}{\mu(G_{p,n})\vee\nu(G_{p,n})}~\beta(P_{p,n})~\|\mu-\nu\|_{\tiny\rm tv}
$$
and
$$
\sup_{\mu,\nu}{\|\Phi_{p,n}(\mu)-\Phi_{p,n}(\nu)\|_{\tiny\rm tv}}=\beta(P_{p,n})
$$
\end{prop}

\subsection{Stability properties}\label{quantitative-contraction}

\subsubsection{Regularity conditions}\label{regularity-conditions}

In this section we present one of the simplest quantitative contraction estimate we known for the normalized
Feynman-Kac semigroups $\Phi_{p,n}$. We consider the following regularity conditions.
\begin{description}

\item ${\bf H_m(G,M)}$ There exists some integer $m\geq 1$, such that for any $n\geq 0$, and any
 $((x,x^{\prime}),A)\in \left(E_{n}^2\times\Ea_n\right)$ and any $n\geq 0$ we have
$$
 M_{n,n+m}(x,A)\leq \chi_m ~M_{n,n+m}(x^{\prime},A)\quad \mbox{and}\quad g=\sup_{n\geq 0}{g_n}<\infty ~
$$
for some finite
parameters $\chi_m,g<\infty$, and some integer $m\geq 1$.
\item ${\bf H_0(G,M)}$ 
\begin{equation}\label{GpnPpnMixbeta}
 \rho:=\sup_{n\geq 0}{\left(g_{n}\beta(M_{n+1})\right)}<1\quad\mbox{and}\quad g=\sup_ng_n<\infty
\end{equation}

\end{description}

Both conditions are related to the stability properties of  the reference Markov chain model 
$X_n$ with probability transition $M_n$. They implies that the chain $X_n$  tends to merges exponentially
fast the random states starting from {\em any} two different locations.

One natural strategy to obtain some useful quantitative contraction estimate
for the Markov transitions $P_{p,n}$ is to write this transition in terms of the composition of Markov
transitions.  

\begin{lemma}\label{lemmPRbeta}
For any $0\leq p\leq q\leq n$, we have
$$
P_{p,n}=R^{(n)}_{p,q}P_{q,n}\quad\mbox{and}\quad P_{p,n}=R^{(n)}_{p+1}R^{(n)}_{p+2}\ldots R^{(n)}_{n-1}R^{(n)}_{n}
$$
with the  triangular array of  Markov transitions $R^{(n)}_{p,q}$ and $(R^{(n)}_{q})_{1\leq q\leq n}$ defined by
$$
R^{(n)}_{p,q}(f):=\frac{Q_{p,q}\left(G_{q,n}f\right)}{Q_{p,q}\left(G_{q,n}\right)}=
\frac{P_{p,q}\left(G_{q,n}f\right)}{P_{p,q}\left(G_{q,n}\right)}
$$
and
$$
R^{(n)}_{p}(f)=\frac{Q_{p}(  G_{p,n} f)}{Q_{p}(  G_{p,n} )}=\frac{M_{p}(  G_{p,n} f)}{M_{p}(  G_{p,n} )}
$$
In addition, for any $0\leq p\leq q\leq n$ we have 
\begin{equation}\label{rpn-gpn}
\beta\left(R^{(n)}_{p,q}\right)\leq g_{q,n}~\beta\left(P_{p,q}\right)
~\mbox{
and}
~~
\log{g_{p,n}}\leq 
\sum_{p\leq q<n} \left(g_q-1\right)~\beta(P_{p,q})
\end{equation}

\end{lemma}
\preuve
Using the decomposition
$$
Q_{p,n}(f)=Q_{p,q}(Q_{q,n}(f))=Q_{p,q}\left(Q_{q,n}(1)~P_{q,n}(f)\right)
$$
we easily check the first assertion. The l.h.s. inequality in (\ref{rpn-gpn}) is a direct consequence of (\ref{trickch4}).
Using (\ref{multiplicative}), the proof of the r.h.s.  inequality in (\ref{rpn-gpn}) is based on the fact that
\begin{eqnarray*}
\frac{
G_{p,n}(x)}{G_{p,n}(y)}&=&
\frac{
\prod_{p\leq q<n} \Phi_{p,q}(\delta_x)(G_q)
}{
\prod_{p\leq q<n} \Phi_{p,q}(\delta_y)(G_q)}\\
&=&
\exp{\left\{
\sum_{p\leq q<n}\left(\log{\Phi_{p,q}(\delta_x)(G_q)}-\log{\Phi_{p,q}(\delta_y)(G_q)}\right)
\right\}}
\end{eqnarray*}
Using the fact that
$$
\log{y}-\log{x}=\int_0^1 \frac{(y-x)}{x+t(y-x)}~dt
$$
for any positive numbers $x,y$, we prove that
$$
\begin{array}{l}
\displaystyle\frac{
G_{p,n}(x)}{G_{p,n}(y)}
\\
\\=
\exp{\left\{
\sum_{p\leq q<n}
\displaystyle\int_0^1 \frac{(\Phi_{p,q}(\delta_x)(G_q)-\Phi_{p,q}(\delta_y)(G_q))}{\Phi_{p,q}(\delta_y)(G_q)+t(\Phi_{p,q}(\delta_x)(G_q)-\Phi_{p,q}(\delta_y)(G_q))}~dt
\right\}}\\
\\
\leq \displaystyle 
\exp{\left\{
\sum_{p\leq q<n}~\displaystyle\widetilde{g}_q\times\left(
\Phi_{p,q}(\delta_x)(\widetilde{G}_q)-\Phi_{p,q}(\delta_y)(\widetilde{G}_q)\right)
\right\}}
\end{array}
$$
with $$\widetilde{G}_q:=G_q/\mbox{\rm osc}(G_q)\quad \mbox{\rm and}\quad 
\widetilde{g}_q:={\mbox{\rm osc}(G_q)}/{\inf G_q}\leq g_q-1$$ We end the proof of the desired estimates
using (\ref{trickch4}), and proposition~\ref{propcontract-de}. This completes the proof of the lemma.
\cqfd

\subsubsection{Quantitative contraction theorems}\label{quantitative-contract-sec}

This section is mainly concerned with the proof of two contraction theorems that can be derived under the
couple of regularity conditions presented in section~\ref{regularity-conditions}.

\begin{theorem}\label{theocontractref}
We assume that condition
$
{\bf H_m(G,M)}
$
is satisfied for some finite
parameters $\chi_m,g<\infty$, and some integer $m\geq 1$. In this situation, we have the uniform
estimates
\begin{equation}\label{contract-ref}
\sup_{0\leq p\leq n}{g_{p,n}}\leq \chi_mg^m \quad\mbox{\rm and}\quad \sup_{p\geq 0}{\beta(P_{p,p+km})}\leq \left( 1-g^{-(m-1)} \chi_m^{-2} \right)^k
\end{equation}
In addition, for any couple of measures $\nu,\mu\in\Pa(E_p)$, and
for any $f\in \mbox{\rm Osc}(E_n)$ we have the decomposition 
\begin{equation}\label{weak-contract-ref}
\left|[\Phi_{p,n}(\mu)-\Phi_{p,n}(\nu)](f)\right|
\leq \rho_m  \left( 1-\kappa_m\right)^{(n-p)/m}~
~\left|(\mu-\nu)D_{p,n,\mu}(f)\right|
\end{equation}
for some function $D_{p,n,\mu}(f)\in \mbox{\rm Osc}(E_p)$ whose values only depends on the parameters $(p,n,\mu)$, and some
parameters $\rho_m<\infty$ and $\kappa_m\in ]0,1]$ such that
\begin{equation}\label{rhokappa}
\rho_m\leq \chi_mg^m  \left( 1-g^{-(m-1)} \chi_m^{-2} \right)^{-1}\quad
\mbox{and}
\quad
\kappa_m\geq g^{-(m-1)} \chi_m^{-2} 
\end{equation}
\end{theorem}
\preuve
For any non negative function $f$, we notice that
\begin{eqnarray*}
R^{(n)}_{p,p+m}(f)(x)&=&\frac{Q_{p,p+m}(  G_{p+m,n} f)(x)}{Q_{p,p+m}(  G_{p+m,n} )(x)}\\
&\geq & g^{-(m-1)} \chi_m^{-2} ~
\frac{M_{p,p+m}(  G_{p+m,n} f)(x^{\prime})}{M_{p,p+m}(  G_{p+m,n} )(x^{\prime})}
\end{eqnarray*}
and for any $p+m\leq n$
$$
\frac{G_{p,n}(x)}{G_{p,n}(x^{\prime})}=\frac{Q_{p,p+m}(G_{p+m,n})(x)}{Q_{p,p+m}(G_{p+m,n})(x^{\prime})}
\leq g^m\frac{
M_{p,p+m}(G_{p+m,n})(x)}{M_{p,p+m}(G_{p+m,n})(x^{\prime})}\leq \chi_mg^m 
$$
For $p\leq n\leq p+m$, this upper bound remains valid. We conclude that 
$$
G_{p,n}(x)\leq \chi_mg^m ~G_{p,n}(x^{\prime})\quad \mbox{\rm and}\quad
\beta\left(R^{(n)}_{p,p+m}\right)\leq 1-g^{-(m-1)} \chi_m^{-2} 
$$
In the same way as above, we have
$$
n=km\Rightarrow
P_{p,p+km}=R^{(n)}_{p,p+m}R^{(n)}_{p+m,p+2m}\ldots R^{(n)}_{p+(k-1)m,p+km}
$$
and
$$
\beta(P_{p,p+km})\leq \prod_{1\leq l\leq k}\beta(R^{(n)}_{p+(l-1)m,p+lm})\leq \left( 1-g^{-(m-1)} \chi_m^{-2} \right)^k
$$
This ends the proof of (\ref{contract-ref}).

The proof of (\ref{weak-contract-ref}) is based on  the decomposition
$$
(\mu-\nu)S_{G_{p,n},\mu}P_{p,n}(f)=\beta(S_{G_{p,n},\mu}P_{p,n})\times (\mu-\nu)D_{p,n,\mu}(f)
$$
with
$$
D_{p,n,\mu}(f):=S_{G_{p,n},\mu}P_{p,n}(f)/\beta(S_{G_{p,n},\mu}P_{p,n})
$$
On the other hand, we have
$$
S_{G_{p,n},\mu}(x,y)\geq (1-\|G_{p,n}\|)\Rightarrow \beta(S_{G_{p,n},\mu})\leq \|G_{p,n}\|\
$$
and
$$
\beta(S_{G_{p,n},\mu})/\nu(G_{p,n})
\leq g_{p,n}
\leq \chi_mg^m
$$
Finally, we observe that
$$
 \beta(P_{p,n})\leq\beta(P_{p,p+\lfloor(n-p)/m\rfloor})
$$
from which we conclude that
$$
\frac{1}{\nu(G_{p,n})}~\beta(S_{G_{p,n},\mu}P_{p,n})\leq 
 \chi_mg^m
\beta(P_{p,p+\lfloor(n-p)/m\rfloor})
$$
The end of the proof is now a direct consequence of the contraction estimate (\ref{contract-ref}). This ends the proof of the theorem.
\cqfd

\begin{theorem}\label{theoHGM}
We assume that condition ${\bf H_0(G,M)}$  is satisfied for some $\rho<1$.
In this situation, for 
any couple of measures $\nu,\mu\in\Pa(E_p)$, and
for any $f\in \mbox{\rm Osc}(E_n)$ we have the decomposition 
$$
\left|[\Phi_{p,n}(\mu)-\Phi_{p,n}(\nu)](f)\right|
\leq \rho^{(n-p)}~
~\left|(\mu-\nu)D_{p,n,\mu}(f)\right|
$$
for some function $D_{p,n,\mu}(f)\in \mbox{\rm Osc}(E_p)$, whose values only depends on the parameters $(p,n,\mu)$. In addition, for any $0\leq p\leq n$, we have the estimates
$$
\beta(P_{p,n})\leq \rho^{n-p}
\quad
\mbox{and}
\quad g_{p,n}\leq \exp{\left((g-1)~{(1-\rho^{n-p})}/{(1-\rho)}\right)}
$$
\end{theorem}
\proof
Using proposition~\ref{propcontract-de}, and recalling that
$\beta(S_{G_{n-1},\mu})\leq \|G_{n-1}\|$, we readily prove that
\begin{eqnarray*}
\left|\left[\Phi_{n}(\mu)-\Phi_{n}(\nu)\right](f)\right|&\leq& g_{n-1}\beta(M_n)~\left|(\mu-\nu)D_{n,\mu}(f)\right|\\
&\leq& \rho~\left|(\mu-\nu)D_{n,\mu}(f)\right|
\end{eqnarray*}
with  the function
$$
D_{n,\mu}(f)=S_{G_{n-1},\mu}M_{n}(f)/\beta(S_{G_{n-1},\mu}M_n)\in \mbox{\rm Osc($E_{n-1}$)}
$$
Now, we can prove the theorem by induction on the parameter $n\geq p$. For $n=p$, the desired
result follows from the above discussion. Suppose we have
$$
\left|[\Phi_{p,n-1}(\mu)-\Phi_{p,n-1}(\nu)](f)\right|
\leq \rho^{(n-p-1)}~
~\left|(\mu-\nu)D_{p,n-1,\mu}(f)\right|
$$
for any $f\in \mbox{\rm Osc}(E_{n-1})$, and some functions $D_{p,n-1,\mu}(f)\in \mbox{\rm Osc}(E_{p})$. In this case, we have
$$
\begin{array}{l}
\left|[\Phi_n\left(\Phi_{p,n-1}(\mu)\right)-\Phi_n\left(\Phi_{p,n-1}(\nu)\right)](f)\right|
\\
\\
\leq g_{n-1}\beta(M_n)~\left|(\Phi_{p,n-1}(\mu)-\Phi_{p,n-1}(\nu))D_{n,\Phi_{p,n-1}(\mu)}(f)\right|
\end{array}$$
for any $f\in \mbox{\rm Osc}(E_{n})$, with $D_{n,\Phi_{p,n-1}(\mu)}(f)\in \mbox{\rm Osc}(E_{n-1})$.

Under our assumptions, we conclude that
$$
\begin{array}{l}
\left|[\Phi_n\left(\Phi_{p,n-1}(\mu)\right)-\Phi_n\left(\Phi_{p,n-1}(\nu)\right)](f)\right|
\leq \rho^{(n-p)}~
~\left|(\mu-\nu)D_{p,n,\mu}(f)\right|
\end{array}$$
with the function $$D_{p,n,\mu}(f):=D_{p,n-1,\mu}\left(D_{n,\Phi_{p,n-1}(\mu)}(f)\right)\in \mbox{\rm Osc}(E_p)$$
The proof of the second assertion is a direct consequence of proposition~\ref{propcontract-de}, and lemma~\ref{lemmPRbeta}.
This ends the proof of the theorem.
\cqfd

\begin{cor}\label{control-tau-kappa}
Under any of the conditions ${\bf H_m(G,M)}$, with $m\geq 0$, the functions $\tau_{k,l}$
and $\kappa$
defined in (\ref{def-tau}) are uniformly bounded; that is, for any $k,l\geq 1$ we have that
$$
\overline{\tau}_{k,l}(m):=\sup_{n\geq 0}\sup_{0\leq p\leq n}{\tau_{k,l}(n)}<\infty\quad\mbox{\rm and}
\quad \overline{\kappa}(m):=\sup_{n\geq 0}\kappa(n)<\infty
$$
 In addition, for any $m\geq 1$, we have
\begin{eqnarray*}
 \overline{\kappa}(m)&\in &\left[1,\chi_m g^m\right]\\
\overline{\tau}_{k,l}(m)&\leq& m~
 (\chi_mg^{m})^k~/\left(1- \left( 1-g^{-(m-1)} \chi_m^{-2} \right)^{l}\right)
\end{eqnarray*}
and for $m=0$, we have the estimates
\begin{eqnarray*}
 \overline{\kappa}(0)&\leq& \exp{\left((g-1)/(1-\rho)\right)}\\
\overline{\tau}_{k,l}(0)&\leq&  \exp{\left(k(g-1)/(1-\rho)\right)}/(1-\rho^l)
\end{eqnarray*}

\end{cor}
\subsubsection{Some illustrations}

We illustrate the regularity conditions presented in section~\ref{regularity-conditions} 
with  three different types of Feynman-Kac models, related respectively to time discretization techniques,
simulated annealing type schemes, and path space models.
 
Of course, a complete analysis of the regularity properties of the 20 
Feynman-Kac application models presented in section~\ref{FK-applications-sec} 
would lead to a too long discussion.  

In some instances, the regularity conditions stated in section~\ref{regularity-conditions}
can be directly translated into regularity properties of the reference Markov chain model and the adaptation
potential function.

 In other instances, the regularity properties of the Feynman-Kac semigroup depend
on some important tuning parameters, including discretization time steps, and cooling schedule in simulated annealing time
models. In section~\ref{Mh-Gh-ref-sec}, we illustrate the regularity property ${\bf H_0(G,M)}$ 
stated in (\ref{GpnPpnMixbeta}) in the context of time discretization model with geometric style clocks introduced in (\ref{h-ref-gm}). In section~\ref{sec-appli-saa}, we present some tools to tune the cooling parameters of the 
annealing model discussed in (\ref{i-MHM}), so that the resulting semigroups are exponentially stable.

For degenerate indicator style functions, we can use a one step integration technique to transform the
model into a Feynman-Kac model on smaller state spaces with positive potential functions. In terms of particle absorption
models, this technique allows to turn a hard obstacle model into a soft obstacle particle model. Further details on this integration technique can be found in~\cite{Delm04,soft}.

Last, but not least, in some important applications, including Feynman-Kac models on path spaces,
the limiting  semigroups are unstable, in the sense that they don't forget their initial conditions.  Nevertheless, in some situations
it is still possible to control uniformly in time the quantities $g_{p,n}$ discussed in section~\ref{sec-fk-sg}

\subsubsection{Time discretization models}\label{Mh-Gh-ref-sec}

We consider the potential functions $G_n$ and Markov transitions 
$M_n$ are given by (\ref{h-ref-gm}),
for some non negative function $V_n$, some positive parameter $\lambda_n$
and some Markov transition
$K_n$ s.t. 
$$
\beta(K_n)\leq \kappa_n<1\quad
h\leq  h_n=(1-\kappa_n)/[v_{n-1}+\alpha]\quad\mbox{\rm and}\quad
\lambda_n\in ]0,1/h]
$$ 
with $v_{n-1}:=\mbox{\rm osc}(V_{n-1})$, and 
for some $\alpha>0$. We also assume that $v=\sup_n v_n<\infty$.

In this situation, for any
$\lambda_n\in \left[\frac{1}{h_n},\frac{1}{h}\right]$,
 we have
\begin{eqnarray*}
g_{n-1}\beta(M_n)&\leq& e^{v_{n-1}h}
\left(1-\lambda_n h~(1-\kappa_n)
\right)\\
&\leq& e^{-h\left(\lambda_n(1-\kappa_n)-v_{n-1}\right)}\leq e^{-\alpha h}
\end{eqnarray*}
from which we conclude that
${\bf H_0(G,M)}$ is met with 
$$
g=\sup_n g_n\leq e^{hv}\quad\mbox{\rm and}\quad \rho\leq e^{-\alpha h}
$$
\subsubsection{Interacting simulated annealing model}\label{sec-appli-saa}
We consider the Feynman-Kac annealing model discussed in (\ref{i-MHM}).
We further assume that
 $$K^{k_n}_n(x,dy)\geq \epsilon_n~\nu_n(y)$$ for some
$k_n\geq 1$, some $\epsilon_n>0$, and some measure $\nu_n$. 

In this situation, we have
$$
\Ma_{n,\beta_{n}}^{k_n}(x,dy)\geq K^{k_n}_n(x,dy)~e^{-\beta_n k_n v}\geq \epsilon_n~e^{-\beta_n k_n v}~\nu_n(dy)
$$
with $v:=\mbox{\rm osc}(V)$.  if we choose $m_n=k_nl_n$, this implies that
$$
\beta(M_n)=
\beta\left(\Ma_{n,\beta_{n}}^{m_n}\right)\leq\beta\left(\Ma_{n,\beta_{n}}^{k_n}\right)^{l_n}\leq  \left(1-\epsilon_n~e^{-\beta_n k_n v}\right)^{l_n}
$$
Therefore, for any given $\rho^{\prime}\in ]0,1[$ we can chose $l_n$ such that
$$
l_n\geq \frac{\log{(1/\rho^{\prime})}+v(\beta_n-\beta_{n-1})}{
\log{1/(1-\epsilon_n~e^{-\beta_n k_n v})}
}
$$
so that
$$
g_{n-1}\beta(M_n)\leq e^{v(\beta_n-\beta_{n-1})}\times  \left(1-\epsilon_n~e^{-\beta_n k_n v}\right)^{l_n}\leq \rho^{\prime}
\Rightarrow \rho\leq \rho^{\prime}$$
For any function $\beta~:~x\in[0,\infty[\mapsto \beta(x)$, with a decreasing derivative  $\beta^{\prime}(x)$ s.t.
$\lim_{x\rightarrow\infty}\beta^{\prime}(x)=0$ and $\beta^{\prime}(0)<\infty$,
we also notice that 
$$
 g=\sup_{n\geq 0} g_n\leq \sup_{n\geq 0}e^{v\beta^{\prime}_{n}}\leq e^{v\beta^{\prime}(0)}
$$
\subsubsection{Historical processes}\label{historical-sg-control}

We return to the historical Feynman-Kac models introduced in section~\ref{historical-sec}.
Using the equivalence principle (\ref{equivalence-principle}), we have proved that
the $n$-time marginal models associated with a Feynman-Kac model  on path space
coincide with the original Feynman-Kac measure  (\ref{defQn}). 

We write ${\bf Q}_{p,n}$ and ${\bf P}_{p,n}$ the Feynman-Kac semigrousp defined as 
 $Q_{p,n}$ and $P_{p,n}$, by replacing  $(G_n,M_n)$ by $({\bf G}_n,{\bf M}_n)$. By construction, we have
  $$
 {\bf Q}_{p,n}(\un)({\bf x}_p)=Q_{p,n}(\un)(x_p)
 $$
 for any ${\bf x}_p=(x_0,\ldots,x_p)\in E_p$. Therefore, if we set 
 $${\bf G}_{p,n}({\bf x}_p):= {\bf Q}_{p,n}(\un)({\bf x}_p)
 $$
 then we find that
 $$
 {\bf g}_{p,n}:=\sup_{{\bf x}_p,{\bf y}_p}\frac{{\bf G}_{p,n}({\bf x}_p)}{{\bf G}_{p,n}({\bf y}_p)}=
 \sup_{{x}_p,{y}_p}\frac{{G}_{p,n}({x}_p)}{{G}_{p,n}({y}_p)}=g_{p,n}
 $$
 On the other hand, we cannot expect the Dobrushin's ergodic coefficient of the historical process semigroup to decrease but we always have
 have $\beta({\bf P}_{p,n})\leq 1$. 
 
In summary,  when the reference Markov chain $X_n$ satisfies  the condition ${\bf H_m(G,M)}$ stated in
the beginning of section~\ref{regularity-conditions},
 for some $m\geq 1$, we always have the estimates
 \begin{equation}\label{gbeta-historical}
  {\bf g}_{p,n}\leq \chi_m g^m\quad\mbox{\rm and}\quad \beta({\bf P}_{p,n})\leq 1
 \end{equation}
 and 
  \begin{equation}\label{gbeta-tau-historical}
 \tau_{k,l}(n)\leq (n+1) (\chi_m g^m)^k\quad\mbox{\rm and}\quad \kappa(n)\leq \chi_m g^m
 \end{equation}
with the functions $ \tau_{k,l}$, and $\kappa$ introduced in (\ref{def-tau}) 

We end this section with some Markov chain Monte Carlo technique often used in practice to stabilize the
genealogical tree based approximation model. To describe with some precision this stochastic method,
we consider the Feynman-Kac measures $\eta_n\in\Pa({\bf E_n})$ associated with the potential
function ${\bf G}_n$ and the Markov transitions ${\bf M}_n$ of the historical process defined respectively in (\ref{eqGb-ref})
and in (\ref{eqMb-ref}). We notice that $\eta_n$ satisfy the updating-prediction equation
$$
\eta_n=\Psi_{{\bf G}_n}(\eta_n){\bf M}_n
$$
This equation on the set of measures on path spaces is unstable, in the sense that its initial condition
is always kept in memory by the historical Markov transitions ${\bf M_n}$. One idea to stabilize this system is to
incorporate an additional Markov chain Monte Carlo move at every time step. More formally, let us suppose that we have
a dedicated Markov chain Monte Carlo transition $K_n$ from the set ${\bf E_n}$ into itself, and
such that
$$
\eta_n=\eta_nK_n
$$
In this situation, we also have that
\begin{equation}\label{genea-regularized}
\eta_n=\Psi_{{\bf G}_n}(\eta_n){\bf M}^{\prime}_n\quad \mbox{\rm with}\quad {\bf M^{\prime}_n}:={\bf M}_nK_n
\end{equation}

By construction, the mean field particle approximation of the equation (\ref{genea-regularized}) is a genealogical tree type evolution
model with path space particles on the state spaces  ${\bf E_n}$.
The updating-selection transitions are related to the potential function ${\bf G_n}$ on the state spaces ${\bf E_n}$,
and the mutation-exploration mechanisms from ${\bf E_n}$ into ${\bf E_{n+1}}$ are dictated by the Markov transitions
${\bf M^{\prime}_{n+1}}$.  

Notice that this mutation transition is decomposed into two different stages. Firstly, we extend the selected path-valued particles with an elementary move according to the Markov transition $M_n$. Then, from every of these extended paths,  
we perform 
a Markov chain Monte Carlo sample according to the Markov transition $K_n$.
\subsection{Mean field particle models}\label{mean-field-sec}
\subsubsection{Interacting particle systems}

With the exception of some very special cases, the  measures $\eta_n$
 cannot be represented in a closed form, even on in finite dimensional state-spaces.
Their numerical estimation using deterministic type grid  approximations 
requires extensive calculations, and their rarely cope with high
dimensional problems. In the same vein, harmonic type approximation schemes, or related
linearization style techniques such as the extended Kalman filter often provide poor estimations
result for highly nonlinear models. In contrast with these conventional techniques, mean field particle models can be thought 
as a stochastic adaptive grid approximation scheme. These advanced Monte Carlo methods take advantage of
 the nonlinearities 
of the model, so that to design an interacting selection-recycling mechanism.

Formally speaking, discrete generation mean field  particle
models are based on the fact that the flow of probability measures
$\eta_n$ satisfy a
non linear evolution equation of the following form 
\begin{equation}
\eta _{n+1}(dy)=\int \eta _{n}(dx)K_{n+1,\eta _{n}}(x,dy)  \label{nonlin1}
\end{equation}%
for some collection of Markov transitions $K_{n+1,\eta }$, indexed by the
time parameter $n\geq 0$ and the set of probability measures $\mathcal{P}%
(E_{n})$.

The choice of the McKean transitions $K_{n+1,\eta _{n}}$ is not unique. For instance, we can choose
$$
K_{n+1,\eta _{n}}(x,dy) =\Phi_{n+1}(\eta_n)(dy)
$$
and more generally
$$
\begin{array}{l}
K_{n+1,\eta _{n}}(x,dy) \\
\\=\epsilon(\eta_n) G_n(x)~M_{n+1}(x,dy)+(1-\epsilon(\eta_n) G_n(x))~\Phi_{n+1}(\eta_n)(dy)
\end{array}$$ 
for any $\epsilon(\eta_n)$ s.t. $\epsilon(\eta_n) G_n(x)\in [0,1]$.
Note that we can define sequentially a Markov chain sequence $(\overline{X}_n)_{n\geq 0}$ such that
$$
\PP\left(\overline{X}_{n+1}\in dx~|~\overline{X}_n\right)=K_{n+1,\eta _{n}}\left(\overline{X}_n,dx\right)\quad\mbox{\rm with}\quad
\mbox{\rm Law}(\overline{X}_n)=\eta_n
$$ 
From the practical point of view, this Markov chain can be seen as a perfect sampler of the 
flow of the distributions (\ref{nonlin1}) of the random states $\overline{X}_n$. For a more thorough discussion on these
nonlinear Markov chain models, we refer the reader to section 2.5 in the book~\cite{Delm04}.  

The mean field particle interpretation of this nonlinear measure
valued model is the $E_{n}^{N}$-valued Markov chain 
\begin{equation*}
\xi _{n}=\left( \xi _{n}^{1},\xi _{n}^{2},\ldots ,\xi _{n}^{N}\right) \in
E_{n}^{N}
\end{equation*}%
with elementary transitions defined as 
\begin{equation}
\mathbb{P}\left( \xi _{n+1}\in dx~|~\xi _{n}\right)
=\prod_{i=1}^{N}~K_{n+1,\eta _{n}^{N}}(\xi _{n}^{i},dx^{i}) 
\label{meanfield}
\end{equation}
with
$$
\eta _{n}^{N}:=\frac{1}{N}\sum_{j=1}^{N}~\delta _{\xi _{n}^{j}}
$$
In the above displayed formula, $dx$ stands for an infinitesimal
neighborhood of the point $x=(x^{1},\ldots ,x^{N})\in E_{n+1}^{N}$. The
initial system $\xi _{0}$ consists of $N$ independent and identically
distributed random variables with common law $\eta _{0}$.

We let $\Ga_{n}^{N}:=\sigma \left( \xi _{0},\ldots ,\xi _{n}\right) $ be the natural
filtration associated with the $N$-particle approximation model defined
above. 

The particle model associated with the parameter $\epsilon(\eta_n) =1$ coincides with the genetic type stochastic
algorithm presented in section~\ref{sec-ips-fk-intro}.

Furthermore, using the equivalence principles (\ref{equivalence-principle}) presented in section~\ref{historical-sec},
we can check that the genealogical tree model discussed above coincides with the
mean field $N$-particle interpretation of the Feynman-Kac measures $(\gamma_n,\eta_n)$  associated with the pair 
 $({\bf G}_n,{\bf M}_n)$ on the path spaces $ {\bf E_n}$. In this context, we recall that
 $\eta_n=\QQ_n$, and
  the $N$-particle approximation measures are given by
\begin{equation}\label{empirique-arbre}
 \eta^N_n:=\frac{1}{N}\sum_{i=1}^N\delta_{\left(\xi^{i}_{0,n},\xi^{i}_{1,n},\ldots,\xi^{i}_{n,n}\right)}\in \Pa({\bf E_n})=\Pa(E_0\times\ldots\times E_n)
\end{equation}

\subsection{Local sampling errors}\label{loc-sampling-sec}

The local sampling errors induced by the mean field particle model (\ref{meanfield}) are
expressed in terms of the empirical random field sequence $V^N_n$
defined by  
$$
V_{n+1}^{N}=\sqrt{N}~\left[\eta_{n+1}^{N}-\Phi_{n+1}\left(  \eta_{n}^{N}\right)  \right]
$$
Notice that $V_{n+1}^{N}$ is alternatively defined by the following stochastic perturbation
formulae
\begin{equation}\label{def-VNn}
\eta_{n+1}^{N}=\Phi_{n+1}\left(  \eta_{n}^{N}\right)  +\frac{1}{\sqrt{N}%
}~V_{n+1}^{N}%
\end{equation}
For $n=0$, we also set
$$
V_{0}^{N}=\sqrt{N}~\left[\eta_{0}^{N}-\eta_0  \right]\Leftrightarrow
\eta_{0}^{N}=\eta_0+\frac{1}{\sqrt{N}}~V_{0}^{N}
$$

In this interpretation, the $N$-particle model can also
be interpreted as a stochastic perturbation of the limiting system
$$
\eta_{n+1}=\Phi_{n+1}\left(  \eta_{n}\right)
$$
It is rather elementary to check that
\begin{eqnarray*}
\label{unbiasloc}\mathbb{E}\left( V^{N}_{n+1}(f)\left| ~\Ga^N_{n}\right. \right)
&=&0\\
\mathbb{E}\left(V^{N}_{n+1}(f)^2\left| ~\Ga^N_{n}\right.
\right)   & =& \eta^N_n\left[K_{n+1,\eta^N_n}\left(f-K_{n+1,\eta^N_n}(f)\right)^2\right]
\end{eqnarray*}

\begin{definition}
We denote by $\sigma^2_n$ the uniform local variance parameter given by
\begin{equation}\label{defsigma-ref}
\sigma^{2}_n:=\sup{\mu\left(
K_{n,\mu}\left[ f_{n}-K_{n,\mu}(f_{n})\right]^{2}\right) }\leq 1
\end{equation}
In the above displayed formula the supremum is taken over all  functions
$f_n\in\mbox{\rm Osc}(E_n)$, and all  probability measures $\mu$ on $E_n$, with $n\geq 1$. For $n=0$, we set
$$
\sigma^{2}_0=\sup_{f_0\in\mbox{\rm Osc}(E_0)}\eta_0\left([f_0-\eta_0(f_0)]^2\right)\leq 1
$$
\end{definition}

We close this section with a brief discussion on these uniform local variance parameters in the context
of continuous time discretization models.
When the discrete time model $K_{n,\mu}=K_{n,\mu}^{(h)}$ comes from a discretization
of the continuous time model with time step $\Delta t=h(\leq 1)$, we often have that
\begin{equation}\label{discret-continu}
K_{n,\mu}=Id+ h L_{n,\mu}+\mbox{\rm O}(h^2)
\end{equation}
for some  infinitesimal generator $L_{n,\mu}$.
In this situation, we also have that
$$
L_{K_{n,\mu}}:=K_{n,\mu}-Id=h L_{n,\mu}+\mbox{\rm O}(h^2)
$$
For any Markov transition $K$, we notice that
\begin{eqnarray*}
K([f-K(f)]^2)&=&K(f^2)-K(f)^2\\
&=&L_K(f^2)-2fL_K(f)-(L_K(f))^2\\
&=&\Gamma_{L_K}(f,f)-(L_K(f))^2
\end{eqnarray*}
with the carr\'e du champ\textquotedblright\ function
$\Gamma_{L_{K}}(f,f)$ defined for any $x\in E$ by
\begin{eqnarray*}
\Gamma_{L_{K}}(f,f)(x)&=&L_{K}\left(  \left[  f-L_{K}(f)(x)\right]
^{2}\right)  (x)\\&=&L_{K}(f^{2})(x)-2f(x)~L_{K}(f)(x)
\end{eqnarray*}

When $K=K_{n,\mu}$ and $L_{K_{n,\mu}}=h L_{n,\mu}+\mbox{\rm O}(h^2)$, we find that
\begin{eqnarray*}
\mu\left[K_{n,\mu}\left[ f_{n}-K_{n,\mu}(f_{n})\right]^{2}\right]&=&\mu\Gamma_{L_{K_{n,\mu}}}(f,f)~h-h^2~\mu(L_{K_{n,\mu}}(f)^2)\\
&=&h~\mu\left(\Gamma_{L_{n,\mu}}(f,f)\right)+\mbox{\rm O}(h^2)
\end{eqnarray*}

\section{Empirical processes}\label{chapter-empirique}

\subsection{Introduction}\label{intro-proc-empirique}

The aim of this chapter is to review some more or less well known stochastic  techniques for analyzing  
the concentration properties of empirical processes associated with independent random sequences. 
The discussion at the start of this section provides some basic definitions on empirical processes associated with
sequences of independent random variables on general measurable state spaces. In section~\ref{statements-sec},
we state and comment the main results of this section. 
Section~\ref{finite-marg-mod-sec} is concerned with finite marginal models. In section~\ref{sect-empirique-state}, we extend
these results at the level of the empirical processes. Besides the fact that the concentration inequalities for empirical processes
holds for supremum of empirical processes over infinite collection of functions, these inequalities are more crude with greater constants than the ones for marginal models. These two sections also contains two new perturbation theorems that apply to nonlinear functional of 
empirical processes. The proofs of these theorems combine  Orlicz's norm techniques, Kintchine's type inequalities, maximal inequalities, as well as Laplace-Cram\`er-Chernov estimation methods. These
4 complementary methodologies are presented respectively in section~\ref{sec-Orlicz-proofs}, section~\ref{marginal-sec}, section~\ref{maximal-proof-sec}, and section~\ref{Laplace-estimate-sec}.

Let $(\mu^i)_{i\geq 1}$ 
be a sequence of probability measures on a given measurable state space
$(E,\Ea)$. 
During the further development of this section, we fix an integer $N
\geq 1$. 
To clarify the presentation, we slightly abuse the  notation
and we denote respectively 
by \index{$m(X)$}
$$
m(X)=\frac{1}{N}\sum_{i=1}^N\delta_{X^i}\quad\mbox{\rm and}\quad
\mu=\frac{1}{N}\sum_{i=1}^N\mu^i
$$ the
$N$-empirical
measure associated with a collection of independent 
random variables $X=(X^i)_{i\geq 1}$, with respective distributions 
$(\mu^i)_{i\geq 1}$, and the $N$-averaged measure associated with the
sequence of measures $(\mu^i)_{i\geq 1}$. We also consider the empirical random field sequences
$$
V(X)=\sqrt{N}~\left(m(X)-\mu\right)
$$
We also set
\begin{equation}\label{def-sigma-f}
\sigma(f)^2:=\EE\left(V(X)(f)^2\right)=\frac{1}{N}\sum_{i=1}^N\mu^i([f-\mu^i(f)]^2)
\end{equation}

\begin{remark}
The rather abstract models presented above can be used to analyze the local sampling random fields
models associated with a mean field particle model discussed in section~\ref{loc-sampling-sec}.

To be more precise, given the information on the $N$-particle model at time $(n-1)$, the 
sequence of random variables $\xi^i_n$ are independent random sequences
with a distribution that depends on the current state $\xi^i_{n-1}$. That is, at any given fixed time horizon
$n$ and given  $\Ga^N_{n-1}$, 
we have
\begin{equation}\label{appli-xi-X}
X^i=\xi^i_n\in E=E_n\quad\mbox{\rm and}\quad\mu^i(dx):=K_{n,\eta^N_{n-1}}(\xi^i_{n-1},dx)
\end{equation}
In this case, we find that
$$
m(X)=\eta^N_n\quad\mbox{\rm and}\quad V(X)=V_n^N
$$
and
\begin{eqnarray*}
\sigma(f)^2&=&\mathbb{E}\left(V^{N}_{n+1}(f)^2\left| ~\Ga^N_{n}\right.
\right)\\
&=& \eta^N_n\left[K_{n+1,\eta^N_n}\left(f-K_{n+1,\eta^N_n}(f)\right)^2\right]\
\end{eqnarray*}
\end{remark}
Let $\Fa$ be a given collection  of 
measurable functions $f:E\rightarrow\RR$ such that $\|f\|\leq 1$.
We associate with $\Fa$
the Zolotarev seminorm on $\Pa(E)$
defined by\index{Zolotarev seminorm}\index{$\|\point\|_{\Fa}$}
$$
{\|\mu - \nu\|}_{\Fa}=\sup \bigl \{|\mu(f)-\nu(f)|; \, f\in \Fa \bigr \},
$$
(see for instance~\cite{rachev}). 
No generality is lost and much convenience is
gained
by supposing that the unit and the null functions
 $f=\un$ and $f=0\in\Fa$. Furthermore, to avoid some unnecessary
technical measurability questions,
we shall also suppose that $\Fa$ is separable in the sense that it contains
a countable and dense subset.

We measure
the size of a given class $\Fa$ in terms of the
covering numbers $N(\e,\Fa,\LL_2(\mu))$ defined as the minimal number of
$\LL_2(\mu)$-balls of radius $\e >0$ needed to cover $\Fa$. We shall also use the following
uniform covering numbers and entropies.

We end this section with the last of the notation to be used in this chapter
dedicated to empirical processes concentration inequalities.

\begin{definition}
By
$\Na(\e,\Fa)$, $\e >0$, and by $I(\Fa)$ we denote the uniform covering numbers and
entropy integral \index{Entropy integral}given by
\begin{eqnarray*}
\Na(\e,\Fa)\!&=&\!\sup \bigl \{\Na(\e,\Fa,\LL_2(\eta)); \,\eta\in\Pa(E)  \bigr \}
\\
I(\Fa)\!&=&\!\int_{0}^2\!\sqrt{\log{ (1+\Na(\e,\Fa))   }}\,d\e
\end{eqnarray*}
\end{definition}
The concentration inequalities stated in this section are expressed in terms of the inverse of 
the couple of functions defined below.
\begin{definition}
We let $(\epsilon_0,\epsilon_1)$ be the functions on $\RR_+$ defined by
$$
\epsilon_0(\lambda)=\frac{1}{2}\left(\lambda-\log{(1+\lambda)}\right)
 \quad\mbox{
and}\quad
\epsilon_1(\lambda)=(1+\lambda)\log{(1+\lambda)}-\lambda
$$
\end{definition}
Rather crude estimates can be derived using the following upper bounds
$$
\epsilon_0^{-1}(x)\leq 2 (x+\sqrt{x})
\quad\mbox{\rm and}\quad
\epsilon_1^{-1}(x)\leq 
\frac{x}{3}+\sqrt{2x}
$$
A proof of these elementary inequalities and refined estimates can be found in the recent article~\cite{rio-2009}.
\subsection{Statement of the main results}\label{statements-sec}

\subsubsection{Finite marginal models}\label{finite-marg-mod-sec}

The main result of this section is a quantitative concentration inequality for the finite marginal models
$$
f\mapsto V(X)(f)=\sqrt{N}\left(m(X)-\mu(f)\right)
$$

In the following theorem, we provide Kintchine's type mean error bounds, and related Orlicz norm
estimates.  The detailed proofs of these results are housed in section~\ref{marginal-sec}. The last quantitative concentration inequality is 
a direct consequence of  (\ref{simga2N}), and it is proved in remark~\ref{cor-last-assertion}.

\begin{theorem}\label{theo-kintchine}
For any integer $m\geq1$, and any measurable function $f$ we have
the $\LL_m$-mean error estimates
\begin{equation}\label{vxLp}
\EE(\left|V(X)(f)\right|^m)^{1/m}\leq b(m)~\left(\mbox{\rm osc}(f)\wedge\left[2~\mu(|f|^{m^{\prime}})^{1/m^{\prime}}\right]\right)
\end{equation}
and
$$
\mathbb{E}\left( \left| V(X)(f)\right| ^{m} \right) ^{\frac{1}{m}}\leq~ 6 b(m)^2
~
\max{\left(
\sqrt{2}\sigma(f), 
\left[\frac{
2\sigma(f)^2}{N^{\frac{m^{\prime}}{2}-1}}
\right]^{1/m^{\prime}}
\right)}
$$
with the smallest even integer $m^{\prime}\geq m$,
and the collection of
constants $b(m)$ defined in (\ref{collec}).

In particular, for any $f\in\mbox{\rm Osc(E)}$, we have
\begin{equation}\label{orlicz-vxf}
\pi_{\psi}(V(X)(f))\leq \sqrt{3/8}
\end{equation}
and for any $N$ s.t. $2\sigma^2(f)N\geq 1$ we have
\begin{equation}\label{simga2N}
 \mathbb{E}\left( \left| V(X)(f)\right| ^{m} \right) ^{\frac{1}{m}}\leq~ 6\sqrt{2}~ b(m)^2\sigma(f)
\end{equation}
In addition,  the  probability of the event
$$
\left|V(X)(f)\right|\leq 6\sqrt{2}~\sigma(f)\left[1+ \epsilon_0^{-1}(x)\right)
$$
is greater than $1-e^{-x}$, for any $x\geq 0$. 
\end{theorem}

In section~\ref{sec-concentration-orlicz} dedicated to concentration properties
of random variables $Y$ with finite Orlicz norms $\pi_{\psi}(Y)<\infty$, we shall prove that
 the probability of the event
$$
Y\leq \pi_{\psi}(Y)~\sqrt{y+\log{2}}
$$
is greater than $1-e^{-y}$, for any $y\geq 0$ (cf. lemma~\ref{lem-Orlicz-L2m}).
 This implies that
 the probability of the events
$$
\left|V(X)(f)\right|\leq \frac{1}{2}~\sqrt{3(x+\log{2})/2}
$$ 
is greater than $1-e^{-x}$, for any $x\geq 0$.

Ou next objective is to derive concentration inequalities for nonlinear functional of the empirical random
field $V(X)$. To introduce precisely these objects, we need another around of notation.

For any measure $\nu$, and any sequence of measurable functions $f=(f_1,\ldots,f_d)$, we write
$$
\nu(f):=\left[\nu(f_1),\ldots,\nu(f_d)\right]
$$

\begin{definition}
We associate with the
second order smooth function $F$ on $\RR^d$, for some $d\geq 1$,  the random
functionals defined by
\begin{equation}\label{F(m(X))}
\begin{array}{l}
f=(f_i)_{1\leq i\leq d}\in\mbox{\rm Osc}(E)^d\\
\\
\mapsto F(m(X)(f))=F(m(X)(f_1),\ldots,m(X)(f_d))\in\RR
\end{array}
\end{equation}
Given 
a probability measure
$\nu$,  and a collection of functions
$(f_i)_{1\leq i\leq d}\in\mbox{\rm Osc}(E)^d$, 
we set
\begin{equation}\label{DnuFf}
D_{\nu}(F)(f)=\nabla F(\nu(f))~f^{\top}
\end{equation}
\end{definition}
Notice that
$$
\mbox{\rm osc}\left(D_{\nu}(F)(f)\right)\leq \left\|\nabla F(\nu(f))\right\|_1:=
\sum_{i=1}^d\left|\frac{\partial F}{\partial u^i}(\nu(f))\right|
$$
We also introduce the following constants
\begin{equation}\label{DnuFf2}
\left\|\nabla^2F_f\right\|_1
:=
\sum_{i,j=1}^d~\sup_{}
{\left| \frac{\partial^2 F}{\partial u^i\partial u^j}(\nu(f))\right|}
\end{equation}
In the r.h.s. display, the supremum is taken over all probability measures $\nu\in\Pa(E)$.

The next theorem extend the exponential inequalities stated in theorem~\ref{theo-kintchine} to this class of nonlinear functionals. It also provide 
 more precise concentration properties in terms of  the variance functional 
 $\sigma$ defined in (\ref{def-sigma-f}). The proof of this theorem is postponed to section~\ref{perturbation-emp-sec}.

\begin{theorem}\label{theo-perturb-marginal}
Let $F$ be 
a second order smooth function on $\RR^d$, for some $d\geq 1$.
For  any collection of functions
$(f_i)_{1\leq i\leq d}\in\mbox{\rm Osc}(E)^d$, and any $N\geq 1$, the probability of the events
$$
\begin{array}{l}
\left[F(m(X)(f))-F(\mu(f))\right]\\
\\
\leq
\displaystyle\frac{1}{2N}~\left\|\nabla^2F_f\right\|_1~\left[3/2+\epsilon_0^{-1}(x)
\right]\\
\\\hskip2cm+\left\|\nabla F(\mu(f))\right\|_1^{-1}~\sigma^2(D_{\mu}(F)(f)) ~\epsilon_1^{-1}
\left(\displaystyle\frac{x \left\|\nabla F(\mu(f))\right\|_1^2}{N\sigma^2(D_{\mu}(F)(f))}\right)\end{array}
$$
is greater than $1-e^{-x}$, for any $x\geq 0$. In the above display, $D_{\mu}(F)(f)$ stands for the first order function
defined in  (\ref{DnuFf}). 
\end{theorem}

\subsubsection{Empirical processes}\label{sect-empirique-state}

The objective of this section is to extend the  quantitative concentration theorems, theorem~\ref{theo-kintchine}
and theorem~\ref{theo-perturb-marginal} at the level of the empirical process associated with a class of function
$\Fa$. These processes are given by the mapping
$$
f\in\Fa\mapsto V(X)(f)=\sqrt{N}\left(m(X)-\mu(f)\right)
$$
Our main result in this direction is the following theorem, whose proof is postponed to
 section~\ref{maximal-proof-sec}.
\begin{theorem}\label{theo-proc-emp}
For any class of functions $\Fa$,  with $I(\Fa)<\infty$, we have
$$
\pi_{\psi}\left(\left\|V(X)\right\|_{\Fa}\right)\leq 12^2~ \int_{0}^{2}\sqrt{\log{(8+\Na(\Fa,\epsilon)^2)}}~d\epsilon
$$
\end{theorem}
\begin{remark}
Using the fact that $\log{(8+x^2)}\leq 4\log{x}$, for any $x\geq 2$, we obtain the rather crude estimate
$$
 \int_{0}^{2}\sqrt{\log{(8+\Na(\Fa,\epsilon)^2)}}~d\epsilon\leq 2\int_{0}^{2}\sqrt{\log{\Na(\Fa,\epsilon)}}~d\epsilon
$$
We check the first observation, using the fact that $\theta(x)=4\log{x}-\log{(8+x^2)}$ is a non decreasing function
on $\RR_+$,
and $\theta(2)=\log{(4\times 4)}-\log{(4\times 3)}\geq 0$. 
\end{remark}
Various examples of classes of functions with finite covering
and entropy integral are given in the book of Van der Vaart and Wellner~\cite{wellner}
(see for instance p. 86,  p. 129, p. 135, and exercise 4 on p.150, and p. 155). The estimation
of the quantities introduced above often depends on several deep 
results on combinatorics
that are not discussed here.

To illustrate these mathematical objects,
we mention that, for the set of indicator
functions 
\begin{equation}\label{d-cells}
\Fa=\left\{1_{\prod_{i=1}^d(-\infty,x_i]}~;~(x_i)_{1\leq i\leq d}\in \RR^d\right\}
\end{equation}
of cells in $E=\RR^d$, we have
$$
\Na(\e,\Fa)\leq c ~(d+1) (4e)^{d+1}~\epsilon^{-2d}
$$
for some universal constant $c<\infty$.
This implies that
$$
\sqrt{\log{\Na(\e,\Fa)}}\leq \sqrt{\log{[c (d+1) (4e)^{d+1}]}}+ \sqrt{(2d)}~\sqrt{\log{(1/\epsilon)}}
$$
An elementary calculation gives $$\int_0^2~\sqrt{\log{(1/\epsilon)}}\leq 2~\int_0^{\infty}x^2e^{-x^2}dx=\sqrt{\pi/4}\leq 1$$ from which we conclude that
\begin{equation}\label{d-cells-analysis}
\int_{0}^{2}\sqrt{\log{\Na(\Fa,\epsilon)}}~d\epsilon\leq 2\sqrt{\log{[c (d+1) (4e)^{d+1}]}}+\sqrt{(2d)}\leq c^{\prime}\sqrt{d}
\end{equation}
for some universal constant $c<\infty$.

For $d=1$, we also have that $\Na(\e,\Fa)\leq 2/\epsilon^{2}$ (cf. p. 129 in~\cite{wellner}) and therefore
$$
\int_{0}^{2}\sqrt{\log{\Na(\Fa,\epsilon)}}~d\epsilon\leq 3\sqrt{2}
$$
\begin{remark}
In this chapter, we have assumed that the 
class of functions $\Fa$ is such that $\sup_{f\in \Fa}{\|f\|}\leq 1$. When $\sup_{f\in \Fa}{\|f\|}\leq c_{\Fa}$, for some
finite constant $c_{\Fa}$, using theorem~\ref{theo-proc-emp}, it is also readily checked that
\begin{equation}\label{theo-F-etendue}
\pi_{\psi}\left(\left\|V(X)\right\|_{\Fa}\right)\leq 12^2~ \int_{0}^{2c_{\Fa}}\sqrt{\log{(8+\Na(\Fa,\epsilon)^2)}}~d\epsilon
\end{equation}
\end{remark}

We mention that the uniform entropy condition $I(\Fa)<\infty$ is 
required in Glivenko-Cantelli and Donsker theorems for empirical processes associated with 
non necessarily independent random sequences~\cite{dm-ledoux}.

Arguing as above, we prove that
 the probability of the events
$$
\left\|V(X)\right\|_{\Fa}\leq I_1(\Fa)~\sqrt{x+\log{2}}
$$ 
is greater than $1-e^{-x}$, for any $x\geq 0$, with some constant
$$
I_1(\Fa)\leq 12^2~ \int_{0}^{2}\sqrt{\log{(8+\Na(\Fa,\epsilon)^2)}}~d\epsilon
$$

As for marginal models~\ref{F(m(X))}, our next objective is to extend
theorem~\ref{maximal-proof-sec} to empirical processes associated with some 
classes of functions. Here, we consider the empirical processes $$f\in \Fa_i\mapsto m(X)(f)\in \RR
$$associated with $d$ classes of functions $\Fa_i$, $1\leq i\leq d$, defined in section~\ref{intro-proc-empirique}. We further assume that $\|f_i\|\vee\mbox{\rm osc}(f_i)\leq 1$, for any $f_i\in \Fa_i$, and we set
$$
\Fa:=\prod_{1\leq i\leq d}\Fa_i\quad\mbox{\rm and}\quad
\pi_{\psi}(\|V(X)\|_{\Fa}):=
\sup_{1\leq i\leq d}\pi_{\psi}(\|V(X)\|_{\Fa_i})
$$

Using theorem~\ref{theo-proc-emp}, we mention that
\begin{eqnarray*}
\pi_{\psi}(\|V(X)\|_{\Fa})&\leq&
12^2~ \int_{0}^{2}\sqrt{\log{(8+\Na(\Fa,\epsilon)^2)}}~d\epsilon
\end{eqnarray*}
with
$$
\Na(\Fa,\epsilon):=\sup_{1\leq i\leq d}\Na(\Fa_i,\epsilon)
$$

We set
$$
\left\|\nabla F_{\mu}\right\|_{\infty}
:=\sup
{\left| \frac{\partial F}{\partial u^i}(\mu(f))\right|}
\quad\mbox{\rm
and}
\quad
\left\|\nabla^2F\right\|_{\infty}=
\sup
{\left| \frac{\partial^2 F}{\partial u^i\partial u^j}(\nu(f))\right|}
$$
The supremum in the l.h.s. is taken over all  ${1\leq i\leq d}$ and all ${f\in\Fa}$; and the 
supremum in the r.h.s. is taken over all  ${1\leq i,j\leq d}$, $\nu\in\Pa(E)$, and all ${f\in\Fa}$.

We are now in position to state the final main result of this section.
The proof of the next theorem is housed in the end of section~\ref{perturbation-emp-sec}.
\begin{theorem}\label{theo-i-empirique}
Let $F$ be a second order smooth function on $\RR^d$, for some $d\geq 1$. For any  
classes of functions $\Fa_i$, $1\leq i\leq d$, and for any $x\geq 0$, the probability of the following events
$$
\begin{array}{l}
\sup_{f\in\Fa}\left|F(m(X)(f))-F(\mu(f))\right|\\
\\
\leq 
\displaystyle\frac{d}{\sqrt{N}}~ \pi_{\psi}\left(\left\|V(X)\right\|_{\Fa}\right)~
\left\|\nabla F_{\mu}\right\|_{\infty}
\left(1+2\sqrt{x}\right)\\
\\\hskip2cm+\displaystyle\frac{1}{2N}\left\|\nabla^2F\right\|_{\infty} \left(d~\pi_{\psi}(\|V(X)\|_{\Fa})\right)^2\left(1+\epsilon_0^{-1}\left(\frac{x}{2}\right)
\right)
\end{array}
$$
is greater than $1-e^{-x}$.
\end{theorem}

\subsection{A reminder on Orlicz' norms}\label{sec-Orlicz-proofs}
In this section, we have collected some important properties of Orlicz' norms.
The first section, section~\ref{comparison-sec} is concerned with rather elementary 
comparison properties.
In section~\ref{sec-concentration-orlicz}, we present a natural way to obtain
Laplace estimates, and related concentration inequalities, using simple Orlicz' norm upper
bounds.

\subsubsection{Comparison properties}\label{comparison-sec}

This short section is mainly concerned with the proof of the following 
three comparison properties.
\begin{lemma}\label{lemme-0}
For any non negative variables $(Y_1,Y_2)$ we have
$$
Y_1\leq Y_2\Longrightarrow \pi_{\psi}(Y_1)\leq \pi_{\psi}(Y_2)
$$
as well as
\begin{equation}\label{lem-Orlicz-L2m-compare}
\left(\forall m\geq 0\quad\EE\left(Y_1^{2m}\right)\leq \EE\left(Y_2^{2m}\right)\right)
\Rightarrow
\pi_{\psi}(Y_1)\leq\pi_{\psi}(Y_2)
\end{equation}
In addition, for any pair of independent random variables $(X,Y)$ on some measurable state space,
and any measurable function $f$, we have
\begin{equation}\label{maj-key}
\left(\pi_{\psi}(f(x,Y))\leq c\quad\mbox{\rm for $\PP$-a.e. x}\right)\Longrightarrow\pi_{\psi}(f(X,Y))\leq c
\end{equation}

\end{lemma}
\preuve
The first assertion is immediate, ad the second assertion comes from the fact that
$$
\EE\left(\exp{\left(\frac{Y_1}{\pi_{\psi}(Y_2)}\right)^2}-1\right)\leq 
\sum_{m\geq 1}\frac{1}{m!}~\frac{\EE(Y_2^{2m})}{\pi_{\psi}(Y_2)^{2m}}=\EE\left(\frac{Y_2}{\pi_{\psi}(Y_2)}\right)~\leq 1
$$
The last assertion comes from the fact that
$$
\EE\left(\EE\left(\psi(f(X,Y)/c)\left|X\right.\right)\right)\leq 1\Rightarrow
\pi_{\psi}(f(X,Y))\leq c
$$
This ends the proof of the lemma.
\cqfd
\subsubsection{Concentration properties}\label{sec-concentration-orlicz}

The following lemma provides a simple way to transfer a control on Orlicz' norm into  moment or Laplace estimates, which in turn can be used to derive quantitative concentration inequalities

\begin{lemma}\label{lem-Orlicz-L2m}
For any non negative random variable $Y$, and any integer $m\geq 0$, we have 
\begin{equation}\label{L2-pi-psi}
\EE\left(Y^{2m}\right)\leq m!~\pi_{\psi}(Y)^{2m}\quad\mbox{\rm and}\quad
\EE\left(Y^{2m+1}\right)\leq (m+1)!~\pi_{\psi}(Y)^{2m+1}
\end{equation}
In addition, for any $t\geq 0$ we have the Laplace estimates
$$
\EE\left(e^{tY}\right)\leq \min{\left(2~e^{\frac{1}{4}(t\pi_{\psi}(Y))^2}~,~(1+t\pi_{\psi}(Y))~e^{(t\pi_{\psi}(Y))^2}\right)}
$$
In particular, for any $x\geq 0$ the probability of the event
\begin{equation}\label{expo-event}
Y\leq \pi_{\psi}(Y)~\sqrt{x+\log{2}}
\end{equation}
is greater than $1-e^{-x}$.
\end{lemma}

\begin{remark}
For a Gaussian and centred random variable $Y$, s.t. $E(Y^2)=1$, we recall that $\pi_{\psi}(Y)=\sqrt{8/3}$. In this situation, letting $y=\sqrt{8(x+\log{2})/3}$ in (\ref{expo-event}), we find that
$$
\PP\left(|Y|\geq y\right)\leq 2~e^{-\frac{1}{2}~\frac{3}{4}~y^2}
$$
Working directly with the Laplace Gaussian function
$\EE\left(e^{tY}\right)=e^{t^2/2}$, we remove the factor $3/4$.
In this sense, we loose a factor $3/4$ using the Orlicz's concentration
property (\ref{expo-event}).

In this situation, the l.h.s. moment estimate in (\ref{L2-pi-psi}) takes the form
\begin{equation}\label{ineg-b2m-pi}
b(2m)^{2m}=\frac{(2m)!}{m!}~2^{-m}\leq m!~(8/3)^m
\end{equation}
while using  Stirling's approximation of the factorials we obtain the estimate
$$
\frac{(2m)!}{m!^2}\simeq~\sqrt{2/m}~4^m \left(\leq (8/3)^m\right)
$$

\end{remark}
\begin{remark}
Given a sequence of independent Gaussian and centred random variables $Y_i$, s.t. $E(Y_i^2)=1$, for $i\geq 1$, and any sequence of non negative numbers $a_i$, we have
$$
\pi_{\psi}\left(\sum_{i=1}^na_iY_i\right)=\sqrt{8/3}~\sqrt{\sum_{i=1}^na_i^2}:=\sqrt{8/3}~\|a\|_2~
$$
while
$$
\sum_{i=1}^na_i\pi_{\psi}\left(Y_i\right)=\sqrt{8/3}~\sum_{i=1}^na_i:=\sqrt{8/3}~\|a\|_1
$$
Notice that
$$
\|a\|_2\leq \|a\|_1\leq \sqrt{n}~\|a\|_2
$$
When the coefficients $a_i$ are almost equal, we can loose a factor $\sqrt{n}$
using the triangle inequality, instead of estimating directly with the Orlicz norm of the Gaussian
mixture. In this sense, it is always preferable to avoid the use of the triangle inequality, and to estimate directly the Orlicz 
norms of linear combinations of "almost Gaussian" random variables.
\end{remark}

Now, we come to the proof of the lemma.

{\bf Proof of lemma~\ref{lem-Orlicz-L2m}:}

For any $m\geq 1$, we have
$$
x^{2m}\leq m!~\sum_{n\geq 1}\frac{x^{2n}}{n!}=m!~\psi(x)$$
$$
\Downarrow
$$
$$
\EE\left(\left[\frac{Y}{\pi_{\psi}(Y)}\right]^{2m}\right)\leq m!~\EE\left({\psi}\left(\frac{Y}{\pi_{\psi}(Y)}\right)\right)\leq m!
$$
For odd integers, we simply use Cauchy-Schwartz inequality to check that
$$
\EE\left(Y^{2m+1}\right)^2\leq\EE\left(Y^{2m}\right)\EE\left(Y^{2(m+1)}\right)
\leq (m+1)!^2~\pi_{\psi}(Y)^{2(2m+1)}
$$
This ends the proof of the first assertion.

We use Cauchy-Schwartz's inequality
to  check that
$$
\EE\left(Y^{2m+1}\right)^2\leq \EE\left(Y^{2m}\right)~\EE\left(Y^{2(m+1)}\right)\leq (m+1)!^2~\pi_{\psi}(Y)^{2(2m+1)}
$$
for any non negative random variable $Y$,
so that
$$
\EE\left(Y^{2m+1}\right)\leq (m+1)!~\pi_{\psi}(Y)^{(2m+1)}
$$
Recalling that $(2m)!\geq m!^2$ and $(m+1)\leq (2m+1)$, we find that
\begin{eqnarray*}
\EE\left(e^{tY}\right)&=&\sum_{m\geq 0}\frac{t^{2m}}{(2m)!}~\EE\left(Y^{2m}\right)+
\sum_{m\geq 0}\frac{t^{2m+1}}{(2m+1)!}~\EE\left(Y^{2m+1}\right)\\
&\leq &\sum_{m\geq 0}\frac{t^{2m}}{m!}~\pi_{\psi}(Y)^{2m}+
\sum_{m\geq 0}\frac{t^{2m+1}}{m!}~\pi_{\psi}(Y)^{(2m+1)}\\
&=&\left(1+t\pi_{\psi}(Y)\right)~\exp{(t\pi_{\psi}(Y))^2}
\end{eqnarray*}
On the other hand, using the estimate
$$
tY=\left(\frac{t\pi_{\psi}(Y)}{\sqrt{2}}\right)~\left(\frac{\sqrt{2}~Y}{\pi_{\psi}(Y)}\right)\leq
\frac{(t\pi_{\psi}(Y))^2}{4}+\left(\frac{Y}{\pi_{\psi}(Y)}\right)^2
$$
we prove that
$$
\EE\left(e^{tY}\right)\leq 2~\exp{\left(\frac{(t\pi_{\psi}(Y))^2}{4}\right)}
$$
The end of the proof of the Laplace estimates is now completed.
To prove the last assertion, we use the fact that for any $y\geq 0$
\begin{eqnarray*}
\PP\left(Y\geq y\right)&\leq& 2~\exp{\left(-\sup_{t\geq 0}\left(ty-(t\pi_{\psi}(Y))^2/4\right)\right)}\\
&=&2~\exp{\left[-\left(y/\pi_{\psi}(Y)\right)^2\right]}
\end{eqnarray*}
This implies that
$$
\PP\left(Y\geq \pi_{\psi}(Y)\sqrt{x+\log{2}}\right)\leq 2\exp{\left[-\left(x+\log{2}\right)\right]}=e^{-x}
$$
This ends the proof of the lemma.
\cqfd

\subsubsection{Maximal inequalities}

Let us now put together the Orlicz's norm properties derived in section~\ref{sec-Orlicz-proofs}
to establish a series of more or less well known maximal inequalities.
More general results can be found in the books~\cite{pollard,wellner}, or in the lecture notes~\cite{gine}.

We emphasize that in the literature on empirical processes, maximal inequalities are often presented in terms
 of universal constant $c$ without further information on their magnitude. In the present section, we shall try to estimate some of these universal constants explicitly.

To begin with, we consider a couple of maximal inequalities over finite sets.

\begin{lemma}\label{lem-Orlicz-1}
For any finite collection of non negative random variables $(Y_i)_{i\in I}$, and any collection of non negative numbers $(a_i)_{i\in I}$, we have
$$
\sup_{i\in I}\EE(\psi(Y_i/a_i))\leq 1\Rightarrow
\EE\left(\max_{i\in I}{Y_i}\right)\leq \psi^{-1}(|I|)\times \max_{i\in I}a_i
$$
\end{lemma}
\preuve
We check this claim using the following estimates
\begin{eqnarray*}
\psi\left(\frac{\EE\left(\max_{i\in I}{Y_i}\right)}{\max_{i\in I}a_i}
\right)&\leq& \psi\left(\EE\left(\max_{i\in I}{
({Y_i}/{a_i})}
\right)
\right)\\
&\leq& \EE\left(\psi\left(\max_{i\in I}{
({Y_i}/{a_i})}
\right)\right)\\
&\leq &\EE\left(\sum_{i\in I}\psi\left({
({Y_i}/{a_i})}
\right)\right)\leq |I|
\end{eqnarray*}
This ends the proof of the lemma.
\cqfd
Working a little harder, we  prove the following lemma.
\begin{lemma}\label{lem-Orlicz-2}
For any finite collection of non negative random variables $(Y_i)_{i\in I}$, we have
$$
\pi_{\psi}\left(
\max_{i\in I}{Y_i}
\right)\leq \sqrt{6\log{(8+|I|)}}~\max_{i\in I}{\pi_{\psi}(Y_i)}
$$
\end{lemma}
\preuve
Without lost of generality, we assume that $\max_{i\in I}{\pi_{\psi}(Y_i)}\leq 1$, and
$I=\{1,\ldots,|I|\}$. In this situation, it suffices to check that
$$
\psi\left(
\frac{\max_{1\leq i\leq |I|}Y_i}{ \sqrt{6\log{(8+|I|)}}}
\right)\leq \psi\left(\max_{1\leq i\leq |I|}
\frac{Y_i}{ \sqrt{6\log{(8+i)}}}
\right)\leq 1
$$
Firstly, we notice that 
for any $i\geq 1$ and $x\geq 3/2$ we have
$$
\frac{1}{\log{(8+i)}}+\frac{1}{\log{x}}\leq
\frac{1}{\log{9}}+\frac{1}{\log{(3/2)}}\leq 3
$$
and therefore
$$
3\log{(8+i)}\log{(x)}\geq \log{(x(8+i))}
$$
We check the first estimate  using the fact that
$$
\log(3)\leq 5\log(3/2)\Rightarrow\log(3)+\log(3/2)\leq 6\log(3/2)\leq 3\log(3/2)\log(9)
$$
Using these observations, we have
$$
\begin{array}{l}
\PP\left(\max_{1\leq i\leq |I|}\left(\frac{Y_i}{\sqrt{6\log{(8+i)}}}\right)^2>\log{x}\right)\\
\\=
\PP\left(\max_{1\leq i\leq |I|}\left(\frac{Y_i}{\sqrt{6\log{(x)}\log{(8+i)}}}\right)^2>1\right)\\
\\
\leq \PP\left(\max_{1\leq i\leq |I|}\frac{Y_i}{\sqrt{2\log{(x(8+i))}}}>1\right)
\\
\\
\leq \sum_{i=1}^{|I|}\PP\left(
Y_i>\sqrt{2\log{(x(8+i))}}
\right)\leq \sum_{i=1}^{|I|}~e^{-2\log{(x(8+i))}}~\EE\left(e^{Y_i^2}\right)\\
\end{array}
$$
This implies that
$$
\begin{array}{l}
\PP\left(\max_{1\leq i\leq |I|}\left(\frac{Y_i}{\sqrt{6\log{(8+i)}}}\right)^2>\log{x}\right)\\
\\
\leq \displaystyle\frac{2}{x^2}\sum_{i=1}^{|I|}~\frac{1}{(8+i)^2}\leq \frac{2}{x^2}~\int_{8}^{\infty}\frac{1}{u^2}~du=\frac{1}{(2x)^2}
\end{array}
$$
If we set
$$
Z_I:=\exp{\left\{\left(\max_{1\leq i\leq |I|}\frac{Y_i}{\sqrt{6\log{(8+i)}}}\right)^2\right\}}
$$
then we have
\begin{eqnarray*}
\EE\left(Z_I\right)&=&\int_0^{\infty}\PP\left(Z_I>x\right)~dx\\
&\leq& \frac{3}{2}+\int_{\frac{3}{2}}^{\infty}\frac{1}{(2x)^2}~dx=\frac{3}{2}\left(1+\frac{1}{4}\right)=\frac{15}{8}\leq 2
\end{eqnarray*}
and therefore
$$
\psi\left(\max_{1\leq i\leq |I|}\left(\frac{Y_i}{\sqrt{6\log{(8+i)}}}\right)\right)\leq 1
$$
This ends the proof of the lemma.
\cqfd

The following technical lemma is pivotal in the analysis of maximal inequalities for sequences of random variables
indexed by
 infinite but separable subsets equipped with a pseudo-metric, under some Lipschitz regularity conditions w.r.t. 
 the Orlicz's norm. 

\begin{lemma}\label{lem-Orlicz-3}
We assume that the index set $(I,d)$ is a separable, and totally bounded 
pseudo-metric space, with finite diameter $$
d(I):=\sup_{(i,i)\in I^2}d(i,j)<\infty
$$

We let $(Y_i)_{i\in I}$ be a separable and $\RR$-valued stochastic process indexed by $I$ and such that
$$
\pi_{\psi}(Y_i-Y_j)\leq c~d(i,j)
$$
for some finite constant $c<\infty$.
We also assume that $Y_{i_0}=0$, for some $i_0\in I$. Then, we have
$$
\pi_{\psi}\left(\sup_{i\in I}Y_i\right)\leq 12~c~ \int_{0}^{d(I)}\sqrt{6\log{(8+\Na(I,d,\epsilon)^2)}}~d\epsilon
$$

\end{lemma}
\preuve
Replacing $Y_i$ by $Y_i/d(I)$, and $d$ by $d/d(I)$, 
 there is no loss of generality to assume that $d(I)\leq 1$. In the same way, 
 Replacing $Y_i$ by $Y_i/c$, we can also assume that $c\leq 1$.
For a given finite subset $J\subset I$, with $i_0\in J$, we let
$J_k=\{i^k_1,\ldots,i^k_{n_k}\}\subset J$, be the centers of $n_k=\Na(J,d,2^{-k})$ balls of radius
at most $2^{-k}$ covering $J$. For $k=0$, we set $J_0=\{i_0\}$. We also consider the mapping $\theta_k~:~i\in J\mapsto \theta_k(i)\in J_k$ s.t.
$$
\sup_{i\in J}d(\theta_k(i),i)\leq 2^{-k}
$$
The set $J$ being finite, there exist some sufficiently integer $k^{\star}_J$ s.t. $d(\theta_k(i),i)=0$, for any $k\geq k^{\star}_J$; and therefore $Y_i=Y_{\theta_k(i)}$, for any $i\in J$, and any $k\geq k^{\star}_J$. This implies that
$$
Y_i=\sum_{k=1}^{k^{\star}_J}~\left[Y_{\theta_k(i)}-Y_{\theta_{k-1}(i)}\right]
$$
We also notice that
$$
d(\theta_k(i),\theta_{k-1}(i))\leq d(\theta_k(i),i)+d(i,\theta_{k-1}(i))\leq 2^{-k}+2^{-(k-1)}=3\times 2^{-k}
$$
and
$$
\sup_{(i,j)\in (J_k\times J_{k-1})~:~d(i,j)\leq 3\times 2^{-k}
 }\pi_{\psi}\left(Y_i-Y_j\right)\leq 3\times 2^{-k}
$$
Using lemma~\ref{lem-Orlicz-2} we prove that
\begin{eqnarray*}
\pi_{\psi}\left(\sup_{i\in J}Y_i\right)&\leq& 
\sum_{k=1}^{k^{\star}_J}~\pi_{\psi}\left(\sup_{i\in J}\left[Y_{\theta_k(i)}-Y_{\theta_{k-1}(i)}\right)\right]\\
&\leq &3~\sum_{k=1}^{k^{\star}_J}~
\sqrt{6\log{(8+\Na(J,d,2^{-k})^2)}}~2^{-k}
\end{eqnarray*}
On the other hand, we have $$2\left(2^{-k}-2^{-(k+1)}\right)=2^{-k}$$ and
$$
\sqrt{6\log{(8+\Na(J,d,2^{-k})^2)}}~2^{-k}\leq 2~ \int_{2^{-(k+1)}}^{2^{-k}}\sqrt{6\log{(8+\Na(J,d,\epsilon)^2)}}~d\epsilon
$$
from which we conclude that
$$
\pi_{\psi}\left(\sup_{i\in J}Y_i\right)\leq 6~ \int_{0}^{1/2}\sqrt{6\log{(8+\Na(J,d,\epsilon)^2)}}~d\epsilon
$$
Using the fact that the $\epsilon$-balls with center in $I$ and intersecting $J$ are necessarily contained in an $(2\epsilon)$-ball with center in $J$, we also have
$$\Na(J,d,2\epsilon)\leq \Na(I,d,\epsilon)$$
This implies that
$$
\pi_{\psi}\left(\sup_{i\in J}Y_i\right)\leq 12~ \int_{0}^{1}\sqrt{6\log{(8+\Na(I,d,\epsilon)^2)}}~d\epsilon
$$
The end of the proof is now a direct consequence of the monotone convergence theorem with increasing series of finite subsets exhausting $I$.
This ends the proof of the lemma.
\cqfd

\subsection{Marginal inequalities}\label{marginal-sec}\index{Kintchine's inequalities}

This section is mainly concerned with the
 proof of the theorem~\ref{theo-kintchine}. This result is a more or less direct consequence of the following
technical lemma of separate interest.

\begin{lemma}\label{kintchine}
\label{lemP} Let $M _{n}:=\sum_{0\leq p\leq n}\Delta_{p}$ be a real valued
martingale with symmetric and independent increments $(\Delta_{n})_{n\geq0}$.
For any integer $m\geq1$, and any $n\geq0$, we have
\begin{equation}
\mathbb{E}\left( \left| M _{n}\right| ^{m}\right) ^{\frac{1}{m}%
}\leq~b(m)~ \mathbb{E}\left( \left[ M \right] _{n}^{m^{\prime}/2}\right)
^{\frac{1}{m^{\prime}}}\label{bdg}
\end{equation}
with the smallest even integer $m^{\prime}\geq m$, the bracket process $$\left[
M \right] _{n}:=\sum_{0\leq p\leq n}\Delta^{2}_{p}$$
and the collection of
constants $b(m)$ defined in (\ref{collec}). 
In addition, for any $m\geq 2$, we have
\begin{equation}
\mathbb{E}\left( \left| M _{n}\right| ^{m}\right) ^{\frac{1}{m}}
\leq  b(m)~\sqrt{(n+1)}~\left(\frac{1}{n+1}\sum_{0\leq p\leq n}\EE\left(|\Delta_p|^{m^{\prime}}\right)\right)^{\frac{1}{m^{\prime}}}\label{bdg2}
\end{equation}
\end{lemma}

{\bf Proof of theorem~\ref{theo-kintchine}:}
We consider a
collection of independent copies $X^{\prime}=(X^{\prime i})_{i\geq 1}$
of 
the random variables $X=(X^{i})_{i\geq 1}$. We consider the martingale sequence $M=(M_i)_{1\leq i\leq N}$ with symmetric and independent increments defined for any $1\leq j\leq N$ by
the following formula
$$
M_j:=\frac{1}{\sqrt{N}}\sum_{i=1}^j~\left[f(X^i)-f(X^{\prime i})\right]
$$
By construction, we have
\begin{eqnarray*}
V(X)(f)&=&\frac{1}{\sqrt{N}}\sum_{i=1}^N~\left(f(X^i)-\mu^i(f)\right)=\EE\left(M_N\left|X\right.\right)
\end{eqnarray*}
Comibing this conditioning property with the estimates provided in lemma~\ref{kintchine},
 the proof of the first assertion is now easily completed. 
 
 The Orlicz norm estimate (\ref{orlicz-vxf}) come from the fact that
 for any $f\in\mbox{\rm Osc(E)}$, we have
  $$
\EE(\left|V(X)(f)\right|^{2m})\leq b(2m)^{2m}=\EE\left(U^{2m}\right) 
 $$
for a Gaussian and centred random variable $U$, s.t. $E(U^2)=1$. Using the comparison lemma, lemma~\ref{lem-Orlicz-L2m-compare}, we find that
$$
\pi_{\psi}(V(X)(f))\leq\pi_{\psi}(U)=\sqrt{8/3}
$$

Applying Kintchine's inequalities (\ref{bdg}), we prove that
$$
\mathbb{E}\left( \left| V(X)(f)\right| ^{m} \right) ^{\frac{1}{m}}\leq~b(m)~\EE\left(\left[\frac{1}{N}\sum_{i=1}^N\left[ f(X^{j})-f(X^{\prime j})\right]^2\right]^{m^{\prime}/2}\right)^{1/m^{\prime}}
$$
By construction, we notice that  for any $f\in\mbox{\rm osc(E)}$, and any $p\geq 2$, we have
$$
\frac{1}{N}\sum_{j=1}^N
\EE
\left(
\left[ f(X^{j})-f(X^{\prime j})\right]^{p}
\right)\leq 2\sigma(f)^2
$$

By the Rosenthal type inequality stated in theorem 2.5 
in~\cite{johnson},  for any sequence of nonnegative, independent
and bounded random variables $(Y_i)_{i\geq 1}$,
 we have the rough estimate
$$
\EE\left[\sum_{i=1}^N Y_i^p\right]^{1/p}\leq 2p\max{\left(
\sum_{i=1}^N \EE(Y_i), \left[\sum_{i=1}^N \EE(Y_i^p)\right]^{1/p}
\right)}
$$
for any $p\geq 1$. If we take $p=m^{\prime}/2$, and 
$$
Y_i=\frac{1}{N}~\left[ f(X^{i})-f(X^{\prime i})\right]^2
$$ 
 we prove that
$$
\begin{array}{l}
\EE\left(\left[\frac{1}{N}\sum_{i=1}^N\left[ f(X^{i})-f(X^{\prime i})\right]^2\right]^{m^{\prime}/2}\right)^{2/m^{\prime}}
\\
\\\leq 
4m\max{\left(
2\sigma(f)^2, \frac{1}{N^{1-\frac{2}{m^{\prime}}}}
\left[
2\sigma(f)^2
\right]^{2/m^{\prime}}
\right)}
\end{array}
$$
for any $f\in\mbox{\rm osc(E)}$. Using Stirling's approximation of factorials
$$
\sqrt{2\pi n}~n^n~e^{-n}~\leq n!\leq e~\sqrt{2\pi n}~n^n~e^{-n}
$$
for any $p\geq 1$ we have 
$$
(2p)^p/b(2p)^{2p}=2^{2p}p^p{p!}/{(2p)!}\leq e^{p+1}\leq 3^{2p}
$$
and
$$
(2p+1)^{p+1/2}/b(2p+1)^{2p+1}=(2p+1)^{p+1}2^p p!/(2p+1)!\leq 
e^{p+2}\leq 3^{2p+1}
$$
This implies that
$$
m^{m/2}/b(m)^m\leq 3^m\Rightarrow \sqrt{m}~b(m)\leq 3 b(m)^2
$$
for any $m\geq 1$. This ends the proof of the theorem.
\cqfd

Now, we come to the proof of the lemma.

\noindent\mbox{\bf Proof of lemma~\ref{kintchine}:}\newline We prove the lemma by induction on the
parameter $n$. The result is clearly satisfied for $n=0$. Suppose the estimate
(\ref{bdg}) is true at rank $(n-1)$. To prove the result at rank $n$, we use
the binomial decomposition
\[
\left(  M_{n-1}+\Delta_{n}\right)  ^{2m}=\sum_{p=0}^{2m}\left(
\begin{array}
[c]{c}%
2m\\
p
\end{array}
\right)  M_{n-1}^{2m-p}~\left(  \Delta_{n}\right)  ^{p}%
\]
Using the symmetry condition, all the odd moments of $\Delta_{n}$ are null.
Consequently, we find that
\[
\mathbb{E}\left(  \left(  M_{n-1}+\Delta_{n}\right)  ^{2m}\right)  =\sum
_{p=0}^{m}\left(
\begin{array}
[c]{c}%
2m\\
2p
\end{array}
\right)  ~\mathbb{E}\left(  M_{n-1}^{2(m-p)}\right)  ~\mathbb{E}\left(
\Delta_{n}^{2p}\right)
\]
Using the induction hypothesis, we prove that the above expression is upper
bounded by the quantity
\[%
\begin{array}
[c]{c}%
\sum_{p=0}^{m}\left(
\begin{array}
[c]{c}%
2m\\
2p
\end{array}
\right)  ~2^{-(m-p)}~(2(m-p))_{(m-p)}~\mathbb{E}\left(  \left[  M\right]
_{n-1}^{m-p}\right)  ~\mathbb{E}\left(  \Delta_{n}^{2p}\right)
\end{array}
\]
To take the final step, we use the fact that
\[
\left(
\begin{array}
[c]{c}%
2m\\
2p
\end{array}
\right)  ~2^{-(m-p)}~(2(m-p))_{(m-p)}~=\frac{2^{-m}~(2m)_{m}}{2^{-p}~(2p)_{p}%
}~\left(
\begin{array}
[c]{c}%
m\\
p
\end{array}
\right) 
\]
and $(2p)_{p}\geq2^{p}$,
to conclude that
\begin{align*}
\mathbb{E}\left(  \left(  M_{n-1}+\Delta_{n}\right)  ^{2m}\right)   &
\leq2^{-m}~(2m)_{m}~\sum_{p=0}^{m}\left(
\begin{array}
[c]{c}%
m\\
p
\end{array}
\right)  ~\mathbb{E}\left(  \left[  M\right]  _{n-1}^{m-p}\right)
~\mathbb{E}\left(  \Delta_{n}^{2p}\right) \\
& =2^{-m}~(2m)_{m}~\mathbb{E}\left(  \left[  M\right]  _{n}^{m}\right)
\end{align*}
For odd integers we use twice the Cauchy-Schwarz inequality to deduce that
\begin{align*}
\mathbb{E}(\left\vert M_{n}\right\vert ^{2m+1})^{2}  & \leq\mathbb{E}%
(M_{n}^{2m})~\mathbb{E}(M_{n}^{2(m+1)})\\
& \leq 2^{-(2m+1)}~(2m)_{m}~(2(m+1))_{(m+1)}~\mathbb{E}\left(  \left[
M\right]  _{n}^{m+1}\right)  ^{\frac{2m+1}{m+1}}%
\end{align*}
We conclude that
\[
\mathbb{E}(\left\vert M_{n}\right\vert ^{2m+1})\leq2^{-(m+1/2)}~\frac
{(2m+1)_{(m+1)}}{\sqrt{m+1/2}}~\mathbb{E}\left(  \left[  M\right]  _{n}%
^{m+1}\right)  ^{1-\frac{1}{2(m+1)}}%
\]
The proof of (\ref{bdg}) is now completed. Now, we come to the proof of (\ref{bdg2}). For any $m^{\prime}\geq 2$
we have

$$
 \left[ \frac{1}{n+1} \sum_{0\leq p\leq n}\Delta_p^2\right]^{m^{\prime}/2}\leq 
 \frac{1}{n+1} \sum_{0\leq p\leq n}\mathbb{E}\left( |\Delta_p|^{m^{\prime}}\right)
$$
and therefore
$$
\mathbb{E}\left( \left[ M\right]_n^{m^{\prime}/2}\right)^{\frac{1}{m^{\prime}}}\leq
(n+1)^{1/2}~\left(\frac{1}{n+1} \sum_{0\leq p\leq n}\mathbb{E}\left( |\Delta_p|^{m^{\prime}}\right)\right)^{\frac{1}{m^{\prime}}}
$$

This ends the proof of the lemma. \hfill\hbox{\vrule height 5pt width 5pt
depth 0pt}\medskip\newline

\subsection{Maximal inequalities}\label{maximal-proof-sec}

The main goal of this section is to prove theorem~\ref{theo-proc-emp}.
We begin with the basic symmetrization technique. We consider a
collection of independent copies $X^{\prime}=(X^{\prime i})_{i\geq 1}$ of 
the random variables $X=(X^{i})_{i\geq 1}$.
Let $\e=(\e_i)_{i\geq 1}$ constitute
a sequence that is independent and identically distributed with
$$P(\e_1=+1)=P(\e_1=-1)=1/2
$$ 
We also consider the empirical random field sequences
$$
V_{\epsilon}(X):=\sqrt{N}~
m_{\e}(X)
$$

We also assume that
$(\epsilon,X,X^{\prime})$ are independent. We associate with the pairs
 $(\epsilon,X)$ and $(\epsilon,X^{\prime})$  the 
random measures 
$
m_{\epsilon}(X)=\frac{1}{N}\sum_{i=1}^N\,\epsilon_i~\delta_{X^i}
$ and $m_{\epsilon}(X^{\prime})=\frac{1}{N}\sum_{i=1}^N\,\epsilon_i~\delta_{X^{\prime i}}$.

We notice that 
\begin{eqnarray*}
\|m(X)-\mu\|^p_{\Fa}&&=\sup_{f\in\Fa}|m(X)(f)-\EE(m(X^{\prime})(f))|^p
\\
&\leq& \EE(\|m(X)-m(X^{\prime})\|^p_{\Fa}~|X)
\end{eqnarray*}
and in view of the symmetry of the random variables
$(f(X^i)-f(X^{\prime i}))_{i\geq 1}$ we have
$$
\EE(\|m(X)-m(X^{\prime})\|^p_{\Fa})=
\EE(\|m_{\epsilon}(X)-m_{\epsilon}(X^{\prime})
\|^p_{\Fa})
$$
from which we conclude that
\begin{equation}\label{VV}
E\left(\|V(X)\|_{\Fa}^p\right)\leq 2^{p}\;
E\left(\|V_{\epsilon}(X)\|_{\Fa}^p\right)
\end{equation}
By using the Chernov-Hoeffding inequality 
for any $x^1,\ldots,x^N\in E$, the empirical process
$$
f\longrightarrow V_{\epsilon}(x)(f):=\sqrt{N}~
m_{\e}(x)(f)
$$
is sub-Gaussian for the norm 
$
\|f\|_{L_2(m(x))}
=m(x)(f^2)^{1/2}
$.
Namely, for any couple of functions $f,g$ and any $\d>0$ we have
$$
\EE\left(\left[V_{\epsilon}(x)(f)-V_{\epsilon}(x)(g)\right]^2\right)=
\|f-g\|^2_{\LL_2(m(x))}
$$
and by Hoeffding's inequality
$$
P\left(\left|V_{\epsilon}(x)(f)-V_{\epsilon}(x)(g)\right|\geq \d\right)\leq 2\;e^{-\frac{1}{2}{\d^2}/{\|f-g\|^2_{\LL_2(m(x))}}}
$$
If we set $Z=\left(\frac{V_{\epsilon}(x)(f)}{\sqrt{6}\|f\|_{\LL_2(m(x))}}\right)^2$, then we find that
\begin{eqnarray*}
\EE\left(e^{Z}\right)-1&=&\int_0^{\infty} e^t~\PP\left(Z\geq t\right)~dt\\
&=&\int_0^{\infty} e^t~\PP\left(\left|V_{\epsilon}(x)(f)\right|\geq \sqrt{6t}~\|f\|_{\LL_2(m(x))}\right)~dt\\
&\leq &2~\int_0^{\infty} e^t~e^{-3t}~dt=1
\end{eqnarray*}
from which we conclude that
$$
\pi_{\psi}\left(V_{\epsilon}(x)(f)-V_{\epsilon}(x)(g)\right)\leq \sqrt{6}\|f-g\|_{\LL_2(m(x))}
$$

Combining the maximal inequalities stated in lemma~\ref{lem-Orlicz-3} and the conditioning property (\ref{maj-key}) we find that
$$
\pi_{\psi}\left(\left\|V_{\epsilon}(X)\right\|_{\Fa}\right)\leq J(\Fa)
$$
with
$$
J(\Fa)\leq 2~6^2~ \int_{0}^{2}\sqrt{\log{(8+\Na(\Fa,\epsilon)^2)}}~d\epsilon\leq 
c~I(\Fa)<\infty
$$
for some finite universal constant $c<\infty$.
Combining (\ref{VV}) with (\ref{lem-Orlicz-L2m-compare}), this implies that
$$
\pi_{\psi}\left(\left\|V(X)\right\|_{\Fa}\right)\leq 2~J(\Fa)$$

This ends the proof of the theorem.
\cqfd

\subsection{Cram\'er-Chernov inequalities}\label{Laplace-estimate-sec}
\subsubsection{Some preliminary convex analysis}\label{sconvex}

In this section, we present some basic Cram\'er-Chernov tools to derive quantitative concentration inequalities.
We begin by recalling some preliminary convex analysis on  Legendre-Fenchel transforms.
We associate with any convex function $$L~:~t\in \mbox{\rm Dom(L)}\mapsto L(t)\in \RR_+$$ defined in some
domain $\mbox{\rm Dom(L)}\subset \RR_+$, with $L(0)=0$, the Legendre-Fenchel transform 
$L^{\star}$ defined by the variational formula
$$
\forall \lambda\geq 0\qquad
L^{\star}(\lambda):=\sup_{t\in \mbox{\rm\small Dom(L)}}{\left(\lambda t-L(t)\right)}
$$
Note that $L^{\star}$ is a convex increasing function with $L^{\star}(0)=0$ and its inverse
  $\left(L^{\star}\right)^{-1}$ is a concave increasing function.

 We let $L_A$ be  the log-Laplace transform of a random variable
$A$ defined on some
domain $
\mbox{\rm Dom($L_A$)}
\subset \RR_+
$  by the formula
$$
L_A(t):=\log{\EE(e^{t A})}
$$

H\"older's inequality implies that $L_A$ is convex.
Using the Cram\'er-Chernov-Chebychev inequality, we find that
$$
\log{
\PP\left(A\geq \lambda\right)
}\leq -L_A^{\star}(\lambda)
\quad\mbox{\rm and}\quad 
\PP\left(
A\geq \left(L_A^{\star}\right)^{-1}(x)
\right)\leq e^{-x}
$$
for any $\lambda\geq 0$ and any $x\geq 0$.

The next lemma provides some key properties of Legendre-Fenchel transforms that will
be used in several places in the further development of the lecture notes. 

\begin{lemma}\label{lem-el-L}
\begin{itemize}
\item For any convex functions $(L_1,L_2)$, such that
$$
\forall t\in \mbox{\rm Dom($L_2$)}\quad L_1(t)\leq L_2(t)
\quad\mbox{\rm and}\quad
\mbox{\rm Dom($L_2$)}\subset \mbox{\rm Dom($L_1$)}
$$
we have
$$
L^{\star}_2\leq L^{\star}_1\quad \mbox{\rm and}\quad
(L^{\star}_1)^{-1}\leq (L^{\star}_2)^{-1}
$$
\item  if we have
$$
\forall  t\in v^{-1}\mbox{\rm Dom($L_2$)}=\mbox{\rm Dom($L_1$)}\qquad 
L_1(t)=u~L_2(v~t)
$$
for some positive numbers $(u,v)\in\RR_+^2$, then we have
$$
 L^{\star}_1(\lambda)=u~L^{\star}_2\left(\frac{\lambda}{uv}\right)
 \quad \mbox{\rm and}\quad
( L^{\star}_1)^{-1} (x)=uv~( L^{\star}_2)^{-1} \left(\frac{x}{u}\right)
$$
for any $\lambda\geq 0$, and any $\forall x\geq 0$. 

\item Let $A$ be a random variable
with a finite log-Laplace transform. For any $a\in\RR$, we have
$$
L_A(t)=-at+L_{A+a}(t)
$$
as well as
$$
L_A^{\star}(\lambda)=L_{A+a}^{\star}(\lambda+a)\quad\mbox{\rm and}\quad
\left(L_A^{\star}\right)^{-1}(x)=-a+\left(L_{A+a}^{\star}\right)^{-1}(x)
$$
\end{itemize}
\end{lemma}

We illustrate this technical lemma with the detailed analysis of three
convex increasing functions of current use in the further development of these notes
\begin{itemize}
\item $L(t)={t^2}/{(1-t)}$, $t\in [0,1[$
\item $L_0(t):=-t-\frac{1}{2}\log{(1-2t)}$, $t\in [0,1/2[$.\label{L_0def}
\item $L_1(t):=e^t-1-t$
\end{itemize}

In the first situation, we readily check that
$$
L^{\prime}(t)=\frac{1}{(1-t)^2}-1\quad\mbox{\rm and}\quad
L^{\prime\prime}(t)=\frac{2}{(1-t)^3}
$$
An elementary manipulation yields that
$$
L^{\star}(\lambda)=\left(\sqrt{\lambda+1}-1\right)^2
$$
and
$$
\left(L^{\star}\right)^{-1}(x)=\left(1+\sqrt{x}\right)^2-1=x+2\sqrt{x}
$$

In the second situation, we have
$$
L^{\prime}_0(t)=\frac{1}{1-2t}-1\quad\mbox{\rm and}\quad
L^{\prime\prime}_0(t)=\frac{2}{(1-2t)^2}
$$
from which we find that
$$
L^{\star}_0(\lambda)=\frac{1}{2}\left(\lambda-\log{(1+\lambda)}\right)
$$
We also notice that
$$
L_0(t)=t^2~\sum_{p\geq 0}\frac{2}{2+p}~(2t)^p
\leq \overline{L}_0(t):=\frac{t^2}{1-2t}=\frac{1}{4}~L(2t)
$$
for every $t\in [0,1/2[$
Using lemma~\ref{lem-el-L}, we prove that
\begin{eqnarray}
  \overline{L}_0^{\star}(\lambda)&=&
\frac{1}{4}~L^{\star}(2\lambda)\leq L^{\star}_0(\lambda)\nonumber\\
\left(L^{\star}_0\right)^{-1}(x)&\leq& \left(\overline{L}^{\star}_0\right)^{-1}(x)=\frac{1}{2}\left(L^{\star}\right)^{-1}(4x)=2(x+\sqrt{x})\label{borne-L0star-inv}
\end{eqnarray}

In the third situation, we have
$$
L^{\prime}_1(t)=e^t-1\quad\mbox{\rm and}\quad
L^{\prime\prime}_1(t)=e^t
$$
from which we conclude that
$$
L_1^{\star}(\lambda)=(1+\lambda)\log{(1+\lambda)}-\lambda
$$
On the other hand, using the fact that $2\times 3^p\leq (p+2)!$, for any $p\geq 0$, we prove that
we have
$$
L_1(t)=\frac{t^2}{2}\sum_{p\geq 0}\frac{2\times 3^p}{(p+2)!}~\left(\frac{t}{3}\right)^p
\leq\overline{L}_1(t):=\frac{t^2}{2(1-t/3)}=\frac{9}{2}~L\left(\frac{t}{3}\right)
$$
for every $t\in [0,1/3[$. This implies that
$$  \overline{L}_1^{\star}(\lambda)=
\frac{9}{2}~L^{\star}\left(\frac{2\lambda}{3}\right)\leq L^{\star}_1(\lambda)
$$
and therefore
\begin{equation}\label{borne-L1star-inv}
\left(L^{\star}_1\right)^{-1}(x)\leq \left(\overline{L}^{\star}_1\right)^{-1}(x)
=\frac{3}{2}\left(L^{\star}\right)^{-1}\left(\frac{2x}{9}\right)=\left(\frac{x}{3}+\sqrt{2x}\right)
\end{equation}

Another crucial ingredient in the concentration analysis of the sum
of two random variables is a deep technical lemma of
J. Bretagnolle and E. Rio~\cite{rio}. In the further development of this chapter,
we use this argument to obtain a large family of concentration inequalities that
 are asymptotically "almost sharp" in a wide variety of situations.

\begin{lemma}[J. Bretagnolle \& E. Rio~\cite{rio}]
For any  pair of random variables $A$ and $B$ 
with  finite log-Laplace transform in a 
neighborhood of $0$, we have

\begin{equation}\label{breta}
\forall x\geq 0\quad (L_{A+B}^\star)^{-1}(x)  \leq (L_A^\star)^{-1} (x) + (L_B^\star)^{-1} (x)
\end{equation}

\end{lemma}

We also quote the following reverse type formulae that allows to turn most
of the concentration inequalities developed in these notes into Bernstein style 
exponential inequalities.

\begin{lemma}\label{bernstein-lem}
For any $(u,v)\in\RR_{+}$, we have 
$$
u~\left(L^{\star}_0\right)^{-1}(x)+v~\left(L^{\star}_1\right)^{-1}(x)\leq \left(L^{\star}_{a(u,v),b(u,v)}\right)^{-1}(x)
$$
with the functions
$$
a(u,v):=\left(2u+\frac{v}{3}\right)\quad\mbox{\rm and}\quad b(u,v):=\left(\sqrt{2}~u+v\right)^2
$$
and the Laplace function
$$
L_{a,b}(t)=\frac{b}{2a^2}~L\left( at\right)
\quad
\mbox{\rm with}
\quad
L_{a,b}^{\star}(\lambda)\geq \frac{\lambda^2}{2(b+\lambda a)} 
$$
\end{lemma}
\proof
Using the estimates (\ref{borne-L0star-inv}) and (\ref{borne-L1star-inv}) we prove that
\begin{eqnarray*}
u~\left(L^{\star}_0\right)^{-1}(x)+v~\left(L^{\star}_1\right)^{-1}(x)
&\leq& 2~u~(x+\sqrt{x})+v~\left(\frac{x}{3}+\sqrt{2x}\right)\\
&=& a(u,v)~ x+~\sqrt{2x ~b(u,v)}
\end{eqnarray*}
with
$$
a(u,v):=\left(2u+\frac{v}{3}\right)\quad\mbox{\rm and}\quad b(u,v):=\left(\sqrt{2}~u+v\right)^2
$$
Now, using lemma~\ref{lem-el-L}, we observe that
\begin{equation}\label{Lab}
a~x~+~\sqrt{2x b}~=\left(L^{\star}_{a,b}\right)^{-1}(x)
\quad\mbox{\rm
with}
\quad
L_{a,b}(t)=\frac{b}{2a^2}~L\left( at\right)
\end{equation}
Finally, we have

$$
L^{\star}(\lambda)=\left(\sqrt{\lambda+1}-1\right)^2\geq \frac{(\lambda/2)^2}{(1+\lambda/2)}
$$
The r.h.s. inequality can be easily checked using the fact that
\begin{eqnarray*}
\left(\sqrt{1+2\lambda}-1\right)^2&=&2\left(
\frac{(1+\lambda)^2-(1+2\lambda)}{(1+\lambda)+\sqrt{1+2\lambda}}\right)\\
&\geq&\frac{\lambda^2}{(1+\lambda)}\qquad\left(\Leftarrow~\sqrt{1+2\lambda}\leq (1+\lambda)\right)
\end{eqnarray*}
This implies that
$$
L_{a,b}^{\star}(\lambda)=\frac{b}{2a^2}~L^{\star}\left(\frac{2a}{b}~\lambda\right)\geq \frac{\lambda^2}{2(b+\lambda a)} 
$$
This ends the proof of the lemma.
\cqfd
\subsubsection{Concentration inequalities}

In this section, we investigate some elementary concentration inequalities for 
bounded  and chi-square type random variables. We also apply these results to empirical processes
associated with independent random variables.

\begin{prop}\label{prop-1-ex}
Let $A$ be a centred random variable such that
$A\leq 1$. If we set $\sigma_A=\EE(A^2)^{1/2}$, then for
any $t\geq 0$, we have
\begin{equation}\label{ex-LA}
L_A(t)\leq \sigma^2_A~L_1(t)
\end{equation}
In addition, the probability of the following events
\begin{eqnarray*}
A&\leq & \sigma^2_A~(L_1^{\star})^{-1}\left(\frac{x}{\sigma^2_A}\right)\leq \frac{x}{3}+\sigma_A~\sqrt{2x}
\end{eqnarray*}
is greater than $1-e^{-x}$,
for any $x\geq 0$.

\end{prop}

\preuve
To prove (\ref{ex-LA}) we use the fact the decomposition
$$
\EE\left(e^{tA}-1-A\right)=\EE\left(L_1(tA)1_{X<0}\right)+\EE\left(L_1(tA)1_{X\in [0,1]}\right)
$$
Since we have
$$
\forall x\leq 0\quad L_1(tx)\leq (tx)^2/2$$
and
$$
\forall x\in [0,1]\quad L_1(tx)=x^2\sum_{n\geq 2} x^{n-2}t^n/n!\leq x^2L_1(t)
$$
 we conclude that
\begin{eqnarray*}
\EE\left(e^{tA}\right)&\leq& 1+\frac{t^2}{2}~\EE(A^21_{A<0})+L_1(t)\EE\left(A^21_{A\in [0,1]}\right)\\
&\leq& 1+L_1(t)\sigma^2_A\leq e^{L_1(t)\sigma^2_A}
\end{eqnarray*}

Using lemma~\ref{lem-el-L}, we readily prove that
$$
 (L_A^{\star})^{-1}(x)\leq \sigma^2_A~(L_1^{\star})^{-1}\left(\frac{x}{\sigma^2_A}\right)\leq \frac{x}{3}+\sigma_A~\sqrt{2x}
$$
This ends the proof of the proposition.
\cqfd

\begin{prop}\label{prop-LV(X)}
For any measurable function $f$, with $0<\mbox{\rm osc(f)}\leq a$,  any $N\geq 1$,
and any $t\geq 0$, we have
\begin{equation}\label{LV(X)}
L_{\sqrt{N}V(X)(f)}(t)\leq N~\sigma^2(f/a)~L_1(at)
\end{equation}

In addition, the probability of the following events
\begin{eqnarray}
V(X)(f)&\leq & a^{-1} \sigma^2(f)\sqrt{N}~\left(L_1^{\star}\right)^{-1}\left(\frac{x a^2}{N\sigma^2(f)}\right)\nonumber\\
&\leq& \frac{x a}{3\sqrt{N}}+\sqrt{2x\sigma(f)^2}\label{compare-2}
\end{eqnarray}
is greater than $1-e^{-x}$, for any $x\geq 0$.
\end{prop}
\preuve
Replacing $f$ by $f/a$, there is no loss of generality to assume that $a=1$.
Using the same arguments as the ones we used in the proof of proposition~\ref{prop-1-ex}, we find that
$$
\log{\EE\left(e^{t (f(X^i)-\mu^i(f)) }\right)}\leq ~\mu^i\left(\left[f-\mu^i(f)\right]^2\right)~L_1(t)
$$
from which we conclude that
\begin{eqnarray*}
L_N(t)&:=&\log{\EE\left(e^{t \sqrt{N}~V(X)(f)}\right)}\\
&=&\sum_{i=1}^N
\log{\EE\left(e^{t (f(X^i)-\mu^i(f)) }\right)}\leq \overline{L}_N(t):=N~\sigma^2(f)~L_1(t)
\end{eqnarray*}
By lemma~\ref{lem-el-L}, we have
$$
\left(L_N^{\star}\right)^{-1}(x)\leq \left(\overline{L}_N^{\star}\right)^{-1}(x)=
N \sigma^2(f)~\left(L_1^{\star}\right)^{-1}\left(\frac{x }{N\sigma^2(f)}\right)
$$
This ends the proof of the proposition.
\cqfd
\begin{prop}\label{prop-LB}
For any random variable  $B$  such that
$$
\EE\left(|B|^m\right)^{1/m}\leq b(2m)^{2}~c\quad\mbox{\rm with}\quad c<\infty
$$
for any $m\geq 1$, with the finite constants $b(m)$ defined in (\ref{collec}), we have
\begin{equation}\label{LB}
L_{B}(t) \leq ct+L_0(ct)
\end{equation}
for any $0\leq ct<1/2$. In addition,  the probability of the following events
\begin{eqnarray*}
B\leq c\left[1+ \left(L^{\star}_0\right)^{-1}(x)\right]\leq c\left[1+ 2(x+\sqrt{x})\right]
\end{eqnarray*}
is greater than $1-e^{-x}$, for any $x\geq 0$.
\end{prop}
\preuve
Replacing $B$ by $B/c$, there is no loss of generality to assume that $c=1$.
We recall that $b(2m)^{2m}=\mathbb{E}(U^{2m})$ for every centred Gaussian random
variable with $\mathbb{E}(U^{2})=1$ and
\[
\forall t\in\lbrack0,1/2)\quad\sum_{m\geq0}\frac{t^{m}}{m!}~b(2m)^{2m}=\frac{1}{\sqrt{1-2t}}=\mathbb{E}(\exp{\left\{  tU^{2}\right\}  }%
)
\]
This implies that
$$
\mathbb{E}(\exp{\left\{  tB\right\}  }
)\leq  \sum_{m\geq 0}\frac{t^{m}}{m!}~b(2m)^{2m}=\frac{1}{\sqrt{1-2t}}
$$
for any $0\leq t<1/2$. In other words, we have
$$
L_{B-1}(t):=\log{\mathbb{E}(\exp{\left\{  t(B-1)\right\}  }}\leq L_0(t)$$
and
$$
L_B(t)=t+L_{B-1}(t)\leq t+L_0(t)
$$
from which we conclude that
$$
L_B^{\star}(\lambda)=L^{\star}_{B-1}(\lambda-1)\Rightarrow
\left(L_B^{\star}\right)^{-1}(x)=1+\left(L_{B-1}^{\star}\right)^{-1}(x)\leq 1+\left(L_0^{\star}\right)^{-1}(x)
$$
This ends the proof of the proposition.
\cqfd

 \begin{remark}\label{cor-last-assertion}
 We end this section with some comments on the estimate
 (\ref{simga2N}).  Using the fact that $b(m)\leq b(2m)$ (see for instance (\ref{b2mbm})) we readily 
 deduce from (\ref{simga2N}) that
$$
\mathbb{E}\left( \left| V(X)(f)\right| ^{m} \right) ^{\frac{1}{m}}\leq~ 6\sqrt{2}~ b(2m)^2\sigma(f)
$$
for any $m\geq 1$, and for any $N$ s.t. $2\sigma^2(f)N\geq 1$. Thus, if we set
$$
B=\left|V(X)(f)\right|\quad\mbox{\rm and}\quad c=6\sqrt{2}~\sigma(f)
$$
in proposition~\ref{prop-LB}, we prove that
for any $N$ s.t. $2\sigma^2(f)N\geq 1$, and for any $0\leq t<1/(12\sqrt{2}~\sigma(f))$
$$
L_{\left|V(X)(f)\right|}(t) \leq 6\sqrt{2}~\sigma(f)~t+L_0(6\sqrt{2}~\sigma(f)~t)
$$
In addition, the probability of the following events
\begin{eqnarray*}
\left|V(X)(f)\right|&\leq&6\sqrt{2}~\sigma(f)\left[1+ \left(L^{\star}_0\right)^{-1}(x)\right]\\
&\leq & 6\sqrt{2}~\sigma(f)\left[1+ 2(x+\sqrt{x})\right]
\end{eqnarray*}
is greater than $1-e^{-x}$, for any $x\geq 0$.

When $N$ is chosen so that $2\sigma^2(f)N\geq 1$, using (\ref{compare-2}) we  improve the above inequality. 
Indeed, using this concentration inequality implies that  
for any $f\in\mbox{\rm Osc}(E)$,
the probability of the following events
\begin{eqnarray}
V(X)(f)
&\leq& \sqrt{2}~\sigma(f)~\left(
\frac{x }{3}+\sqrt{x}\right)
\end{eqnarray}
is greater than $1-e^{-x}$, for any $x\geq 0$. 
\end{remark}

\subsection{Perturbation analysis}\label{perturbation-emp-sec}

This section is mainly concerned with the proof of theorem~\ref{theo-perturb-marginal},
and theorem~\ref{theo-i-empirique}. 
We recall that for any
second order smooth function $F$ on $\RR^d$, for some $d\geq 1$,  $F(m(X)(f))$ stands for the random
functionals defined by
$$
\begin{array}{l}
f=(f_i)_{1\leq i\leq d}\in\mbox{\rm Osc}(E)^d\\
\\
\mapsto F(m(X)(f))=F(m(X)(f_1),\ldots,m(X)(f_d))\in\RR
\end{array}
$$

Both results rely on the following 
second order decomposition of independent interest.

\begin{prop}\label{prop-dl-Fmx}
For any $N\geq 1$, we have the decomposition
$$
\sqrt{N}\left[F(m(X)(f))-F(\mu(f))\right]=
V(X)\left[D_{\mu}(F)(f)\right]+\frac{1}{\sqrt{N}}~R(X)(f)
$$
with a first order functional $D_{\mu}(F)(f)$ defined in (\ref{DnuFf}), and a second order term
$R(X)(f)$ such that
$$
\EE\left(\left|R(X)(f)\right|^m\right)^{1/m}
\leq ~\frac{1}{2}~b(2m)^{2}\left\|\nabla^2F_f\right\|_1
$$
 for any $m\geq 1$, with the parameter defined in (\ref{DnuFf2}).
\end{prop}

\preuve
Using a Taylor first order expansion, we have
$$
\sqrt{N}\left[F(m(X)(f))-F(\mu(f))\right]=
\nabla F(\mu(f))~V(X)(f)^{\top}+\frac{1}{\sqrt{N}}~R(X)(f)
$$
with the second order remainder term
$$
\begin{array}{l}
R(X)(f)\\
\\
:=\int_0^1(1-t)
V(X)(f)~\nabla^2F\left(t m(X)(f)+(1-t)\mu(f)\right)~V(X)(f)^{\top}~dt
\end{array}
$$
We notice that
$$
\nabla F(\mu(f))~V(X)(f)^{\top}=V(X)\left[\nabla F(\mu(f))~f^{\top}\right]
$$
and
$$
\mbox{\rm osc}\left(\nabla F(\mu(f))~f^{\top}\right)\leq \sum_{i=1}^d\left|\frac{\partial F}{\partial u^i}(\mu(f))\right|
$$
It is also easily checked that
$$
\begin{array}{l}
\EE\left(\left|R(X)(f)\right|^m\right)^{1/m}\\
\\
\leq \frac{1}{2}\sum_{i,j=1}^d~\sup_{\nu\in\Pa(E)}
{\left| \frac{\partial^2 F}{\partial u^i\partial u^j}(\nu(f))\right|}~\EE\left(\left|V(X)(f_i)V(X)(f_j)\right|^m\right)^{1/m}
\end{array}
$$
and for any $1\leq i,j\leq d$, we have
$$
\begin{array}{l}
\EE\left(\left|V(X)(f_i)V(X)(f_j)\right|^m\right)^{1/m}\\
\\
\leq \EE\left(V(X)(f_j)^{2m}\right)^{1/(2m)}
\EE\left(V(X)(f_j)^{2m}\right)^{1/(2m)}\leq b(2m)^{2}
\end{array}
$$
This ends the proof of the proposition.
\cqfd

We are now in position to prove theorem~\ref{theo-perturb-marginal}.

{\bf Proof of theorem~\ref{theo-perturb-marginal}:}

We set
$$
N\left[F(m(X)(f))-F(\mu(f))\right]=
A+B
$$
with
$$
A=\sqrt{N}~V(X)\left[D_{\mu}(F)(f)\right]\quad\mbox{\rm and}\quad
B=R(X)(f)
$$
Combining proposition~\ref{prop-LV(X)} with proposition~\ref{prop-LB},
if we set
$$
g=D_{\mu}(F)(f)
\quad
a=\left\|\nabla F(\mu(f))\right\|_1
\quad\mbox{\rm and}\quad
c=\frac{1}{2}\left\|\nabla^2F_f\right\|_1
$$
then we have
\begin{eqnarray*}
L_A(t)&\leq& N\sigma^2(g/a) L_1(at)\\
L_B(t)&=&ct+L_{B-c}(t)\quad\mbox{\rm with}\quad L_{B-c}(t)\leq L_0(ct)
\end{eqnarray*}
On the other hand, we have
$$
\left(L_A^{\star}\right)^{-1}(x)\leq N~a~\sigma^2(g/a) \left(L_1^{\star}\right)^{-1}
\left(\frac{x}{N\sigma^2(g/a)}\right)\\
$$
and using the fact that
$$
L_B(t)=ct+L_{B-c}(t)
$$
we prove that
$$
L^{\star}_B(\lambda)=L_{B-c}^{\star}(\lambda-c)
\Rightarrow\begin{array}[t]{rcl}
\left(L^{\star}_B\right)^{-1}(x)&=&c+\left(L^{\star}_{B-c}\right)^{-1}(x)\\
&\leq & c\left(1+\left(L_0^{\star}\right)^{-1}(x)\right)
\end{array}
$$
Using Bretagnolle-Rio's lemma, we find that
\begin{eqnarray*}
(L_{A+B}^\star)^{-1}(x)& \leq& (L_A^\star)^{-1}(x) + (L_B^\star)^{-1} (x)\\
&\leq &N~a^{-1}~\sigma^2(g) \left(L_1^{\star}\right)^{-1}
\left(\frac{x a^2}{N\sigma^2(g)}\right)+c\left(1+\left(L_0^{\star}\right)^{-1}(x)\right)
\end{eqnarray*}
This ends the proof of the theorem.
\cqfd

Now, we come to the proof of theorem~\ref{theo-i-empirique}.

{\bf Proof of theorem~\ref{theo-i-empirique}:}

We consider the empirical processes $$f\in \Fa_i\mapsto m(X)(f)\in \RR
$$ associated with $d$ classes of functions $\Fa_i$, $1\leq i\leq d$, defined in section~\ref{intro-proc-empirique}. We further assume that $\|f_i\|\vee\mbox{\rm osc}(f_i)\leq 1$, for any $f_i\in \Fa_i$, and we set
\begin{eqnarray*}
\pi_{\psi}(\|V(X)\|_{\Fa})&:=&
\sup_{1\leq i\leq d}\pi_{\psi}(\|V(X)\|_{\Fa_i})
\end{eqnarray*}

Using theorem~\ref{theo-proc-emp}, we have that
\begin{eqnarray*}
\pi_{\psi}(\|V(X)\|_{\Fa})&\leq&
12^2~ \int_{0}^{2}\sqrt{\log{(8+\Na(\Fa,\epsilon)^2)}}~d\epsilon
\end{eqnarray*}
with
$$
\Na(\Fa,\epsilon):=\sup_{1\leq i\leq d}\Na(\Fa_i,\epsilon)
$$

Using proposition~\ref{prop-dl-Fmx}, for any collection of functions 
$$f=(f_i)_{1\leq i\leq d}\in \Fa:=\prod_{i=1}^d\Fa_i$$
we have
$$
\begin{array}{l}
\sqrt{N}\sup_{f\in\Fa}\left|F(m(X)(f))-F(\mu(f))\right|\\
\\
\leq 
\left\|\nabla F_{\mu}\right\|_{\infty}
\sum_{i=1}^d \left\|V(X)\right\|_{\Fa_i}+\frac{d}{2\sqrt{N}}\left\|\nabla^2F\right\|_{\infty}~\sum_{i=1}^d \left\|V(X)\right\|_{\Fa_i}^2
\end{array}
$$
If we set
$$
A:=\left\|\nabla F_{\mu}\right\|_{\infty}
\sum_{i=1}^d \left\|V(X)\right\|_{\Fa_i}
$$
then we find that
$$
\pi_{\psi}(A)\leq \left\|\nabla F_{\mu}\right\|_{\infty}
\sum_{i=1}^d \pi_{\psi}\left(\left\|V(X)\right\|_{\Fa_i}\right)
$$
By lemma~\ref{lem-Orlicz-L2m}, this implies that
$$
\EE\left(e^{tA}\right)\leq (1+t\pi_{\psi}(A))~e^{(t\pi_{\psi}(A))^2}\leq
e^{at+\frac{1}{2}t^2b}
$$
with $b=2a^2$ and
$$
a=\pi_{\psi}(A)\leq 
\left\|\nabla F_{\mu}\right\|_{\infty}
\sum_{i=1}^d \pi_{\psi}\left(\left\|V(X)\right\|_{\Fa_i}\right)
$$
Notice that
$$
L_{A-a}(t)\leq L(t)=\frac{1}{2}t^2b
$$
Recalling that
$$
L^{\star}(\lambda)=\frac{\lambda^2}{2b}\quad\quad\mbox{\rm and}
\quad \left(L^{\star}\right)^{-1}(x)=\sqrt{2bx}
$$
we conclude that
\begin{eqnarray*}
\left(L_{A}^{\star}\right)^{-1}(x)&=&a+
\left(L_{A-a}^{\star}\right)^{-1}(x)\\
&\leq& a+\sqrt{2bx}=\pi_{\psi}(A)\left(1+2\sqrt{x}\right)
\end{eqnarray*}

Now, we come to the analysis of the second order term defined by
$$
B=\frac{d}{2\sqrt{N}}\left\|\nabla^2F\right\|_{\infty}~\sum_{i=1}^d \left\|V(X)\right\|_{\Fa_i}^2
$$
Using  the inequality
 $$\left(\sum_{i=1}^d a_i\right)^m\leq d^{m-1}\sum_{i=1}^d a_i^m$$ which is valid for any $d\geq$,  any $m\geq 1$, and any 
sequence of real numbers $(a_i)_{1\leq i\leq d}\in\RR^d_+$, we prove  that
$$
\EE(B^m)\leq \beta^m~d^{m-1}\sum_{i=1}^d \EE\left(\left\|V(X)\right\|_{\Fa_i}^{2m}\right)
$$
with
$$
\beta:=\frac{d}{2\sqrt{N}}\left\|\nabla^2F\right\|_{\infty}~
$$
Combining lemma~\ref{lem-Orlicz-L2m} with theorem~\ref{theo-proc-emp}, we conclude that
$$
\EE(B^m)\leq m!~\left(\beta~d~ \pi_{\psi}(\|V(X)\|_{\Fa})^{2}\right)^m
$$
If we set 
$$
b:=\beta~d~ \pi_{\psi}(\|V(X)\|_{\Fa})^{2}
$$
then we have that
$$
\EE\left(e^{tB}\right)\leq \sum_{m\geq 0} (bt)^m=\frac{1}{1-bt}=e^{bt}\times e^{2L_0(bt/2)}
$$
for any $0\leq t<1/b$ with the convex increasing function $L_0$ introduced on page~\pageref{L_0def}, so that
$$
2L_0(bt/2)=-bt-\log{(1-bt)}
$$
Using lemma~\ref{lem-el-L}, we prove that
$$
L_{B-b}(t)\leq 2L_0(bt/2)$$
and
\begin{eqnarray*} 
\left(L_{B}^{\star}\right)^{-1}(x)&=&b+\left(L_{B-b}^{\star}\right)^{-1}(x)\\
&\leq& b
\left(1+\left(L_0^{\star}\right)^{-1}\left(\frac{x}{2}\right)
\right)\\
&=&\frac{1}{2\sqrt{N}}\left\|\nabla^2F\right\|_{\infty} \left(d~\pi_{\psi}(\|V(X)\|_{\Fa})\right)^{2}\left(1+\left(L_0^{\star}\right)^{-1}\left(\frac{x}{2}\right)
\right)
\end{eqnarray*}
Finally, using the Bretagnolle-Rio's lemma, we prove that
$$
\begin{array}{l}
\left(L_{A+B}^{\star}\right)^{-1}(x)\\
\\
\leq d \pi_{\psi}\left(\left\|V(X)\right\|_{\Fa}\right)~
\left[
\left\|\nabla F_{\mu}\right\|_{\infty}
\left(1+2\sqrt{x}\right)\right.\\
\\
\qquad+\left.
\frac{1}{2\sqrt{N}}\left\|\nabla^2F\right\|_{\infty} \left(d~\pi_{\psi}(\|V(X)\|_{\Fa})\right)\left(1+\left(L_0^{\star}\right)^{-1}\left(\frac{x}{2}\right)
\right)\right]
\end{array}$$
This ends the proof of the theorem~\ref{theo-i-empirique}.
\cqfd

\section{Interacting empirical processes}\label{ips-emp-sec}
\subsection{Introduction}\label{intro-ips-empirique}
This short chapter is concerned with the concentration analysis of  
sequences of empirical processes associated with
conditionally independent random variables. 

In preparation for the work in
chapter~\ref{fk-particle-chap} on the collection of Feynman-Kac particle models introduced in section~\ref{sec-ips-fk-intro}, we consider a  general class of 
interaction particle processes
with non necessarily mean field type dependency. 

Firstly, we analyze the concentration properties of integrals
of local sampling error sequences, with general random but predictable test functions.  These results will be used to analyze
the concentration properties of the first order fluctuation terms of the particle models.

We also present a stochastic perturbation technique to analyze the
second order type decompositions. We consider finite marginal models and empirical processes. We close the chapter with the analysis of the covering numbers and the entropy parameters of linear transformation of classes of functions.
 
 We end this introductory section, with the precise description of the main mathematical objects
we shall analyze in the further development of the chapter. 

We let $X_n^{(N)}=(X^{(N,i)}_n)_{1\leq i\leq N}$ 
be a Markov chain on some product state spaces $E_n^N$, for some $N\geq 1$.  
We also let $\Ga_n^N$
be the increasing $\sigma$-field generated by the random sequence
$(X_p^{(N)})_{0\leq p\leq n}$. We further assume that $(X^{(N,i)}_n)_{1\leq i\leq N}$
are conditionally independent, given $\Ga_{n-1}^N$.

As traditionally, when there is no possible confusion, 
we simplify notation and suppress the index $(\point)^{(N)}$ so 
that we write $(X_n,X^i_n,\Ga_n)$ instead of 
$(X^{(N)}_n,X^{(N,i)}_n,\Ga^N_n)$.
 
In this simplified notation, we also denote by $\mu^i_n$
the conditional distribution of the random state $X^i_n$ given the 
$\Ga_{n-1}$; that is, we have that
$$
\mu^i_n=\mbox{\rm Law}(X^i_n~|~\Ga_{n-1})
$$
Notice that the conditional distributions
$$
\mu_n:=\frac{1}{N}\sum_{i=1}^N \mu^i_n
$$
represent the local conditional mean  of  the occupation measures $$
m(X_n):=\frac{1}{N}\sum_{i=1}^N\delta_{X^i_n}$$ 

At this level of generality, we cannot obtain
any kind of concentration properties for the deviations of  the occupation measures $
m(X_n)$ around some deterministic limiting value.

In chapter~\ref{fk-particle-chap}, dedicated to  particle approximations of Feynman-Kac measures
$\eta_n$, we shall deal with
mean field type random measures $\mu_n^i$, in the sense that the randomness only depends
on the location of the random state $X^i_{n-1}$ and on the current occupation measure
$m(X_{n-1})$. In this situation, the fluctuation of $m(X_n)$ around the limiting deterministic measures $\eta_n$ will be expressed in terms of second order Taylor's type expansions w.r.t. the local sampling errors
$$
V(X_p)=\sqrt{N}\left(m(X_p)-\mu_p\right)
$$
from the origin $p=0$, up to the current time $p=n$.

The first order terms will be expressed in terms of integral formulae of predictable functions $f_p$
w.r.t. the local sampling error measures $V(X_p)$. These stochastic first order expansions are defined below.
\begin{definition}
For any sequence of $\Ga_{n-1}$-measurable random function $f_n\in\mbox{\rm Osc}(E_n)$, and any numbers
$
a_n\in\RR_+$,
we set
\begin{equation}\label{defVnX}
V_{n}(X)(f)=\sum_{p=0}^n a_p~V(X_p)(f_p)
\end{equation}
\end{definition}
For any $\Ga_{n-1}$-measurable random function $f_n\in\mbox{\rm Osc}(E_n)$, we have
\begin{eqnarray*}
\EE\left(V(X_n)(f_n)\left|\Ga_{n-1}\right.\right)&=&0\\
\EE\left(V(X_n)(f_n)^2\left|\Ga_{n-1}\right.\right)&=&\sigma_n^N(f_n)^2:=\frac{1}{N}\sum_{i=1}^N \mu^i_n\left(\left[f_n-\mu_n^i(f_n)\right]^2\right)
\end{eqnarray*}
We also assume that we have an almost sure estimate
\begin{equation}\label{Condition-sigma}
\sup_{N\geq 1}{\sigma_n^N(f_n)^2}\leq \sigma_n^2\quad\mbox{\rm for some positive constant $\sigma_n^2\leq 1$.}
\end{equation}

\subsection{Finite marginal models}

We will now derive a quantitative contraction inequality for the general random fields models
of the following form
\begin{equation}\label{defWVR}
W_{n}(X)(f)=V_n(X)(f)+\frac{1}{\sqrt{N}}~R_n(X)(f)
\end{equation}
with $V_n(X)(f)$ defined in (\ref{defVnX}), and a second order term such that
$$
\EE\left(\left|R_n(X)(f)\right|^m\right)^{1/m}\leq b(2m)^2~c_n
$$
for any $m\geq 1$, for some finite constant $c_n<\infty$ whose values only depend on the parameter $n$.

For a null remainder term $R_n(X)(f)=0$, these concentration properties are
easily derived using proposition~\ref{prop-LV(X)}.

\begin{prop}\label{prop-inter-ips}
We let $V_{n}(X)(f)$ be the random field sequence defined in (\ref{defVnX}).
For any  $t\geq 0$, we have that
$$
L_{\sqrt{N}V_{n}(X)(f)}(t)\leq N~\overline{\sigma}^2_n~
L_1(ta^{\star}_n)
$$
with the parameters
$$
\overline{\sigma}^2_n:=\sum_{0\leq p\leq n}\sigma^2_p\quad\mbox{\rm and}
\quad
a^{\star}_n:=\max_{0\leq p\leq n} a_p
$$
In addition,  the probability of the following events
\begin{eqnarray*}
V_{n}(X)(f)&\leq& 
\sqrt{N}~a^{\star}_n~\overline{\sigma}^2_n~\left(L_{1}^{\star}\right)^{-1}\left(\frac{x}{N\overline{\sigma}^2_n}\right)\\
&\leq & a^{\star}_n\left(\frac{x}{3\sqrt{N}}+\sqrt{2\overline{\sigma}^2_n}~x\right)
\end{eqnarray*}
is greater than $1-e^{-x}$, for any $x\geq 0$.
\end{prop}

\preuve
By proposition~\ref{prop-LV(X)}, we have 
$$
\EE\left(e^{t \sqrt{N}V_{n}(X)(f)}\left|\Ga_{n-1}\right.\right)=e^{t \sqrt{N}V_{n-1}(X)(f)}~
\EE\left(e^{(t a_n) \sqrt{N}V(X_n)(f_n)  }\left| \Ga_{n-1}\right.\right)
$$
with
$$
\log{\EE\left(e^{(t a_n) \sqrt{N}V(X_n)(f_n)}\left|\Ga_{n-1}\right.\right)}\leq N~\sigma^2_n~L_1(ta^{\star}_n)
$$
This clearly implies that
$$
L_{\sqrt{N}V_{n}(X)(f)}(t)\leq \overline{L}_1(t):=N~\overline{\sigma}^2_n~
L_1(ta^{\star}_n)
$$
Using lemma~\ref{lem-el-L}, we conclude that
\begin{eqnarray*}
\left(L_{\sqrt{N}V_{n}(X)(f)}^{\star}\right)^{-1}(x)&\leq& \left(\overline{L}_{1}^{\star}\right)^{-1}(x)\\
&=&
N~a^{\star}_n~\overline{\sigma}^2_n~\left(L_{1}^{\star}\right)^{-1}\left(\frac{x}{N\overline{\sigma}^2_n}\right)
\end{eqnarray*}
The last assertion is a direct consequence of (\ref{borne-L1star-inv}). This ends the proof of the proposition.
\cqfd

\begin{theorem}\label{theo-ref-marginal-ips}
We let $W_{n}(X)(f)$ be the random field sequence defined in (\ref{defWVR})

In this situation, the probability of the events
$$
\sqrt{N}~W_{n}(X)(f)\leq c_n~\left(1+\left(L_0^{\star}\right)^{-1}(x)\right)+
N~a^{\star}_n~\overline{\sigma}^2_n~\left(L_{1}^{\star}\right)^{-1}\left(\frac{x}{N\overline{\sigma}^2_n}\right)
$$
is greater than $1-e^{-x}$, for any $x\geq 0$. In the above display, $\overline{\sigma}_n$ stands for the
variance parameter definition in proposition~\ref{prop-inter-ips}
\end{theorem}

\preuve
We set
$
\sqrt{N}~W_{n}(X)(f)=A_n+B_n
$,
with
$$
A_n=\sqrt{N}V_n(X)(f)\quad\mbox{\rm and}\quad
B_n=R_n(X)(f)
$$
By  proposition~\ref{prop-LB} and proposition~\ref{prop-inter-ips}, we have
$$
L_{A_n}(t)\leq  \overline{L}_{A_n}(t):=N~\overline{\sigma}^2_n(f)~
L_1(ta^{\star}_n)$$
and
$$
L_{B_n-c_n}(t) \leq \overline{L}_{B_n-c_n}(t):=L_0(c_nt)
$$
We recall that
$$
L_{B_n}(t)=c_nt+L_{B_n-c_n}(t)\Rightarrow L_{B_n}^{\star}(\lambda)=L_{B_n-c_n}^{\star}(\lambda-c_n)
$$
Using lemma~\ref{lem-el-L}, we also have that
\begin{eqnarray*}
\left(\overline{L}_{B_n}^{\star}\right)^{-1}(x)&=&c_n+\left(L_{B_n-c_n}^{\star}\right)^{-1}(x)
\\&\leq& c_n+\left(\overline{L}_{B_n-c_n}^{\star}\right)^{-1}(x)=c_n\left(1+\left(L_0^{\star}\right)^{-1}(x)\right)
\end{eqnarray*}
In the same vein, arguing as in the end of the proof of proposition~\ref{prop-inter-ips}, we have
\begin{eqnarray*}
\left(L_{\sqrt{N}V_{n}(X)(f)}^{\star}\right)^{-1}(x)&\leq& N~a^{\star}_n~\overline{\sigma}^2_n~\left(L_{1}^{\star}\right)^{-1}\left(\frac{x}{N\overline{\sigma}^2_n}\right)
\end{eqnarray*}
The end of the proof is now a direct consequence of the Bretagnolle-Rio's lemma. This ends the proof of the theorem.
\cqfd

\subsection{Empirical processes}

We let $V_{n}(X)$ be the random field sequence defined in (\ref{defVnX}), and we
consider a sequence  of classes of $\Ga_{n-1}$-measurable random 
functions $\Fa_n$,  such that $\|f_n\|\vee\mbox{\rm osc}(f_n)\leq 1$, for any $f_n\in \Fa_n$. 

\begin{definition}
For any $f=(f_n)_{n\geq 0}\in \Fa:=(\Fa_n)_{n\geq 0}$, and any sequence 
of numbers $a=(a_n)_{n\geq 0}\in\RR_+^{\NN}$,
we set
$$
V_{n}(X)(f)=\sum_{p=0}^n a_p~V(X_p)(f_p)\quad\mbox{\rm and}\quad
\left\|V_{n}(X)\right\|_{\Fa}=\sup_{f\in\Fa}{\left|V_{n}(X)(f)\right|}
$$
\end{definition}
We further assume that
for any $n\geq 0$, and any $\epsilon>0$, we have an almost sure estimate
\begin{equation}\label{entropy-cond-random}
\Na\left(\Fa_n,\epsilon\right)\leq \Na_n(\epsilon)
\end{equation}
for some non increasing function $\Na_n(\epsilon)$ such that
$$
b_n:=12^2~ \int_{0}^{2}\sqrt{\log{(8+\Na_n(\epsilon)^2)}}~d\epsilon<\infty
$$
 In this situation, we have
$$
\pi_{\psi}\left(\left\|V_{n}(X)\right\|_{\Fa}\right)\leq \sum_{p=0}^n~a_p~\pi_{\psi}\left(\left\|V(X_p)\right\|_{\Fa_p}\right)
$$
Using theorem~\ref{theo-proc-emp}, given $\Ga_{n-1}$ we have the almost sure upper bound 
\begin{eqnarray*}
\pi_{\psi}\left(\left\|V(X_p)\right\|_{\Fa_p}\right)&\leq& 12^2~ \int_{0}^{2}\sqrt{\log{(8+\Na(\Fa_p,\epsilon)^2)}}~d\epsilon\leq b_p\end{eqnarray*}
Combining lemma~\ref{lemme-0} and lemma~\ref{lem-Orlicz-L2m}, we readily prove the following theorem.
\begin{theorem}\label{theo-class-random}
For any classes of $\Ga_{n-1}$-measurable random 
functions $\Fa_n$ satisfying the entropy condition (\ref{entropy-cond-random}), we have
$$
\pi_{\psi}\left(\left\|V_{n}(X)\right\|_{\Fa}\right)\leq c_n:=\sum_{p=0}^n~a_p b_p
$$
In particular, the probability of the events
$$
\left\|V_{n}(X)\right\|_{\Fa}\leq c_n~\sqrt{x+\log{2}}
$$
is greater than $1-e^{-x}$, for any $x\geq 0$.
\end{theorem}

Next, we consider classes of non random 
functions $\Fa=(\Fa_n)_{n\geq 0}$.
We further assume that $\|f_n\|\vee\mbox{\rm osc}(f_n)\leq 1$, for any $f_n\in \Fa_n$,
and
 $$
 I_1(\Fa):=12^2~ \int_{0}^{2}\sqrt{\log{(8+\Na(\Fa,\epsilon)^2)}}~d\epsilon<\infty
$$
with $$
\Na(\Fa,\epsilon)=\sup_{n\geq 0}\Na(\Fa_n,\epsilon)<\infty
$$

\begin{theorem}
We let $W_{n}(X)(f)$, $f\in\Fa$, be the random field sequence defined by
$$
W_{n}(X)(f)=V_n(X)(f)+\frac{1}{\sqrt{N}}~R_n(X)(f)
$$
with a second order term such that
$$
\EE\left(\sup_{f\in\Fa}{\left|R_{n}(X)(f)\right|}^m\right)\leq m!~c_n^m
$$
for any $m\geq 1$, for some finite constant $c_n<\infty$ whose values only depend on the parameter $n$.
In this situation, the probability of the events
$$
\left\|W_{n}(X)\right\|_{\Fa}\leq 
\left[\sum_{p=0}^n a_p\right]~I_1(\Fa)\left(1+2\sqrt{x}\right)+\frac{c_n}{\sqrt{N}}
\left(1+\left(L_0^{\star}\right)^{-1}\left(\frac{x}{2}\right)
\right)
$$
is greater than $1-e^{-x}$, for any $x\geq 0$.

\end{theorem}
\preuve
We set
$
\sqrt{N}~\left\|W_{n}(X)\right\|_{\Fa}\leq A_n+B_n
$,
with
$$
A_n=\sqrt{N}~{\left\|V_{n}(X)\right\|_{\Fa}}\quad\mbox{\rm and}\quad
B_n=\sup_{f\in\Fa}{\left|R_{n}(X)(f)\right|}
$$
Using the fact that
$$
\sup_{f\in\Fa}{\left|V_{n}(X)(f)\right|}\leq 
\sum_{p=0}^n a_p~\left\|V(X_p)\right\|_{\Fa_p}
$$
by lemma~\ref{lemme-0}, we have
$$
\pi_{\psi}\left(\left\|V_{n}(X)\right\|_{\Fa}\right)\leq \sum_{p=0}^n a_p~\pi_{\psi}\left(
\left\|V(X_p)\right\|_{\Fa_p}\right)
$$
Using theorem~\ref{theo-proc-emp}, we also have that
$$
\pi_{\psi}\left(
\left\|V(X_p)\right\|_{\Fa_p}\right)
\leq I_1(\Fa):=12^2~ \int_{0}^{2}\sqrt{\log{(8+\Na(\Fa,\epsilon)^2)}}~d\epsilon
$$
with
$$
\Na(\Fa,\epsilon)=\sup_{n\geq 0}\Na(\Fa_n,\epsilon)
$$
This implies that
$$
\pi_{\psi}\left(A_n\right)\leq \overline{a}_n~\sqrt{N}~I_1(\Fa)
\quad
\mbox{\rm with} 
\quad
\overline{a}_n:=\sum_{p=0}^n a_p
$$
By lemma~\ref{lem-Orlicz-L2m}, we have
$$
\EE\left(e^{tA_n}\right)\leq (1+t\pi_{\psi}(A_n))~e^{(t\pi_{\psi}(A_n))^2}\leq
e^{\alpha_n t+\frac{1}{2}t^2\beta_n}
$$
with $$\beta_n=2\alpha_n^2\quad\mbox{\rm and}\quad
\alpha_n=\pi_{\psi}(A_n)
$$
Notice that
$$
L_{A_n-\alpha_n}(t)\leq L_n(t):=\frac{1}{2}t^2\beta_n
$$
Recalling that
$$
L^{\star}_n(\lambda)=\frac{\lambda^2}{2\beta_n}\quad\quad\mbox{\rm and}
\quad \left(L^{\star}_n\right)^{-1}(x)=\sqrt{2\beta_nx}
$$
we conclude that
\begin{eqnarray*}
\left(L_{A_n}^{\star}\right)^{-1}(x)&=&\alpha_n+
\left(L_{A_n-\alpha_n}^{\star}\right)^{-1}(x)\\
&\leq& \alpha_n+\sqrt{2\beta_nx}=\pi_{\psi}(A_n)\left(1+2\sqrt{x}\right)
\end{eqnarray*}

On the other hand, under our assumption, we also have that
$$
\EE\left(e^{tB_n}\right)\leq \sum_{m\geq 0} (c_nt)^m=\frac{1}{1-c_nt}=e^{c_nt}\times e^{2L_0(c_nt/2)}
$$
for any $0\leq t<1/c_n$ with the convex increasing function $L_0$ introduced on page~\pageref{L_0def}, so that
$$
2L_0(c_nt/2)=-c_nt-\log{(1-c_nt)}
$$
Using lemma~\ref{lem-el-L}, we conclude that
$$
L_{B_n-c_n}(t)\leq 2L_0(c_nt/2)$$
and
\begin{eqnarray*} 
\left(L_{B_n}^{\star}\right)^{-1}(x)&=&c_n+\left(L_{B_n-c_n}^{\star}\right)^{-1}(x)\\
&\leq& c_n
\left(1+\left(L_0^{\star}\right)^{-1}\left(\frac{x}{2}\right)
\right)
\end{eqnarray*}
The end of the proof is now a direct consequence of the
Bretagnolle-Rio's lemma.
\cqfd
\subsection{Covering numbers and entropy methods}

In this final section, we derive some properties of covering numbers for some classes of functions.
These two results are the key to derive uniform concentration inequalities w.r.t. the time
parameter for Feynman-Kac particle models. This subject is investigated in chapter~\ref{fk-particle-chap}.

We let $(E_n,\Ea_n)_{n=0,1}$ be a pair of measurable state spaces,
and $\Fa$ be a separable collection  of 
measurable functions $f:E_1\rightarrow\RR$ such that $\|f\|\leq 1$ and $\mbox{\rm osc}(f)\leq 1$.

We consider a Markov transition $M(x_0,dx_1)$ from $E_0$ into $E_1$, a probability measure $\mu$ on $E_0$,
and a function $G$ from $E_0$ into $[0,1]$ . We associate with these objects the class of functions
$$
G\cdot M(\Fa)=\left\{G~M(f)~:~f\in \Fa\right\}$$
and
$$
G\cdot (M-\mu M)(\Fa)=\left\{G~\left[M(f)-\mu M(f)\right]~:~f\in \Fa\right\}
$$
\begin{lemma}
For any $\epsilon>0$, we have
$$
 \Na\left[G\cdot M(\Fa),\epsilon\right]\leq \Na(\Fa,\epsilon)
$$
\end{lemma}
\preuve
For any probability measure $\eta$ on $E_0$, we let $\left\{f_1,\ldots,f_{n_{\epsilon}}\right\}$ be the centers of $n_{\epsilon}=
\Na(\Fa,\LL_2(\eta),\epsilon)$ 

$\LL_2(\eta)$-balls of radius at most $\epsilon$ covering $\Fa$. For any $f\in \Fa$, there exists some $1\leq i\leq n_{\epsilon}$
such that
$$
\eta\left(\left[G(f-f_i)\right]^2\right)^{1/2}\leq \eta\left(\left[(f-f_i)\right]^2\right)^{1/2}\leq \epsilon
$$
 This implies that
 $$
 \Na\left(G\cdot\Fa,\LL_2(\eta),\epsilon\right)\leq \Na(\Fa,\LL_2(\eta),\epsilon)
 $$
 In much the same way, we let  $\left\{f_1,\ldots,f_{n_{\epsilon}}\right\}$  be the
  $n_{\epsilon}=\Na(\Fa,\LL_2(\eta M),\epsilon)$ centers of
$\LL_2(\eta M)$-balls of radius at most $\epsilon$ covering $\Fa$. In this situation, for any $f\in \Fa$, there exists some $1\leq i\leq n_{\epsilon}$
such that
 $$
 \eta\left(\left[(M(f)-M(f_i))\right]^2\right)^{1/2}\leq \eta M\left(\left[(f-f_i)\right]^2\right)^{1/2}\leq \epsilon
 $$ 
 This implies that
  $$
 \Na\left(M(\Fa),\LL_2(\eta),\epsilon\right)\leq \Na(\Fa,\LL_2(\eta M),\epsilon)
 $$
This ends the proof of the lemma. 
\cqfd

One of the simplest way to control the covering numbers of the second class of functions is
to assume that $M$ satisfies the following condition
 $M(x,dy)\geq \delta \nu(dy)$, for any $x\in E_0$, and for some measure $\nu$, and some $\delta\in ]0,1[$.
 Indeed, in this situation we observe that 
 $$
 M_{\delta}(x,dy)=\frac{M(x,dy)-\delta \nu(dy)}{1-\delta}
 $$
 is a Markov transition and 
 $$
 (1-\delta)~\left[M_{\delta}(f)(x)-M_{\delta}(f)(y)\right]=\left[M(f)(x)-M(f)(y)\right]
 $$
This implies that
$$
 (1-\delta)~\left[M_{\delta}(f)(x)-\mu M_{\delta}(f)\right]=\left[M(f)(x)-\mu M(f)\right]
$$
and
$$
\eta\left[\left(M(f)(x)-\mu M(f)\right)^2\right]\leq  2(1-\delta)~\eta M_{\delta,\mu} (|f|^2)
$$
with the Markov transition
$$
M_{\delta,\mu}(x,dy)=\frac{1}{2}~\left[M_{\delta}(x,dy)+\mu M_{\delta}(x,dy)\right]
$$
We let  $\left\{f_1,\ldots,f_{n_{\epsilon}}\right\}$  be the
  $n_{\epsilon}=\Na(\Fa,\LL_2\left(\eta M_{\delta,\mu}\right),\epsilon/2)$ centers of
$\LL_2\left(\eta M_{\delta,\mu}\right)$-balls of radius at most $\epsilon$ covering $\Fa$. 
If we set
$$
\overline{f}=M(f)-\mu M(f)
\quad
\mbox{\rm and}
\quad
\overline{f}_i=M(f_i)-\mu M(f_i)
$$
then we find that
$$
\overline{f}-\overline{f}_i=M(f-f_i)-\mu M(f-f_i)
$$
from which we prove that
$$
\eta\left[\left(\overline{f}-\overline{f}_i\right)^2\right]^{1/2}\leq 2(1-\delta)~\left[\eta M_{\delta,\mu} (|f-f_i|^2) \right]^{1/2}
$$
We conclude that
$$
 \Na\left((M-\mu M)(\Fa),\LL_2(\eta),2 \epsilon (1-\delta) \right)\leq \Na(\Fa,\LL_2\left(\eta M_{\delta,\mu}\right),\epsilon)
$$
and therefore
$$
 \Na\left((M-\mu M)(\Fa),2 \epsilon (1-\delta) \right)\leq \Na(\Fa,\epsilon)
 $$
or equivalently
$$
  \Na\left(\frac{1}{1-\delta}~(M-\mu M)(\Fa),\epsilon  \right)\leq \Na(\Fa,\epsilon/2)
$$
In more general situations, we quote the following result.
\begin{lemma}\label{estimation-Nepsilon}
For any $\epsilon>0$, we have
$$
 \Na\left[G\cdot (M-\mu M)(\Fa),2\epsilon \beta(M)\right]\leq \Na(\Fa,\epsilon)
$$
\end{lemma}
\preuve
We consider a Hahn-Jordan orthogonal decomposition
$$
(M(x,dy)-\mu M(dy))=M_{\mu}^+(x,dy)-M_{\mu}^-(x,dy)
$$
with $$
M_{\mu}^+(x,dy)=(M(x,\point)-\mu M)^+\quad\mbox{\rm and}\quad
M_{\mu}^-(x,dy)=(M(x,\point)-\mu M)^-
$$
with
$$
\left\|M(x,\point)-\mu M\right\|_{\rm\tiny tv}=M_{\mu}^+(x,E_1)=M_{\mu}^-(x,E_1)\leq \beta(M)
$$
By construction, we have
\begin{eqnarray*}
M(f)(x)-\mu M(f)=M_{\mu}^+(x,E_1)~\left( \overline{M}^+_{\mu}(f)(x)-\overline{M}^-_{\mu} (f)(x) \right)
\end{eqnarray*}
with
$$
\overline{M}^+_{\mu}(x,dy):=\frac{M_{\mu}^+(x,dy)}{M_{\mu}^+(x,E_1)}\quad\mbox{\rm and}
\quad
\overline{M}^-_{\mu}(x,dy):= \frac{M_{\mu}^-(x,dy)}{M_{\mu}^-(x,E_1)}
$$
This implies that
$$
\left|M(f)(x)-\mu M(f)\right|\leq 2\beta(M)~\overline{M}_{\mu} (|f|)(x) 
$$
with
$$
\overline{M}_{\mu}(x,dy)=\frac{1}{2}\left( \overline{M}^+_{\mu}(x,dy)+ \overline{M}^-_{\mu}(x,dy)\right)
$$
One concludes that
$$
\eta\left[\left(M(f)(x)-\mu M(f)\right)^2\right]^{1/2}\leq 2\beta(M)~\left[\eta \overline{M}_{\mu} (|f|^2) \right]^{1/2}
$$

We let  $\left\{f_1,\ldots,f_{n_{\epsilon}}\right\}$  be the
  $n_{\epsilon}=\Na(\Fa,\LL_2\left(\eta \overline{M}_{\mu}\right),\epsilon/2)$ centers of
$\LL_2\left(\eta \overline{M}_{\mu}\right)$-balls of radius at most $\epsilon$ covering $\Fa$. 

If we set
$$
\overline{f}=M(f)-\mu M(f)
\quad
\mbox{\rm and}
\quad
\overline{f}_i=M(f_i)-\mu M(f_i)
$$
then we find that
$$
\overline{f}-\overline{f}_i=M(f-f_i)-\mu M(f-f_i)
$$
from which we prove that
$$
\eta\left[\left(\overline{f}-\overline{f}_i\right)^2\right]^{1/2}\leq 2\beta(M)~\left[\eta \overline{M}_{\mu} (|f-f_i|^2) \right]^{1/2}
$$
In this situation, for any $\overline{f}\in (M-\mu M)(\Fa)$, there exists some $1\leq i\leq n_{\epsilon}$ such that
$$
\eta\left[\left(\overline{f}-\overline{f}_i\right)^2\right]^{1/2}\leq \beta(M)~\epsilon
$$
We conclude that
$$
 \Na\left((M-\mu M)(\Fa),\LL_2(\eta),\epsilon \beta(M)\right)\leq \Na(\Fa,\LL_2\left(\eta \overline{M}_{\mu}\right),\epsilon/2)
$$
and therefore
\begin{eqnarray*}
 \Na\left(G\cdot(M-\mu M)(\Fa),\epsilon \beta(M)\right)&\leq& \Na\left((M-\mu M)(\Fa),\epsilon \beta(M)\right)\\
 &\leq& \Na(\Fa,\epsilon/2)
\end{eqnarray*}
This ends the proof of the lemma.
\cqfd

\section{Feynman-Kac particle processes}\label{fk-particle-chap}
\subsection{Introduction}

In this chapter, we investigate the concentration properties of the collection of
Feynman-Kac particle measures introduced in section~\ref{sec-ips-fk-intro}, in terms of the contraction parameters $\tau_{k,l}(n)$ and $\overline{\tau}_{k,l}(m)$ introduced in definition~\ref{deftaukl}, and in corollary~\ref{control-tau-kappa}.

In the first section, section~\ref{first-order-BBG-ref}, 
we present some basic  first order decompositions of the Boltzmann-Gibbs transformation
associated with some regular potential function.

In section~\ref{stoch-perturbation-sec}, we combine the semigroup techniques developed in chapter~\ref{fk-sg-chap},
with a stochastic perturbation analysis to derive first order integral expansions in terms of local random fields and 
and Feynman-Kac transport operators.

In section~\ref{concentration-ips-fk}, we combine these key formulae with the concentration analysis of interacting empirical processes developed in
chapter~\ref{ips-emp-sec}. We derive quantitative concentration estimates for finite marginal models, as well as for empirical processes w.r.t. some classes of functions. The final two sections, section~\ref{section-pfe}, and section~\ref{section-backward-ref},  are devoted respectively to particle free energy models, and backward particle Markov models.

\subsection{First order expansions}\label{first-order-BBG-ref}
For any positive potential
function $G$, any measures $\mu$ and $\nu$, and any function $f$ on $E$, we have
\begin{equation}\label{ref-decomp-BBG}
\begin{array}{l}
\left[\Psi_G(\mu)-\Psi_G(\nu)\right](f)\\
\\
=\frac{1}{\mu(G_{\nu})}~(\mu-\nu)(d_{\nu}\Psi_G(f))\\
\\
=\left(1-\frac{1}{\mu(G_{\nu})}~(\mu-\nu)(G_{\nu})\right)~(\mu-\nu)(d_{\nu}\Psi_G(f))
\end{array}
\end{equation}
with the functions
\begin{equation}\label{deffdPsi}
d_{\nu}\Psi_G(f):=G_{\nu}~(f-\Psi_G(\nu)(f))\quad\mbox{\rm and}\quad
G_{\nu}:=G/\nu(G)
\end{equation}
Notice that
$$
\left|\left[\Psi_G(\mu)-\Psi_G(\nu)\right](f)\right|\leq g~\left|(\mu-\nu)(d_{\nu}\Psi_G(f))\right|
$$
and
$$
\|d_{\nu}\Psi_G(f)\|\leq g~\mbox{\rm osc}(f)\quad\mbox{\rm 
with} \quad
g:=\sup_{x,y}\left({G(x)}/{G(y)}\right)
$$ 
It is also important to observe that
\begin{eqnarray*}
\left|\left[\Psi_G(\mu)-\Psi_G(\nu)\right](f)\right|&\leq& \frac{1}{\mu(G^{\prime})}~\left|(\mu-\nu)(d_{\nu}^{\prime}\Psi_G(f))\right|\\
&\leq & g~\left|(\mu-\nu)(d_{\nu}^{\prime}\Psi_G(f))\right|
\end{eqnarray*}
with the integral operator $d_{\nu}^{\prime}\Psi_G$ from $\mbox{\rm Osc(E)}$ into itself defined by
$$
d_{\nu}^{\prime}\Psi_G(f):=G^{\prime}~(f-\Psi_G(\nu)(f))\quad\mbox{\rm and}\quad
G^{\prime}:={G}/{\|G\|}
$$
Using lemma~\ref{estimation-Nepsilon}, we readily prove the following lemma.
\begin{lemma}
We let   $\Fa$ be  separable collection of 
measurable functions $f:E^{\prime}\rightarrow\RR$ on some  possibly different state space $E^{\prime}$, 
and such that $\|f\|\leq 1$ and $\mbox{\rm osc}(f)\leq 1$.
For any Markov transition
$M$ from $E$ into $E^{\prime}$, we set
$$
d_{\nu}^{\prime}\Psi_GM(\Fa):=\left\{d_{\nu}^{\prime}\Psi_G(M(f))~:~f\in \Fa\right\}
$$
In this situation, we have the uniform estimate
\begin{equation}\label{estimation-dpsiM}
\sup_{\nu\in\Pa(E)}{ \Na\left[d_{\nu}^{\prime}\Psi_GM(\Fa),2\epsilon \beta(M)\right]}\leq \Na(\Fa,\epsilon)
\end{equation}

\end{lemma}

\subsection{A stochastic perturbation analysis}\label{stoch-perturbation-sec}
Mean field particle models can be thought as a stochastic perturbation technique
for solving nonlinear measure valued equations of the form
$$
\eta_n=\Phi_n\left(\eta_{n-1}\right)
$$
The random perturbation term is encapsulated into
the sequence of local random sampling  errors $(V_n^N)_{n\geq 0}$ given by the local perturbation equations
$$
\eta^N_n=\Phi_n\left(\eta^N_{n-1}\right)+\frac{1}{\sqrt{N}}~V^N_n
$$
One natural 
way to control the
fluctuations and the concentration properties of the particle measures
$(\eta^N_n,\gamma^N_n)$ around their limiting values $(\eta_n,\gamma_n)$
is to express the random fields $(W^{\gamma,N}_n,W^{\eta,N}_n)$
defined by
$$
\gamma^N_n=\gamma_n+\frac{1}{\sqrt{N}}~W^{\gamma,N}_n
\qquad
\eta^N_n=\eta_n+\frac{1}{\sqrt{N}}~W^{\eta,N}_n
$$
in terms of the empirical random fields $(V_n^N)_{n\geq 0}$.

As shown in (\ref{appli-xi-X}), it is important to recall that the local sampling random fields models 
$V_n^N$
belong to the class of empirical processes we analyzed in section~\ref{intro-proc-empirique}. The 
stochastic analysis developed in chapter~\ref{chapter-empirique} applies directly to these models.
For instance, using theorem~\ref{theo-kintchine}
we have the quantitative almost sure estimate of the amplitude of the stochastic perturbations
\begin{equation}
\label{reflVm-ref}
\mathbb{E}\left( \left| V^{N}_{n}(f)\right| ^{m} \left|\Ga^N_{n-1}\right. \right) ^{1/m}\leq ~b(m)
\end{equation}
 for any $m\geq1$ and any test function $f\in\mbox{\rm Osc}(E_n)$

The first order expansions presented in the further development of this section will be
expressed in terms of the random functions $d^N_{p,n}(f)$ and $G^N_{p,n}$, and the first order functions $d_{p,n}(f)$ defined below.
\begin{definition}\label{defidN-pn}
For any $0\leq p\leq n$, and any function $f$ on $E_n$, we denote by $d_{p,n}(f)$
the function on $E_p$ defined by
$$
d_{p,n}(f)=d_{\eta_p}\Psi_{G_{p,n}}(P_{p,n}(f))
$$
For any $N\geq 1$, and any $0\leq p\leq n$,
 we also denote by $G^N_{p,n}$, $d^N_{p,n}(f)$, and $d^{\prime N}_{p,n}(f)$
 the $\Ga^N_{p-1}$-measurable random functions on $E_p$ given by
$$
G^N_{p,n}:=\frac{G_{p,n}}{\Phi_{p}(\eta^N_{p-1})(G_{p,n})}
\qquad
d^N_{p,n}(f):=d_{\Phi_{p}(\eta_{p-1}^{N})}\Psi_{G_{p,n}}(P_{p,n}(f))
$$
and
$$
d^{\prime N}_{p,n}(f):=d^{\prime}_{\Phi_{p}(\eta_{p-1}^{N})}\Psi_{G_{p,n}}(P_{p,n}(f))
$$

\end{definition}

Notice that
$$
\|G^N_{p,n}\|\leq g_{p,n}\quad\mbox{\rm and}\quad
\|d_{p,n}(f)\|\vee \|d^N_{p,n}(f)\|\leq g_{p,n}~\beta(P_{p,n})
$$
as well as
$$
\left\|d^{\prime N}_{p,n}(f)\right\|\leq \beta(P_{p,n})\quad\mbox{\rm and}\quad
\mbox{\rm osc}\left(d^{\prime N}_{p,n}(f)\right)\leq 2\beta(P_{p,n})
$$

As promised, the next theorem presents some key first order decompositions which are the progenitors for our
other results. Further details
on these expansions and their use in the bias and the 
fluctuation analysis of Feynman-Kac particle models can be found
~\cite{dm96,dmsp,Delm04,dd-2012}. 

\begin{theorem}\label{theo-key-decomp-ref-1}
For any $0\leq p\leq n$, and any function $f$ on $E_n$, 
we have the decomposition 
\begin{eqnarray}
W^{\eta,N}_n(f)
&=&\sum_{p=0}^{n}
\frac{1}{\eta_{p}^{N}(G^N_{p,n})}~V^N_p(d^N_{p,n}(f))\label{martingaledec21}
\end{eqnarray}
and the $\LL_m$-mean error estimates
\begin{equation}\label{Lm-bounds}
\EE\left(\left|W^{\eta,N}_n(f)\right|^m\right)^{1/m}
\leq 2~b(m)~
\tau_{1,1}(n)
\end{equation}
with the parameter $\tau_{1,1}(n)$ defined in (\ref{def-tau}).

In addition, we have
\begin{equation}\label{heq-ok-form}
\begin{array}{l}
W^{\eta,N}_n(f)
=\displaystyle\sum_{p=0}^{n}V^N_p\left[d_{p,n}(f)\right] +\displaystyle\frac{1}{\sqrt{N}}~R^N_n(f)
\end{array}
\end{equation}
with a second order remainder term
$$
R^N_n(f):=-\sum_{p=0}^{n-1}~\displaystyle\frac{1}{\eta^N_p(\overline{G}_p)}~W^{\eta,N}_p(\overline{G}_p)~W^{\eta,N}_p\left[d_{p,n}(f)\right]
$$
such that
\begin{equation}\label{estimation-reste}
\sup_{f\in\mbox{\rm Osc}(E_n)}{\EE\left[\left|R^{N}_n(f)\right|^m\right]^{1/m}}\leq 4~b(2m)^2
\tau_{2,1}(n)
\end{equation}
with the parameter $\tau_{2,1}(n)$ defined in (\ref{def-tau}).
\end{theorem}

\proof

The proof of (\ref{martingaledec21}) is based on the telescoping sum decomposition
$$
 \eta_{n}^{N}-\eta_{n}
 =\sum_{p=0}^{n}~\left[  \Phi_{p,n}(\eta_{p}^{N})-\Phi_{p,n}\left(
\Phi_{p}(\eta_{p-1}^{N})\right)  \right] 
$$
Recalling that
$$
\Phi_{p,n}(\mu)=\Psi_{G_{p,n}}(\mu)P_{p,n}
$$
we prove that
$$
\Phi_{p,n}(\eta_{p}^{N})-\Phi_{p,n}\left(
\Phi_{p}(\eta_{p-1}^{N})\right) =\left[\Psi_{G_{p,n}}\left(\eta_{p}^{N}\right)-\Psi_{G_{p,n}}\left(\Phi_{p}(\eta_{p-1}^{N})\right)\right]P_{p,n}
$$
Using (\ref{ref-decomp-BBG}), we have
$$
\begin{array}{l}
\sqrt{N}\left[\Psi_{G_{p,n}}\left(\eta_{p}^{N}\right)-\Psi_{G_{p,n}}\left(\Phi_{p}(\eta_{p-1}^{N})\right)\right](f)
\\
\\
=V^N_p\left[d_{\Phi_{p}(\eta_{p-1}^{N})}\Psi_{G_{p,n}}(f)\right]\\
\\
\hskip3cm-\displaystyle\frac{1}{\sqrt{N}}
\frac{1}{\eta^N_p(G^N_{p,n})}~V^N_p(G^N_{p,n})~V^N_p\left[d_{\Phi_{p}(\eta_{p-1}^{N})}\Psi_{G_{p,n}}(f)\right]
\end{array}
$$
The proof of (\ref{heq-ok-form}) is based on the  telescoping sum decomposition
$$
\eta_{n}^{N}-\eta_{n}
 =\sum_{p=0}^{n}~\left[\eta_p^N\overline{Q}_{p,n}-\eta_{p-1}^N\overline{Q}_{p-1,n}\right]
$$
with the convention $\eta_{-1}^N\overline{Q}_{-1,n}=\eta_0\overline{Q}_{0,n}=\eta_n$, for $p=0$. Using the fact that
$$
\eta_{p-1}^N\overline{Q}_{p-1,n}(f)=\eta_{p-1}^N(\overline{G}_{p-1})\times \Phi_{p}\left(\eta^N_{p-1}\right)\overline{Q}_{p,n}(f)
$$
we prove that
$$
\left[\eta_{n}^{N}-\eta_{n}\right](f)=\sum_{p=0}^n \left[\eta_p^N-\Phi_{p}\left(\eta^N_{p-1}\right)\right]\overline{Q}_{p,n}(f)+R_n^N(f)
$$
with the second order remainder term
$$
R^N_n(f):=\sum_{p=1}^n\left(1-\eta^N_{p-1}(\overline{G}_{p-1})\right)\times \Phi_{p}\left(\eta^N_{p-1}\right)\overline{Q}_{p,n}(f)
$$
Replacing $f$ by the centred function $(f-\eta_n(f))$, and using the fact that
$$
1=\eta_{p-1}(\overline{G}_{p-1})\quad\mbox{\rm and}\quad
\eta_{p}\left[d_{p,n}(f-\eta_n(f))\right]=0
$$
we conclude that
$$
\left[\eta_{n}^{N}-\eta_{n}\right](f)=\sum_{p=0}^n \left[\eta_p^N-\Phi_{p}\left(\eta^N_{p-1}\right)\right](d_{p,n}(f))+\overline{R}_n^N(f)
$$
with the second order remainder term
\begin{eqnarray*}
\overline{R}^N_n(f)&:=&\sum_{p=1}^n\left[\eta_{p-1}-\eta^N_{p-1}\right](\overline{G}_{p-1})\\
&&\qquad\times \left[\Psi_{G_{p-1}}\left(\eta^N_{p-1}\right)-\Psi_{G_{p-1}}\left(\eta_{p-1}\right)\right]\left(M_{p}\left(d_{p,n}(f)\right)\right)\\
&=&-\frac{1}{N}\sum_{p=1}^nW^{\eta,N}_{p-1}(\overline{G}_{p-1})\\
&&\quad\times 
\frac{1}{\eta^N_{p-1}(\overline{G}_{p-1})}~
W^{\eta,N}_{p-1}
\left(d_{\eta_{p-1}}\Psi_{G_{p-1}}\left(M_{p}
\left(d_{p,n}(f)\right)\right)\right)
\end{eqnarray*}
Finally, we observe that
\begin{eqnarray*}
d_{\eta_{p-1}}\Psi_{G_{p-1}}\left(M_{p}
\left(d_{p,n}(f)\right)\right)&=&\frac{G_{p-1}}{\eta_{p-1}(G_{p-1})}\left(
M_{p}(d_{p,n}(f))-\eta_p(d_{p,n}(f))
\right)\\
&=&\frac{G_{p-1}}{\eta_{p-1}(G_{p-1})}~
M_{p}(d_{p,n}(f))\\
&=&\overline{Q}_{p-1,p}(d_{p,n}(f))=
\overline{Q}_{p-1,n}(f-\eta_n(f))\\
&=&d_{p-1,n}(f)
\end{eqnarray*}
This ends the proof of (\ref{heq-ok-form}).

Combining (\ref{martingaledec21}) with the almost sure estimates
$$
\EE\left(\left|V^N_p(d^N_{p,n}(f_n))\right|^m\left|\Ga^N_{p-1}\right.\right)^{1/m}\leq 2 b(m)\|G_{p,n}^N\|~\beta(P_{p,n})
$$
for any $m\geq 1$, we easily prove (\ref{Lm-bounds}).
Using the fact that
$$
\displaystyle\frac{1}{\inf_xG^N_{p,n}(x)}~
\EE\left[\left|V^N_p(G^N_{p,n})\right|^m\left|\Ga^N_{p-1}\right.\right]^{1/m}\leq 2~b(m)~g_{p,n}
$$
and
$$
\EE\left[\left|V^N_p\left[d^N_{p,n}(f)\right]\right|^m\left|\Ga^N_{p-1}\right.\right]^{1/m}\leq 2~b(m)~g_{p,n}~\beta(P_{p,n})
$$
 we prove (\ref{estimation-reste}). This ends the proof of the theorem.
\cqfd

\subsection{Concentration inequalities}\label{concentration-ips-fk}

\subsubsection{Finite marginal models}\label{finite-marg-ips-concentration}

In section~\ref{loc-sampling-sec}, dedicated to the variance analysis of the local sampling
models, we have seen that the empirical random fields $V^N_n$ satisfy the regularity  
conditions of the general 
interacting empirical process models $V(X_n)$ presented in section~\ref{intro-ips-empirique}.

 To be more precise, we have the formulae
$$
\sum_{p=0}^{n}V^N_p\left[d_{p,n}(f)\right] =\sum_{p=0}^n a_p~V^N_p\left[\delta_{p,n}(f)\right] 
$$
with the functions
$$
\delta_{p,n}(f)=d_{p,n}(f)/a_p\in \mbox{\rm Osc}(E_p)\cap\mathcal{B}_1(E_p)
$$
for  any finite constants
$$
a_p\geq 2\sup_{0\leq p\leq n}{\left(g_{p,n}\beta(P_{p,n})\right)}
$$

Therefore, if we fix the final time horizon $n$, with  some
slightly abusive notation we have the  formulae
$$
\sum_{p=0}^{n}V^N_p\left[d_{p,n}(f)\right] =\sum_{p=0}^{n}~\alpha_p~V(X_p)\left[f_p\right] 
$$
with 
$$
X_p=\xi_p=\left(\xi^i_p\right)_{1\leq i\leq N}\quad \alpha_p=\sup_{0\leq q\leq n}{a_q}=a^{\star}_n
\quad\mbox{\rm and}\quad
f_p=\delta_{p,n}(f)/a^{\star}_n
$$ 
We also notice that
\begin{eqnarray*}
\EE\left(V^N_p\left[\delta_{p,n}(f)\right]^2 \left|\Ga^N_{p-1}\right.\right)&\leq& 
\frac{1}{(a^{\star}_n)^2}~\sigma_p^2~\mbox{\rm osc}(d_{p,n}(f))^2\\
&\leq &\frac{4}{(a^{\star}_n)^2}~\sigma_p^2~\left\|d_{p,n}(f)\right\|^2\end{eqnarray*}
and therefore
$$
\EE\left(V^N_p\left[\delta_{p,n}(f)\right]^2 \left|\Ga^N_{p-1}\right.\right)\leq \frac{4}{(a^{\star}_n)^2}~\sigma_p^2~g_{p,n}^2\beta(P_{p,n})^2
$$
with the uniform local variance parameters $\sigma_p^2$ defined in (\ref{defsigma-ref}). 

This shows that the regularity
 condition stated in
(\ref{Condition-sigma}) is met by replacing the parameters $\sigma_p$ in the variance formula (\ref{Condition-sigma})  by
the constants
$2\sigma_pg_{p,n}\beta(P_{p,n})/a^{\star}_n$, with the uniform local variance parameters $\sigma_p$ defined in (\ref{defsigma-ref}).

 Using 
theorem~\ref{theo-ref-marginal-ips}, we easily prove the following exponential concentration property.
\begin{theorem}[\cite{rio-2009}]
For any $n\geq 0$, any $f\in\mbox{\rm Osc}(E_n)$, and any $N\geq 1$, the probability of the event
$$
\left[\eta^N_{n}-\eta_n\right](f)\leq \frac{4\tau_{2,1}(n)}{N}~\left(1+\left(L_0^{\star}\right)^{-1}(x)\right)+
2b_n~\overline{\sigma}^2_n~\left(L_{1}^{\star}\right)^{-1}\left(\frac{x}{N\overline{\sigma}^2_n}\right)
$$
is greater than $1-e^{-x}$, for any $x\geq 0$, with
$$
\overline{\sigma}^2_n:=
\frac{1}{b_n^2}~\sum_{0\leq p\leq n} ~g_{p,n}^2~\beta(P_{p,n})^2~\sigma_p^2
$$
 for 
  any choice of $
b_n\geq \kappa(n)
$.
In the above display,  $\tau_{2,1}(n)$ and $\kappa(n)$ stands for the parameter defined in (\ref{def-tau}), and
 $\sigma_n$ is the uniform local variance parameter defined in (\ref{defsigma-ref}).
\end{theorem}

We illustrate the impact of this theorem with two applications. The first one is concerned with
regular and stable Feynman-Kac models satisfying the regularity conditions presented 
in section~\ref{quantitative-contraction}. The second one is concerned with 
the concentration properties of the genealogical tree based models developed in section~\ref{genealogy-sec}.

In the first situation, combining corollary~\ref{control-tau-kappa}, with the estimates (\ref{borne-L0star-inv}) and (\ref{borne-L1star-inv})
we prove the following uniform concentration inequalities w.r.t. the time horizon.

\begin{cor}\label{cor-eta-1-ref}
We assume that one of the regularity conditions ${\bf H_m(G,M)}$ stated in section~\ref{regularity-conditions}
is met for some $m\geq 0$, and we set  $$
p_m(x)=
4\overline{\tau}_{2,1}(m)~\left(1+2(x+\sqrt{x})\right)+\frac{2}{3}~\overline{\kappa}(m)~x
$$
and
$$
q_m(x)=
\sqrt{8\sigma^2\overline{\tau}_{2,2}(m)~x}
\quad\mbox{\rm with}\quad
\sigma^2=\sup_{n\geq 0}{\sigma_n^2}$$ 
In the above displayed formula, $\overline{\tau}_{2,2}(m)$ and $\overline{\kappa}(m)$ stands for the parameters  defined in corollary~\ref{control-tau-kappa}, and  $\sigma_n$ is the uniform local variance parameter defined in (\ref{defsigma-ref}).

In this situation, 
for any $n\geq 0$, any $f\in\mbox{\rm Osc}(E_n)$,  any $N\geq 1$,  and for any $x\geq 0$, the probability of the event
$$
\left[\eta^N_{n}-\eta_n\right](f)\leq \frac{1}{N}~p_m(x)+\frac{1}{\sqrt{N}}~q_m(x)
$$
is greater than $1-e^{-x}$.
\end{cor}

In the same vein, using the estimates (\ref{borne-L0star-inv}) and (\ref{borne-L1star-inv}),
concentration inequalities for genealogical tree models can be derived easily using the estimates
(\ref{gbeta-tau-historical}). 

\begin{cor}\label{ref-cor-genealogical-tree-intro}
We let $\eta^N_n$ be the occupation measure of the genealogical tree model presented in
(\ref{empirique-arbre}). We also set $\sigma^2=\sup_{n\geq 0}{\sigma_n^2}$, the supremum of is the uniform local variance parameters $\sigma_n^2$  defined in (\ref{defsigma-ref}), and
$$
p_{n,m}(x)=4(\chi_mg^m)^2\left(1+2(x+\sqrt{x})\right)
+\frac{2}{3}~\frac{\chi_mg^m}{(n+1)}~x$$
and
$$ q_m(x)=(\chi_mg^m) 
~\sqrt{8\sigma^2 x}
$$
In this situation, for any $n\geq 0$, any ${\bf f_n}\in\mbox{\rm Osc}({\bf E_n})$, and any $N\geq 1$, the probability of the event
$$
\left[\eta^N_{n}-\QQ_n\right]({\bf f_n})\leq \frac{n+1}{N}~p_{n,m}(x)+\sqrt{\frac{n+1}{N}}~q_m(x)
$$
is greater than $1-e^{-x}$, for any $x\geq 0$.

\end{cor}

\subsubsection{Empirical processes}\label{interacting-ips-sec}

The main aim of this section is to derive concentration inequalities for particle empirical processes. Several consequences of this general theorem are also discussed, including uniform estimates w.r.t. the time parameter, and concentration properties of genealogical particle processes.

\begin{theorem}
We let  $\Fa_n$ be a separable collection  of 
measurable functions $f_n$ on $E_n$, such that $\|f_n\|\leq 1$, $\mbox{\rm osc}(f_n)\leq 1$, with finite entropy $I(\Fa_n)<\infty$.
$$
\pi_{\psi}\left(\left\|W^{\eta,N}_n\right\|_{\Fa_n}\right)\leq c_{\Fa_n}~\tau_{1,1}(n)
$$
with the parameter $\tau_{1,1}(n)$ defined in (\ref{def-tau}) and
\begin{equation}\label{def-cFa}
 c_{\Fa_n}\leq 24^2~ \int_{0}^{1}\sqrt{\log{(8+\Na(\Fa_n,\epsilon)^2)}}~d\epsilon
\end{equation}
In particular, for any $n\geq 0$,  and any $N\geq 1$,  the probability of the following event
$$
\sup_{f\in\Fa_n}{\left|\eta^N_n(f)-\eta_n(f)\right|}\leq \frac{c_{\Fa_n}}{\sqrt{N}}~\tau_{1,1}(n)~\sqrt{x+\log{2}}
$$
is greater than $1-e^{-x}$, for any $x\geq 0$.
\end{theorem}
\preuve
Using (\ref{martingaledec21}), for any function $f_n\in\mbox{\rm Osc}(E_n)$ we have
the estimate
$$
\left|W^{\eta,N}_n(f_n)\right|
\leq 2~\sum_{p=0}^{n} g_{p,n}\beta(P_{p,n})
\left|V^N_p(\delta_{p,n}^N(f_n))\right|
$$
with the $\Ga^N_{p-1}$-measurable random functions $\delta_{p,n}^N(f_n)$ on $E_p$ defined by
$$
\delta_{p,n}^N(f_n)=\frac{1}{2\beta(P_{p,n})}~d^{\prime N}_{p,n}(f_n)
$$
By construction, we have
 $$
\left\|\delta_{p,n}^N(f_n)\right\|\leq 1/2\quad\mbox{\rm and}\quad  \delta_{p,n}^N(f_n)\in \mbox{\rm Osc}(E_p)
$$
 Using the uniform estimate (\ref{estimation-dpsiM}),
if we set
$$
\Ga^N_{p,n}:=\delta_{p,n}^N(\Fa_n)=\left\{  \delta_{p,n}^N(f)~:~f\in\Fa_n\right\}
$$
then we also prove the almost sure upper bound
$$
\sup_{N\geq 1}{ \Na\left[\Ga^N_{p,n},\epsilon\right]}\leq \Na(\Fa_n,\epsilon/2)
$$
The end of the proof is now a direct consequence of theorem~\ref{theo-class-random}.
This ends the proof of the theorem.
\cqfd

\begin{cor}\label{cor-empirique-ips-1}
We consider time homogeneous Feynman-Kac models on some common measurable state space $E_n=E$.
We also let  $\Fa$ be a separable collection  of 
measurable functions $f$ on $E$, such that $\|f\|\leq 1$, $\mbox{\rm osc}(f)\leq 1$, with finite entropy $I(\Fa)<\infty$.

We also assume that one of the regularity conditions ${\bf H_m(G,M)}$ stated in section~\ref{regularity-conditions}
is met for some $m\geq 0$.
In this situation, for any $n\geq 0$,  and any $N\geq 1$,  the probability of the following event
$$
\sup_{f\in\Fa_n}{\left|\eta^N_n(f)-\eta_n(f)\right|}\leq \frac{c_{\Fa}}{\sqrt{N}}~\overline{\tau}_{1,1}(m)~\sqrt{x+\log{2}}
$$
is greater than $1-e^{-x}$, for any $x\geq 0$.
\end{cor}

In the same vein,  using the estimates
(\ref{gbeta-tau-historical}), we easily prove the following corollary.

\begin{cor}\label{cor-2-gen-emp-ips}
We also assume that one of the regularity conditions ${\bf H_m(G,M)}$ stated in section~\ref{regularity-conditions}
is met for some $m\geq 0$. 

We let  $\Fa_n$ be a separable collection  of 
measurable functions ${\bf f_n}$ on the path space ${\bf E_n}$, such that $\|{\bf f_n}\|\leq 1$, $\mbox{\rm osc}({\bf f_n})\leq 1$, with finite entropy $I(\Fa_n)<\infty$.

We also let $\eta^N_n$ be the occupation measure of the genealogical tree model presented in
(\ref{empirique-arbre}). In this situation, for any $n\geq 0$,  and any $N\geq 1$,  the probability of the following event
$$
\sup_{{\bf f_n}\in\Fa_n}{\left|\eta^N_n({\bf f_n})-\QQ_n({\bf f_n})\right|}\leq \frac{c_{\Fa_n}}{\sqrt{N}}~(n+1)~\chi_mg^m~\sqrt{x+\log{2}}
$$
is greater than $1-e^{-x}$, for any $x\geq 0$ with the constant $c_{\Fa_n}$ defined in (\ref{def-cFa}).

\end{cor}

The following corollaries are a direct consequence of (\ref{d-cells-analysis}).

\begin{cor}\label{cor-cells-cv}
We assume that the conditions stated in corollary~\ref{cor-empirique-ips-1} are satisfied.
When $\Fa$ stands for the indicator functions (\ref{d-cells}) of cells in $E=\RR^d$,  for some $d\geq 1$,
 the probability of the following event
$$
\sup_{f\in\Fa}{\left|\eta^N_n(f)-\eta_n(f)\right|}\leq c~~\overline{\tau}_{1,1}(m)~\sqrt{\frac{d}{N}~(x+1)}
$$
is greater than $1-e^{-x}$, for any $x\geq 0$, for some universal constant $c<\infty$ that doesn't depend on the dimension. 

\end{cor}

\begin{cor}\label{cor-cells-cv-gen}
We assume that the conditions stated in corollary~\ref{cor-2-gen-emp-ips} are satisfied.
When $\Fa_n$ stands for product functions of  indicator of cells (\ref{d-cells})  in the path space
${\bf E_n}=\left(\RR^{d_0}\times\ldots,\times\RR^{d_n}\right)$,  for some $d_p\geq 1$, $p\geq 0$,
 the probability of the following event
$$
\sup_{{\bf f_n}\in\Fa_n}{\left|\eta^N_n({\bf f_n})-\QQ_n({\bf f_n})\right|}\leq c~(n+1)~\chi_mg^m~\sqrt{\frac{\sum_{0\leq p\leq n}d_p}{N}~(x+1)}
$$
is greater than $1-e^{-x}$, for any $x\geq 0$, for some universal constant $c<\infty$ that doesn't depend on the dimension. 

\end{cor}

\subsection{Particle free energy models}\label{section-pfe}

\subsubsection{introduction}
The main aim of this section is to analyze the concentration properties of the particle free energy models
introduced in section~\ref{sec-intro-pfree}.   More formally, the unnormalized particle random field
models discussed in this section are defined below. 

 \begin{definition}
 We denote by $\overline{\gamma}^N_n$  
 the normalized models defined by the following formulae
 $$
 \overline{\gamma}^N_n(f)={\gamma^N_n(f)}/{\gamma_n(\un)}=\eta^N_n(f)~
 \prod_{0\leq p<n}\eta_p^N(\overline{G}_p)
 $$
 with the normalized potential functions $$\overline{G}_n:=G_n/\eta_n(G_n)$$
 We also let ${W}^{\gamma,N}_n$ and $\overline{W}^{\gamma,N}_n$ be the random field particle models
 defined by
 $$
 {W}^{\gamma,N}_n=\sqrt{N}~\left[\gamma^N_n-\gamma_n\right]
 \quad\mbox{\rm and}\quad \overline{W}^{\gamma,N}_n:= {W}^{\gamma,N}_n(f)/\gamma_n(\un)
 $$
 \end{definition}
 
  These unnormalized particle 
models $\gamma^N_n$  have a particularly simple form. They are defined in terms
of product of empirical mean values $\eta^N_p(G_p)$ of the potential functions $G_p$ w.r.t. the flow of normalized particle
 measures $\eta^N_p$ after the $p$-th mutation stages, with $p<n$.

 Thus, the concentration properties of $\gamma^N_n$ should be related in some way to the one of the interacting processes $\eta^N_n$
 developed in section~\ref{concentration-ips-fk}. 
 
To begin with, we mention that
  $$
\overline{\gamma}^N_n(\un):={\gamma^N_n(\un)}/{\gamma_n(\un)}=
\prod_{0\leq p<n}\eta_p^N(\overline{G}_p)=1+\frac{1}{\sqrt{N}}~\overline{W}^{\gamma,N}_n(\un)
$$
 For more general functions we also observe that
for any function $f$ on $E_n$, s.t.
$\eta_n(f)=1$, we have the decompositions
$$
\begin{array}{l}
\overline{W}^{\gamma,N}_n(f)\\
\\
=\sqrt{N}~
\left[
\left(1+\frac{1}{\sqrt{N}}~\overline{W}^{\gamma,N}_n(\un)\right)\left(\eta_n(f)+\frac{1}{\sqrt{N}}~W^{\eta,N}_n(f)\right)-\eta_n(f)
\right]\\
\\
=\sqrt{N}~
\left[
\left(1+\frac{1}{\sqrt{N}}~\overline{W}^{\gamma,N}_n(\un)\right)\left(1+\frac{1}{\sqrt{N}}~W^{\eta,N}_n(f)\right)-1
\right]
\end{array}
$$
We readily  deduce the following second order decompositions of the fluctuation errors
$$
\overline{W}^{\gamma,N}_n(f)=\left[\overline{W}^{\gamma,N}_n(\un)+W^{\eta,N}_n(f)\right]+\frac{1}{\sqrt{N}}~\left(\overline{W}^{\gamma,N}_n(\un)~W^{\eta,N}_n(f)\right)
$$

This decomposition allows to reduce the concentration properties of $\overline{W}^{\gamma,N}_n(f)$
to the ones of $W^{\eta,N}_n(f)$ and $\overline{W}^{\gamma,N}_n(\un)$.

 In the first part of this section, we provide some
key decompositions of $ {W}^{\gamma,N}_n$ in terms of the local sampling errors $V^N_n$, as well as a pivotal exponential formula connecting the fluctuations of
 the particle free energies in terms of the fluctuations of the  potential empirical mean values.
 
 In the second part of the section, we derive first order expansions, and logarithmic concentration inequalities for particle free energy ratios $\overline{\gamma}^N_n(\un)={\gamma}^N_n(\un)/{\gamma}_n(\un)$.

 \subsubsection{Some key decomposition formulae}
 
 This section is mainly concerned with the proof of the following decomposition theorem.
  
\begin{theorem}\label{theo-key-decomp-ref-gamma}
For any $0\leq p\leq n$, and any function $f$ on $E_n$, 
we have the decompositions 
\begin{eqnarray}
W^{\gamma,N}_n(f)&=&\sum_{p=0}^{n}
\gamma^N_{p}(\un)~V^N_p(Q_{p,n}(f))  \label{martingaledec1}\\
\overline{W}^{\gamma,N}_n(f)&=&\sum_{p=0}^{n}
\overline{\gamma}^N_{p}(\un)~V^N_p(\overline{Q}_{p,n}(f))  \label{martingaledec2}
\end{eqnarray}
with the normalized Feynman-Kac semigroup
$$
\overline{Q}_{p,n}(f)=
{Q_{p,n}(f)}/{\eta_pQ_{p,n}(\un)}
$$
In addition, we have the exponential formulae
\begin{equation}
\begin{array}{l}
\overline{W}^{\gamma,N}_n(\un)\\
\\
=
\displaystyle\sqrt{N}\left(
\exp{
\left\{
\frac{1}{\sqrt{N}}
\displaystyle\int_0^1
\sum_{0\leq p<n}\frac{W^{\eta,N}_p(\overline{G}_p)}{1+\frac{t}{\sqrt{N}}~W^{\eta,N}_p(\overline{G}_p)}
~dt\right\}
}-1
\right)
 \label{log-form}
 \end{array}
 \end{equation}
 
 \end{theorem}
\preuve
We use the telescoping sum decomposition
$$
\gamma_{n}^N-\gamma_{n}=\sum_{p=0}^n\left(\gamma^N_pQ_{p,n}-\gamma^N_{p-1}Q_{p-1,n}\right)
$$
with the conventions $Q_{n,n}=Id$, for $p=n$; and $\gamma^N_{-1}Q_{-1,n}=\gamma_0Q_{0,n}$, for $p=0$.
Using the fact that
$$
\gamma_p^N(\un)=\gamma^N_{p-1}(G_{p-1})\quad\mbox{\rm and}
\quad
\gamma^N_{p-1}Q_{p-1,n}(f)=\gamma^N_{p-1}(G_{p-1}M_p(Q_{p,n}(f))
$$
we prove that
$$
\gamma^N_{p-1}Q_{p-1,n}=\gamma_p^N(\un)~\Phi_{p}\left(\eta^N_{p-1}\right)Q_{p,n}
$$
The end of the proof of the first decomposition is now easily completed.  We prove
 (\ref{martingaledec2}) using the following formulae
\begin{eqnarray*}
\overline{Q}_{p,n}(f)(x)&=&\frac{\gamma_p(\un)}{\gamma_n(\un)}~Q_{p,n}(f)(x)\\
&=&
Q_{p,n}(f)(x)~\prod_{p\leq q<n}\eta_q(G_q)^{-1}=
\frac{Q_{p,n}(f)(x)}{\eta_pQ_{p,n}(\un)}
\end{eqnarray*}

The proof of (\ref{log-form}) is based on the fact that
$$
\log{y}-\log{x}=\int_0^1 \frac{(y-x)}{x+t(y-x)}~dt
$$
for any positive numbers $x,y$. Indeed, we have the formula
\begin{eqnarray*}
\log{\left(\gamma^N_n(\un)/\gamma_n(\un)\right)}&=&\log{\left(1+\frac{1}{\sqrt{N}} \frac{W^{\gamma,N}_n(\un)}{\gamma_n(\un)}\right)}\\
&=&\sum_{0\leq p<n}\left(
\log{\eta^N_p(G_p)}-\log{\eta_p(G_p)}\right)\\
&=&\frac{1}{\sqrt{N}}\sum_{0\leq p<n}\displaystyle\int_0^1 \frac{W^{\eta,N}_p(G_p)}{\eta_p(G_p)+\frac{t}{\sqrt{N}}~W^{\eta,N}_p(G_p)}~dt
\end{eqnarray*}
This ends the proof of the theorem.
\cqfd
 \subsubsection{Concentration inequalities}

 Combining the exponential formulae (\ref{log-form}) with the expansions
  (\ref{martingaledec2}), we derive first order decompositions for the random
  sequence $ \sqrt{N}~\log{\overline{\gamma}^N_n(\un)}$. These expansions will be expressed in terms
  of the random predictable functions defined below.

\begin{definition}

We let $h^N_{q,n}$ be the random $\Ga^N_{q-1}$-measurable functions 
given by
$$
h^N_{q,n}:=\sum_{q\leq p<n}d^N_{q,p}(\overline{G}_p)
$$
with the functions $d^N_{q,p}(\overline{G}_p)$ given in definition~\ref{defidN-pn}.
\end{definition}

\begin{lemma}\label{lem-ref-gamma-conc}
For any $n\geq 0$ and any $N\geq 1$, we have
\begin{equation}\label{decomp-log-gamma}
 \sqrt{N}~\log{\overline{\gamma}^N_n(\un)}=\sum_{0\leq q<n}
V^N_q\left(h^N_{q,n}\right)+\frac{1}{\sqrt{N}}~R^{N}_n
\end{equation}
with   a second remainder order term $R^{N}_n$ such that
$$
\EE\left(\left|R^{N}_n\right|^m\right)^{1/m}\leq b(2m)^2 r(n)
$$
for any $m\geq 1$, 
with  some constant 
$$
r(n)\leq 8\sum_{0\leq p<n}~g_p~
\left(
2g_p~\tau_{1,1}(p)^2+\tau_{3,1}(p)
\right)
$$
\end{lemma}

\preuve

 Using the exponential formulae  (\ref{log-form}), we have
 \begin{eqnarray*}
 \sqrt{N}~\log{\overline{\gamma}^N_n(\un)}&=& \sqrt{N} \log{\left(1+\frac{1}{\sqrt{N}}\overline{W}^{\gamma,N}_n(\un)\right)}\\&=&
\sum_{0\leq p<n}\int_0^1
\frac{W^{\eta,N}_p(\overline{G}_p)}{1+\frac{t}{\sqrt{N}}~W^{\eta,N}_p(\overline{G}_p)}
~dt
\end{eqnarray*}
This implies that
$$
 \sqrt{N}~\log{\overline{\gamma}^N_n(\un)}=\sum_{0\leq p<n}W^{\eta,N}_p(\overline{G}_p)+\frac{1}{\sqrt{N}}~R^{N,1}_n
$$
with the (negative) second order remainder term
$$
R^{N,1}_n=-\sum_{0\leq p<n}~\int_0^1t~
\frac{W^{\eta,N}_p(\overline{G}_p)^2}{1+\frac{t}{\sqrt{N}}~W^{\eta,N}_p(\overline{G}_p)}
~dt
$$
On the other hand, using (\ref{martingaledec2}) we have
$$
\sum_{0\leq p<n}W^{\eta,N}_p(\overline{G}_p)=\sum_{0\leq q<n}
V^N_q\left(h^N_{q,n}\right)+\displaystyle\frac{1}{\sqrt{N}}~
R^{N,2}_n
$$
with the second order remainder order term
$$
R^{N,2}_n:=-\sum_{0\leq q\leq p<n}\displaystyle\frac{1}{\eta^N_q(G^N_{q,p})}~V^N_q(G^N_{q,p})~V^N_q\left[d^N_{q,p}(\overline{G}_p)\right]
$$
This gives the decomposition (\ref{decomp-log-gamma}),
with  the second remainder order term
$$
R^N_n:=R^{N,1}_n+R^{N,2}_n
$$
Using the fact that
$$
1+\frac{t}{\sqrt{N}}~W^{\eta,N}_p(\overline{G}_p)=t~\eta^N_p(\overline{G}_p)+(1-t)\geq t~g^-_p
$$
for any $t\in ]0,1]$, with $g^-_p:= \inf_x\overline{G}_p(x)$, we find that
$$
\left|R^{N,1}_n\right|\leq \sum_{0\leq p<n}~\frac{1}{g^-_p~}
W^{\eta,N}_p(\overline{G}_p)^2$$
Using (\ref{Lm-bounds}), we prove  that
$$
\EE\left(\left|R^{N,1}_n\right|^r\right)^{1/r}\leq 
 4 b(2r)^2\sum_{0\leq p<n}~\frac{1}{g^-_p~}\mbox{\rm osc}\left(\overline{G}_p\right)^2~\tau_{1,1}(p)^2
$$
from which we conclude that
$$
\EE\left(\left|R^{N,1}_n\right|^r\right)^{1/r}\leq 
 (4 b(2r))^2\sum_{0\leq p<n}~g_p^2~\tau_{1,1}(p)^2
$$
In much the same way, we have
$$
\begin{array}{l}
\EE\left(\left|R^{N,2}_n\right|^r\right)^{1/r}\\
\\
\leq 
\sum_{0\leq q\leq p<n}g_{q,p} ~
\EE\left(\left|V^N_q(G^N_{q,p})\right|^{2r}\right)^{1/2r}~\EE\left(\left|V^N_q\left[d^N_{q,p}(\overline{G}_p)\right]
\right|^{2r}\right)^{1/2r}
\end{array}$$
and using (\ref{reflVm-ref}), we prove that
$$
\begin{array}{l}
\EE\left(\left|R^{N,2}_n\right|^r\right)^{1/r}
\leq 8b(2r)^2 
\sum_{0\leq q\leq p<n}~g_p~g_{q,p}^3 ~ \beta(P_{q,p})
\end{array}
$$
This ends the proof of the lemma.
\cqfd

We are now in position to state and to prove the following concentration theorem.

\begin{theorem}\label{theo-concentration-gamma}
 For any $N\geq 1$, $\epsilon\in \{+1,-1\}$, $n\geq 0$, and for any
 $$
\varsigma^{\star}_n\geq \sup_{0\leq q\leq n}{\varsigma_{q,n}}\quad\mbox{\rm with}\quad
\varsigma_{q,n}:=\frac{4}{n}\sum_{q\leq p<n}g_{q,p}g_p~\beta(P_{q,p})
$$
 the probability of the following events
$$
\frac{\epsilon}{n}\log{\overline{\gamma}^N_n(\un)}
\leq \frac{1}{N}~\overline{r}(n)~\left(1+\left(L_0^{\star}\right)^{-1}(x)\right)+
\varsigma^{\star}_n~\overline{\sigma}^2_n~\left(L_{1}^{\star}\right)^{-1}\left(\frac{x}{N\overline{\sigma}^2_n}\right)
$$
is greater than $1-e^{-x}$, for any $x\geq 0$, with the parameters
$$
\overline{\sigma}^2_n:=\sum_{0\leq q<n}\sigma_q^2~(\varsigma_{q,n}/\varsigma^{\star}_n)^2\quad\mbox{\rm and}\quad\overline{r}(n)=r(n)/n
$$

\end{theorem}

Before getting into the  proof of the theorem, we present simple arguments 
to  derive exponential concentration inequalities for the quantities $\left| \overline{\gamma}^N_n(\un)-1\right|$.
Suppose that for any $\epsilon\in \{+1,-1\}$, the probability of events
$$
\frac{\epsilon}{n} \log{\overline{\gamma}^N_n(\un)} \leq \rho_n^N(x)
$$
is greater than $1-e^{-x}$, for any $x\geq 0$, for some function $\rho_n^N$ such that
$$\rho_n^N(x)\rightarrow_{N\rightarrow\infty}0$$
 In this case, 
the probability of event
$$
-\left(1-e^{-n\rho_n^N(x)}\right)\leq \overline{\gamma}^N_n(\un)-1\leq e^{n\rho_n^N(x)}-1
$$
is greater than $1-2e^{-x}$, for any $x\geq 0$. Choosing $N$ large enough so that $\rho_n^N(x)\leq 1/n$ we 
have 
$$
-2n\rho_n^N(x)\leq -\left(1-e^{-n\rho_n^N(x)}\right)\quad\mbox{\rm and}\quad e^{n\rho_n^N(x)}-1\leq 2n\rho_n^N(x)
$$
from which we conclude that
the probability of event
$$
\PP\left(\left| \overline{\gamma}^N_n(\un)-1\right|\leq 2n~\rho_n^N(x)\right)
\geq 1-2e^{-x}
$$

Now, we come to the proof of the theorem.

{\bf Proof of theorem~\ref{theo-concentration-gamma}:}

We use the same line of arguments as the ones we used in section~\ref{finite-marg-ips-concentration}.
Firstly, we observe that
\begin{eqnarray*}
\left\|h^N_{q,n}\right\|
&\leq& 
\sum_{q\leq p<n}\left\|d^N_{q,p}(\overline{G}_p)\right\|\\
&\leq& 
\sum_{q\leq p<n}g_{q,p}~\mbox{\rm osc}(P_{q,p}(\overline{G}_p))\leq 2
\sum_{q\leq p<n}g_{q,p}g_p~\beta(P_{q,p})=c_{q,n}/2
\end{eqnarray*}
and $\mbox{\rm osc}(h^N_{q,n})\leq c_{q,n}$. Now, we use the following decompositions

$$
\sum_{0\leq q<n}
V^N_q\left(h^N_{q,n}\right)=a^{\star}_{n}~\sum_{0\leq q<n}~
V^N_q\left(\delta^N_{q,n}\right)
$$
with the $\Ga^N_{q-1}$-measurable functions
$$
\delta^N_{q,n}=h^N_{q,n}/a^{\star}_{n}\in  \mbox{\rm Osc}(E_q)\cap\mathcal{B}_1(E_q)
$$
and for any constant 
$
a^{\star}_{n}\geq \sup_{0\leq q\leq n}{c_{q,n}}$.

On the other hand, we have the almost sure variance
estimate
\begin{eqnarray*}
\EE\left(V^N_q\left[\delta^N_{q,n}\right]^2 \left|\Ga^N_{q-1}\right.\right)&\leq& 
\sigma_q^2~\mbox{\rm osc}(h^N_{q,n})^2/a^{\star 2}_n\leq \sigma_q^2c_{q,n}^2/a^{\star 2}_n\end{eqnarray*}
from which we conclude that
$$
\EE\left(V^N_q\left[\delta^N_{q,n}\right]^2 \right)\leq \sigma_q^2c_{q,n}^2/a^{\star 2}_n
$$
This shows that the regularity
 condition stated in
(\ref{Condition-sigma}) is met by replacing the parameters $\sigma_q$ in the variance formula (\ref{Condition-sigma})  by
the constants
$\sigma_qc_{q,n}/a^{\star}_n$, with the uniform local variance parameters $\sigma_p$ defined in (\ref{defsigma-ref}).

The end of the proof is now a direct consequence of
theorem~\ref{theo-ref-marginal-ips}. This ends the proof of the theorem.
\cqfd
\begin{cor}\label{cor-gammaN-ref-intro}
 We assume that one of the regularity conditions ${\bf H_m(G,M)}$ stated in section~\ref{regularity-conditions}
is met for some $m\geq 0$, and we set 
$$
p_m(x):=c_1(m)~\left(1+2(x+\sqrt{x})\right)+c_2(m)~x
\quad\mbox{\rm
and}\quad
q_m(x)=c_3(m)\sqrt{x}
$$
with the parameters
\begin{eqnarray*}
 c_1(m)&=&(4g
\overline{\tau}_{1,1}(m))^2
+ 8g\overline{\tau}_{3,1}(m)\\
c_2(m)&=&4(\beta_m g^{m+1})/3\quad\mbox{\rm and}\quad c_3(m)=
4g
\sqrt{2\overline{\tau}_{2,2}(m)\sigma^2}
\end{eqnarray*}
In the above displayed formula, $\overline{\tau}_{2,2}(m)$ and $\overline{\kappa}(m)$ stands for the parameters  defined in corollary~\ref{control-tau-kappa}, and  $\sigma_n$ is the uniform local variance parameter defined in (\ref{defsigma-ref}).

 In this situation,
for any $N\geq 1$,  and any $\epsilon\in \{+1,-1\}$, the probability of each of the following events
$$
\frac{\epsilon}{n}\log{\overline{\gamma}^N_n(\un)}
\leq \frac{1}{N}~p_m(x)
+\frac{1}{\sqrt{N}}~
q_m(x)
$$
is greater than $1-e^{-x}$, for any $x\geq 0$, 

\end{cor}
\preuve

 Under condition ${\bf H_m(G,M)}$,  we have
 $$
r(n)/n\leq (4g
\overline{\tau}_{1,1}(m))^2
+ 8g\overline{\tau}_{3,1}(m)
$$
and
for any $p<n$
\begin{eqnarray*}
\varsigma_{p,n}^2&=&\left(\frac{4g}{n} \right)^2  (n-p)^2\left(\frac{1}{n-p}\sum_{p\leq q<n} g_{p,q}~\beta(P_{p,q})\right)^2
\\
&\leq& \frac{(4g)^2}{n}
\frac{(n-p)}{n}\sum_{p\leq q<n} g_{p,q}^2~\beta(P_{p,q})^2
\end{eqnarray*}
This implies that
$$
\sum_{0\leq p<n} \varsigma_{p,n}^2\leq \frac{(4g)^2}{n}
\sum_{0\leq q<n}\tau_{2,2}(q)\leq (4g)^2\overline{\tau}_{2,2}(m)
$$
In much the same way, we prove that $\varsigma^{\star}_n\leq 4\beta_m g^{m+1}
 $. The end of the proof is now a consequence of the 
 estimates (\ref{borne-L0star-inv}) and (\ref{borne-L1star-inv}). This ends the proof of the corollary.
\cqfd

\subsection{Backward particle Markov models}\label{section-backward-ref}

This section is concerned with the concentration properties of the backward Markov
particle measures defined in (\ref{backNnU}). Without further mention, we assume that the 
Markov transitions $M_n$ satisfy the regularity
condition (\ref{reg-H}), and we consider the random fields defined below.

\begin{definition}
We let $W_{n}^{\Gamma ,N}$ and $W_{n}^{\QQ ,N}$ be
random field models defined by
$$
W_{n}^{\Gamma ,N}=\sqrt{N}~\left(\Gamma^N_n-\Gamma_n\right)\quad\mbox{\rm and}\quad
W_{n}^{\QQ ,N}=\sqrt{N}~\left(\QQ^N_n-\QQ_n\right)
$$
\end{definition}

The analysis of the fluctuation  random fields of backward particle models is a little more involved than the 
one of the genealogical tree particle models. The main difficulty is to deal with the nonlinear dependency
of these backward particle Markov chain models with the flow of particle measures $\eta^N_n$. 

In section~\ref{conditioning-back-sec}, we provide some preliminary key backward conditioning principles. We also introduce some predictable integral operators involved in the first order expansions of the fluctuation random
fields discussed in section~\ref{stoch-perturbation-sec-back}. In section~\ref{part-approx-sec}, we illustrate these models in the context of additive functional models. In section~\ref{concentration-back-fin}, we put together  the semigroup techniques developed in earlier sections to derive a series of quantitative
concentration inequalities.

\subsubsection{Some preliminary conditioning principles}\label{conditioning-back-sec}
By definition of the unnormalized Feynman-Kac measures $\Gamma_n$,
we have
$$
\Gamma_n(d(x_0,\ldots,x_n))=\Gamma_p(d(x_0,\ldots,x_p))~\Gamma_{n|p}(x_p,d(x_{p+1},\ldots,x_n))
$$
with
$$
\Gamma_{n|p}(x_p,d(x_{p+1},\ldots,x_n))=~\prod_{p<q\leq n}Q_q(x_{q-1},dx_q)
$$
This implies that
$$
\QQ_n(d(x_0,\ldots,x_n))=\QQ_{n,p}(d(x_0,\ldots,x_p))\times \QQ_{n|p}(x_p,d(x_{p+1},\ldots,x_n))
$$
with the $\QQ_n$-distribution of the random states $(X_0,\ldots,X_p)$
$$
\QQ_{n,p}(d(x_0,\ldots,x_p)):=\frac{1}{\eta_p(G_{p,n})}~\QQ_p(d(x_0,\ldots,x_p))~G_{p,n}(x_p)
$$
and the  $\QQ_n$-conditional distribution of $(X_{p+1},\ldots,X_n)$ given the random state $X_p=x_p$
defined by
$$
 \QQ_{n|p}(x_p,d(x_{p+1},\ldots,x_n))=\frac{1}{\Gamma_{n|p}(\un)(x_p)}~\Gamma_{n|p}(x_p,d(x_{p+1},\ldots,x_n))
$$

Now, we discuss some backward conditioning principles. Using the backward Markov chain formulation (\ref{backward}), we have
$$
\mathbb{Q}_n(d(x_0,\ldots,x_n))=\eta_n(dx_n)~\QQ_{n|n}(x_n,d(x_0,\ldots,x_{n-1}))
$$
with the $\QQ_n$-conditional distribution of $(X_{0},\ldots,X_{n-1})$ given the terminal 
random state $X_n=x_n$
defined by the backward Markov transition
$$
\QQ_{n|n}(x_n,d(x_0,\ldots,x_{n-1})):=
\prod_{q=1}^{n}\MM_{q,%
\eta_{q-1}}(x_q,dx_{q-1})
$$

By construction, the $\QQ^N_n$-conditional distribution of $(X_{0},\ldots,X_{n-1})$ given the terminal 
random state $X_n=x_n$
is also defined by the particle backward Markov transition given by
$$
\QQ^N_{n|n}(x_n,d(x_0,\ldots,x_{n-1})):=
\prod_{q=1}^{n}\MM_{q,%
\eta^N_{q-1}}(x_q,dx_{q-1})
$$
We check this claim, using the fact that
$$
\mathbb{Q}^N_n(d(x_0,\ldots,x_n))=\eta_n^N(dx_n)~\QQ^N_{n|n}(x_n,d(x_0,\ldots,x_{n-1}))
$$

\begin{definition}
For any $0\leq p\leq n$ and $N\geq 1$, we denote by
$D^N_{p,n}$ and $L^N_{p,n}$ the $\Ga^N_{p-1}$-measurable integral operators
defined by
\begin{equation}
\begin{array}{l}
D^N_{p,n}(x_p,d(y_0,\ldots,y_n))
\\
\\
:=\QQ^N_{p|p}(x_p,d(y_0,\ldots,y_{p-1}))~\delta_{x_p}(dy_p)~
\Gamma_{n|p}(x_{p},d(y_{p+1},\ldots ,y_{n}))
\end{array}
\label{defiDpnM}
\end{equation}
and
\begin{equation}
\begin{array}{l}
L^N_{p,n}(x_p,d(y_0,\ldots,y_n))\\
\\
:=\QQ^N_{p|p}(x_p,d(y_0,\ldots,y_{p-1}))~\delta_{x_p}(dy_p)~
\QQ_{n|p}(x_{p},d(y_{p+1},\ldots ,y_{n}))\end{array}\label{defiLpn}
\end{equation}
\end{definition}
For $p\in\{0,n\}$, we use the convention 
\begin{eqnarray*}
D^N_{n,n}(x_n,d(y_0,\ldots,y_n))&=&L^N_{n,n}(x_n,d(y_0,\ldots,y_n))\\
&=&\QQ^N_{n|n}(x_n,d(y_0,\ldots,y_{n-1}))~\delta_{x_n}(dy_n)
\end{eqnarray*}
and 
\begin{eqnarray*}
D^N_{0,n}(x_0,d(y_0,\ldots,y_n))&=&\delta_{x_0}(dy_0)~\Gamma_{n|0}(x_{0},d(y_{1},\ldots ,y_{n}))\\
L^N_{0,n}(x_0,d(y_0,\ldots,y_n))&=&\delta_{x_0}(dy_0)~\QQ_{n|0}(x_{0},d(y_{1},\ldots ,y_{n}))
\end{eqnarray*}

The main reason for introducing these integral operators comes from the following
integral transport properties.

\begin{lemma}\label{lemderefDN}
For any $0\leq p\leq n$, and any $N\geq 1$, and any function 
${\bf f}_n$ on the path space ${\bf E}_{n}$, we have the almost sure formulae
\begin{equation}
\eta^N_{p}D^N_{p,n}= \eta^N_{p}(G_p)\times \Phi_{p+1}(\eta^N_p)D^N_{p+1,n}
\label{ref-decomp-k}
\end{equation}
and
\begin{eqnarray}
\frac{\eta^N_{p}D^N_{p,n}({\bf f}_n)}{\eta^N_{p}D^N_{p,n}(\un)}
&=& \Psi_{G_{p,n}}\left(\eta _{p}^{N}\right)L_{p,n}^{N}({\bf f}_n)\label{ref-decomp-k3}\\
&=&
\Psi_{G_{p+1,n}}\left(\Phi _{p+1}(\eta _{p}^{N})\right)L_{p+1,n}^{N}({\bf f}_n)
\label{ref-decomp-k2}
\end{eqnarray}
\end{lemma}
\proof

We check (\ref{ref-decomp-k}) using the fact that
\begin{equation}\label{ref-preuve-back-1}
\begin{array}{l}
\Gamma_{n|p}(x_{p},d(y_{p+1},\ldots ,y_{n}))\\
\\=
Q_{p+1}(x_p,dy_{p+1})\Gamma_{n|p+1}(x_{p+1},d(y_{p+2},\ldots ,y_{n}))
\end{array}
\end{equation}
and
\begin{equation}\label{ref-preuve-back-2}
\eta^N_p(dx_p)Q_{p+1}(x_p,dy_{p+1})=\eta^N_pQ_{p+1}(dy_{p+1})\times \MM_{p+1,\eta^N_p}(y_{p+1},dx_p)
\end{equation}
More precisely, we have
$$
\begin{array}{l}
\eta^N_p(dx_p)D^N_{p,n}(x_p,d(y_0,\ldots,y_n))
\\
\\
:=\eta^N_p(dx_p)Q_{p+1}(x_p,dy_{p+1})
\QQ^N_{p|p}(x_p,d(y_0,\ldots,y_{p-1}))\\
\\\hskip4cm\times~\delta_{x_p}(dy_p)~
~\Gamma_{n|p+1}(y_{p+1},d(y_{p+2},\ldots ,y_{n}))
\end{array}
$$
Using (\ref{ref-preuve-back-2}), this implies that
$$
\begin{array}{l}
\eta^N_p(dx_p)D^N_{p,n}(x_p,d(y_0,\ldots,y_n))
\\
\\
:=\eta^N_pQ_{p+1}(dy_{p+1})~\MM_{p+1,\eta^N_p}(y_{p+1},dx_p)
\QQ^N_{p|p}(x_p,d(y_0,\ldots,y_{p-1}))\\
\\\hskip4cm\times~\delta_{x_p}(dy_p)~
~\Gamma_{n|p+1}(y_{p+1},d(y_{p+2},\ldots ,y_{n}))
\end{array}
$$
from which we conclude that
$$
\begin{array}{l}
\eta^N_pD^N_{p,n}({\bf f_n})
\\
\\
:=\int~\eta^N_pQ_{p+1}(dy_{p+1})~
\QQ^N_{p+1|p+1}(y_{p+1},d(y_0,\ldots,y_{p}))\\
\\\hskip3cm\times~
~\Gamma_{n|p+1}(y_{p+1},d(y_{p+2},\ldots ,y_{n}))
~{\bf f_n}\left(y_0,\ldots,y_n\right)\\
\\
=\left(\eta^N_pQ_{p+1}\right)D^N_{p+1,n}({\bf f_n})
\end{array}
$$
This ends the proof of the first assertion.
Now, using (\ref{ref-decomp-k}) we have
$$
\frac{\eta^N_{p}D^N_{p,n}({\bf f_n})}{\eta^N_{p}D^N_{p,n}(\un)}=\frac{\Phi_{p+1}\left(\eta^N_{p}\right)D^N_{p+1,n}({\bf f_n})}{\Phi_{p+1}\left(\eta^N_{p}\right)D^N_{p+1,n}(\un)}
$$
Recalling that
$$
D_{p,n}^{N}(\un)=Q_{p,n}(\un)=G_{p,n}
$$
we readily prove  (\ref{ref-decomp-k3}) and  (\ref{ref-decomp-k2}). 
This ends the proof of the lemma.
\cqfd
\subsubsection{Additive functional models}\label{part-approx-sec}\index{Additive functional}

In this section we provide a brief discussion on the action of the operators $D^N_{p,n}$ and $L^N_{p,n}$ on 
additive linear functionals 
\begin{equation}\label{additivefunctions}
{\bf f_n}(x_0,\ldots,x_n)=\sum_{p=0}^nf_p(x_p)
\end{equation}
associated with some collection of functions $f_n\in\mbox{\rm Osc}(E_n)
$.

\begin{equation*}
D_{p,n}^{N}({\bf f_n})=Q_{p,n}(1)\left[~\sum_{0\leq q<p}\left[ \MM_{p,\eta
_{p-1}^{N}}\ldots \MM_{q+1,\eta _{q}^{N}}\right] (f_{q})+\sum_{p\leq q\leq
n}R^{(n)}_{p,q}(f_q)\right]
\end{equation*}
with  triangular array of  Markov transitions $R^{(n)}_{p,q}$ introduced in definition~\ref{lemmPRbeta}.
By definition of $L^N_{p,n}$, we also have that 
\begin{equation*}
L_{p,n}^{N}({\bf f_n})=\sum_{0\leq q<p}\left[ \MM_{p,\eta _{p-1}^{N}}\ldots
\MM_{q+1,\eta _{q}^{N}}\right] (f_{q})+\sum_{p\leq q\leq n}R^{(n)}_{p,q}(f_q)
\end{equation*}%
using the estimates (\ref{rpn-gpn}) , we prove the 
following upper bounds
\begin{equation*}
\begin{array}{l}
\mbox{\rm osc}(L_{p,n}^{N}({\bf f_{n})}) \\ 
\\ 
\leq \sum_{0\leq q<p}\beta \left( \MM_{p,\eta _{p-1}^{N}}\ldots \MM_{q+1,\eta
_{q}^{N}}\right) +\sum_{p\leq q\leq
n}~g_{q,n}\times
\beta\left(P_{p,q}\right)
\end{array}%
\end{equation*}%

There are many ways to control the Dobrushin operator norm
of the product of the random matrices defined in (\ref{random-matrix-ref}). For instance, we can use
the multiplicative formulae
$$
\beta \left( \MM_{p,\eta _{p-1}^{N}}\ldots \MM_{q+1,\eta
_{q}^{N}}\right)\leq \prod_{p<k\leq q}\beta \left( \MM_{k+1,\eta
_{k}^{N}}\right)
$$
One of the simplest way to proceed, is to assume that
\begin{equation}\label{contract-dualsg}
H_n(x,y)\leq \tau~H_n(x,y^{\prime})
\end{equation}
for any $x,y,y^{\prime}$, and for some finite constant $\tau<\infty$. In this situation,
we find that
$$
\MM_{k+1,\eta
_{k}^{N}}(y,dx)\leq \tau^2~ \MM_{k+1,\eta
_{k}^{N}}(y^{\prime},dx)$$
from which we conclude that
$$
\beta \left( \MM_{k+1,\eta
_{k}^{N}}\right)
\leq 1-\tau^{-2}
$$
We further assume that the condition ${\bf H_m(G,M)}$ stated in section~\ref{regularity-conditions}
is met for some $m\geq 1$. In this situation, we have
$$
\begin{array}{l}
\mbox{\rm osc}\left( L_{p,n}^N({\bf f_n})\right)\\
\\\leq 
\sum_{0\leq q<p}\left(1-\tau^{-2}\right)^{(p-q) }
+
\chi_mg^m\sum_{p\leq q\leq n}
\left(1-g^{-(m-1)}\chi_m^{-2}\right)^{\lfloor (q-p)/m\rfloor}
\end{array}$$
from which we prove the following uniform estimates
\begin{equation}\label{refdtet}
\sup_{0\leq p\leq n}\mbox{\rm osc}\left( L_{p,n}^N({\bf f_n})\right)\leq \tau^2+
m~g^{2m-1}\chi_m^{3}
\end{equation}

\subsubsection{A stochastic perturbation analysis}\label{stoch-perturbation-sec-back}

As in section~\ref{stoch-perturbation-sec}, we develop a stochastic perturbation analysis
that allows to express $W_{n}^{\Gamma ,N}$ and $W_{n}^{\QQ ,N}$
in terms of the local sampling random fields $(V^N_p)_{0\leq p\leq n}$.

These  first order expansions presented  will be
expressed in terms of the first order functions  $d_{\nu}\Psi_{G}(f)$ introduced in
(\ref{deffdPsi}), and the random $\Ga^N_{p-1}$-measurable functions
$ G^N_{p,n}$  introduced in 
definition~\ref{defidN-pn}.

\begin{definition}\label{def-dpn-bf}
For any $N\geq 1$, any $0\leq p\leq n$, and any function ${\bf f_n}$ on the path space 
${\bf E_n}$, we 
let $ {\bf d_{p,n}^{N}}({\bf f_n})$ be the  $\Ga_{p-1}^N$-measurable functions
 $$
 {\bf d_{p,n}^{N}}({\bf f_n})= d_{\Phi _{p}(\eta _{p-1}^{N})}\Psi_{G_{p,n}}\left(L_{p,n}^{N}({\bf f_n})\right)
 $$
\end{definition}

We are now in position to state and to prove the following decomposition theorem.

\begin{theorem}
\label{theo-nonb-seca} For any $0\leq p\leq n$, and any function ${\bf f_n}$ on the path
space $E^{n+1}$, we have 
\begin{equation}\label{unbias-Gamma}
\mathbb{E}\left( \Gamma _{n}^{N}({\bf f_n})\left\vert \Ga^N_p\right. \right) =\gamma _{p}^{N}\left( D_{p,n}^{N}({\bf f_n})\right)
\end{equation}
In addition, we have
\begin{eqnarray}
W_{n}^{\Gamma ,N}({\bf f_n})&=&\sum_{p=0}^{n}\gamma
_{p}^{N}(1)~V_{p}^{N}\left( D_{p,n}^{N}({\bf f_n})\right)\label{decomposition1}\\
W_{n}^{\QQ ,N}({\bf f_n})&=&\sum_{p=0}^{n}\displaystyle\frac{1}{\eta _{p}^{N}(G^N_{p,n})}
V^N_p\left({\bf d_{p,n}^{N}}({\bf f_n})\right)\label{decomposition2-empirique}
\\
&=&\sum_{p=0}^{n}
V^N_p\left({\bf d_{p,n}^{N}}({\bf f_n})\right)\nonumber\\
&&\quad-
\sum_{p=0}^{n}\displaystyle\frac{1}{\eta _{p}^{N}(G^N_{p,n})}~\frac{1}{\sqrt{N}}~V_p^N(G_{p,n}^N)\times V^N_p\left({\bf d_{p,n}^{N}}({\bf f_n})\right)\nonumber
\\&&\label{decomposition2}
\end{eqnarray}
\end{theorem}

\preuve

To prove the first assertion, we use a backward induction on the parameter $p
$. For $p=n$, the result is immediate since we have 
\begin{equation*}
\Gamma _{n}^{N}({\bf f_n})=\gamma _{n}^{N}(\un)~\eta _{n}^{N}\left(
D_{n,n}^{N}({\bf f_n})\right) 
\end{equation*}%
We suppose that the formula is valid at a given rank $p\leq n$. In this
situation, using the fact that $D_{p,n}^{N}({\bf f_n})$ is a  $\Ga^{N}_{p-1}$-measurable function,
we prove that
\begin{equation}\label{refrhss-b}
\begin{array}{l}
\mathbb{E}\left( \Gamma _{n}^{N}({\bf f_n})\left\vert ~\Ga^{N}_{p-1}\right. \right)  \\
\\
=\mathbb{E}\left( \gamma
_{p}^{N}\left( D_{p,n}^{N}({\bf f_n})\right) \left\vert ~\Ga^{N}_{p-1}\right. \right)   
=(\gamma _{p-1}^{N}Q_{p})D_{p,n}^{N}({\bf f_n})
\end{array}
\end{equation}
Applying (\ref{ref-decomp-k}), we also have that
$$
 \gamma^N_{p-1}Q_{p}D^N_{p,n}=\gamma^N_{p-1}D^N_{p-1,n}
$$
from which we conclude that the desired formula is satisfied at rank $(p-1)$.
This ends the proof of first assertion.

Now, combining lemma~\ref{lemderefDN} and (\ref{unbias-Gamma}), the proof of the second
assertion is simply based on the following decomposition 
$$
\begin{array}{l}
\left( \Gamma _{n}^{N}-\Gamma _{n}\right) ({\bf f_n}) \\
\\
=\sum_{p=0}^{n}\left[ 
\mathbb{E}\left( \Gamma _{n}^{N}({\bf f_n})\left\vert ~\Ga^N_{p}\right. \right) -\mathbb{E}\left( \Gamma _{n}^{N}({\bf f_n})\left\vert ~%
\Ga^N_{p-1}\right. \right) \right]  \\
\\
=\sum_{p=0}^{n}\gamma _{p}^{N}(1)~\left( \eta _{p}^{N}\left(
D_{p,n}^{N}({\bf f_n})\right) -\displaystyle\frac{1}{\eta _{p-1}^{N}(G_{p-1})}~\eta
_{p-1}^{N}\left( D_{p-1,n}^{N}({\bf f_n})\right) \right) 
\end{array}
$$
To prove the final decomposition, we use the fact that
\begin{equation*}
\lbrack \mathbb{Q}_{n}^{N}-\mathbb{Q}_{n}]({\bf f_n})=\sum_{0\leq p\leq n}\left( 
\frac{\eta _{p}^{N}D_{p,n}^{N}({\bf f_n})}{\eta _{p}^{N}D_{p,n}^{N}(1)}-\frac{%
\eta _{p-1}^{N}D_{p-1,n}^{N}({\bf f_n})}{\eta _{p-1}^{N}D_{p-1,n}^{N}(1)}\right) 
\end{equation*}%
with the conventions $\eta _{-1}^{N}D_{-1,n}^{N}=\eta _{0}\Gamma_{n|0}$%
, for $p=0$.  

Finally, we use (\ref{ref-decomp-k3}) and  (\ref{ref-decomp-k2}) to check that
$$
\lbrack \mathbb{Q}_{n}^{N}-\mathbb{Q}_{n}]=\sum_{0\leq p\leq n}\left( 
\Psi_{G_{p,n}}\left(\eta _{p}^{N}\right)-\Psi_{G_{p,n}}\left(\Phi _{p}(\eta _{p-1}^{N})\right)\right)L_{p,n}^{N}
$$
We end the proof using the first order expansions of Boltzmann-Gibbs transformation developed in section~\ref{ref-decomp-BBG}
$$
\begin{array}{l}
\sqrt{N}~\left(\Psi_{G_{p,n}}\left(\eta _{p}^{N}\right)-\Psi_{G_{p,n}}\left(\Phi _{p}(\eta _{p-1}^{N})\right)\right)L_{p,n}^{N}({\bf f_n})\\
\\
=\displaystyle\frac{1}{\eta _{p}^{N}(G^N_{p,n})}~V^N_p\left( {\bf d_{p,n}^{N}}({\bf f_n})\right)\\
\\
=V^N_p\left( {\bf d_{p,n}^{N}}({\bf f_n})\right)-\displaystyle\frac{1}{\eta _{p}^{N}(G^N_{p,n})}~\frac{1}{\sqrt{N}}~V_p^N(G_{p,n}^N)\times V^N_p\left( {\bf d_{p,n}^{N}}({\bf f_n})\right)
\end{array}
$$
This ends the proof of the theorem.
\cqfd
\subsubsection{Concentration inequalities}\label{concentration-back-fin}

{\bf  Finite marginal models}

Given a bounded function ${\bf f_n}$ on the path space ${\bf E_n}$,
we further assume that we have some almost sure estimate
\begin{equation}\label{def-lpn-ref}
\sup_{N\geq 1}{\mbox{\rm osc}\left( L_{p,n}^N({\bf f_n})
\right)}\leq l_{p,n}({\bf f_n})
\end{equation}
for some finite constant  $l_{p,n}({\bf f_n})~\left(\leq \|{\bf f_n}\|\right)$. For instance, for additive functionals of the form
(\ref{additivefunctions}), 
we have proved in section~\ref{part-approx-sec} the following uniform estimates
$$
\mbox{\rm osc}\left( L_{p,n}^N({\bf f_n})\right)\leq \tau^2+
m~g^{2m-1}\chi_m^{3}
$$
which are valid for any $N\geq 1$ and any $0\leq p\leq n$; as soon as
the mixing condition ${\bf H_m(G,M)}$ stated in section~\ref{regularity-conditions}
is met for some $m\geq 1$, and the regularity (\ref{contract-dualsg}) is satisfied for some finite 
$\tau$.

For any additive functional ${\bf f_n}$ of the form
(\ref{additivefunctions}), we denote by $\overline{\bf f}_n={\bf f_n}/(n+1)$ the normalized additive functional.

\begin{lemma}
For any $N\geq 1$, $n\geq 0$, and any bounded function ${\bf f_n}$ on the path space ${\bf E_n}$,
we have the first order decomposition
\begin{equation}\label{decomp-fin-ref}
W_{n}^{\QQ ,N}({\bf f_n})=\sum_{p=0}^{n}
V^N_p\left({\bf d_{p,n}^{N}}({\bf f_n})\right)+\frac{1}{\sqrt{N}}~R^N_n({\bf f_n})
\end{equation}
with a second order remain term $R^N_n({\bf f_n})$ such that
$$
\EE\left(\left|R^N_{n}({\bf f_n})\right|^m\right)^{1/m}\leq b(2m)^2~r_n({\bf f_n})
$$
for any $m\geq 1$, 
with some finite constant
\begin{equation}\label{def-rn-ref}
r_n({\bf f_n})\leq 4\sum_{0\leq p\leq n}g^2_{p,n}~l_{p,n}({\bf f_n})~
\end{equation}
\end{lemma}
\proof
Firstly, we notice that
$$
\left\| {\bf d_{p,n}^{N}}({\bf f_n})\right\|\leq \frac{\|G_{p,n}\|}{\Phi _{p}(\eta _{p-1}^{N})(G_{p,n})}~\mbox{\rm osc}\left(L_{p,n}^{N}({\bf f_n})\right)\leq  g_{p,n}~\mbox{\rm osc}\left(L_{p,n}^{N}({\bf f_n})\right)
$$
Using (\ref{decomposition2}), we find the decomposition (\ref{decomp-fin-ref})
with the second order remainder term
$$
R^N_n({\bf f_n})=-
\sum_{p=0}^{n}R^N_{p,n}({\bf f_n})
$$
with
$$
R^N_{p,n}({\bf f_n}):=\frac{1}{\eta _{p}^{N}(G^N_{p,n})}~V_p^N(G_{p,n}^N)~V^N_p\left({\bf d_{p,n}^{N}}({\bf f_n})\right)
$$
On the other hand, we have
$$
\begin{array}{l}
\EE\left(\left|R^N_{p,n}({\bf f_n})\right|^m\left|\Ga^N_{p-1}\right.\right)^{1/m}\\
\\
\leq
g_{p,n}~~
 \EE\left(\left|V_p^N(G_{p,n}/\|G_{p,n}\|)\right|^{2m}\left|\Ga^N_{p-1}\right.\right)^{1/(2m)}
\\
\\\hskip5cm\times ~~\EE\left(\left|V^N_p\left(d_{p,n}^{N}({\bf f_n})\right)\right|^{2m}\left|\Ga^N_{p-1}\right.\right)^{1/(2m)} 
 \end{array}
$$
Using (\ref{reflVm-ref}), we prove that
$$
\EE\left(\left|R^N_{p,n}({\bf f_n})\right|^m\left|\Ga^N_{p-1}\right.\right)^{1/m}\leq 4b(2m)^2~g^2_{p,n}~~l_{p,n}({\bf f_n})~
$$
The end of the proof  is now clear. This ends the proof of the lemma.
\cqfd

\begin{theorem}\label{le-theo-fin}
For any $N\geq 1$, $n\geq 0$, and any bounded function ${\bf f_n}$ on the path space ${\bf E_n}$,
the probability of the events
$$
\left[\QQ ^N_n-\QQ_n\right]({\bf f_n})\leq \frac{r_n({\bf f_n})}{N}~\left(1+\left(L_0^{\star}\right)^{-1}(x)\right)+
2b_n~\overline{\sigma}^2_n~\left(L_{1}^{\star}\right)^{-1}\left(\frac{x}{N\overline{\sigma}^2_n}\right)
$$
is greater than $1-e^{-x}$, for any $x\geq 0$, with
$$
\overline{\sigma}^2_n:=\frac{1}{b_n^2}\sum_{0\leq p\leq n}~
~g_{p,n}^2~l_{p,n}({\bf f_n})^2\sigma_p^2
$$
and for any choice of $b_n\geq \sup_{0\leq p\leq n}{g_{p,n}~l_{p,n}({\bf f_n})}$. In the above displayed
formulae, $\sigma_n$ are the uniform local variance parameters  defined in (\ref{defsigma-ref}), 
$l_{p,n}({\bf f_n})$ and $r_n({\bf f_n})$ are the parameters defined respectively in (\ref{def-lpn-ref}) and (\ref{def-rn-ref}).

\end{theorem}

Before getting into the  proof of the theorem, we present some direct consequences of these concentration inequalities
for normalized additive functional (we use the estimates (\ref{borne-L0star-inv}) and (\ref{borne-L1star-inv})).

\begin{cor}\label{cor-fin-1}
We assume that the mixing condition ${\bf H_m(G,M)}$ stated in section~\ref{regularity-conditions}
is met for some $m\geq 1$, and the regularity (\ref{contract-dualsg}) is satisfied for some finite 
$\tau$. We also suppose that the parameters  $\sigma_n$ defined in (\ref{defsigma-ref}) are uniformly bounded
$\sigma=\sup_{n\geq 0}\sigma_n<\infty$ and we set
$$
c_1(m):= 2g^m\chi_m~\left(\tau^2+
m~g^{2m-1}\chi_m^{3}\right)
\quad\mbox{\rm and}\quad
c_2(m):= 2 (g^m\chi_m)c_1(m)
$$
In this notation,
for any $N\geq 1$, $n\geq 0$, and any normalized additive functional $\overline{\bf f}_n$ 
on the path space ${\bf E_n}$,
the probability of the events
$$
\begin{array}{l}
\left[\QQ ^N_n-\QQ_n\right](\overline{\bf f}_n)\\
\\
\leq \displaystyle\frac{c_2(m)}{N}~\left(1+\left(L_0^{\star}\right)^{-1}(x)\right)+
c_1(m)\sigma^2~\left(L_{1}^{\star}\right)^{-1}\left(\displaystyle\frac{x}{N(n+1)\sigma^2}\right)
\end{array}
$$
is greater than $1-e^{-x}$, for any $x\geq 0$.
\end{cor}
\begin{cor}\label{ref-cor-final-back-intro}
We assume that the assumptions of corollary~\ref{cor-fin-1} are satisfied, 
and we set
$$
p_{m,n}(x)=c_2(m)(1+2(x+\sqrt{x}))+\frac{c_1(m)}{3(n+1)}~x$$
and
$$
q_{m,n}(x)=c_1(m)~\sqrt{\frac{2x\sigma^2}{(n+1)}}
$$
with the constants $c_1(m)$ and $c_2(m)$ defined in corollary~\ref{cor-fin-1}.

In this situation, the probability of the events
$$
\left[\QQ ^N_n-\QQ_n\right](\overline{\bf f}_n)\\
\\
\leq \displaystyle\frac{1}{N}~p_{m,n}(x)+\frac{1}{\sqrt{N}}~q_{m,n}(x)
$$
is greater than $1-e^{-x}$, for any $x\geq 0$.
\end{cor}

{\bf Proof of theorem~\ref{le-theo-fin}:}

We use the same line of arguments as the ones we used in section~\ref{finite-marg-ips-concentration}.
Firstly, we notice that
$$
\left\| {\bf d_{p,n}^{N}}({\bf f_n})\right\|\leq  g_{p,n}~l_{p,n}({\bf f_n})
$$
This yields the decompositions
$$
\sum_{p=0}^{n}
V^N_p\left({\bf d_{p,n}^{N}}({\bf f_n})\right)=\sum_{p=0}^{n} a_p
V^N_p\left({\bf \delta_{p,n}^{N}}({\bf f_n})\right)
$$
with the functions
$$
{\bf \delta_{p,n}^{N}}({\bf f_n})={\bf d_{p,n}^{N}}({\bf f_n})/a_p\in \mbox{\rm Osc}(E_p)\cap\Ba_1(E_p)
$$
and for any finite constants
$$
a_p\geq 2\sup_{0\leq p\leq n}{g_{p,n}~l_{p,n}({\bf f_n})}
$$
On the other hand, we also have that
\begin{eqnarray*}
\EE\left(V^N_p\left({\bf \delta_{p,n}^{N}}({\bf f_n})\right)^2\left|\Ga^N_{p-1}\right.\right)&\leq &\frac{4}{(a^{\star}_n)^2}~\sigma_p^2~
~g^2_{p,n}~l_{p,n}({\bf f_n})^2
\end{eqnarray*}
with $a^{\star}:=\sup_{0\leq p\leq n}{a_p}$, and the uniform local variance parameters $\sigma_p^2$ defined in (\ref{defsigma-ref}).

This shows that the regularity
 condition stated in
(\ref{Condition-sigma}) is met by replacing the parameters $\sigma_p$ in the variance formula (\ref{Condition-sigma})  by
the constants
$2\sigma_pg_{p,n}l_{p,n}({\bf f_n})
/a^{\star}_n$, with the uniform local variance parameters $\sigma_p$ defined in (\ref{defsigma-ref}).
 
 Using 
theorem~\ref{theo-ref-marginal-ips}, we easily prove the desired concentration property.
This ends the proof of the theorem.\cqfd
{\bf Empirical processes}

Using the same line of arguments as the ones we used in section~\ref{interacting-ips-sec}, we prove the following
concentration inequality.

\begin{theorem}
We let  $\Fa_n$ be a separable collection  of 
measurable functions ${\bf f_n}$ on ${\bf E_n}$, such that $\|{\bf f_n}\|\leq 1$, $\mbox{\rm osc}({\bf f_n})\leq 1$, with finite entropy $I(\Fa_n)<\infty$.
$$
\pi_{\psi}\left(\left\|W^{\QQ,N}_n\right\|_{\Fa_n}\right)\leq c_{\Fa_n}~\sum_{p=0}^{n} g_{p,n}~\|l_{p,n}\|_{\Fa_n}
$$
with the functional  ${\bf f_n}\in\Fa_n\mapsto l_{p,n}({\bf f_n})$ defined in (\ref{def-lpn-ref}) and
$$
 c_{\Fa_n}\leq 24^2~ \int_{0}^{1}\sqrt{\log{(8+\Na(\Fa_n,\epsilon)^2)}}~d\epsilon
$$
In particular, for any $n\geq 0$,  and any $N\geq 1$,  the probability of the following event
$$
\sup_{{\bf f_n}\in\Fa_n}{\left|\QQ^N_n({\bf f_n})-\QQ_n({\bf f_n})\right|}\leq \frac{c_{\Fa_n}}{\sqrt{N}}~\sum_{p=0}^{n} g_{p,n}~\|l_{p,n}\|_{\Fa_n}~\sqrt{x+\log{2}}
$$
is greater than $1-e^{-x}$, for any $x\geq 0$.
\end{theorem}

\preuve
Using (\ref{decomposition2-empirique}), for any function ${\bf f_n}\in\mbox{\rm Osc}({\bf E_n})$ we have
the estimate
$$
\left|W^{\QQ,N}_n({\bf f_n})\right|
\leq 2~\sum_{p=0}^{n} g_{p,n}~l_{p,n}({\bf f_n})
\left|V^N_p(\delta_{p,n}^N(f_n))\right|
$$
with the $\Ga^N_{p-1}$-measurable random functions $\delta_{p,n}^N({\bf f_n})$ on $E_p$ defined by
$$
 \delta_{p,n}^N({\bf f_n})=\frac{1}{2l_{p,n}({\bf f_n})}~\frac{G_{p,n}}{\|G_{p,n}\|}~\left[
L_{p,n}^N({\bf f_n})-\Psi_{G_{p,n}}\left(\Phi_p\left(\eta^N_{p-1}\right)\right)L_{p,n}^N({\bf f_n})
\right]
$$
By construction, we have
 $$
\left\|\delta_{p,n}^N({\bf f_n})\right\|\leq 1/2\quad\mbox{\rm and}\quad  \delta_{p,n}^N({\bf f_n})\in \mbox{\rm Osc}(E_p)
$$
 Using the uniform estimate (\ref{estimation-dpsiM}),
if we set
$$
\Ga^N_{p,n}:=\delta_{p,n}^N(\Fa_n)=\left\{  \delta_{p,n}^N({\bf f_n})~:~{\bf f_n}\in\Fa_n\right\}
$$
then we also prove the almost sure upper bound
$$
\sup_{N\geq 1}{ \Na\left[\Ga^N_{p,n},\epsilon\right]}\leq \Na(\Fa_n,\epsilon/2)
$$
The end of the proof is now a direct consequence of theorem~\ref{theo-class-random}, and theorem~\ref{lem-Orlicz-L2m}.
This ends the proof of the theorem.
\cqfd

\begin{definition}
We let  $\Fa=\left(\Fa_n\right)_{n\geq 0}$, be a sequence of separable collections $\Fa_n$   of 
measurable functions $f_n$ on $E_n$, such that $\|f_n\|\leq 1$, $\mbox{\rm osc}(f_n)\leq 1$, 
and finite entropy $I(\Fa_n)<\infty$. 

For any $n\geq 0$, we set
 $$
 J_{n}(\Fa):=24^2~\sup_{0\leq q\leq n}{\int_{0}^{1}\sqrt{\log{(8+\Na(\Fa_q,\epsilon)^2)}}~d\epsilon}
 $$
 
We also denote by $\Sigma_n(\Fa)$ the collection of additive functionnals 
defined by
$$
\begin{array}{l}
\Sigma_n(\Fa)
=\left\{{\bf f}_n\in \Ba({\bf E}_n)~\mbox{\rm such that}\quad \forall~ {\bf x_n}=(x_0,\ldots,x_n)\in{\bf E_n}\right.
\\
\\
\left.\hskip1cm~{\bf f_n}({\bf x_n})=\sum_{p=0}^nf_p(x_p)\quad  \mbox{\rm with}\quad
f_p\in\Fa_p,~\mbox{\rm for}~0\leq p\leq n\right\}
\end{array}
$$
and 
$$
\overline{\Sigma}_n(\Fa)
=\left\{{\bf f}_n/(n+1)~:~{\bf f}_n\in \Sigma_n(\Fa)\right\}
$$
\end{definition}

\begin{theorem}
For any $N\geq 1$, and any $n\geq 0$, we have
$$
\pi_{\psi} \left(\left\|W_{n}^{\QQ ,N}\right\|_{\Sigma_n(\Fa)}\right)\leq a_n~ J_n(\Fa)~
\sum_{p=0}^{n}~g_{p,n}
$$
 with some constant
 $$
 a_n\leq  \sum_{0\leq q<p} ~ \beta_{p,q}
+\sum_{p\leq q\leq n} \beta\left(R^{(n)}_{p,q}\right)
 $$
and for any   a collection of $[0,1]$-valued parameters $\beta_{p,q}$ such that
$$
\forall 0\leq q\leq p\qquad\sup_{N\geq 1}{\beta\left(\MM_{p,\eta _{p-1}^{N}}\ldots
\MM_{q+1,\eta _{q}^{N}}\right)}\leq \beta_{p,q}
$$
\end{theorem}
\proof
By definition of the operator $ {\bf d_{p,n}^{N}}$ given in definition~\ref{def-dpn-bf}, for additive functionals ${\bf f_n}$ of the form (\ref{additivefunctions}), we find that
\begin{eqnarray*}
 {\bf d_{p,n}^{\prime N}}({\bf f_n})&:=& \Phi _{p}(\eta _{p-1}^{N})(G_{p,n})\times {\bf d_{p,n}^{N}}({\bf f_n})\\
 &=& G_{p,n}~\left[Id
-\Psi_{G_{p,n}}\left(\Phi _{p}(\eta _{p-1}^{N})\right)\right] L_{p,n}^{N}({\bf f_n})\\
 &=&\sum_{0\leq q<p}  {\bf d_{p,q,n}^{(N,1)}}(f_q)+ \sum_{p\leq q\leq n}{\bf d_{p,q,n}^{(N,2)}}(f_q)
 \end{eqnarray*}
 with
 \begin{eqnarray*}
 {\bf d_{p,q,n}^{(N,1)}}(f_q)&=& G_{p,n}~\left[Id
-\Psi_{G_{p,n}}\left(\Phi _{p}(\eta _{p-1}^{N})\right)\right] \left( \MM_{p,q}^{(N)}(f_{q})\right)\\
 {\bf d_{p,q,n}^{(N,2)}}( f_q)&=&G_{p,n}~\left[Id
-\Psi_{G_{p,n}}\left(\Phi _{p}(\eta _{p-1}^{N})\right)\right]\left(R^{(n)}_{p,q}(f_q)\right)
 \end{eqnarray*}
and
$$
  \MM_{p,q}^{(N)}:=
\MM_{p,\eta _{p-1}^{N}}\ldots
\MM_{q+1,\eta _{q}^{N}}
$$
This implies that
 \begin{equation}\label{back-emp-pi-psi}
\begin{array}{l}
 \left|V^N_p\left({\bf d_{p,n}^{\prime N}}({\bf f_n})\right)\right|
 \\
 \\\leq 2\|G_{p,n}\|~\left(
 \sum_{0\leq q<p}  \beta_{p,q}~
 V^N_p\left(   {\bf \delta_{p,q,n}^{(N,1)}}(f_q) \right)\right.\\
 \\
 \hskip3.5cm+\left.\sum_{p\leq q\leq n} \beta\left(R^{(n)}_{p,q}\right)~
 V^N_p\left( {\bf \delta_{p,q,n}^{(N,2)}}(f_q)\right)\right)
 \end{array}
 \end{equation}
 with
 $$
 \begin{array}{l}
  {\bf \delta_{p,q,n}^{(N,1)}}(f_q)\\
  \\= \frac{1}{2
\beta\left(  \MM_{p,q}^{(N)}\right)
  }~\frac{G_{p,n}}{\|G_{p,n}\|}~G_{p,n}~\left[Id
-\Psi_{G_{p,n}}\left(\Phi _{p}(\eta _{p-1}^{N})\right)\right]
  \MM_{p,q}^{(N)}(f_{q})
  \end{array}
  $$
  and
 $$
 \begin{array}{l}
  {\bf \delta_{p,q,n}^{(N,2)}}(f_q)\\
  \\= \frac{1}{2
\beta\left( R_{p,q}^{(n)}\right)
  }~\frac{G_{p,n}}{\|G_{p,n}\|}~G_{p,n}~\left[Id
-\Psi_{G_{p,n}}\left(\Phi _{p}(\eta _{p-1}^{N})\right)\right]
  R_{p,q}^{(n)}(f_{q})
  \end{array}
  $$ 
By construction, we have that
$$
\left\| {\bf \delta_{p,q,n}^{(N,i)}}(f_q)\right\|\leq 1/2\quad\mbox{\rm and}\quad
\mbox{\rm osc}\left( {\bf \delta_{p,q,n}^{(N,i)}}(f_q)\right)\leq 1
$$  
for any $i\in\{1,2\}$. We set
$$
\Ga^{(N,i)}_{p,q,n}:=  {\bf \delta_{p,q,n}^{(N,i)}}(\Fa_q)=\left\{   {\bf \delta_{p,q,n}^{(N,i)}}(f_q)~:~ f_q\in\Fa_q\right\}
$$
Using the uniform estimate (\ref{estimation-dpsiM}),
we also prove the almost sure upper bound
$$
\sup_{N\geq 1}{ \Na\left[\Ga^{(N,i)}_{p,q,n},\epsilon\right]}\leq \Na(\Fa_q,\epsilon/2)
$$
Using theorem~\ref{theo-class-random}, we prove that
$$
\begin{array}{l}
\pi_{\psi} \left(\sup_{{\bf f}_n\in\Sigma_n(\Fa)}\left|V^N_p\left({\bf d_{p,n}^{\prime N}(f_n)}\right)\right|\right)
 \leq a_n~\|G_{p,n}\|~ J_n(\Fa)
  \end{array}
 $$
The end of the proof of the theorem is now a direct consequence of the decomposition (\ref{decomposition2-empirique}).
This ends the proof of the theorem.
\cqfd

\begin{cor}
We further assume that the condition ${\bf H_m(G,M)}$ stated in section~\ref{regularity-conditions}
is met for some $m\geq 1$,  the condition (\ref{contract-dualsg}) is satisfied for some $\tau$, and we have
$J(\Fa):=\sup_{n\geq 0}J_n(\Fa)<\infty$. In this situation, we have the uniform estimates
$$
\sup_{n\geq 0}\pi_{\psi} \left(\left\|W_{n}^{\QQ ,N}\right\|_{\overline{\Sigma}_n(\Fa)}\right)\leq c_{\Fa}(m)
$$
with some constant
$$
c_{\Fa}(m)\leq  \chi_m g^m~\left(\tau^2+
m~g^{2m-1}\chi_m^{3}
\right)~ J(\Fa)
$$
In particular, for any time horizon $n\geq 0$,  and any $N\geq 1$,  the probability of the following event
$$
\left\|\QQ^N_n-\QQ_n\right\|_{\overline{\Sigma}_n(\Fa)}\leq \frac{1}{\sqrt{N}}~c_{\Fa}(m)
~\sqrt{x+\log{2}}
$$
is greater than $1-e^{-x}$, for any $x\geq 0$.
\end{cor}

\proof
When the conditions ${\bf H_m(G,M)}$ and (\ref{contract-dualsg}) are satisfied, we proved in   section~\ref{part-approx-sec} that
$$
 \sum_{0\leq q<p} ~ \beta_{p,q}
+\sum_{p\leq q\leq n} \beta\left(R^{(n)}_{p,q}\right)\leq  \tau^2+
m~g^{2m-1}\chi_m^{3}
$$
We end the proof of the theorem, recalling that $g_{p,n}\leq \chi_m g^m$.  This ends the proof of the theorem.
\cqfd

We end this section with a direct consequence of (\ref{d-cells-analysis}).

\begin{cor}\label{cor-fin-lecture}
We further assume that the condition ${\bf H_m(G,M)}$ stated in section~\ref{regularity-conditions}
is met for some $m\geq 1$,  the condition (\ref{contract-dualsg}) is satisfied for some $\tau$.
We let $\Fa_n$ be the set of product functions of  indicator of cells  in the path space
${E_n}=\RR^{d}$,  for some $d\geq 1$, $p\geq 0$. 

In this situation, for any time horizon $n\geq 0$,  and any $N\geq 1$,  the probability of the following event
$$
\left\|\QQ^N_n-\QQ_n\right\|_{\overline{\Sigma}_n(\Fa)}\leq ~c(m)~\sqrt{\frac{d}{{N}}(x+1)}
$$
is greater than $1-e^{-x}$, for any $x\geq 0$, with some constant
$$
c(m)\leq  c~\chi_m g^m~\left(\tau^2+
m~g^{2m-1}\chi_m^{3}
\right)~
$$
In the above display, $c$ stands for some finite universal constant.
\end{cor}





\begin{thebibliography}{99}

\bibitem{adamczak}
R. Adamczak. A tail inequality for suprema of unbounded empirical processes with applications to Markov chains.
{\em EJP}, Vol. 13, no. 34, pp. 1000?1034 (2008).

\bibitem{andrieu-doucet}
C. Andrieu, A. Doucet, and R. Holenstein
Particle Markov chain Monte Carlo methods. {\em Journal Royal Statistical Society B}, vol. 72, no. 3, pp. 269-342 (2010).

\bibitem{MR2001g:82084}
N.~Bellomo and M.~Pulvirenti.
\newblock Generalized kinetic models in applied sciences.
\newblock In {\em Modeling in Applied Sciences}. Modeling and Simulation
in  Science,  Engineering, and
  Technology, 1--19. Birkh\"auser, Boston (2000).

\bibitem{MR2001a:82003}
N.~Bellomo and M.~Pulvirenti, editors.
\newblock {\em Modeling in Applied Sciences}.
\newblock Modeling and Simulation in Science, Engineering, and Technology.
  Birkh\"auser, Boston (2000).



\bibitem{benj2}
  {E. Canc\`es, B. Jourdain and T. Leli\`evre},
{Quantum Monte Carlo
simulations of fermions. A mathematical analysis of the fixed-node
approximation},
   {\em ESAIM: M2AN},
        vol. {16},
     no. {9},
     pp. {1403-1449},
     {(2006)}.
     
     

\bibitem{carmona-fouque}
R. Carmona, J.-P. Fouque, and D. Vestal: Interacting Particle Systems for the Computation of Rare Credit Portfolio Losses, Finance and Stochastics, 13(4), 2009 (p. 613-633).
     
     
\bibitem{caron1}     
Fr. Caron, P. Del Moral, A. Doucet, and M. Pace 
Particle approximations of a class of branching distribution flows arising in multi-target tracking
{\em HAL-INRIA RR-7233 (2010) [29p]}, To appear in {\em SIAM Journal on Control and Optimization} (2011). 


\bibitem{caron2}     
Fr. Caron, P. Del Moral, A. Doucet, and M. Pace.
On the Conditional Distributions of Spatial Point Processes
{\em HAL-INRIA RR-7232}. To appear in {\em Advances in Applied Probability} (2011). 



\bibitem{caron3}     
Fr. Caron, P. Del Moral, M. Pace, and B.N. Vo.
On the Stability and the Approximation of Branching Distribution Flows, with Applications to Nonlinear Multiple Target Filtering
{\em HAL-INRIA RR-7376 [50p]}. To appear in {\em Stoch. Analysis and Applications} (2011). 



\bibitem{caron4}
F. Caron, P. Del Moral, A. Tantar, E. Tantar. Simulation particulaire, Dimensionnement en conditions extr\^emes. {\em Research contract no. 2010-IFREMER-01} with IFREMER, (92p) September (2010).


\bibitem{cerou}
F. Cerou, P. Del Moral, Fr. Le Gland, and P. Lezaud.
Genealogical Models in Entrance Times Rare Event Analysis 
{\em Alea}, Latin American Journal of Probability And Mathematical Statistics (2006). 

\bibitem{cerou2}
F. Cerou, P. Del Moral, A. Guyader, F. LeGland, P. Lezaud, H. Topart
Some recent improvements to importance splitting (preliminary version) 
{\em Proceedings of 6th International Workshop on Rare Event Simulation}, Bamberg, Germany (2006).


\bibitem{cerou3}
F. Cerou, P. Del Moral, T. Furon et A. Guyader. Sequential Monte Carlo for Rare Event Estimation.
{\em Statistics and Computing}, to appear  (2011).

\bibitem{cgrv}
Fr. Cerou, A. Guyader, R. Rubinstein, R. Vaismana
Smoothed Splitting Method for Counting. {\em Stochastic models}, to appear (2011).

\bibitem{dawson}
D.A. Dawson. Stochastic evolution equations and related measure processes. {\em J. Multivariate Anal.}, vol. 3, pp. 1-52, (1975).

\bibitem{dawson2}
D.A. Dawson.  Measure-valued Markov processes. {\em \'Ecole d'\'Et\'e de Probabilit\'es de Saint-Flour 1991}.
 Lecture Notes Math. 1541, pp.1-260. Springer, Berlin (1993).

\bibitem{dawson-dm} 
D.A. Dawson, and P. Del Moral.
Large deviations for interacting processes in the strong topology.
{\em Statistical Modeling and Analysis for Complex Data Problem},
P. Duchesne and B. Rmillard Editors, pp. 179--209, Springer (2005). 

\bibitem{dawson3}
D.A. Dawson, K.J. Hochberg, V. Vinogradov. On weak convergence of branching particle systems undergoing spatial motion. {\em Stochastic analysis: Random fields and measure-valued processes}, pp.65-79. {\em Israel Math. Conf. Proc. 10. Bar-Ilan Univ., Ramat Gan.} (1996).

\bibitem{dm96}
P. Del Moral.
Non Linear Filtering: Interacting Particle Solution.
{\em Markov Processes and Related Fields}, Volume 2 Number 4, 555--580 (1996). 

\bibitem{dm98}
P. Del Moral.
Measure Valued Processes and Interacting Particle Systems. Application to Non Linear Filtering Problems.
{\em Annals of Applied Probability}, vol. 8 , no. 2, 438-495 (1998). 

\bibitem{djjp}
     {P. Del Moral, J. Jacod. and Ph. Protter},
 {The Monte-Carlo Method for filtering with discrete-time observations},
   {\em Probability Theory and Related Fields},
    vol. {120},
    pp. {346--368},
    {(2001)}.
 
 \bibitem{lezaud}  
P. Del Moral, P. Lezaud.
Branching and interacting particle interpretation of rare event probabilities. 
{\em Stochastic Hybrid Systems: Theory and Safety Critical Applications}, eds. H. Blom and J. Lygeros. Springer-Verlag, Heidelberg (2006).    
    
\bibitem{dmsp}
P. Del Moral, and L. Miclo.
 \emph{Branching and Interacting Particle Systems Approximations of Feynman-Kac Formulae with Applications to Non-Linear Filtering}.
Sminaire de Probabilits XXXIV, Ed. J. Azma and M. Emery and M. Ledoux and M. Yor, Lecture Notes in Mathematics, Springer-Verlag Berlin, Vol. 1729, 1-145 (2000). 

\bibitem{dm3}
    {P. Del Moral, L. Miclo},
    {Particle Approximations of Lyapunov Exponents Connected to Schroedinger Operators and Feynman-Kac Semigroups},
    {\em ESAIM Probability and Statistics},
    vol. {7},
    pp. {169-207},
    {(2003)}.

\bibitem{Delm04}
P. Del Moral. \emph{Feynman-Kac Formulae: Genealogical and
Interacting Particle Systems with Applications}, New York: Springer-Verlag (2004).

\bibitem{dd-2012}
P. Del Moral, A. Doucet.
{\em Sequential Monte Carlo
\& Genetic particle models. Theory and Practice}, Chapman \& Hall, Green series, Statistics and Applied Probability, to appear  (2012).

\bibitem{delmoraldoucetjasra2006}
P. Del Moral, A. Doucet and A. Jasra.
Sequential Monte Carlo samplers. \emph{J. Royal Statist. Soc.} B, vol. 68, pp.
411--436 (2006).

\bibitem{ddj-abc-2011}
P. Del Moral, A. Doucet and A. Jasra.
An Adaptive Sequential Monte Carlo Method for Approximate Bayesian Computation. {\em Research report, Imperial College} (2008).

\bibitem{soft} P. Del Moral and A. Doucet. Particle Motions in Absorbing
Medium with Hard and Soft Obstacles. \emph{Stochastic Analysis and
Applications}, vol. 22, no. 5, pp. 1175-1207 (2004).

\bibitem{dm-metropolis}
P. Del Moral, A. Doucet.
On a class of genealogical and interacting Metropolis Models. 
{\em S\'eminaire de Probabilit\'es XXXVII, Ed. J. Az\'ema and M. Emery and M. Ledoux and M. Yor, Lecture Notes in Mathematics 1832, Springer-Verlag Berlin}, pp. 415--446 (2003). 

\bibitem{sss}
P. Del Moral, A. Doucet, and S. S. Singh.
A Backward Particle Interpretation of Feynman-Kac Formulae. {\em ESAIM M2AN},
vol 44, no. 5, pp. 947--976  (sept. 2010).

\bibitem{sss2}
P. Del Moral, A. Doucet, and S. S. Singh.
Forward Smoothing Using Sequential Monte Carlo 
{\em Technical Report 638. Cambridge University Engineering Department.} (2010).

\bibitem{sss3}
P. Del Moral, A. Doucet, and S. S. Singh.
Uniform stability of a particle approximation of the optimal filter derivative.
{\em  Technical Report CUED/F-INFENG/TR 668. Cambridge University Engineering Department.} (2011)
 

\bibitem{dm-garnier1}
P. Del Moral, J. Garnier.
Genealogical Particle Analysis of Rare events.
{\em Annals of Applied Probability}, vol. 15, no. 4, 2496--2534 (2005). 

\bibitem{dm-garnier2}
P. Del Moral, J. Garnier.
Simulations of rare events in fiber optics by interacting particle systems (preliminary version) 
{\em Optics Communications}, Elsevier eds  (2006) 



\bibitem{DeGu:la}
P.~{Del Moral} and A.~Guionnet.
\newblock Large deviations for interacting particle systems. {A}pplications to
  nonlinear filtering problems.
\newblock {\em Stochastic Processes Appl.}, 78:69--95, 1998.

\bibitem{DeGu:ce}
P.~{Del Moral} and A.~Guionnet.
\newblock A central limit theorem for nonlinear filtering using interacting
  particle systems.
\newblock {\em Ann. Appl. Probab.}, 9(2):275--297  (1999).

\bibitem{dmg-cras}
P.~{Del Moral} and A.~Guionnet.
On the stability of Measure Valued Processes with Applications to filtering.
{\em C.R. Acad. Sci. Paris}, t. 329, Serie I, pp. 429-434 (1999).


\bibitem{DeGu:ont}
P.~{Del Moral} and A.~Guionnet.
\newblock On the stability of interacting processes with applications to
  filtering and genetic algorithms.
\newblock {\em Ann. Inst. Henri Poincar\'e}, 37(2):155--194 (2001).

\bibitem{liming-shulan}
P. Del Moral, S. Hu, L.M. Hu.
Moderate deviations for mean field particle models, in preparation (2011).


\bibitem{dm-jacod}
P.~{Del Moral} and J.~Jacod.
\newblock The {M}onte-{C}arlo method for filtering with discrete-time
  observations: Central limit theorems.
\newblock In {\em Numerical Methods and stochastics (Toronto, ON, 1999)},
  volume~34 of Fields Inst. Commun., pages 29--53. American Mathematical 
Society,
  Providence, RI (2002).


\bibitem{dkm}
P. Del Moral, M.A. Kouritzin, and L. Miclo.
On a class of discrete generation interacting particle systems. 
{\em Electronic Journal of Probability} , no. 16, 26 pp. (electronic). (2001). 

\bibitem{dm-ledoux}
P. Del Moral, and M. Ledoux.
On the Convergence and the Applications of Empirical Processes for Interacting Particle Systems and Nonlinear Filtering. 
{\em Journal of Theoretical Probability}, Vol. 13, No. 1, 225-257 (2000).

\bibitem{rio-2009}P.~Del Moral, and E. Rio. \newblock Concentration
inequalities for mean field particle models. Technical report HAL-INRIA
RR-6901, 2009. To appear in the Annals of Applied Probability (2011). 


\bibitem{dm-zajic}
P. Del Moral, and T. Zajic.
A note on Laplace-Varadhan's integral lemma. 
{\em Bernoulli}, vol.9, no. 1, pp. 49--65 (2003). 

\bibitem{dimasi}
   {G.B. Di Masi, M. Pratelli and W.G. Runggaldier}.
 {An
approximation for the nonlinear filtering problem with error bounds.}
   {\em Stochastics},
    vol. {14},
    no. {4},
    pp. {247-271},
    {(1985)}.
    
    
    
\bibitem{douc}
R. Douc, A. Guillin and J. Najim.
Moderate Deviations for Particle Filtering.
{\em The Annals of Applied Probability}, Vol. 15, No. 1B,  pp. 587-614 (2005).


\bibitem{doucet} A.~Doucet, N. D{e Freitas} and N.~Gordon, editors. %
\newblock {\em Sequential Monte Carlo Methods in Pratice}. %
\newblock{Statistics for engineering and Information Science}. Springer, New
York (2001).

\bibitem{dynkin}
E.B. Dynkin.
 An introduction to branching measure-valued processes. {\em CRM Monograph Series}. Vol. 6. Amer. Math. Soc., Providence (1994).

\bibitem{benjamin}
    {M. El Makrini, B. Jourdain and T. Leli\`evre},
    {Diffusion
Monte Carlo method: Numerical analysis in a simple case.},
    {\em ESAIM: M2AN},
    vol. {41},
    no. {2},
    pp. {189--213},
     {(2007)}.

\bibitem{kurtz}
S.N. Ethier, T.G. Kurtz. {\em Markov processes: Characterization and convergence}. Wiley, New York (1995).

\bibitem{everett}
C.J. Everett, S. Ulam. {\em USAEC, Los Alamos report, LADC-533, LADC-534}, declassified (1948).


\bibitem{everett2}
C.J. Everett, S. Ulam. {\em Proc. Nat. Acad. Sciences}, no. 33, vol. 403 (1948).

\bibitem{feller}
W. Feller. Diffusion processes in genetics. {\em Proc. Second Berkeley Symp. Math. Statist. Prob., University of California Press, Berkeley}, pp. 227- 246 (1951).

\bibitem{ferrari}
P. A. Ferrari, H. Kesten, S. Mart\'{\i}nez, P. Picco. Existence of quasi-stationary distributions. A renewal dynamical approach.
{\em Ann. Probab.}, vol. 23, pp. 501?521 (1995).



\bibitem{gine}
E. Gin\'e. Lectures on some aspects of the bootstrap. {\em Lecture Notes in Maths}, vol. 1665, Springer (1997).

\bibitem{gosselin}
F. Gosselin. Asymptotic behavior of absorbing Markov chains conditional on nonabsorption for applications in conservation biology. {\em Ann. Appl. Probab.}, vol. 11, pp. 261?284 (2001).

\bibitem{GrMe:sto}
C.~Graham and S.~M\'el\'eard.
\newblock Stochastic particle approximations for generalized {Boltzmann} models
  and convergence estimates.
\newblock {\em Ann. Probab.}, 25(1):115--132 (1997).

\bibitem{MR2003a:82062}
C. Graham and S. M{\'e}l{\'e}ard.
\newblock Probabilistic tools and {M}onte-{C}arlo approximations for some
  {B}oltzmann equations.
\newblock In {\em CEMRACS 1999 (Orsay)}, volume~10 of  ESAIM Proceedings, 
pages 77--126 (electronic). Soci\'et\'e de Math\'ematiques Appliqu\'ees 
et Industrielles, Paris (1999).


\bibitem{ghm2011}
A. Guyader, N. Hengartner, E. Matzner-L\o ber.
Simulation and Estimation of Extreme Quantiles and
Extreme Probabilities.
 {\em Applied Mathematics \& Optimization}, to appear (2011).


\bibitem{holland}
J.H. Holland. Adaptation in Natural and Artificial Systems.University of Michigan Press, Ann Arbor, 1975.

\bibitem{johansen}
A. M. Johansen, P. Del Moral, and A. Doucet
Sequential Monte Carlo Samplers for Rare Events (preliminary version) 
{\em Proceedings of 6th International Workshop on Rare Event Simulation}, Bamberg, Germany (2006). 


\bibitem{johnson}
W.B. Johnson, G. Schechtman, J. Zinn.
Best Constants in Moment Inequalities for Linear Combinations of Independent and Exchangeable Random Variables. {\em Ann. Probab.},  vol. 13, no. 1, pp. 234-253 (1985).

\bibitem{klein}
Th. Klein, and E. Rio. Concentration around the mean for maxima of empirical processes. 
{\em Ann. Probab.}, vol. 33 , no. 3, 1060?1077 (2005).

\bibitem{kac}
   {M. Kac},
 {On distributions of certain Wiener functionals},
    {\em Trans. American Math. Soc.},
        vol. {61},
     no. {1},
     pp. {1-13},
    {(1949)}.


 \bibitem{korez2}
    {H. Korezlioglu and W.J. Runggaldier}.
    { Filtering for
nonlinear systems driven by nonwhite noises : An approximating scheme.},
    {\em Stochastics and Stochastics Rep.},
    vol. {44},
    no. {1-2},
    pp. {65-102},
    {(1983)}.
    
    
    \bibitem{ledoux-c}
    M. Ledoux. The concentration of measure phenomenon.
   {\em AMS Monographs, Providence} (2001).
   
   
   
  \bibitem{legall}
J.F. Le Gall, Spatial branching processes, random snakes, and partial different equations. {\em lectures given in the Nachdiplomsvorlesung, ETH Z\"urich}. Published in the series Lectures in Mathematics ETH Z\"urich, Birkhauser (1999).




\bibitem{maler0}
    {R. Mahler}.
   {A theoretical foundation for the Stein-Winter Probability Hypothesis Density (PHD) multi-target tracking approach}.
     {\em Proc. MSS Nat'l Symp. on Sensor and Data Fusion},
     vol. {1},
      {(2000)}.
  



\bibitem{maler}
  {R. Mahler}.
    {Multi-target Bayes filtering via first-order multi-target moments}.
     {\em IEEE Trans. Aerospace \& Electronic Systems},
     vol. {39},
      no. {4},
    pp. {1152-1178},
     {(2003)}.

\bibitem{Me:asympt}
S.~M\'el\'eard.
\newblock Asymptotic behaviour of some interacting particle systems;
  {McKean}-{Vlasov} and {Boltzmann} models.
\newblock In D.~Talay and L.~Tubaro, editors, {\em Probabilistic Models for
  Nonlinear Partial Differential Equations, Montecatini Terme, 1995}, Lecture
  Notes in Mathematics 1627. Springer-Verlag, Berlin (1996).

\bibitem{MR99g:60103}
S.~M{\'e}l{\'e}ard.
\newblock Convergence of the fluctuations for interacting diffusions with jumps
  associated with {B}oltzmann equations.
\newblock {\em Stochastics Stochastics Rep.}, 63(3--4):195--225 (1998).

\bibitem{MR2002a:60164}
S.~M{\'e}l{\'e}ard.
\newblock Probabilistic interpretation and approximations of some {B}oltzmann
  equations.
\newblock In {\em Stochastic models (Spanish) (Guanajuato, 1998)}, volume~14 of
  {\em Aportaciones Mat. Investig.}, pages 1--64. Soc. Mat. Mexicana, M\'exico,
  (1998).

\bibitem{MR98m:82064}
S.~M{\'e}l{\'e}ard.
\newblock Stochastic approximations of the solution of a full {B}oltzmann
  equation with small initial data.
\newblock {\em ESAIM Probab. Stat.}, 2:23--40 (1998).


\bibitem{MR36:4647}
H.P. McKean, Jr.
\newblock A class of {M}arkov processes associated with nonlinear parabolic
  equations.
\newblock {\em Proc. Natl. Acad. Sci. U.S.A.}, 56:1907--1911 (1966).


\bibitem{MR38:1759}
H.P. McKean, Jr.
\newblock Propagation of chaos for a class of non-linear parabolic equations.
\newblock In {\em Stochastic Differential Equations (Lecture Series in
  Differential Equations, Session 7, Catholic University, 
1967)}, pages 41--57. Air
  Force Office of Scientific Research, Arlington, VA (1967).
  
  
     
   \bibitem{ulam}
   N. Metropolis and 
S. Ulam. The Monte Carlo method.
{\em Journal of the American Statistical Association}, No.  247, Vol. 44, pp.  335-341 (1949).

  
  \bibitem{pollard}
  D. Pollard. {\em Empirical processes : theory and applications}. NSF-CBMS Regional Conference Series in Probability and Statistics. Publication of the Institute of Mathematical Statistics, Hayward, CA, vol. 2 (1990).
  
  \bibitem{picard}
    {J. Picard}.
    {Approximation of the non linear filtering
problems and order of convergence.}
    {\em Lecture Notes in Control and Inf. Sc},
    vol. {61},
    {(1984)}.


\bibitem{MR98k:35214}
M.~Pulvirenti.
\newblock Kinetic limits for stochastic particle systems.
\newblock In {\em Probabilistic Models for Nonlinear Partial Differential
  Equations (Montecatini Terme, 1995)}, volume 1627 of {Lecture Notes in
  Mathematics}, pages 96--126. Springer, Berlin (1996).

\bibitem{rachev}
S.T. Rachev.
\newblock {\em Probability Metrics and the Stability of Stochastic Models}.
\newblock Wiley, New York (1991).


\bibitem{rio}
Rio E., {\em Local invariance principles and their applications to density estimation}.  Probability and related Fields, vol. 98,  21-45 (1994).


\bibitem{rousset}
    {M. Rousset},
    {On the control of an interacting particle
approximation of Schroedinger ground states},
    {\em SIAM J. Math. Anal.},
    vol. {38},
    no. {3},
    pp. {824--844},
    {(2006)}.
  
  \bibitem{talagrand-1}
  M. Talagrand. Concentration of measure and isoperimetric inequalities in product spaces. {\em Inst. Hautes \'Etudes Sci. Publ. Math.}, No. 81, pp. 73-205  (1995).
    \bibitem{talagrand-2}
M. Talagrand. A new look at independence. {\em Ann. Probab.}, vol. 24, pp. 1-34  (1996).
    
  \bibitem{talagrand-3}
    M. Talagrand. New concentration inequalities in product spaces. {\em Invent. Math.}, vol. 125, pp. 505-563  (1996).    
    
\bibitem{island2}
Z. Skolicki and K. De Jong. The influence of migration sizes
and intervals on island models. 
{\em In Proceedings of the Genetic
and Evolutionary Computation Conference (GECCO-2005).
ACM Press}, (2005).

\bibitem{david}
D. Steinsaltz and S. N. Evans. Markov mortality models: implications of quasistationarity and varying initial distributions. {\em Theor. Pop. Biol.}, vol. 65, pp. 319-337 (2004).


\bibitem{stettner}
W. Runggaldier, L. Stettner
Approximations of Discrete Time Partially Observed Control Problems
Applied Mathematics Monographs CNR, Giardini Editori, Pisa (1994).

\bibitem{Sz:topi}
A.S. Sznitman.
\newblock Topics in propagation of chaos.
\newblock In P.L. Hennequin, editor, {\em Ecole d'Et\'e de Probabilit\'es de
  Saint-Flour XIX-1989}, Lecture Notes in Mathematics 1464. Springer-Verlag,
 Berlin (1991).

\bibitem{wellner}
A.N.~Van der Vaart and J.A. Wellner.
\newblock {\em Weak Convergence and Empirical Processes with Applications to
  Statistics}.
\newblock Springer Series in Statistics. Springer, New York (1996).



\bibitem{vanneste}
J. Vanneste.
Estimating generalized Lyapunov exponents for products of random matrices
{\em Phys. Rev. E}, vol. 81, 036701 (2010).

\bibitem{vila}
J.-P. Vila, V. Rossi. Nonlinear filtering in discrete time: A particle convolution
approach, Biostatistic group of Montpellier, {\em Technical Report 04-03,
(available at http://vrossi.free.fr/recherche.html)}, (2004).

\bibitem{island1}
D. Whitley, S. Rana, and R. B. Heckendorn. 
The island model
genetic algorithm: On separability, population size and
convergence. 
{\em Journal of Computing and Information
Technology,} 7(1):33?47, (1999).

\end{thebibliography}
\end{document}